\pdfoutput=1
\documentclass[11pt]{article}

\usepackage[T1]{fontenc}
\usepackage[utf8]{inputenc}
\usepackage{lmodern}
\usepackage[margin=1.05in]{geometry}
\usepackage{amsmath,amssymb,amsthm,mathtools}
\usepackage{microtype}
\usepackage{booktabs,longtable,array}
\usepackage{enumitem}
\usepackage{natbib}
\usepackage{xcolor}
\definecolor{linkblue}{HTML}{214D72}
\usepackage{hyperref}
\hypersetup{colorlinks=true,linkcolor=linkblue,citecolor=linkblue,urlcolor=linkblue,
 pdftitle={Exact Simes inequalities for two-sided Gaussian tests},
 pdfauthor={Xianyang Zhang},pdfsubject={Global-null two-sided Gaussian Simes inequality}}
\usepackage{bookmark}
\setlist{nosep}
\allowdisplaybreaks[2]
\newcommand{\Prob}{\mathbb P}
\newcommand{\Expect}{\mathbb E}
\newcommand{\Corr}{\mathcal C}
\newcommand{\barPhi}{\overline\Phi}
\theoremstyle{remark}
\theoremstyle{plain}

\numberwithin{equation}{section}
\title{Exact Simes inequalities for two-sided Gaussian tests}
\author{Xianyang Zhang\thanks{Department of Statistics, Texas A\&M University. Email: \href{mailto:zhangxiany@stat.tamu.edu}{zhangxiany@stat.tamu.edu}.}}
\date{}
\newtheorem{theorem}{Theorem}[section]
\theoremstyle{plain}
\newtheorem{lemma}[theorem]{Lemma}
\theoremstyle{plain}
\newtheorem{proposition}[theorem]{Proposition}
\theoremstyle{plain}
\newtheorem{corollary}[theorem]{Corollary}
\theoremstyle{remark}

\theoremstyle{plain}
\newtheorem{lmtheorem}[theorem]{Theorem}
\theoremstyle{plain}

\theoremstyle{plain}
\newtheorem{sltheorem}[theorem]{Theorem}
\theoremstyle{plain}
\newtheorem{sllemma}[theorem]{Lemma}
\theoremstyle{remark}
\newtheorem{slremark}[theorem]{Remark}
\theoremstyle{plain}
\newtheorem{tritheorem}[theorem]{Theorem}
\theoremstyle{plain}
\newtheorem{trilemma}[theorem]{Lemma}
\theoremstyle{plain}
\newtheorem{triproposition}[theorem]{Proposition}
\theoremstyle{plain}

\theoremstyle{plain}

\theoremstyle{remark}

\theoremstyle{plain}
\begin{document}
\maketitle
\begin{abstract}
The Simes test combines marginal p-values to test a global null
hypothesis. Whether it controls the type I error for two-sided Gaussian
tests with arbitrary covariance remains open. We prove that, for every
fixed number $m$ of tests, the rejection probability does not exceed
the nominal level when that level is sufficiently small. The same range
of levels works for every covariance matrix. Equality holds if and only
if the p-values form independent groups with identical values within
each group. For $m=3$, a computer-assisted proof establishes validity
at $0.05$, and an analytic argument extends it to an explicit interval
around $0.05$.

As the number $m$ of tests increases, we prove validity in two settings: each test statistic
is correlated with at most a fixed number of others; or the statistics
form independent blocks whose sizes grow sufficiently slowly, at a
specified rate of order $\log\log m$. Both settings allow arbitrary
correlations within these restrictions. For any fixed $0<a<1$, validity
holds for every $m$ above a threshold $M$ and every level
$0<\alpha\le a$; the same $M$ works for all allowed covariances.
Further results quantify the difference between the nominal level and
the rejection probability, and permit more correlated neighbors when
their correlations are bounded away from $\pm1$. We do not obtain
numerical values for the general small-level cutoffs or the required
number of tests.
\end{abstract}

\section{Introduction}\label{sec:introduction}

In multiple hypothesis testing, a natural first question is whether all
the null hypotheses are true. The Simes test evaluates this global null
hypothesis by comparing the ordered marginal p-values with a sequence of
linearly increasing thresholds. For independent uniform p-values, the
probability of rejecting a true global null equals the nominal level
\citep{simes1986}. Several positive-dependence conditions also ensure
that this type I error probability does not exceed the nominal level
\citep{sarkar1998,by2001}. For two-sided Gaussian tests, however, it
remains unknown whether the same conclusion holds for every covariance
matrix. We establish it in several settings that allow both positive and
negative correlations.

Let $Z=(Z_1,\ldots,Z_m)\sim N_m(0,\Sigma)$, where $m$ is the number of
tests and $\Sigma$ is a positive semidefinite correlation matrix. The
two-sided p-values are $P_i=2\barPhi(|Z_i|)$, where $\barPhi=1-\Phi$
is the standard normal survival function. Write $P_{(1)}\le\cdots\le P_{(m)}$
for their ordered values. The rejection probability at level $\alpha$ is
\begin{equation}\label{eq:target}
 F^{(m)}_\Sigma(\alpha)
 =\Prob\left\{P_{(k)}\le\frac{k\alpha}{m}
                  \text{ for some }1\le k\le m\right\},
 \qquad 0<\alpha<1.
\end{equation}
When the dimension is clear, we also write $F_\Sigma=F_\Sigma^{(\dim\Sigma)}$.
We ask whether the test controls type I error at the nominal level,
that is, whether $F^{(m)}_\Sigma(\alpha)\le\alpha$ for every
$m$, $\alpha$, and $\Sigma$. We call the test valid at a given level
and covariance when this inequality holds. All null hypotheses are true in this formulation; our
results do not assert false discovery rate (FDR) control when some null
hypotheses are false. Singular covariance matrices are allowed.

Our first result concerns a fixed number of tests. For each $m$, there
is a number $a_m>0$ such that the Simes test is valid at every level
$0<\alpha\le a_m$, for every covariance matrix
(Theorem~\ref{main:uniform}). Thus the range of valid levels can be chosen
without knowing the covariance. We also characterize equality:
$F^{(m)}_\Sigma(\alpha)=\alpha$ if and only if the p-values can be divided
into independent groups, with identical p-values within each group.
For Gaussian statistics, the latter occurs exactly when the statistics in a
group agree up to sign. Independent p-values and repeated copies of one
p-value are two examples. Every other covariance matrix gives a rejection
probability strictly below $\alpha$ in this range.

The theorem does not give a numerical value for $a_m$, so it does not
by itself establish validity at a conventional level such as $0.05$.
For $m=3$, we address that level directly: a computer-assisted argument
proves validity for every covariance matrix at $0.05$, and an analytic
argument extends it to an explicit interval around $0.05$
(Theorem~\ref{main:tri-theorem} and Corollary~\ref{main:tri-interval}).
The equality cases are the same. We do not know whether this interval
overlaps the small-level interval, so we do not establish validity at every level
in $(0,0.05]$.

Our remaining results allow the number of tests to increase. We
consider two settings. In the first, each statistic is correlated with
at most $D$ other statistics, where $D$ is fixed as $m$ increases
(Theorem~\ref{main:fixed_degree}). Correlations may have either sign and
may be arbitrarily close to, or equal to, $\pm1$. This condition allows
a chain of correlations to connect arbitrarily many statistics; it does
not require small independent blocks.

In the second setting, the statistics form mutually independent blocks,
with unrestricted covariance within each block. The maximum block size
may grow at a rate of order $\log\log m$, with a sufficiently small
coefficient (Theorem~\ref{main:growing_blocks}
and~\eqref{main:growing_simple}). In either
setting, fix any $0<a<1$ and the stated restrictions on dependence.
There is then a threshold $M$ such that validity holds for every $m\ge M$
and every $0<\alpha\le a$. The same $M$ works for all these levels and
all allowed covariances.

We also study the difference between $\alpha$ and the rejection probability.
When all pairwise correlations have magnitude at most a fixed
$r<1$, we approximate this difference and control the error
even when the correlations are very small. This permits the number of
correlated neighbors to grow at a specified rate of order
$\log m/\log\log m$ (Theorem~\ref{main:graph_theorem}). The full theorem
also permits groups of exact copies of bounded size, defined in
Section~\ref{main:physical_weights}. Theorem~\ref{main:mixed_theorem}
combines this setting with additional mutually independent blocks of
bounded size; those blocks are also independent of the first set of
statistics. The growing-block result applies to the separate setting
consisting entirely of independent blocks.

These results include covariance matrices that fail positive regression
dependence on a subset (PRDS), a standard sufficient condition for the
Simes inequality. Section~\ref{main:related_work} explains this condition,
gives examples where it fails, and compares our results with earlier
work that also establishes validity beyond PRDS. We begin with the
small-level theorem. The section
\hyperref[comparison:guide]{``How to read the proofs''} before
Appendix~\ref{sl:section} summarizes the appendices and explains how they fit together.

\section{Uniform validity at small levels}
\label{main:smalllevels}

For a fixed number $m$ of tests, we seek a range of levels that ensures
validity for every covariance matrix. Let $\Corr_m$ be the compact set
of all $m\times m$ positive semidefinite correlation matrices. We first
describe the matrices that give equality. Let $\mathcal E_m$ consist of
those matrices for which the statistics form independent groups, with
each statistic in a group equal to the same standard normal variable or
its negative. The two-sided p-values are identical within each group.
Equivalently, after reordering the statistics and changing selected
signs, the covariance matrix is
$\operatorname{diag}(J_{k_1},\ldots,J_{k_r})$, where $J_k$ is the
$k\times k$ all-ones matrix and the positive integers $k_1,\ldots,k_r$
sum to $m$. The theorem gives a common interval of levels over
$\Corr_m$ and shows that $\mathcal E_m$ contains all the equality cases.

\begin{theorem}[Covariance-uniform small-level validity]
\label{main:uniform}
For every fixed integer $m\ge2$, there exists $a_m\in(0,1)$ such that
\begin{equation}\label{main:uniformeq}
 \sup_{\Sigma\in\Corr_m}F^{(m)}_\Sigma(\alpha)=\alpha
 \qquad\text{for every }0<\alpha\le a_m.
\end{equation}
At every such level, $F^{(m)}_\Sigma(\alpha)=\alpha$ if and only if
$\Sigma\in\mathcal E_m$.
\end{theorem}

The cutoff in Theorem~\ref{main:uniform} depends on the number of tests
but not on their covariance. The covariance may therefore change with
$\alpha$ and approach a singular matrix at any rate. Throughout the
common interval, independent groups of signed copies give equality;
every other covariance gives a strictly smaller rejection probability.
This is an inequality at each level, so it determines the sign of the
error that an asymptotic approximation could leave unresolved.
Section~\ref{main:related_small} explains how it strengthens the
first-order result and small-level assertion of \citet{chen2026}.

The cutoff $a_m$ may depend on $m$, so Theorem~\ref{main:uniform} does
not give a common cutoff when the dimension is unrestricted. Such a
cutoff can be obtained for the graph and block classes studied
later, provided the class does not depend on the testing level.
For each such class, the later theorem gives validity at every level
up to a chosen ceiling, in every dimension at or above a common
threshold. Below that threshold, only finitely many dimensions
remain. Taking the minimum of that ceiling and their small-level
cutoffs gives one
positive cutoff for every dimension in the class.

\subsection{Why repeated p-values enter the proof}

The proof compares nearly identical statistics with exact copies of a
single statistic. This reduces the number of distinct Gaussian
variables, but it does not reduce the number of tests. For example,
$(U,U,V)$ contains three p-values, with thresholds $\alpha/3$,
$2\alpha/3$, and $\alpha$, even though only two variables describe it.
The proof must retain these repetitions because they affect which
ranks can cross a threshold.

To record the repetitions, let $X\sim N_q(0,\Gamma)$ with
$\Gamma\in\Corr_q$, and repeat $X_i$
exactly $w_i$ times, where the positive integers $w_1,\ldots,w_q$ sum
to $m$. We call the $q$ variables representatives. The $m$
entries after repetition are the individual tests, or physical
coordinates.
Write $F_w(\Gamma,\alpha)$ for the rejection probability with the
original thresholds $k\alpha/m$. A sign change leaves a two-sided
p-value unchanged, so the notation also covers signed copies.

Independent representatives give exact equality even when their
p-values are repeated:
\begin{equation}\label{main:duplicateequality}
 F_w(I_q,\alpha)=\alpha.
\end{equation}
To see why, let the rejection count be the largest rank at which a
Simes threshold is crossed, or zero if there is no crossing. Assign
each group its fraction of the rejected tests, with share zero when
there is no rejection. These shares sum to the indicator that Simes
rejects. The calculation uses a leave-one-group-out argument: temporarily
set every p-value in a selected group to zero while keeping all $m$
tests in the rule. The resulting rejection count depends only on the
other groups and is independent of the selected representative.
Restoring the selected group's uniform p-value then gives its expected
share $w_i\alpha/m$ under the original distribution. Summing over groups gives
\eqref{main:duplicateequality}. Lemma~\ref{sl:pairidentity} gives the
calculation and the corresponding two-representative inequality.

We can now keep $m$ fixed and induct on the number $q$ of
representatives. The induction proves validity and the equality
classification for every choice of copy counts summing to $m$.
Representatives that already agree up to sign can be merged by adding
their copy counts; this leaves the rejection event unchanged.
For nearly identical representatives, we first align their signs,
then take their weighted mean, divide it by its standard deviation,
and repeat this new
representative with their combined copy count. We call this operation
contraction. It reduces $q$ while preserving $m$, which is why the
induction must cover repeated p-values. The equality class for the
representative covariance is $\mathcal E_q$.

\subsection{Three comparisons that complete the proof}

The difficult case is a covariance matrix close to one that gives
equality. Its rejection probability may be arbitrarily close to
$\alpha$, so an error bound tending to zero need not determine whether
the probability is above or below $\alpha$. The proof uses three
comparisons, according to how close the statistics are to exact copies
relative to the spacing between rejection thresholds. We describe that
scale first, then explain why the three cases rule out a failure of
the theorem. Full proofs are in Appendix~\ref{sl:section}.

Suppose that validity or the stated equality classification fails
arbitrarily close to level zero. For some fixed copy counts, we can then
choose levels $\alpha_n\downarrow0$ and correlation matrices
$\Gamma_n\notin\mathcal E_q$ such that
$F_w(\Gamma_n,\alpha_n)\ge\alpha_n$. The weak inequality includes both
a violation of validity and an additional equality case.
On $\mathcal E_q$, equality already follows from
\eqref{main:duplicateequality} after merging signed copies.

Put $b_n=\Phi^{-1}(1-\alpha_n/2)$. The Gaussian rejection cutoff
at Simes rank $k$ is $\Phi^{-1}(1-k\alpha_n/(2m))$.
For fixed $m$, the gaps between these cutoffs are of order $1/b_n$.
After changing the sign of one statistic if necessary, the difference
between the two statistics of a pair $(i,j)$ has standard deviation
$\sqrt{2(1-|\Gamma_{ij,n}|)}$. We compare this standard deviation,
which describes how much the two statistics typically differ, with
the gap between rejection cutoffs. The quantity
\[
 b_n^2(1-|\Gamma_{ij,n}|)
\]
is one half of the squared ratio of this standard deviation to $1/b_n$.
We call it the rescaled distance. A positive finite limit means that
the typical difference and the cutoff gap have the same order; a zero
limit means that the typical difference becomes negligible compared
with the gap.

We now select a subsequence of the levels and correlation matrices
$(\alpha_n,\Gamma_n)$ chosen above, keeping the same copy counts,
so that every pair's rescaled distance converges to a limit in
$[0,\infty]$. This is possible because there are only finitely many
pairs. We continue to use $n$ to index this subsequence.

Two representatives belong to the same cluster when their limiting
rescaled distance is finite. To see why this rule forms well-defined
groups, consider three representatives $i,j,k$. The square roots of
the rescaled distances satisfy the triangle inequality: the distance
from $i$ to $k$ is at most the sum of the distances from $i$ to $j$
and from $j$ to $k$, allowing sign changes. Thus, if pairs $(i,j)$ and $(j,k)$ have finite limiting rescaled
distances, so does $(i,k)$.
In particular, two clusters cannot overlap without being the same
group. The unit-vector representation of correlations in Appendix~\ref{sl:section}
gives this triangle inequality. Rescaled distances between different
clusters tend to infinity, although the clusters need not be independent.

Within each cluster, choose one representative as an anchor. For
sufficiently large $n$, change the sign of each remaining
representative, if needed, so that its
correlation with the anchor is positive. Every representative is then
close to the same anchor, which makes all pairwise correlations in the
cluster tend to one. Thus one sign is chosen for each statistic,
and these choices work for all pairs in the cluster. There are only finitely many sign choices, so
we can select a further subsequence on which the signs stay fixed.
These changes leave the two-sided p-values unchanged.

\paragraph{1. Every pair has diverging rescaled distance.}
No two representatives remain close on the threshold-spacing scale.
The rejection probability can be written exactly as $\alpha$, minus
one nonnegative term for each pair of representatives, plus an error
(Lemma~\ref{sl:pairidentity}). Each pair term is the decrease in
rejection probability caused by that pair's correlation when all other
p-values are set to one; the sum of the pair terms is the pair
contribution. The error is what remains, and it arises only
when three or more representatives have small p-values at the same
time. We show that the pair contribution dominates the error, giving
strict inequality whenever the representatives are not independent.

This comparison must remain valid even when the pair contribution
is arbitrarily small. If every correlation is at most a fixed small
bound in absolute value, the lower bound
for the pair contribution and the upper bound for the error are both
multiples of the sum of the squared correlations, and this sum cancels
when the two bounds are compared. Otherwise, we compare the error with the decrease caused by
a pair having the largest absolute correlation. Even if this
correlation approaches one in magnitude while its
rescaled distance diverges, the bounds for this decrease and for the
error retain their dependence on the distance from perfect correlation.
In both situations the error is smaller than the pair contribution. Lemma~\ref{sl:separated} makes this comparison precise
and gives the eventual strict inequality.

\paragraph{2. A cluster retains a positive rescaled distance.}
Suppose two members of a cluster have a finite positive limiting
distance. Their differences remain comparable to the spacing between
rejection thresholds. In this case a limiting calculation already
detects a strict decrease: the rejection probability divided by
$\alpha_n$ has a limit below one.

Count every repeated test in the cluster, giving $s$ entries with
covariance matrix $\Sigma_n$. After aligning signs, write
$b_n^2(1-\Sigma_{ij,n})\to L_{ij}<\infty$.
Condition on the first statistic being near a rejection threshold.
The rescaled positions of the other statistics then have a Gaussian
limit described by $Y$, where
\[
 \Expect Y_i=-L_{1i},\qquad
 \operatorname{Cov}(Y_i,Y_j)=L_{1i}+L_{1j}-L_{ij},
 \qquad Y_1=0,
\]
and the cluster's normalized contribution converges to
\begin{equation}\label{main:criticalfunctional}
 \Psi_s(L)
 =\Expect\max_{1\le k\le s}\frac{k}{s}e^{Y_{[k]}}\le1,
\end{equation}
where $Y_{[1]}\ge\cdots\ge Y_{[s]}$. The maximum accounts for a
crossing at any rank within the cluster. Exact copies give the value
one. For the upper bound, the maximum in
\eqref{main:criticalfunctional} is at most
$s^{-1}\sum_i e^{Y_i}$, whose expectation is one because
$\Expect e^{Y_i}=1$ for each $i$. The expected maximum is strictly
smaller than one if some $L_{ij}>0$. The proof establishes this
strictness and justifies taking the limit inside the expectation,
including for singular Gaussian laws.

The cluster containing $s$ tests receives weight $s/m$ in the total
limit. These weights sum to one. Therefore, if one cluster's
normalized contribution is strictly below one and all the others are
at most one, the limiting total rejection probability divided by
$\alpha_n$ is strictly below one. This contradicts the chosen
sequence; Lemmas~\ref{sl:critical} and~\ref{sl:multicritical} provide
the details. The Gaussian extreme-value scale and limiting laws are classical
\citep{husler1989,engelke2015}; here they identify the Simes limit and
its equality cases. If every finite $L_{ij}$ is zero,
however, the limit equals one and a more precise comparison is needed.

\paragraph{3. All within-cluster rescaled distances tend to zero.}
The remaining case has $G<q$ clusters. Within each cluster the
statistics approach exact copies so quickly that the preceding limit
equals one. We therefore compare rejection probabilities before taking
a limit. Replace each cluster by exact copies of one standardized
Gaussian variable. Consider the indices for which at least one
within-cluster pair is not yet an exact copy up to sign, meaning that
its absolute correlation is less than one. Along these indices,
the original vector eventually has a smaller rejection probability
than this contracted vector, with a quantified positive difference.

After aligning signs, let $\delta_n$ be the largest within-cluster
value of one minus a correlation. First take indices with
$\delta_n>0$ and put $\lambda_n=b_n^2\delta_n\to0$.
For a cluster $C$, average all its repeated tests,
$\sum_{i\in C}w_iX_{i,n}/\sum_{i\in C}w_i$, and divide the result by
its standard deviation. Repeat this representative with the cluster's
total copy count. The contracted vector has $G$ representatives,
covariance $\Gamma_n^c$ and copy counts $w'$, but still $m$ tests.
Lemma~\ref{sl:contraction} gives
\begin{equation}\label{main:contraction}
 F_w(\Gamma_n,\alpha_n)
 \le F_{w'}(\Gamma_n^c,\alpha_n)
       -\frac{\alpha_n\sqrt{\lambda_n}}{160m^2}
 \qquad\text{eventually}.
\end{equation}
The proof builds both vectors from the same Gaussian variables and
compares their rejection decisions. Call a rejection lost when the
contracted vector rejects and the original does not, and gained in
the reverse case. A narrow interval around a rejection threshold gives
loss probability bounded below by a positive constant times
$\alpha_n\sqrt{\lambda_n}$.
A change of Gaussian measure bounds the gain by $\alpha_n$ times a
factor exponentially small in $1/\lambda_n$. The loss therefore
dominates. To allow dependence between clusters, the proof combines
the geometry of unit vectors representing the covariance with
conditioning on one Gaussian coordinate. Singular cluster covariances
are included.

The indices with $\delta_n=0$ require a separate observation. All
statistics within each cluster are then already exact signed copies,
and contraction leaves the rejection event unchanged.

\paragraph{Completing the induction.}
There are fewer representatives after contraction. The induction
hypothesis therefore bounds the contracted rejection probability by
$\alpha_n$ for all sufficiently large $n$. The uniform level cutoff
in the induction hypothesis is essential: $\Gamma_n^c$ can change with $n$ and can be
singular. If $\delta_n>0$, subtracting the positive term in
\eqref{main:contraction} gives strict inequality for the original
vector. If $\delta_n=0$, contraction leaves its rejection probability
unchanged, but the contracted covariance cannot be an equality case.
Otherwise expanding its copies would put the original covariance in
$\mathcal E_q$. The equality part of the induction hypothesis gives
strictness in this case too.

Each of the three comparisons contradicts the proposed sequence.
Finally, there are only finitely many positive integer copy-count
vectors with total $m$. Taking the minimum of their positive cutoffs
gives the common cutoff in Theorem~\ref{main:uniform}. The full induction appears in
Section~\ref{sl:inductionsubsection}.

\section{The trivariate inequality at level 0.05}
\label{main:trivariate}

For three tests, we now address the conventional level $1/20$ directly.
We prove validity for every correlation matrix, including singular
matrices, and identify all equality cases. Analytic bounds cover
covariance matrices near the equality cases; a finite computation
covers the remaining covariance matrices. A further analytic argument extends the conclusion and its
equality classification to a common interval around $1/20$.

For three tests, Simes rejects when at least one p-value is at or below
the first threshold, at least two are at or below the second, or all
three are at or below the third. Let
$N(u)=\sum_{i=1}^3\mathbf1\{P_i\le u\}$ count the p-values at or
below $u$. The rejection indicator is
\begin{equation}\label{main:tri-event}
 H=\mathbf1\{N(1/60)\ge1
       \ \text{or}\ N(1/30)\ge2
       \ \text{or}\ N(1/20)=3\}.
\end{equation}
\begin{theorem}[Trivariate validity at level $1/20$]
\label{main:tri-theorem}
Let $Z\sim N_3(0,\Sigma)$ for any positive semidefinite correlation
matrix $\Sigma$, and let $P_i=2\barPhi(|Z_i|)$. Then
\begin{equation}\label{main:tri-bound}
 F^{(3)}_\Sigma(1/20)
 =\Prob\{P_{(1)}\le1/60\ \text{or}\ P_{(2)}\le1/30
                         \ \text{or}\ P_{(3)}\le1/20\}
 \le \frac1{20}.
\end{equation}
Writing $J_k$ for the $k\times k$ all-ones matrix, equality holds
exactly at coordinate-sign changes and permutations of
\begin{equation}\label{main:tri-equality}
 I_3,\qquad J_3,\qquad \operatorname{diag}(J_2,1).
\end{equation}
\end{theorem}

\begin{corollary}[Validity on an interval around $1/20$]
\label{main:tri-interval}
For every
\[
 \alpha\in\left[\frac1{20}-2^{-21},\;
                     \frac1{20}+\frac{2^{-20}}6\right]
\]
and every positive semidefinite $3\times3$ correlation matrix $\Sigma$,
\[
 F^{(3)}_\Sigma(\alpha)\le\alpha.
\]
At every level in this interval, equality holds exactly at the
coordinate-sign changes and permutations of the matrices in
\eqref{main:tri-equality}.
\end{corollary}

The three matrices in~\eqref{main:tri-equality} have a simple
interpretation: three independent p-values, three identical p-values,
or two identical p-values independent of the third. Reordering the
tests and changing the signs of the statistics leave the Simes event
unchanged.

Near an equality case, the difference between the rejection probability
and the nominal level can be arbitrarily small. We control these
neighborhoods analytically. For the remaining covariance matrices,
we divide the region into finitely many cells and bound the probability
throughout every cell. The cells and the bounds recorded for them form a
finite certificate; checking that these bounds are valid and that the
cells cover the remaining region proves the inequality there.
Appendix~\ref{tri:proof} gives the argument, and
Section~\ref{tri:reproducibility} explains the arithmetic assumptions
and what the computation verifies. Section~\ref{tri:level-interval}
proves the interval extension analytically.
Corollary~\ref{main:tri-interval} does not establish validity at every
level below $1/20$ or for four tests. Nor does
\eqref{main:tri-bound} cover three representatives with arbitrary copy
counts: their total number of tests and rejection thresholds would
change.

\subsection{Relation to earlier finite-dimensional results}
\label{main:tri-prior}

The two-sided Gaussian inequality for two tests is classical
\citep{samuelcahn1996}. It also follows from the theorem of
\citet{sarkar1998} for p-values that are multivariate totally positive
of order two (MTP$_2$). A joint density $f$ is MTP$_2$ if
$f(x\vee y)f(x\wedge y)\ge f(x)f(y)$, where $\vee$ and $\wedge$ take
componentwise maxima and minima; this is a strong positive-dependence
condition. When the correlation lies in $(-1,1)$, the absolute values
of two Gaussian statistics, and hence their two-sided p-values, are
MTP$_2$ \citep{karlinrinott1981}. A correlation of $\pm1$ gives
identical p-values, so the rejection probability equals $\alpha$.
Positive-dependence theorems also cover some three-test covariance
models \citep{sarkar1998,sarkar2008}, but their assumptions do not hold
for every correlation matrix.
Validity beyond these assumptions has precedents.
\citet[Lemma~3.4 and Corollary~3.5]{yekutieli2008} treats pairwise
differences of independent Gaussian group means with equal variances,
each difference divided by the same standard error. The result
requires explicit conditions on the rejection thresholds and includes
singular three-test examples that fail PRDS. The table in
\citet{yekutieli2008} shows that, for normal statistics at level $0.05$,
these sufficient conditions hold whenever at most 68 group means are
compared.
Theorem~\ref{main:tri-theorem} covers every
three-dimensional correlation matrix at the specified level and
identifies its equality cases.

Two related results address different aspects of dependence.
\citet[Theorem~2]{ristl2016} consider three standard Gaussian
statistics with an arbitrary correlation matrix. For $\alpha\le1/2$,
they show that the probability that at least two of them exceed the
same one-sided threshold $\Phi^{-1}(1-\alpha)$ is at most $\alpha$.
Their proof studies how this probability changes with the correlations.
We also use covariance differentiation, but apply it to the absolute
values and three thresholds in \eqref{main:tri-event}.
\citet{finnerroters2024} show that, for certain covariance matrices, the
absolute values of a Gaussian vector are positively associated but not
MTP$_2$. Positive association means that any two coordinatewise
nondecreasing functions of the vector of absolute values are
nonnegatively correlated. This property does not by itself imply
\eqref{main:tri-bound} for every covariance matrix.

\citet{reinerbenaim2007} studied the FDR of the
Benjamini--Hochberg (BH) procedure \citep{bh1995} for correlated two-sided normal tests, in
models where some null hypotheses may be false. The paper first treats
two hypotheses, with explicit FDR upper bounds and simulations, and
then extends the analysis to larger families.
\citet[p.~407]{benjamini2010} later cited this work as strong evidence
that the procedure controls the FDR for two-sided tests under any
correlation, adding that a complete proof was still missing. This
assessment predates the counterexamples of \citet{dobriban2026} and
\citet{lei2026}, discussed in Section~\ref{main:related_fdr}; both
include false null hypotheses.
Here all null hypotheses are true: every rejection
is a false discovery, so the FDR equals the rejection probability
in~\eqref{eq:target}. Theorem~\ref{main:tri-theorem} establishes this
global-null inequality for three tests at $1/20$.

\subsection{Analytic estimates and the covariance cover}
\label{main:tri-proof-outline}

To describe the covariance region, write
\[
 \Sigma(a,b,c)=
 \begin{pmatrix}1&a&b\\a&1&c\\b&c&1\end{pmatrix},
 \qquad
 \Delta(a,b,c)=1-a^2-b^2-c^2+2abc.
\]
The parameters $a,b,c$ are the pairwise correlations, and
$\Delta(a,b,c)$ is the determinant. Reordering the statistics and
changing their signs leave the Simes rejection event unchanged. We
can therefore put the correlations in one of two regions, called
closed chambers:
\begin{equation}\label{main:tri-chambers}
 0\le d\le b\le a\le1,\qquad c=\varepsilon d,\qquad
 \varepsilon\in\{-1,1\},\qquad \Delta(a,b,c)\ge0.
\end{equation}
Every correlation matrix in the theorem has a representative in these
regions. We retain their common boundary, so proving the bounds on
both closed chambers covers the full problem.

The first analytic estimate shows how far the rejection probability
falls below the nominal level near independence:
\begin{equation}\label{main:tri-cube}
 \max(|a|,|b|,|c|)\le\frac7{16}
 \quad\Longrightarrow\quad
 F^{(3)}_\Sigma(1/20)
 \le\frac1{20}-\frac{a^2+b^2+c^2}{500}.
\end{equation}
We separate the contribution of each correlated pair from the part
involving all three statistics. A Hermite-polynomial expansion lets
us compute finitely many terms and bound the omitted infinite sum.
The negative contribution and its error bound both retain the factor
$a^2+b^2+c^2$. Their comparison therefore gives strictness even
arbitrarily close to independence.

The second estimate covers every feasible matrix with a nearly
perfect pair:
\begin{equation}\label{main:tri-pair}
 \max(|a|,|b|,|c|)\ge\frac{2047}{2048}
\end{equation}
After aligning signs, replace the pair by two exact copies of its
standardized mean, retaining
both tests in the rejection rule. The resulting vector depends on
only two distinct Gaussian variables, so a bivariate calculation
controls its rejection probability. We then compare the rejections
gained and lost when the copies are separated. The estimates cover
every allowed correlation with the third statistic, including full
duplication, and give precisely the corresponding equality cases
in~\eqref{main:tri-equality}. Lemmas~\ref{tri:cube} and~\ref{tri:pair}
give the proofs and constants.

For the remaining covariance region, divide the parameter space into
finitely many rectangular cells. We need an upper bound throughout
each cell, including its feasible singular matrices. The calculation
uses Price's identity, which differentiates Gaussian expectations
with respect to the covariance \citep{Voigtlaender2021Price}. In each cell requiring a numerical
bound, choose a feasible point $\rho_*=(a_*,b_*,c_*)$. We call this
the reference point: we first bound the rejection probability at the
covariance matrix $\Sigma(\rho_*)$, then extend the bound to the other
feasible covariance matrices in the cell. Starting at independence, where the rejection probability
is known, we move along the straight line to $\rho_*$. Price's identity expresses the derivative of the
probability along this line through changes in the rejection indicator $H$ at its thresholds.
Integrating this derivative, with rigorous bounds on numerical error,
gives an interval containing
$F^{(3)}_{\Sigma(\rho_*)}(1/20)-1/20$.

To extend that bound away from the reference point, we also enclose each
partial derivative with respect to a correlation coordinate in an
interval. Each interval bounds that derivative at every positive definite feasible
covariance in the cell, with the same bounds controlling approaches
to a singular matrix. Thus it controls the change in probability
as the corresponding correlation varies.
Write $\mathcal K_3=\{\rho\in[-1,1]^3:\Sigma(\rho)\succeq0\}$.
For a cell $\mathcal B=\prod_i[\ell_i,r_i]$, let $G_i$ be the interval
bounding the $i$th partial derivative. The resulting upper bound is
\begin{equation}\label{main:tri-cell}
 \sup_{\rho\in\mathcal B\cap\mathcal K_3}
 \left\{F^{(3)}_{\Sigma(\rho)}(1/20)-\frac1{20}\right\}
 \le U_*
       +\sum_{i=1}^3
          \sup\{gv:g\in G_i,\ 
                    v\in[\ell_i-\rho_{*,i},r_i-\rho_{*,i}]\},
\end{equation}
where $U_*$ bounds the probability difference at the reference point.
The sum allows for its largest possible increase elsewhere in the cell.

For each chamber, the cells come from repeatedly halving initial cells
that cover it. Exact arithmetic confirms that the pieces of each halved
cell cover it, so the final cells still cover the whole chamber. Each final cell is of one of three
kinds: it contains no feasible point of the chamber, as exact tests of
the ordering and the determinant confirm; its points in the chamber
lie in one of the two analytic neighborhoods; or its upper bound in~\eqref{main:tri-cell} is
at most $-2^{-20}$. Exact arithmetic also confirms that each reference
point lies in its cell and gives a positive semidefinite correlation
matrix. Each numerical bound is checked by computing its intervals again
with rigorous interval arithmetic: every recorded interval must contain
the recomputed one, which contains the exact value. These checks ensure
that no part of the required region is omitted.
Proposition~\ref{tri:finite-certificate} and its proof give the
verification conditions.

The cover also includes singular covariance matrices. If two distinct
feasible matrices have no perfect pair, the interior of the line
segment between them is positive definite. We can therefore integrate
the derivative bound along the segment. Continuity gives the bound at
its endpoints: represent the Gaussian vectors using one common
standard normal vector and the continuous positive semidefinite
square roots of the
covariance matrices. The rejection indicators then converge except
on marginal threshold surfaces, each of probability zero.
When the chosen covariance matrix $\Sigma(\rho_*)$ is singular,
we bound the integrand over the final integration interval to control
that interval's contribution.
Together with the analytic neighborhoods, these checks give strictness
outside~\eqref{main:tri-equality}. Independence and the exact
two-representative calculation give equality for the three covariance
structures in~\eqref{main:tri-equality}. Changing signs and reordering
the statistics then give all the stated equality cases.

\phantomsection\label{main:tri-replay}
Appendix~\ref{tri:proof} gives the certificate checks and their results.
Sections~\ref{tri:s6} and~\ref{tri:reproducibility} explain how
rounding and integration errors are controlled and state the arithmetic
assumptions.
This computer-assisted verification applies only at $1/20$.
The interval in Corollary~\ref{main:tri-interval} follows analytically
from the certified margin and the two analytic neighborhoods, as shown
in Section~\ref{tri:level-interval}.

\section{Pairwise corrections for covariance graphs}
\label{main:large_m}

We now consider a growing number of tests. A covariance graph joins two
statistics when their covariance is nonzero; a statistic's degree is the
number of its correlated neighbors. Independent statistics give
rejection probability exactly $\alpha$. Under the conditions below, the
amount by which the rejection probability falls below $\alpha$ is approximated
by a sum of contributions from correlated pairs. The error from interactions
among more coordinates is smaller than this sum, even when the nonzero
correlations are arbitrarily small. Correlations may have either sign,
and the degree may grow at the rate stated below. The bound holds
simultaneously over all levels below a fixed upper bound.
Section~\ref{main:mixed_blocks} extends the setting to include arbitrarily
many independent blocks, each containing at most a fixed number of
tests, with unrestricted covariance within each block, and proves
validity in this extended setting.

\subsection{Counting repeated p-values and the pair correction}
\label{main:physical_weights}

We first allow groups of exact copies. Two test statistics have the
same two-sided p-value when they are equal or when one is the negative
of the other. We put such statistics in the same group, whether or not
their signs agree. We describe each group by one Gaussian coordinate,
called its representative, and record the number of statistics in the
group as its multiplicity. Specifically, let $X=(X_1,\ldots,X_n)$ be centered
Gaussian with correlation matrix $\Gamma$, and let
\[
 P_i=2\barPhi(|X_i|),\qquad
 w_i\in\{1,2,\ldots\},\qquad \sum_{i=1}^n w_i=m.
\]
The list tested by Simes contains $w_i$ copies of $P_i$. When a group's
p-value passes a threshold, every copy contributes to the count. The resulting rejection probability is
\begin{equation}
 F_m(\Gamma,w;\alpha)
 =\Prob\left\{\sum_{i=1}^n w_i
       \mathbf 1\{P_i\le k\alpha/m\}\ge k
       \text{ for some }1\le k\le m\right\}.
 \label{main:weighted_simes}
\end{equation}
In Section~\ref{main:smalllevels} this probability is written
$F_w(\Gamma,\alpha)$, with $q$ representatives in place of $n$.
Independent representatives satisfy
\begin{equation}
 F_m(I_n,w;\alpha)=\alpha
 \quad\text{for every integer multiplicity vector }w.
 \label{main:weighted_equality}
\end{equation}
To see this, temporarily set all copies of group $i$ to zero and let $R^{[i]}$
be the largest rank crossing a Simes threshold in this modified vector.
Group $i$ is rejected in the original vector precisely when $P_i\le \alpha R^{[i]}/m$; on this event
the original vector's largest crossing rank is also $R^{[i]}$. The group's rejection share
is the fraction of all rejected tests that belong to it, with value zero
when there is no rejection. Conditioning on all other representatives
shows that its expected share is $\alpha w_i/m$ under independence.
Summing these shares gives
\eqref{main:weighted_equality}; the complete identity, including ties, is
proved in Section~\ref{sec:lm:g:2}.

We next measure the contribution of one correlated pair at a time.
Write
\[
 E(\Gamma)=\{\{i,j\}:i<j,\ \Gamma_{ij}\ne0\}
\]
for the edges of the covariance graph on the representatives. For $e=\{i,j\}\in E(\Gamma)$,
let $\Gamma^e$ retain its actual correlation $\Gamma_{ij}$ and set every
other off-diagonal entry to zero. Keep the copy counts $w_i$ unchanged and define
\begin{equation}
 d_e=\alpha-F_m(\Gamma^e,w;\alpha),\qquad
 \mathcal D_m=\sum_{e\in E(\Gamma)}d_e.
 \label{main:complete_budget}
\end{equation}
Each $\Gamma^e$ is a valid correlation matrix. The quantity $d_e$
measures the decrease from the nominal rejection probability when this
pair is the only correlated pair of representatives. Every other test
remains, with independent representatives and the original copy counts.
We call this decrease the pair deficit. It is computed for the full
testing problem, rather than a separate test of the pair. The next
theorem shows that the sum of these deficits gives the leading correction
to the rejection probability.

\subsection{An error bound relative to the pair correction}
\label{main:graph_expansion}

The theorem fixes upper bounds on the number of copies of each
representative, the magnitude of correlations between different
representatives, and the testing level. Let $\Delta_m$ be a prescribed
upper bound on the number of correlated neighbors of any representative
when there are $m$ tests in total. We call the sequence $(\Delta_m)$ the
\emph{degree envelope}; it is fixed independently of the testing level.
This bound may grow with $m$. Faster growth allows more interactions
among three or more representatives linked by correlations, which the
pair deficits do not capture. The pair correction, which subtracts
$\mathcal D_m$ from $\alpha$ to approximate the rejection probability,
can then become less accurate.
The condition
\begin{equation}
 \Delta_m\log(\Delta_m+1)=o(\log m)
 \label{main:degree_envelope}
\end{equation}
is the simplest growth condition we allow. It is equivalent to
condition~\eqref{main:degree_boundary} in the theorem below with
$\beta=0$: it requires the ratio in~\eqref{main:degree_boundary} to
tend to zero, whereas condition~\eqref{main:degree_boundary} with
$\beta>0$ only requires its limit superior to be at most $\beta$. The
theorem therefore also permits faster growth of $\Delta_m$ and
describes how this growth affects the accuracy of the pair
correction.

\begin{theorem}[Complete pair correction]
\label{main:graph_theorem}
Fix an integer $W\ge1$, $0<r_1<1$, and $0<a_1<1$.
Put
\begin{equation}
 g=\frac{1-r_1}{(1+2r_1)(1+r_1)}.
 \label{main:graph_exponent}
\end{equation}
Fix $0\le\beta<(1-r_1)/(1+r_1)$ and prescribe positive integers $\Delta_m$ such that
\begin{equation}
 \limsup_{m\to\infty}
 \frac{(\Delta_m+1)\log(\Delta_m+1)}{\log m}\le\beta.
 \label{main:degree_boundary}
\end{equation}
For each $m$, consider all positive-semidefinite correlation matrices
$\Gamma$ and multiplicities $w$ satisfying
\begin{equation}
 1\le w_i\le W,\quad \sum_iw_i=m,\quad
 |\Gamma_{ij}|\le r_1\ (i\ne j),\quad
 \max_i\#\{j\ne i:\Gamma_{ij}\ne0\}\le\Delta_m.
 \label{main:graph_class}
\end{equation}
Every nonzero edge has $d_e>0$. Put
\[
 \vartheta=\Theta_{r_1}(\beta):=
 \min\left\{g,\frac{1-\beta}{r_1}-\frac2{1+r_1}\right\}>0.
\]
For every $0<\varepsilon<\vartheta$, there is an integer $M$ such that,
for every $m\ge M$, the following bound holds simultaneously for every
$0<\alpha\le a_1$ and every $(\Gamma,w)$ satisfying
\eqref{main:graph_class} with $E(\Gamma)\ne\varnothing$:
\begin{equation}
 \frac{\left|F_m(\Gamma,w;\alpha)-\alpha+\mathcal D_m\right|}
      {\mathcal D_m}
 \le (\alpha/m)^{\vartheta-\varepsilon}.
 \label{main:graph_relative}
\end{equation}
Consequently, $F_m(\Gamma,w;\alpha)=\alpha-\mathcal D_m+o(\mathcal D_m)$
as $m\to\infty$, uniformly over the same levels, correlation matrices
and multiplicities. If $M$ is also chosen so that
$(a_1/M)^{\vartheta-\varepsilon}\le1/2$, then
\begin{equation}
 F_m(\Gamma,w;\alpha)
 \le\alpha-\tfrac12\mathcal D_m<\alpha
 \label{main:graph_sign}
\end{equation}
for every $m\ge M$ and all the same levels, correlation matrices and
multiplicities.

The integer $M$ depends only on $W,r_1,a_1,\beta,\varepsilon$ and
the prescribed sequence $(\Delta_m)$. The covariance matrix may be singular. Its graph,
signs, nonzero magnitudes and multiplicities may all vary jointly with $m$ and $\alpha$,
and no positive lower bound on a nonzero correlation is required.
\end{theorem}

In~\eqref{main:graph_relative}, the exponent is the power
$\vartheta-\varepsilon$ of $\alpha/m$. Since $0<\alpha/m<1$,
a larger positive power gives a smaller bound on the error relative to
$\mathcal D_m$. Here $\vartheta$ is shorthand for the value
$\Theta_{r_1}(\beta)$ at the chosen degree-growth bound $\beta$;
$g$ is the largest value this function takes. The theorem gives
every fixed positive power below $\vartheta$ by allowing any
$0<\varepsilon<\vartheta$.

At $\beta=0$, condition~\eqref{main:degree_boundary} is equivalent to
\eqref{main:degree_envelope}, and $\vartheta=g$. In fact,
$\vartheta=g$ throughout $0\le\beta\le(1-r_1)/(1+2r_1)$, so the
right side of~\eqref{main:graph_relative} is $(\alpha/m)^{g-\varepsilon}$
throughout this range. Thus the stronger condition~\eqref{main:degree_envelope}
is not needed for this exponent. For example, if
$\Delta_m\sim c\log m/\log\log m$ with
$0<c<(1-r_1)/(1+r_1)$, we may take $\beta=c$. The corresponding
bound is $(\alpha/m)^{\Theta_{r_1}(c)-\varepsilon}$.
This statement does not cover $c=(1-r_1)/(1+r_1)$.
The proof of Theorem~\ref{thm:lm:graph_boundary}, the appendix version
of Theorem~\ref{main:graph_theorem}, obtains the bound by counting
connected sets of correlated representatives and bounding their
contributions. If the correlations in a set are weak, its contribution
is bounded by a multiple of the sum of the squared correlations on its
edges, and this sum is compared with the pair deficits. Otherwise, the
contribution is compared with the deficit of an edge with the largest
absolute correlation in that set.
Theorem~\ref{main:mixed_theorem} permits the same degree range under
its stated correlation cap.

The relative error bound ensures that the approximation error is smaller
than the negative pair correction, giving an exact inequality in all
sufficiently large dimensions. An estimate $F_m=\alpha+o(1)$ could leave
either sign of $F_m-\alpha$ at every finite dimension. Even an absolute error that decreases rapidly with $m$ may exceed
$\mathcal D_m$ if the correlations
shrink still faster. Equation~\eqref{main:graph_relative} handles this
case by comparing the error directly with the pair correction.
The degree-bound sequence $(\Delta_m)$ is fixed before choosing $M$.
Even with $W,r_1,a_1,\beta,\varepsilon$ fixed, different sequences
satisfying~\eqref{main:degree_boundary} may require different values
of $M$. The
correlation cap $r_1<1$ and the upper level $a_1<1$ are also fixed.
The level may approach zero arbitrarily quickly; uniformity as either
of these two fixed bounds approaches one is not asserted.

\subsection{Outline of the proof}
\label{main:graph_roadmap}

The proof shows that each correlated pair, taken alone, lowers the
rejection probability and that all remaining effects are smaller than
the sum $\mathcal D_m$ of the pair deficits. This comparison is
delicate because $\mathcal D_m$ can be arbitrarily small when the
correlations are small. An error bound that merely tends to zero is
therefore not enough: every error must be small relative to
$\mathcal D_m$. The proof has four steps.

\paragraph{1. A lower bound for the pair deficits.}
As in the proof of~\eqref{main:weighted_equality}, we write the
rejection probability as a sum of expected rejection shares. For a
single correlated pair, conditioning on the absolute value of one
member writes its deficit as a sum of nonnegative terms. This
calculation involves only the pair, so it does not require PRDS for
the full vector. Consider the probability that both p-values of the
pair are at most given thresholds, as a function of one threshold with
the other held fixed. A bound on its second derivative then gives, for
all sufficiently large $m$,
\begin{equation}
 \mathcal D_m\ge c\,\log(m/\alpha)(\alpha/m)^2
                    \sum_{e\in E(\Gamma)}|\Gamma_e|^2,
 \label{main:budget_lower}
\end{equation}
where $\Gamma_e=\Gamma_{ij}$ is the correlation on the edge
$e=\{i,j\}$ and $c>0$ depends only on the fixed parameters. These
steps appear in Sections~\ref{sec:lm:g:2}--\ref{sec:lm:g:3}.

\paragraph{2. Very small correlations.}
Suppose first that every correlation is at most $(\alpha/m)^4$ in
absolute value. We expand the joint distribution of the absolute
values $|X_i|$ around independence. Changing the sign of one statistic
reverses the signs of its correlations but leaves the absolute values
unchanged, so the first-order term of the expansion vanishes. Ignoring
terms of third and higher order in the correlations, the joint
distribution equals the independent distribution plus the sum of the
changes caused by single correlated pairs. For the rejection event, these single-pair changes add up to
$-\mathcal D_m$. In this regime, Section~\ref{sec:lm:g:4} bounds the
remaining error by $3(2\sum_e|\Gamma_e|^2)^{3/2}$. Because every
correlation is at most $(\alpha/m)^4$, this third-order error is
smaller than the lower bound~\eqref{main:budget_lower}. This step uses
the whole rejection event in~\eqref{main:weighted_simes}, including
rejections at large ranks. Steps~3 and~4 instead set aside rejections
at ranks above a cutoff and bound their probability by an amount that
does not depend on the correlations. When all correlations are this
small, that amount can exceed $\mathcal D_m$, so nothing is set aside
in this step.

\paragraph{3. Rejections at small ranks.}
Otherwise, some correlation exceeds $(\alpha/m)^4$ in absolute value,
and~\eqref{main:budget_lower} gives
$\mathcal D_m\ge c\log(m/\alpha)(\alpha/m)^{10}$. We first consider
rejections at ranks up to a cutoff $K$. Such a rejection depends only
on the p-values at or below $K\alpha/m$. We write the joint
distribution of these small p-values as a sum of terms, one for each
set of representatives. Each term describes the dependence within its
set that is not already described by its subsets; these terms are
signed measures, so their contributions may be positive or negative.
Single representatives contribute nothing. The term for a set that is
not connected in the covariance graph is a product of the terms for
its connected parts, and such products have smaller order. Sets of two
representatives give the pair correction. It remains to show that
connected sets of three or more representatives contribute less than
$\mathcal D_m$.

If the correlations in a connected set are weak, the bound on its
contribution is a multiple of the sum of the squared correlations on
its edges. Summed over all such sets, the bound is a multiple of
$\sum_e|\Gamma_e|^2$, the same sum that appears
in~\eqref{main:budget_lower}. This sum cancels when the two are
compared, so the comparison remains valid however small the
correlations are.

For stronger correlations, we compare a connected set's contribution
with the deficit of an edge having the largest absolute correlation in
that set. Both are controlled by Gaussian tail probabilities. Up to
logarithmic factors, the probability that three statistics with
pairwise correlations at most $r$ in absolute value all have p-values
at most $\alpha/m$ is at most $(\alpha/m)^{3/(1+2r)}$, whereas the
deficit of an edge with correlation $r$ is at least of order
$(\alpha/m)^{2/(1+r)}$. We call these powers tail exponents. Their
difference is
\[
 \frac{3}{1+2r}-\frac{2}{1+r}
 =\frac{1-r}{(1+2r)(1+r)}.
\]
Its lower bound over the allowed correlations, attained at $r=r_1$, gives $g$ in
\eqref{main:graph_exponent}. Thus the contribution of a connected
triple is at most roughly $(\alpha/m)^g$ times the deficit of its
strongest edge. For larger connected sets, we bound the joint tail
probability of a fixed number $J$ of their coordinates; with
correlations at most $r$ in absolute value, its tail exponent is
$J/[1+(J-1)r]$, which increases to $1/r$ as $J$ grows. The number of
connected sets grows with the degree; under
condition~\eqref{main:degree_boundary}, counting them reduces this
exponent by a fraction $\beta$. Subtracting the pair exponent $2/(1+r)$ and
taking the worst case $r=r_1$, large sets give an exponent close to
$(1-\beta)/r_1-2/(1+r_1)$. The smaller of this and $g$ is the exponent
$\Theta_{r_1}(\beta)$ in Theorem~\ref{main:graph_theorem}.
Under~\eqref{main:degree_envelope}, which is equivalent
to~\eqref{main:degree_boundary} with $\beta=0$, counting the connected
sets does not reduce the exponent for large sets. See
Sections~\ref{sec:lm:g:5}--\ref{sec:lm:g:boundary}.

\paragraph{4. Rejections at larger ranks.}
Finally, we restore the ranks beyond the cutoff $K$. A rejection at
such a rank requires many small p-values. We color the covariance graph
so that no two correlated representatives share a color;
representatives with the same color are then independent. An
exponential-moment (Chernoff) bound, using this independence and
counting every copy, shows that a rejection beyond the cutoff has
probability at most a constant times $(\Delta_m+1)(\alpha/m)^{20}$. The
same bound holds for each of the at most $m\Delta_m/2$ single-pair
covariances $\Gamma^e$ that define $\mathcal D_m$. Even summed over all
these distributions, these probabilities are much smaller than the
lower bound for $\mathcal D_m$ in step~3. All errors are therefore
measured against the complete sum of deficits
in~\eqref{main:complete_budget}.

Gaussian covariance interpolation and decompositions at order-statistic
boundaries are established tools. The inequalities of
\citet{LiShao2002} and \citet{DebickiHashorvaJiLing2017} give useful
comparisons when covariances are ordered. Their general bounds involve
arcsines of correlations, with absolute values or only positive changes
retained. These changes are linear near zero correlation. Such bounds
alone do not capture the cancellation of the linear terms for the
two-sided rejection event or determine the decrease in its full
probability. Our remainder estimate retains that cancellation.

\subsection{Examples and covariance restrictions}\label{sec:graph_examples}

The degree condition permits the sum of absolute correlations in a row
to increase with dimension. For an example, consider a star graph: one
central statistic is correlated with $\Delta$ other statistics, called
its leaves, while the leaves have zero covariance with one another.
Take the covariance matrix to be
\[
 \Sigma_\star=
 \begin{pmatrix}
 1 &(c/\sqrt\Delta)\varepsilon^T\\
 (c/\sqrt\Delta)\varepsilon&I_\Delta
 \end{pmatrix},
 \qquad \varepsilon\in\{-1,1\}^\Delta,\quad 0<c<1.
\]
Its eigenvalues are $1-c$, $1+c$, and $1$ with multiplicity
$\Delta-1$. The absolute correlations in the central statistic's row sum to $c\sqrt\Delta$.
Consider disjoint unions of such stars, each with at most $\Delta_m$
leaves, together with independent leftover coordinates. They satisfy
the conditions of Theorem~\ref{main:graph_theorem} with $W=1$, any
fixed correlation cap $r_1$ with $c\le r_1<1$, and any degree envelope
satisfying~\eqref{main:degree_envelope}. The theorem therefore applies
to them: for every $m\ge M$ and $0<\alpha\le a_1$, the rejection
probability is strictly below $\alpha$, and the decrease is
approximated by $\mathcal D_m$ as in~\eqref{main:graph_relative}. Both
component sizes and absolute row sums may diverge, while the
correlations remain uniformly bounded away from perfect dependence.

\section{Independent blocks with unrestricted covariance}
\label{main:mixed_blocks}

The graph expansion requires correlations between distinct
representatives to stay bounded away from perfect dependence. We now
allow unrestricted covariance within independent blocks of bounded
size. There may be arbitrarily many blocks, and each within-block
covariance may approach singularity at any rate. The next theorem also
allows an additional component satisfying the graph conditions of
Section~\ref{main:large_m}. This component and all the blocks are
mutually independent. We call this component a capped graph because its correlation
magnitudes have a fixed upper bound below one.

\begin{theorem}[A capped graph and arbitrary independent blocks]
\label{main:mixed_theorem}\label{mx:theorem}
Fix integers $W,Q\ge1$, a correlation cap $0<\rho_0<1$, and a level
ceiling $0<a_1<1$. Prescribe a degree envelope satisfying
\eqref{main:degree_boundary} for some
$0\le\beta<(1-\rho_0)/(1+\rho_0)$. Consider a centered Gaussian vector formed
from mutually independent parts $G,A_1,\ldots,A_b$. The part $G$ consists
of representatives with multiplicities at most $W$, correlation cap
$\rho_0$, and covariance-graph degree at most $\Delta_m$. Each $A_j$ has
at most $Q$ physical coordinates and an arbitrary positive-semidefinite
correlation matrix. Either type of part may be absent, and $m$ is the
\emph{total} number of physical coordinates.

There exists $M=M(W,Q,\rho_0,a_1,\beta,(\Delta_m))$ such that
\begin{equation}
 F_\Sigma(\alpha)\le\alpha
 \quad\text{for every }m\ge M\text{ and every }0<\alpha\le a_1,
 \label{main:mixed_validity}
\end{equation}
uniformly over all these vectors. Covariance matrices, multiplicities,
block sizes, and the level may vary jointly with $m$. For $m\ge M$ and
$0<\alpha\le a_1$, equality holds precisely when the representatives in $G$ are independent
and every $A_j$ consists of mutually independent groups of exact signed
copies. For this equality class, equality holds for every dimension and
every $0<\alpha<1$.
\end{theorem}

One dimension threshold works for every level below the chosen upper
bound, even when $\alpha$ is arbitrarily small. The cap $\rho_0$ plays the role of $r_1$ in Theorem~\ref{main:graph_theorem}.

\paragraph{Bounded independent blocks.}\label{main:pure_blocks}\label{mx:pure}
Taking the graph component empty in Theorem~\ref{main:mixed_theorem}
gives one threshold $M(Q,a_1)$ for any number of independent Gaussian
blocks of size at most $Q$, with arbitrary positive semidefinite
covariance within each block. For $m\ge M(Q,a_1)$ and
$0<\alpha\le a_1$, validity and the theorem's exact equality
classification hold uniformly. The restriction is on each block's
size; their total size is unrestricted. Negative correlations,
singularity and arbitrarily fast approaches to exact signed copies
are allowed.

The blocks themselves can also grow, provided they do so sufficiently
slowly as the total number of tests increases. The next theorem gives
an explicit growth condition by keeping track of how the proof's error
bounds depend on block size.

\begin{theorem}[Increasing sizes of arbitrary independent blocks]
\label{main:growing_blocks}
Fix $0<a<1$ and prescribe integers $Q_m\ge2$. Put
\begin{equation}\label{main:growing_envelope}
 h_Q=\frac1{16}(31/10)^{-2(Q-1)},\qquad
 \frac{Q_m+\log\log m}{h_{Q_m}\log m}
 \longrightarrow0.
\end{equation}
Suppose that the $m$ physical coordinates form mutually independent
Gaussian blocks, each of size at most $Q_m$, with an arbitrary
positive-semidefinite correlation matrix within each block. There is an
$M=M(a,(Q_m))$ such that, uniformly over these vectors,
\[
 F_\Sigma(\alpha)\le\alpha
 \qquad(m\ge M,\ 0<\alpha\le a).
\]
For $m\ge M$ and $0<\alpha\le a$, equality holds precisely for laws consisting of mutually
independent groups of exact signed copies. These equality laws attain
$F_\Sigma(\alpha)=\alpha$ at every dimension and every $0<\alpha<1$.
\end{theorem}

The same conclusion, with a dimension threshold that may also depend
on $B$, holds for every fixed $B>B_*:=1+\sqrt{17}/2$
if $h_Q$ is replaced by $h_0(B)B^{-2(Q-1)}$, for the positive constant
$h_0(B)$ specified in Corollary~\ref{gb:general_bases}. Consequently,
a simpler sufficient condition for eventual validity and the same
equality classification is
\begin{equation}\label{main:growing_simple}
 \limsup_{m\to\infty}\frac{Q_m}{\log\log m}
 <\frac1{2\log(1+\sqrt{17}/2)}.
\end{equation}
The strict inequality permits choosing one fixed $B>B_*$ before the
dimension threshold. It does not justify taking $B=B_*$ in the bound for $h_Q$,
or the case of equality in \eqref{main:growing_simple}.
The block covariances, their partition, and the level may vary jointly
with $m$. In particular, no lower bound is imposed on nonzero
correlations or on the distance from an exact signed-copy configuration.
Appendix~\ref{gb:section} tracks the effect of block size on the error
bounds. The stated growth condition keeps these errors below the
positive corrections in the proof. Some unspecified increasing bound
on block size already follows from the bounded-block case of
Theorem~\ref{main:mixed_theorem}, including its eventual equality
classification: increase the permitted size only after the theorem
applies at that size. Theorem~\ref{main:growing_blocks} gives an explicit
growth condition instead. It concerns vectors consisting entirely of
independent blocks. The individual block sizes must still be bounded
in the mixed graph theorem, Theorem~\ref{main:mixed_theorem}.
The dimension threshold depends on the entire sequence $(Q_m)$.

The proof controls each block's contribution when all tests are analyzed
together. On an outcome where Simes rejects, assign each block the
fraction of rejected tests that it contains. These fractions sum to
one. Their expectations therefore sum to the probability that Simes
rejects. More formally, if $R$ denotes the number of
BH rejections at level $\alpha$ and $N_A(u)$ counts
the p-values in a block $A$ that do not exceed $u$, define
\[
 b_A=\Expect\left[\frac{N_A(\alpha R/m)}R\mathbf1\{R>0\}\right].
\]
The integrand is defined to be zero when $R=0$.
Under the global null, $F_\Sigma(\alpha)=\sum_A b_A$. We prove
$b_A\le |A|\alpha/m$ for each original block outside the graph component.
Each proof keeps the joint distribution of all tests outside the
selected block unchanged. For Theorem~\ref{main:mixed_theorem}, we also bound the expected share of the graph
component while retaining the other blocks' distribution; changing it
could change the total rejection count. The nominal share for the
graph component is its number of physical tests times $\alpha/m$. For
each correlated pair of representatives in the component, the
single-edge share deficit is the decrease from this nominal share when
that pair is the only correlated pair of representatives in the
component, while the blocks keep their joint
distribution. Lemma~\ref{mx:lc} gives a relative error bound of the
same form as~\eqref{main:graph_relative} for the component's share,
with $\mathcal D_m$ replaced by the sum of these deficits and with
exponent $\Theta_{\rho_0}(\beta)-\varepsilon$. All bounds apply to contributions under the original joint
distribution, so they can be added regardless of the number of blocks.

An individual test need not satisfy the corresponding share bound. In
the validity ranges above, one coordinate can contribute more than its
nominal share, with other coordinates in its block compensating for the
excess. The next proposition gives such an example, even with independent
Gaussian tests outside the block. This is why we bound the contribution
of the whole block. The proposition uses two forms of PRDS. Call an event increasing if raising any p-value cannot
take the vector out of the event. Pointwise PRDS requires
$\Prob(D\mid P_i=u)$ to be nondecreasing in $u$ for every increasing
event $D$ and every $i$; lower-tail PRDS, a weaker condition, requires
the same of $\Prob(D\mid P_i\le u)$. Section~\ref{main:prds}
discusses both. Write $H_m=\sum_{j=1}^m j^{-1}$.

\begin{proposition}[A coordinate can exceed its nominal share]
\label{main:coordinate_share}
For each $0<a_1<1$, there are explicit constants $c_1,c_2,c_3>0$ depending
only on $a_1$ with the following property. For every $m\ge3$ and every
$0<\alpha\le a_1$, put $t=\alpha/m$, $q=c_1t^4/H_m$, and
$\varepsilon=c_2t^{21}/H_m^3$. There is a positive-definite Gaussian
three-coordinate block with correlation matrix
\begin{equation}
 (1-\varepsilon)
 \begin{pmatrix}
 1&1-2q^2&q\\
 1-2q^2&1&-q\\
 q&-q&1
 \end{pmatrix}+\varepsilon I_3
 \label{main:coordinate_matrix}
\end{equation}
whose third coordinate, when combined with any independent outside
p-value vector of length $m-3$, satisfies
\begin{equation}
 \Expect\left[\frac{\mathbf1\{P_3\le\alpha R/m\}}R
                         \mathbf1\{R>0\}\right]
 \ge t+c_3\frac{t^6}{H_m}>t.
 \label{main:coordinate_excess}
\end{equation}
The outside vector need not have independent coordinates. With independent
uniform outside p-values, the resulting Gaussian example violates both
pointwise PRDS and lower-tail PRDS through an increasing event whose
conditioning thresholds lie in $(0,\alpha]$.
\end{proposition}

Section~\ref{mx:coordinate} gives the constants and proof.
Theorem~\ref{main:mixed_theorem} applies to this example. With
independent Gaussian tests outside the block, the vector consists of
one block of size three and $m-3$ blocks of size one, so the
bounded-block case with $Q=3$ gives $F_\Sigma(\alpha)<\alpha$ for every
$m\ge M(3,a_1)$ and $0<\alpha\le a_1$. The third coordinate takes more
than its nominal share, but the proof bounds the share of the whole
block by $3\alpha/m$, so the other two coordinates of the block take
less than $2\alpha/m$ together. Proposition~\ref{main:coordinate_share}
thus shows excess contribution from one test, not a violation of the
Simes inequality. It also shows that the bounded-block case of
Theorem~\ref{main:mixed_theorem}, even for $Q=3$, includes
positive definite examples failing both PRDS conditions in arbitrarily
large dimensions. The covariance and the event witnessing PRDS failure
may depend on $(m,\alpha)$.

\section{Arbitrary correlations in bounded-degree graphs}
\label{main:uncapped_section}

Independent blocks prevent dependence from extending across block
boundaries. We now allow correlations to connect arbitrarily many tests,
provided each test is correlated with at most a fixed number of others.
There is no restriction on correlation magnitudes, including perfect
correlation. Nearly identical tests may also share correlated neighbors,
and these correlations need not tend to zero.

\begin{theorem}[Arbitrary correlations at fixed covariance degree]
\label{main:fixed_degree}
Fix an integer $D\ge1$ and an upper bound on the level $0<a_1<1$. There exists
$M=M(D,a_1)$ such that
\begin{equation}\label{main:fixed_degree_bound}
 F_\Sigma(\alpha)\le\alpha
 \quad\text{for every }m\ge M\text{ and every }0<\alpha\le a_1,
\end{equation}
uniformly over all positive semidefinite correlation matrices satisfying
\[
 \max_i\#\{j\ne i:\Sigma_{ij}\ne0\}\le D.
\]
For $m\ge M$ and $0<\alpha\le a_1$, equality holds precisely when the coordinates form
independent groups of signed exact copies. For these configurations,
$F_\Sigma(\alpha)=\alpha$ in every dimension and at every
$0<\alpha<1$.
\end{theorem}

Here degree counts nonzero covariances. Their signs and magnitudes are unrestricted, and the
covariance matrix may be singular. Once $D$ and $a_1$ are fixed, the
graph, covariance and level may vary jointly with $m$. Combining this
result with Theorem~\ref{main:uniform} gives a single positive level
cutoff that works in every dimension of the fixed-degree class.

Independent Gaussian blocks of any fixed size at most $Q$, with
arbitrary covariance within each block, can be added by replacing $D$
with $\max(D,Q-1)$. Similarly, representative degree at most $d$ and
multiplicity at most $W$ give physical degree at most $W(d+1)-1$.
These bounds follow directly from counting the nonzero covariances.
Theorem~\ref{main:fixed_degree} does not replace the growing-degree
results of Theorems~\ref{main:graph_theorem}
and~\ref{main:mixed_theorem}. Its thresholds do imply validity along
some sufficiently slow increase in degree: each increase can be delayed
until the theorem applies at that degree. This argument provides no
prescribed growth rate.

The proof in Appendix~\ref{fd:section} first groups statistics that are
nearly equal up to sign. It then chooses between two comparisons,
according to the group's spread relative to the spacing of Gaussian
rejection thresholds. When the spread is large enough on this scale,
we make the group independent of the remaining tests. The resulting
error is smaller than the reduction in rejection probability supplied
by the group. For a smaller spread, we replace the group by copies of
one of its members, retaining that member's dependence on all other
tests. This second comparison handles nearly identical tests that share
correlated neighbors.

Both comparisons require a group to be much closer to exact copies
internally than it is to the other tests. After signs are aligned, a
group's internal spread measures how far its correlations fall below
one, and its gap measures how far its correlations with the other tests
fall below one in absolute value. The proof finds the groups from the
distances between unit vectors whose inner products are the
correlations. Each group's spread is less than a small fixed fraction
of its gap, and every gap is bounded below by a constant that depends
only on $D$.

Each comparison changes the joint distribution of all the tests, so we
compare full rejection probabilities at every step. After the last
step, the shares of all parts refer to the same final distribution, and
their bounds can be added as in Section~\ref{main:mixed_blocks}.

In the next proposition, $\delta_m$ is one minus the correlation of the
nearly identical pair.

\begin{proposition}\label{main:fixed_degree_example}
There are connected positive definite covariance graphs of maximum
degree three with a nearly identical pair having a common neighbor of
fixed nonzero correlation, for which
\[
 \delta_m\longrightarrow0,\qquad
 \delta_m\log(m/\alpha)\longrightarrow\infty
 \quad\text{at every fixed }\alpha\in(0,1).
\]
Such a family can also fail both pointwise and lower-tail PRDS.
\end{proposition}

Each such family has maximum degree three, so
Theorem~\ref{main:fixed_degree} applies with $D=3$: for every
$m\ge M(D,a_1)$ and $0<\alpha\le a_1$, the rejection probability is
strictly below $\alpha$. No earlier increasing-dimension theorem in this
paper covers these families. The pair's correlation tends to one without equaling
one, so the pair is neither a group of exact copies nor within a fixed
correlation cap; this rules out Theorem~\ref{main:graph_theorem} and
the capped graph in Theorem~\ref{main:mixed_theorem}. The graph is
connected, so it does not split into independent blocks of bounded or
slowly growing size. As $\delta_m\to0$, the pair keeps a fixed correlation with its common
neighbor. The condition $\delta_m\log(m/\alpha)\to\infty$ means that
the pair does not approach exact copies on the scale of the spacing
between rejection thresholds. This is the situation described at the start of
this section. When the family fails PRDS, the positive-dependence
results discussed in Section~\ref{main:related_work} do not apply
either.

\section{Related work and covariance examples}
\label{main:related_work}

Under independence, the rejection probability in~\eqref{eq:target}
is exactly $\alpha$ \citep{simes1986}. Classical positive-dependence
theorems extend the inequality to important covariance classes
\citep{sarkar1998,sarkar2008,by2001}.
For example, \citet{sarkar1998} proves it when the p-values are
MTP$_2$. For a centered nonsingular Gaussian vector,
\citet{karlinrinott1981} show that the absolute values are MTP$_2$ if
and only if some diagonal sign matrix $\Lambda$ makes $\Lambda\Sigma^{-1}\Lambda$
nonpositive off the diagonal. This is a condition on the
inverse covariance matrix after coordinate sign changes. Nonnegative
entries in $\Sigma$ alone do not imply this condition. Positive
association of $|Z|$ is also distinct from the dependence conditions
used to prove the Simes inequality.

\citet[p.~229]{bodnar2017} report a conjecture that L\"auter presented
in a 2013 conference talk on the basis of extensive simulations: the
Simes inequality holds for two-sided $F$ or beta tests applied to normal
data with any covariance matrix. For the absolute-Gaussian problem
studied here, \citet[Remark~11]{chi2025} note that their results for
monotone transformations of Gaussian vectors do not cover two-sided
p-values, and they describe the Simes inequality for these p-values as
a well-known conjecture.
\citet{finner2017} study the Simes test for p-values based on
exchangeable test statistics. They give sufficient conditions for asymptotic validity
and possible strictness, and a counterexample showing that
exchangeability alone does not ensure validity. Exchangeable Gaussian
statistics have a common correlation. If this common correlation is
nonnegative and less than one, the criterion above makes their absolute values MTP$_2$, so
finite-dimensional validity follows from \citet{sarkar1998}.
Our increasing-dimension results concern sparse covariance graphs and
independent blocks. We cover a common negative correlation only in
fixed dimension, through Theorem~\ref{main:uniform} and, for three
tests, Theorem~\ref{main:tri-theorem}.
Section~\ref{main:tri-prior} compares our results with
\citet{reinerbenaim2007}. The comparisons below distinguish the
probability being bounded, the dependence assumptions, and an exact
inequality from a limiting approximation.

\subsection{Exact validity at small levels}
\label{main:related_small}

For fixed dimension and a fixed Gaussian covariance with no perfect
pair, \citet[Theorem~2]{chen2026} prove the first-order formula
$F^{(m)}_\Sigma(\alpha)=\alpha+o(\alpha)$.
Their supplement, Section~S2.2, also contains an informal assertion of
validity at sufficiently small levels. The argument there uses
asymptotic equivalence, and the detailed proof in Section~S4.4 concludes
with the first-order formula. That formula does not determine the sign
of the error or give a common cutoff when the covariance varies with
the level. Our small-level theorem establishes the exact inequality
and a cutoff uniform over covariance, strengthening the earlier
small-level assertion in these two respects.

The difficulty is greatest near equality. As $\alpha\downarrow0$, a
covariance sequence may approach independence or exact duplication,
so that the rejection probability divided by $\alpha$ tends to one.
This limit still permits a positive error at every level; continuity
and compactness alone do not determine its sign. Near duplication,
the relevant scale is $c_\alpha^2(1-|\Sigma_{ij}|)$, where
$c_u=\Phi^{-1}(1-u/2)$. Gaussian extreme-value limits at this scale
are described by H\"usler--Reiss theory
\citep{husler1989,engelke2015}. To prove the finite-level inequality,
we also need comparisons that determine how the rejection probability
changes near configurations attaining the limiting equality.

\subsection{Results for fixed and increasing numbers of tests}

The bivariate two-sided problem has classical positive results; see
\citet{samuelcahn1996} and the MTP$_2$ argument in
Section~\ref{main:tri-prior}. For three standard Gaussian statistics
with any correlation matrix and $\alpha\le1/2$,
\citet[Theorem~2]{ristl2016} show that the probability that at least
two of them exceed $\Phi^{-1}(1-\alpha)$ is at most $\alpha$. Their
proof studies how this probability changes with the correlations. We
also differentiate with respect to the correlations, but the Simes event
in~\eqref{eq:target} uses absolute values and compares each of the
three ordered p-values with its own threshold. It rejects if at least
one ordered p-value is at or below its corresponding threshold.
\citet[Theorems~1--3]{finnerroters2024} give covariance matrices for
which the absolute values are positively associated but not MTP$_2$;
their Section~6 discusses consequences such as lower bounds for joint
probabilities by products of marginal probabilities. These results do
not establish $F^{(m)}_\Sigma(\alpha)\le\alpha$ for unrestricted
covariance. \citet[pp.~340--341]{gou2018} restrict attention to
one-sided hypotheses, citing earlier two-sided results that require
additional dependence conditions. Their own procedure and mixture model
concern one-sided p-values and do not address this inequality.

Simes validity is also known for some models that fail PRDS.
\citet[Lemma~3.4 and Corollary~3.5]{yekutieli2008} considers pairwise
differences of independent Gaussian group means with equal variances,
each divided by the same standard error. The paper proves the Simes
inequality under explicit conditions on the testing thresholds.
The proof uses the fact that the standardized difference between
two groups is the sum of the differences through an intermediate group.
Example~2.6 of that paper gives three such standardized pairwise
differences whose absolute values fail PRDS. This earlier result
therefore establishes exact validity beyond PRDS, including singular
covariance models. Our results cover the additional covariance classes
stated above and quantify their deficits below the nominal level.

For increasing dimension, there is extensive theory on Gaussian
extremes, rare-event approximation, and multiple testing under weak
dependence \citep{agg1989,proschan2011,datta2026}.
For example, \citet[Section~2 and Theorem~3.2]{ClarkeHall2009}
derive joint rare-event limits for counts at several testing thresholds,
and \citet[Theorem~1 and Corollary~2]{farcomeni2007} establishes
asymptotic FDR control under weak dependence for BH and for plug-in
procedures, which use an estimate of the proportion of false null
hypotheses. \citet{finner2007} study FDR asymptotics for exchangeable
test statistics, such as equicorrelated normal statistics. These
limiting statements alone do not determine the
sign of the error at a finite dimension.
Theorem~\ref{main:graph_theorem} gives a more detailed expansion for the
covariance graphs studied here: the rejection probability equals
$\alpha$ minus a positive correction, with an error smaller than
that correction, even when the departure from independence is
arbitrarily small.

Independent blocks have also motivated procedures designed for block
dependence. \citet[Result~1]{guosarkar2016} show that a two-stage BH
procedure, which first selects blocks and then tests hypotheses within
the selected blocks, controls the FDR in finite samples under arbitrary
dependence within independent blocks. Our mixed-block theorem instead applies the
ordinary Simes test to all p-values together. With an empty graph
component, it covers centered Gaussian statistics in independent
blocks whose individual sizes have a fixed upper bound. One threshold
for the total number of tests works for every allowed within-block
covariance and every level below a fixed ceiling less than one.
The proof bounds each block's contribution to the full rejection
probability. Validity of separate Simes tests within the blocks would
not by itself bound this combined test.

The growing-block theorem makes the allowed block sizes explicit.
The case of bounded independent blocks in Theorem~\ref{main:mixed_theorem}
already implies, by a diagonal argument, validity for some maximum
block size tending to infinity sufficiently slowly. Theorem~\ref{main:growing_blocks} supplies an explicit growth
condition by tracking how the estimates depend on block size. It applies to pure
independent blocks; the mixed theorem still requires a fixed size
bound for its additional blocks.

Dependency graphs have likewise led to modified testing procedures.
For a nonparametric graph model,
\citet[Theorem~1]{nguyenfithian2025} prove FDR control for a class of
graph-adapted procedures. Their Proposition~11 also bounds the FDR of
ordinary BH, but at an inflated level. Our fixed-degree theorem gives
the exact nominal inequality for ordinary Simes under the centered
Gaussian global null, once the dimension exceeds a common threshold.
That threshold is uniform over all covariances in the fixed-degree
class and all levels below a fixed ceiling. Our graph records nonzero marginal covariances.

The comparison proofs use established Gaussian tools. Price's
identity and its general forms differentiate Gaussian expectations
\citep[Theorem~1]{Voigtlaender2021Price}. Interpolation that preserves
the distributions within blocks also appears in the proof of
\citet[Theorem~1]{royen2014}. Gaussian comparison inequalities bound
how probabilities of rectangles and of order-statistic events change
with the covariance \citep{LiShao2002,DebickiHashorvaJiLing2017}.
Here the estimates must cover every Simes threshold and retain the
dependence of the rejection loss on the original cluster's departure
from exact copies. Comparing the error with that loss on the same
scale keeps the error smaller, even near singular and exact-copy
limits. The final block contributions can then be added under one
joint distribution.

\subsection{Relation to FDR control}\label{main:related_fdr}

Under the global null, every BH rejection is a false
discovery. The procedure's FDR therefore equals the probability of at least one
rejection, namely~\eqref{eq:target}. This identity no longer holds
when false null hypotheses are present.
\citet{dobriban2026} constructs two-sided Gaussian
configurations containing both true and false nulls whose FDR exceeds
the nominal level. \citet[Theorem~2.1 and Propositions~2.2--2.3]{lei2026}
gives stronger lower bounds over dimension, covariance, and mean
vectors. Those constructions use nonzero means, while our problem
concerns the centered global null. Conditional calibration offers
another approach by modifying rejection thresholds under dependence
\citep{fithianlei2020}.

\citet[Theorem~3.1]{ghosh2026} give covariance-specific lower and upper
bounds for the FDR of ordinary BH; the same paper also studies modified
BH procedures that use the covariance. Their upper bound for
two-sided Gaussian tests with positive definite covariance can be
examined using
$G_\tau(u)=\overline\Psi_1(\tau\overline\Psi_1^{-1}(u))$,
where $\overline\Psi_1$ is the chi-squared survival function with one
degree of freedom and $0<\tau\le1$ is a conditional variance.
Concavity gives
$G_\tau(G_\tau^{-1}(\alpha)/m)\ge\alpha/m$.
Under the global null, summing these bounds over all true nulls gives
an FDR upper bound of at least $\alpha$. This bound therefore does not establish
the exact nominal inequality studied here.

\subsection{Covariance models beyond PRDS}
\label{main:prds}

In its pointwise form, PRDS
requires conditional probabilities of increasing events to be
nondecreasing in the conditioning p-value. An event is increasing if raising any p-value
cannot take the vector out of the event. Each main validity theorem
includes models that fail PRDS. To verify this directly, we check
conditional monotonicity; failure of a sufficient MTP$_2$ condition
would not alone prove failure of PRDS. Consider
\begin{equation}\label{main:prds-matrix}
 \Sigma_0=\begin{pmatrix}
 1&-1/2&-1/2\\-1/2&1&-1/2\\-1/2&-1/2&1
 \end{pmatrix},\qquad
 \Sigma_\varepsilon=(1-\varepsilon)\Sigma_0+\varepsilon I_3.
\end{equation}
The first matrix is obtained from the example of pairwise differences in
\citet[Example~2.6]{yekutieli2008} by coordinate sign changes.
For $Z\sim N_3(0,\Sigma_0)$, conditional on $Z_1=x$, write
\[
 Z_2=-x/2+(\sqrt3/2)E,\qquad
 Z_3=-x/2-(\sqrt3/2)E,\qquad E\sim N(0,1).
\]
Fix $a>0$ and put $u_a=2\barPhi(a)$. The event
$D_a=\{P_2\ge u_a\text{ or }P_3\ge u_a\}$ is increasing.
For $0<x<2a$, the same calculation applies conditional on $|Z_1|=x$:
\begin{equation}\label{main:prds-witness}
 \Prob\{D_a\mid P_1=2\barPhi(x)\}
 =2\Phi\left(\frac{a+x/2}{\sqrt3/2}\right)-1.
\end{equation}
The right side increases strictly with $x$, while
$P_1=2\barPhi(x)$ decreases. Thus the conditional probability of
this increasing event decreases as the conditioning p-value increases,
which violates PRDS.

The preceding calculation conditions on an exact p-value. A weaker
condition, which we call \emph{lower-tail PRDS}, conditions on the
p-value being at most a cutoff. It requires
$u\mapsto\Prob(D\mid P_1\le u)$ to be nondecreasing for every
increasing event $D$; \citet[Section~3]{wangramdas2022} use this weaker
condition as their definition of PRDS.
The same covariance also fails this condition.
Take $a=1/4$, let $q(x)$ be the right side
of~\eqref{main:prds-witness}, and put
$H(b)=\Prob(D_a\mid |Z_1|\ge b)$. With $\phi$ the standard normal density,
inclusion--exclusion and the bivariate density bound $1/(\pi\sqrt3)$ give
\[
 H(0)-q(0)\ge
 \phi(0)\left(e^{-1/32}-\frac1{\sqrt3}\right)
       -\frac1{4\pi\sqrt3}>\frac7{96}>0.
\]
The last inequality uses $\phi(0)>1/3$, $e^{-1/32}>31/32$,
$1/\sqrt3<3/5$, and $\pi\sqrt3>5$.
For the density $f_X$ and survival function $\overline F_X$ of $X=|Z_1|$,
$H'(b)=f_X(b)\{H(b)-q(b)\}/\overline F_X(b)>0$ when $b>0$ is sufficiently small.
Since $u=2\barPhi(b)$ decreases with $b$,
$\Prob(D_a\mid P_1\le u)$ decreases strictly as the cutoff $u$
increases. This contradicts the monotonicity required by lower-tail PRDS. These calculations establish
failure somewhere in the full p-value range; they do not identify a
violation below every testing level. Proposition~\ref{main:coordinate_share}
gives a separate construction whose excess rejection share forces
failure of both conventions within $(0,\alpha]$.

The pointwise violation persists for positive definite covariances.
Take $a=1$ and compare $x=1/2$ with $x=3/2$
in~\eqref{main:prds-witness}. The strict gap persists for sufficiently
small $\varepsilon>0$ by continuity of the two conditional Gaussian
probabilities. Their event boundaries have probability zero even at
the singular limit. The eigenvalues of $\Sigma_\varepsilon$ are
$\varepsilon,(3-\varepsilon)/2,(3-\varepsilon)/2$, so the perturbed
matrices are positive definite. Continuity then gives an open set of
positive definite matrices failing pointwise PRDS. For conditioning on
$P_1\le u$, choose two interior cutoffs with a strict violation and apply
continuity to the corresponding joint probabilities. Thus the failure
of lower-tail PRDS also persists under perturbation.

These examples belong to the covariance classes of our main theorems.
For every fixed $m\ge3$, Theorem~\ref{main:uniform} covers
$\Sigma_\varepsilon\oplus I_{m-3}$ at its common small levels.
Corollary~\ref{main:tri-interval} covers $\Sigma_\varepsilon$ on its
displayed interval around $1/20$.
Theorem~\ref{main:mixed_theorem} allows arbitrarily many independent
copies of these trivariate blocks when the individual block-size bound
is at least three, with a common dimension threshold for
$0<\alpha\le a_1<1$. The covariance graphs of these models have
maximum degree two, so Theorem~\ref{main:fixed_degree} applies with
degree bound $D\ge2$. They also satisfy Theorem~\ref{main:growing_blocks} whenever
its block-size envelope is at least three. Finally, PRDS passes to
marginals: an increasing event in a subset of coordinates remains
increasing as an event for the full vector. Adjoining coordinates
therefore cannot remove the displayed failure of PRDS.

We can also obtain connected non-PRDS examples of arbitrarily large
dimension. Fix sufficiently small $0<\varepsilon<1/4$ and
$0<\gamma<\varepsilon$. Start with $n$ independent blocks, each having
covariance $\Sigma_\varepsilon$. Add covariance $\gamma$ between
coordinate three of block $j$ and coordinate one of block $j+1$.
These new edges form a matching, whose symmetric adjacency matrix
has operator norm $\gamma$. The resulting covariance has smallest
eigenvalue at least $\varepsilon-\gamma>0$, a connected graph of
degree at most three, and all off-diagonal magnitudes at most $1/2$.
Each trivariate marginal is unchanged and still fails PRDS. With
unit multiplicities, this family satisfies
Theorem~\ref{main:graph_theorem} at $m=3n$, with correlation cap
$1/2$ and maximum degree at most three.

The fixed-degree theorem also covers the connected common-neighbor
families of Proposition~\ref{main:fixed_degree_example}, including
families that fail both PRDS conventions. These examples illustrate the covariance coverage of our
theorems. Validity beyond PRDS was already known from Yekutieli's
result; our conclusions add the stated uniform inequalities and
comparisons that determine how dependence changes the rejection
probability.

\section{Concluding remarks}
\label{main:discussion}

The general Simes inequality for two-sided Gaussian tests remains open.
A practical next step is to obtain numerical bounds for the level
cutoffs and dimension thresholds, so that our general results can be
assessed at conventional testing levels. Further work could also
determine how quickly block sizes and numbers of correlated neighbors
may grow while retaining validity when correlation magnitudes are
unrestricted.

\providecommand{\tightlist}{\setlength{\itemsep}{0pt}\setlength{\parskip}{0pt}}

\paragraph{Use of AI.}
The results in this paper were developed through many rounds of
discussion and iteration with AI systems, including GPT-6 Astra,
Claude Fable 5.1, and Claude Opus 5.5. The author guided the AI systems
to partition the set of covariance matrices and to pursue an induction
argument; the AI systems filled in the key technical details. These
discussions also helped explore possible strengthenings and alternative
proofs, refine arguments,
examine assumptions and boundary cases, and revise the exposition.
Gemini 3.1 Pro also assisted with wording revisions. The author
is responsible for the final statements,
proofs, and presentation.

\section*{How to read the proofs}
\phantomsection\label{comparison:guide}
\addcontentsline{toc}{section}{How to read the proofs}

The appendices contain the analytic proofs, together with the assumptions and
reported verification checks for the trivariate computation. The arguments share
some probability identities but answer different questions.
The table identifies the proof of each main result and its
prerequisites; each appendix begins with its proof strategy.
Appendices~\ref{sl:section} and~\ref{tri:proof} can be read on their
own; Appendix~\ref{mx:section} builds on Appendix~\ref{sec:lm:graph},
and Appendix~\ref{gb:section} uses both.

\begin{center}
\begin{tabular}{@{}>{\raggedright\arraybackslash}p{.27\textwidth}>{\raggedright\arraybackslash}p{.63\textwidth}@{}}
\toprule
Question & Proof and prerequisites\\
\midrule
Small testing levels (Theorem~\ref{main:uniform}) & Appendix~\ref{sl:section}: compare correlated pairs and nearly identical groups to prove validity and classify equality together, including singular laws.\\[.4em]
Three tests near level $1/20$ (Theorem~\ref{main:tri-theorem}, Corollary~\ref{main:tri-interval}) & Appendix~\ref{tri:proof}: a separate analytic and finite-certificate proof, followed by stability in the level; neither requires the small-level theorem.\\[.4em]
Sparse covariance graphs (Theorem~\ref{main:graph_theorem}) & Appendix~\ref{sec:lm:graph}: pair deficits and a relative error bound. The normal tail estimates are in Lemma~\ref{sl:tails}.\\[.4em]
Independent blocks (Theorems~\ref{main:mixed_theorem} and~\ref{main:growing_blocks}, Proposition~\ref{main:coordinate_share}) & Appendix~\ref{mx:section}: the graph estimates and near-copy comparisons give bounds for whole original blocks; Proposition~\ref{main:coordinate_share} is proved in Section~\ref{mx:coordinate}. Appendix~\ref{gb:section} retains the size dependence needed for growing blocks.\\[.4em]
Arbitrary correlations at fixed degree (Theorem~\ref{main:fixed_degree}, Proposition~\ref{main:fixed_degree_example}) & Appendix~\ref{fd:section}: group nearly identical tests, then make groups independent or replace them by copies. Uses the block bounds in Appendix~\ref{mx:section}. Proposition~\ref{main:fixed_degree_example} is proved in Section~\ref{fd:scope}.\\
\bottomrule
\end{tabular}
\end{center}

All Gaussian test coordinates have mean zero and variance one, and all
null hypotheses are true. Singular covariance matrices are allowed
unless a statement requires positive definiteness. The number $m$
always counts the original tests, so the rejection thresholds remain
$k\alpha/m$. Signed Gaussian copies have identical two-sided p-values
and can be recorded by one representative p-value with an integer
multiplicity.

The two-sided Gaussian threshold is $c_u=\Phi^{-1}(1-u/2)$.
For increasing dimension, we write $t=\alpha/m$ and $L=\log(1/t)$.
A covariance graph joins two coordinates when their covariance is
nonzero. Gaussian subvectors with zero cross-covariance are independent,
even if the covariance within either subvector is singular.

\paragraph{Shares and deficits.}
When there are rejections, the share of a set of coordinates is the
fraction of rejected coordinates that belong to that set. Its share
is zero when there is no rejection. For a partition of the tests,
the expected shares add to the full rejection probability. A uniform
representative with $w$ copies, independent of all remaining p-values,
has expected share $w\alpha/m$.
A deficit is the difference between a stated benchmark and a rejection
probability or expected share. The graph expansion compares the full
rejection probability with $\alpha$. A block of $s$ original tests is
instead assigned the nominal expected share $s\alpha/m$. A pair's
effect on the full probability and its effect on one block's share
are therefore different deficits. The table at the start of
Appendix~\ref{mx:section} lists the main deficits used in
Appendices~\ref{sec:lm:graph}, \ref{mx:section} and~\ref{fd:section},
each with the distribution under which it is computed.

\paragraph{What a comparison changes.}
A reference distribution is a specified comparison model whose
rejection probability or expected share we can bound. For a graph
edge, the reference keeps that pair's correlation and every copy
count; all other representatives are independent of the pair and of
each other. In a block comparison, the selected block may change, but
the joint distribution of the tests outside it is kept fixed, and these
tests remain independent of the block. Different blocks can use different
comparison models. After their shares have been bounded, however,
the shares being summed are those under the original joint
distribution.
The fixed-degree proof follows a different sequence of steps. It first
compares full rejection probabilities through successive isolation and contraction
steps, both defined in the next paragraph. The final distribution is called the terminal law, and all
subsequent share bounds use that same distribution. In each argument,
the deficit must be defined under the distribution used for that
comparison.

Isolation makes a cluster independent of its complement while keeping
the joint distribution within each part unchanged. Contraction
replaces every test in a cluster by a copy of one representative, so
the representative's multiplicity is the cluster's number of tests and
$m$ is unchanged. The representative may be a cluster mean normalized
to variance one or an actual coordinate in the cluster. The proof
must specify which choice it uses, because the choices can have
different covariances with the outside coordinates.
Appendix~\ref{sl:section} uses the normalized cluster mean, and
Appendices~\ref{mx:section}, \ref{gb:section} and~\ref{fd:section} use an
actual coordinate, which keeps its correlations with the
other tests. A rejection is
gained if it occurs only for the original vector and lost if it occurs
only for the reference vector. We bound these two probabilities
separately to determine the direction of the comparison.

\paragraph{Uniformity.}
A uniform estimate uses the same bound or threshold for every
covariance matrix and level allowed by its hypotheses. Each proof
specifies which parameters are fixed and how its constants may depend
on them. If block size grows, that dependence must also be bounded.
Comparing with a reference probability proves validity only after the
reference probability has itself been bounded. These requirements
continue to apply when parameters approach singular or exact-copy
limits.

\paragraph{How the arguments fit together.}
Near covariance matrices at which equality holds, the rejection probability
can differ from $\alpha$ by an arbitrarily small amount. An approximation
$F_\Sigma(\alpha)=\alpha+o(\alpha)$ does not determine the sign of
that difference when the covariance varies with the level. The proofs
therefore retain a positive deficit and show that each remaining
error is smaller than this deficit.

For small levels, the pairwise scale is
$c_\alpha^2(1-|\Gamma_{ij}|)$, where $\Gamma$ is the correlation matrix of the representatives. Representatives separated on this scale
admit a pair expansion with a remainder smaller than the pair
contribution. A positive finite limiting distance inside a cluster
gives a contribution strictly below the limiting equality value.
When the within-cluster distances tend to zero on this scale, we
contract nearly identical coordinates to exact copies and show that the probability of losing a rejection
dominates the probability of gaining one. Induction retains all
multiplicities and proves validity and the equality classification
together. Appendix~\ref{sl:section} gives the argument.

For increasing dimension, the graph proof (Appendix~\ref{sec:lm:graph})
compares interactions among
larger connected sets with the one-edge deficits in the full rejection
probability. Counting those interactions explicitly gives the stated
degree range. The block proof (Appendix~\ref{mx:section}) follows each
original block's expected
share while keeping fixed the joint distribution of all tests outside
that block, which are independent of it.
The resulting error is smaller than the retained loss from departures from exact
copies. Bounds on how the constants depend on block size give the
growing-block result (Appendix~\ref{gb:section}). For unrestricted
correlations at fixed degree (Appendix~\ref{fd:section}), the comparisons
leading to one terminal distribution are organized by clusters of
coordinates, found from distances between unit vectors whose inner
products are the correlations; each cluster is much tighter internally
than its separation from the other tests. The distinction among these
distributions is why the corresponding deficits must be kept separate.

The trivariate proof is separate. Analytic estimates cover
neighborhoods of equality, and a finite certificate gives a strict
margin on the remaining covariance region. Near equality, the lower bound for the deficit and the bound on how the deficit changes with the level
must retain the same factor that tends to zero. Away from equality, the certified margin
allows an absolute bound on that change. Together these bounds give
the level interval using the existing certificate at $1/20$; no new
covariance cover is needed.

\appendix

\providecommand{\slP}{\mathbb P}
\providecommand{\slE}{\mathbb E}
\providecommand{\slC}{\mathcal C}
\providecommand{\slPhi}{\overline\Phi}
\providecommand{\slphi}{\phi}
\providecommand{\slone}{\mathbf 1}
\theoremstyle{remark}
\theoremstyle{plain}

\section{Exact covariance-uniform validity at small levels}
\label{sl:section}

We prove that, in each fixed dimension, the Simes inequality holds at all
sufficiently small levels, with one cutoff for every covariance matrix.
Throughout this section, all Gaussian test statistics
are centered with unit variance, and every null hypothesis
is true.  Let \(\Phi\) and
\(\slphi\) denote the standard normal distribution function and density,
let \(\slPhi=1-\Phi\), and put
\[
 c_u=\Phi^{-1}(1-u/2),\qquad 0<u<1.
\]
For a positive semidefinite correlation matrix
\(\Sigma\in\slC_m\), define \(Z\sim N_m(0,\Sigma)\),
\(P_i=2\slPhi(|Z_i|)\), and
\begin{equation}
 F^{(m)}_\Sigma(\alpha)
 =\slP\left\{P_{(k)}\le\frac{k\alpha}{m}
       \text{ for some }1\le k\le m\right\}.
 \label{sl:event}
\end{equation}
We allow singular Gaussian laws.  Every marginal \(P_i\) remains uniform
in this case.
For each positive integer \(d\), let \(\mathcal E_d\) be the set of
coordinate-sign changes and permutations of
\(\operatorname{diag}(J_{k_1},\ldots,J_{k_r})\), where \(J_k\) is the
\(k\times k\) all-ones matrix and the positive integers
\(k_1,\ldots,k_r\) sum to \(d\).
These are the covariance matrices of independent groups of signed copies.
Theorem~\ref{sl:uniform} restates Theorem~\ref{main:uniform} and adds the
form with repeated p-values used in the induction.  In this form, a
multiplicity vector \(w=(w_1,\ldots,w_q)\) repeats the \(g\)th of \(q\)
Gaussian representatives \(w_g\) times, the representative covariance
is their \(q\times q\) correlation matrix, and the rejection probability
\(F_w\) is defined in Section~\ref{sl:preliminaries}.

\begin{sltheorem}[Exact covariance-uniform small-level validity]
\label{sl:uniform}
For every fixed integer \(m\ge2\), there exists \(a_m\in(0,1)\) such that
\begin{equation}
 \sup_{\Sigma\in\slC_m}F^{(m)}_\Sigma(\alpha)=\alpha,
 \qquad 0<\alpha\le a_m.
 \label{sl:main}
\end{equation}
The cutoff is common to all correlation signs and all singular or
nonsingular covariance matrices.  It is also common to all vectors of
positive integer duplication multiplicities whose total is \(m\).
At each level in this interval, equality holds if and only if
\(\Sigma\in\mathcal E_m\). More generally, for any such multiplicity
vector \(w\) of length \(q\), equality holds if and only if its
representative covariance belongs to \(\mathcal E_q\).
\end{sltheorem}

For \(m=2\), Lemma~\ref{sl:pairidentity} shows that any \(a_2\in(0,1)\)
can be used. For \(m\ge3\), the proof does not give a numerical value
for the cutoff in Theorem~\ref{sl:uniform}. Whether it reaches
\(0.05\) or \(0.1\), or can be chosen independently of \(m\), remains
unresolved.  The conclusion concerns the global null and does not imply
FDR control when some null hypotheses are false.

\paragraph{Proof overview.}
We fix \(m\) and induct on the number \(q\) of Gaussian representatives,
retaining their multiplicities; the induction step argues by
contradiction.  Along a counterexample sequence \(\alpha_n\downarrow0\),
we classify pairs according to the limits of their rescaled distances
\(c_{\alpha_n}^2(1-|\Gamma_{ij,n}|)\), where \(\Gamma_n\) is the
correlation matrix of the representatives.
This scale, which we call the critical scale, compares the differences
between nearly identical Gaussian coordinates with the gaps between
rejection thresholds.  A positive
finite limit means that the two are of comparable size.  An infinite
limit separates the pair at this scale, while a zero limit makes the
pair behave almost as exact copies.
Representatives whose rescaled distances have finite limits form
classes, called clusters or blocks in
Lemmas~\ref{sl:contraction}--\ref{sl:multicritical}.
When all pairs are separated at this scale and the representatives are
not independent, the combined decrease in rejection probability from
pairs exceeds the remainder from interactions involving three or more
coordinates (Lemma~\ref{sl:separated}, which uses Lemma~\ref{sl:weakcube}
near independence). If a class has a
positive finite limiting distance between two of its members, the total
rejection probability divided by the nominal level has a limit strictly
below one (Lemma~\ref{sl:multicritical}).
In the remaining case, every within-class limiting distance is zero
and some class has at least two members, so the members of each class
approach exact duplication faster than this scale.  We treat these classes using Lemma~\ref{sl:contraction}, which
compares the original vector with the vector formed from copies of the
standardized cluster means.  We call a rejection gained if it occurs
only for the original vector, and lost if it occurs only for the
contracted vector.  The probability of losing a rejection is
bounded below by a positive constant times
\(\alpha_n\sqrt{\lambda_n}\), while the
probability of gaining one is
\(O(\alpha_n e^{-c/\lambda_n})\) for a constant \(c>0\); here
\(\lambda_n\to0\) is the rescaled within-cluster spread, defined
precisely in Lemma~\ref{sl:contraction}.  Comparing gains with losses determines the
direction of the inequality.  Contraction lowers the number of
representatives but keeps the denominator \(m\), so the induction
hypothesis, including its equality classification, bounds the
contracted probability.  Exact contraction leaves the rejection
probability unchanged, and the contracted covariance is not an equality
case, since otherwise expanding it would place the original covariance in
the class of independent signed-copy groups; every other branch has a
strict deficit.

\subsection{Normal tails and duplication-preserving identities}
\label{sl:preliminaries}

The proof needs two kinds of preliminary facts. Normal tail estimates
translate the gaps between p-value thresholds into gaps between Gaussian
thresholds. Exact probability identities separate the contribution of each
pair from the remaining dependence. These identities count every copy of a
duplicated p-value, so they still apply after a cluster is contracted.

\begin{sllemma}[Tail and quantile estimates]
\label{sl:tails}
For \(x>0\),
\begin{equation}
 \frac{x}{1+x^2}\slphi(x)\le\slPhi(x)\le\frac{\slphi(x)}x,
 \qquad
 x\le h(x):=\frac{\slphi(x)}{\slPhi(x)}\le x+x^{-1}.
 \label{sl:mills}
\end{equation}
The hazard \(h\) is increasing on the real line.  If \(b=c_\alpha\)
and \(A>0\) is fixed, then, as \(\alpha\downarrow0\),
\begin{align}
 c_{A\alpha}
   &=b-\frac{\log A}{b}+O_A(b^{-3}), \label{sl:quantile}\\
 \frac{\slphi(c_{A\alpha})}{\slphi(b)}
   &\longrightarrow A. \label{sl:densityratio}
\end{align}
The estimates are uniform when \(A\) ranges over any fixed compact subset of \((0,\infty)\).
If \(x_b\to x\in\mathbb R\), then
\begin{equation}
 \frac{\slPhi(b+x_b/b)}{\slPhi(b)}\longrightarrow e^{-x}.
 \label{sl:tailshift}
\end{equation}
\end{sllemma}

\begin{proof}
We obtain the tail bounds by integration by parts.  This gives the
upper bound in~\eqref{sl:mills} and
\[
 \slPhi(x)
 =\frac{\slphi(x)}x-\int_x^\infty\frac{\slphi(t)}{t^2}\,dt
 \ge\frac{\slphi(x)}x-\frac{\slPhi(x)}{x^2},
 \qquad x>0,
 \]
which gives the lower bound.  Also
\(h'(x)=h(x)\{h(x)-x\}>0\): for \(x>0\) this follows from the strict
upper Mills bound, and for \(x\le0\) it follows from \(h(x)>0\).

To compare nearby rejection thresholds, we next determine how a small
shift in a Gaussian threshold changes its tail probability. For bounded
real \(v\), integration of \(h(t)=t+O(t^{-1})\) on the interval with
endpoints \(b\) and \(b+v/b\) gives
\[
 \log\frac{\slPhi(b+v/b)}{\slPhi(b)}
 =-v+O_v(b^{-2}),
 \]
uniformly on bounded \(v\)-intervals.  This proves~\eqref{sl:tailshift}.
Comparing the quantile with the endpoints of fixed intervals of this
form gives
\(c_{A\alpha}-b=O_A(b^{-1})\).  Inserting this bound into
\[
 \log A=-\int_b^{c_{A\alpha}}h(t)\,dt
 \]
gives~\eqref{sl:quantile}.  Squaring the expansion yields
\eqref{sl:densityratio}. On a fixed compact subset of positive $A$,
$\log A$ is uniformly bounded. We can therefore use the same intervals
to locate the quantiles and the same remainder bounds throughout that subset.
\end{proof}

Fix positive integers \(w_1,\ldots,w_q\) with
\(\sum_{g=1}^q w_g=m\).  For \(\Gamma\in\slC_q\), write
\(F_w(\Gamma,\alpha)\) for the rejection probability obtained by
repeating the \(g\)th Gaussian representative p-value \(w_g\) times
and using the original thresholds \(k\alpha/m\).  We also call \(w_g\)
the weight, or copy count, of representative \(g\).  For example, a
representative with weight two supplies
two identical entries to the ordered list of p-values.  We use
``physical coordinates'' to refer to all entries in that original list,
including copies.  Put \(t=\alpha/m\).

For a bivariate standard Gaussian vector with correlation \(s\in[0,1]\),
let
\[
 J_s(x,y)=
 \slP\{2\slPhi(|Z_1|)\le x,\ 2\slPhi(|Z_2|)\le y\}.
 \]
The joint law of the absolute values depends only on the absolute
correlation: changing the sign of either coordinate leaves its p-value
unchanged.  For positive integer
weights \(u,v\) with \(u+v\le m\), define
\begin{equation}
 D_s^{u,v}(t)
 =J_s(ut,(u+v)t)+J_s(vt,(u+v)t)
       -J_s((u+v)t,(u+v)t).
 \label{sl:deficit}
\end{equation}

This pair deficit measures the reduction in rejection probability for
the two active groups from their combined nominal contribution.  The
next lemma separates these pair effects from the remaining dependence.

\begin{sllemma}[Two groups and the exact residual]
\label{sl:pairidentity}
For \(0<\alpha<1\), \(D_s^{u,v}(\alpha/m)\ge0\) for all \(s\in[0,1]\).
It vanishes exactly at \(s=0\) and \(s=1\).
Independent representatives satisfy
\begin{equation}
 F_w(I_q,\alpha)=\alpha.
 \label{sl:independent}
\end{equation}
The same equality holds for every \(\Gamma\in\mathcal E_q\).
For every \(q\) and \(\Gamma\in\slC_q\),
\begin{equation}
 F_w(\Gamma,\alpha)
 =\alpha-\sum_{g<h}D_{|\Gamma_{gh}|}^{w_g,w_h}(t)
       +\slE_\Gamma R_w.
 \label{sl:decomposition}
\end{equation}
Here \(R_w\) is the rejection indicator minus its approximation by one-
and two-group indicators, defined in~\eqref{sl:residualdefinition}; it is
unchanged by changing the sign of any Gaussian coordinate
(that is, it is coordinatewise even) and satisfies
\[
 |R_w|\le C_q^R:=\begin{cases}\binom q2-1,&q\ge3,\\0,&q\le2.\end{cases}
\]
It vanishes unless at least three representative p-values are at
most \(\alpha\).  More precisely, writing
\(N_\alpha=\sum_{g=1}^q\slone\{P_g\le\alpha\}\), we have
\[
 1-\binom{N_\alpha}{2}\le R_w\le\binom{N_\alpha-1}{2}
 \quad(N_\alpha\ge3),\qquad
 |R_w|\le2\binom{N_\alpha}{3}.
\]
Moreover, \(\slE_{I_q}R_w=0\).
The Fr\'echet derivative of its expectation at \(I_q\) vanishes in every
zero-diagonal symmetric direction, that is, along correlation matrices;
consequently the first derivative along every path in \(\slC_q\)
differentiable at its starting point \(I_q\) vanishes.
In particular, \(F_w(\Gamma,\alpha)\le\alpha\) for \(q\le2\).
\end{sllemma}

\begin{proof}
We prove four facts in turn: pair deficits are nonnegative, independent
representatives give equality in the Simes inequality, the remaining terms involve at least
three small p-values and obey the stated bounds, and their expectation and its first
derivative vanish at independence.
Let \(W=(u+v)t<1\), and define
\(f(z)=\slP(P_2\le W\mid P_1=z)\).
Since \(P_1\) is uniform on \([0,1]\), \(J_s(x,W)=\int_0^xf(z)\,dz\) for
\(0\le x\le1\). The three terms of~\eqref{sl:deficit} therefore combine as
\begin{equation}
 D_s^{u,v}(t)
   =\int_0^{ut}\{f(z)-f(z+vt)\}\,dz.
 \label{sl:deficitintegral}
\end{equation}
For \(0<s<1\), put \(B=c_W\), \(\sigma^2=1-s^2\), and
\[
 Q_s(y)=\slP\{|sy+\sigma N|\ge B\},\qquad N\sim N(0,1).
 \]
Then \(f(z)=Q_s(c_z)\), and differentiation gives, for \(y>0\),
\begin{align}
 Q_s'(y)
 &=\frac{s}{\sigma}\slphi\left(\frac{B-sy}{\sigma}\right)
       \left(1-e^{-2sBy/\sigma^2}\right) \nonumber\\
 &=\frac{2s}{\sqrt{2\pi}\sigma}
       e^{-(B^2+s^2y^2)/(2\sigma^2)}
       \sinh\left(\frac{sBy}{\sigma^2}\right)>0.
 \label{sl:foldedderivative}
\end{align}
Thus \(f\) is strictly decreasing on \((0,1)\) when \(0<s<1\).
Since both weights are positive, the integrand in
\eqref{sl:deficitintegral} is strictly positive on \((0,ut)\), and so
the deficit is strictly positive.
At independence \(J_0(x,y)=xy\), and
\[
 D_0^{u,v}(t)
 =\{u(u+v)+v(u+v)-(u+v)^2\}t^2=0.
 \]
At \(s=1\), the two p-values coincide and \(J_1(x,y)=\min(x,y)\),
which again gives zero.  The endpoint cases therefore follow directly
from the corresponding distribution functions.

To prove the equality under independence in~\eqref{sl:independent}, let \(U_g\) be independent uniform
representatives, and let \(R\in\{0,\ldots,m\}\) be the largest rank
at which the step-up inequality holds.  When \(R>0\), precisely the p-values
at most \(Rt\) are rejected: if their number were larger than
\(R\), that number would give a larger crossing.  Hence
\[
 \slone\{R>0\}=\sum_g
       \frac{w_g\slone\{U_g\le Rt\}}{R},
 \]
with summands set to zero if \(R=0\).
Let \(R^{[g]}\) be the step-up count after all copies of group \(g\)
are replaced by zero.  For \(r\ge w_g\),
\[
 \{U_g\le rt,\ R=r\}=\{U_g\le rt,\ R^{[g]}=r\}.
 \]
Indeed, reducing a group already among the first \(r\) positions
does not alter any larger order statistic; conversely, restoring
that group to a value at most \(rt\) preserves the crossing at
\(r\) and cannot create a larger one.  Since \(R^{[g]}\) is independent
of \(U_g\) and \(R^{[g]}\ge w_g\),
\[
 \slP(R>0)
 =\sum_g\sum_{r=w_g}^m\frac{w_g}{r}
       \slP(U_g\le rt,R^{[g]}=r)
 =\sum_g w_gt=\alpha.
 \]
If \(\Gamma\in\mathcal E_q\), align the signs and merge representatives
that coincide almost surely, adding their multiplicities. The resulting
representatives are independent, and the physical p-value list is
unchanged. The same identity therefore gives equality on all of
\(\mathcal E_q\), at every level in \((0,1)\).

We now subtract the contributions that can already be explained by one
or two groups. To isolate a selected set of groups, we set every other
p-value to one while keeping its place in the list. Let \(H\) be the
rejection indicator for the list including all copies, and let \(H_S\)
be the same indicator after every representative outside \(S\)
is replaced by one; its denominator remains \(m\).
Set
\begin{equation}
 R_w=H-\sum_g H_{\{g\}}
           -\sum_{g<h}\bigl(H_{\{g,h\}}-H_{\{g\}}-H_{\{h\}}\bigr).
 \label{sl:residualdefinition}
\end{equation}
One group has rejection probability \(w_gt\).
Two active groups reject on
\[
 \{P_g\le w_gt\}\ \cup\ \{P_h\le w_ht\}\
 \cup\ \{P_g\le(w_g+w_h)t,\ P_h\le(w_g+w_h)t\}.
 \]
The two coordinates of a bivariate standard Gaussian vector are
exchangeable, so \(J_s(x,y)=J_s(y,x)\). By inclusion--exclusion and this
symmetry, the probability of this union is
\((w_g+w_h)t-D_{|\Gamma_{gh}|}^{w_g,w_h}(t)\).
Under independence, this probability is therefore
\((w_g+w_h)t\), which equals \(\alpha\) only when \(w_g+w_h=m\).
Taking expectations in~\eqref{sl:residualdefinition}
proves~\eqref{sl:decomposition}.

The residual can be nonzero only when at least three groups have small
p-values. Indeed, when at most two representative p-values are at most
\(\alpha\), replacing the other p-values by one leaves every crossing
unchanged. Definition~\eqref{sl:residualdefinition} then gives \(R_w=0\).
Its invariance under coordinate sign changes follows directly from
the indicator definition. To bound its size, distinguish groups that
reject on their own from pairs that reject only together. Let
$n_1=\sum_gH_{\{g\}}$ and let $K$ count pairs for which neither group
rejects alone but the two groups reject together. The pair correction is $-1$ when both singletons succeed,
$1$ for a pair counted by $K$, and zero otherwise. Therefore
\[
 R_w=H-n_1+\binom{n_1}2-K,\qquad
 0\le K\le\binom{N_\alpha-n_1}{2}.
\]
If $H=0$, then $n_1=K=0$. If $H=1$ and $N_\alpha\ge3$, substituting the
upper bound on $K$ gives
$R_w\ge1-\binom{N_\alpha}{2}+(N_\alpha-2)n_1$.
For the upper bound, maximize the convex quadratic
$1-n_1+\binom{n_1}2$ over $0\le n_1\le N_\alpha$.
This proves the asymmetric interval in the statement, and hence
\[
 |R_w|\le\mathbf1\{N_\alpha\ge3\}
                \left(\binom{N_\alpha}{2}-1\right)
          \le2\binom{N_\alpha}{3}.
\]
The one-sided bound
$R_w\le\binom{N_\alpha}{3}$ follows as well. These deterministic
identities include ties and arbitrary positive real masses, after the
same total-mass normalization.

Finally, we check the cancellation needed when all correlations are
small. Independence and the vanishing pair deficits give
\(\slE_{I_q}R_w=0\). For a zero-diagonal symmetric direction
\(A\), the density derivative at \(I_q\) is
\[
 \left(\sum_{g<h}A_{gh}z_gz_h\right)
            \prod_{\ell=1}^q\slphi(z_\ell).
 \]
Its integral against the bounded coordinatewise even function \(R_w\) is zero: each
summand changes sign when one of its two coordinates changes sign,
whereas the residual and independent density do not. Differentiation
under the integral is justified by a polynomial times a Gaussian
envelope on a small positive-definite neighborhood of \(I_q\).
\end{proof}

\subsection{Uniform validity under weak correlation}
\label{sl:weaksubsection}

We next use the pair decomposition to obtain validity in a common
neighborhood of independence.  At independence, both the expected residual and
the pair deficits vanish.  We therefore keep their dependence on the squared
correlations and show that the expected residual is smaller at low levels.

\begin{sllemma}[Weak-correlation cube]
\label{sl:weakcube}
Fix \(w\) with \(q\ge3\), and let
\[
 r_0=\frac1{10(q-1)},\qquad M=\frac{11}{10}.
 \]
There exists \(\alpha_w^{\mathrm{weak}}>0\) such that, simultaneously
for \(0<\alpha\le\alpha_w^{\mathrm{weak}}\) and
\(\max_{g<h}|\Gamma_{gh}|\le r_0\),
\[
 F_w(\Gamma,\alpha)\le\alpha.
 \]
Equality holds if and only if \(\Gamma=I_q\).
\end{sllemma}

\begin{proof}
We first bound the expected residual by differentiating the Gaussian density
along a path from independence.  We then obtain a lower bound for the
pair deficits and compare their decay rates.  Write \(A=\Gamma-I_q\),
\(S=\sum_{g<h}\Gamma_{gh}^2=\|A\|_{\mathrm F}^2/2\), and
\(\Gamma_\theta=I_q+\theta A\), \(0\le\theta\le1\).
Every eigenvalue of \(\Gamma_\theta\) lies in \([9/10,11/10]\)
by the row-sum bound \(\|A\|_{\mathrm{op}}\le(q-1)r_0=1/10\).
Let \(f_\theta\) be its Gaussian density and \(K_\theta=\Gamma_\theta^{-1}\).
Direct differentiation gives
\begin{align*}
 \partial_\theta\log f_\theta(z)
    &=-\tfrac12\operatorname{tr}(K_\theta A)
       +\tfrac12z^TK_\theta A K_\theta z,\\
 \partial_\theta^2\log f_\theta(z)
    &=\tfrac12\operatorname{tr}(K_\theta A K_\theta A)
       -z^TK_\theta A K_\theta A K_\theta z.
\end{align*}
Using \(\|K_\theta\|_{\mathrm{op}}\le10/9\),
\(\|A\|_{\mathrm{op}}\le\|A\|_{\mathrm F}=\sqrt{2S}\),
and
\[
 f_\theta(z)\le
 (2\pi)^{-q/2}(10/9)^{q/2}
          e^{-\|z\|^2/(2M)},
 \]
in
\(\partial_\theta^2 f_\theta
 =f_\theta\{(\partial_\theta\log f_\theta)^2+
                 \partial_\theta^2\log f_\theta\}\)
shows that
\begin{equation}
 |\partial_\theta^2f_\theta(z)|
 \le C_q S(1+\|z\|^4)e^{-\|z\|^2/(2M)}.
 \label{sl:densitysecond}
\end{equation}
Here \(C_q<\infty\) depends only on \(q\). We bound each trace by
the Frobenius norm and each quadratic form by the operator norm,
retaining the factor \(S\). Keeping this factor is essential: the
expected residual must vanish together with the pair deficits when correlations
approach zero, however quickly they do so.

Put \(b=c_\alpha\) and assume \(b\ge1\).
For \(j=0,2,4\), integration by parts yields
\[
 \int_{|x|\ge b}|x|^j e^{-x^2/(2M)}\,dx
 \le C_{j,M}b^{j-1}e^{-b^2/(2M)}.
 \]
For \(j=0\), use the upper Mills bound; for \(j=2,4\), use the
recursion
\[
 \int_b^\infty x^j e^{-x^2/(2M)}\,dx
 =Mb^{j-1}e^{-b^2/(2M)}
     +M(j-1)\int_b^\infty x^{j-2}e^{-x^2/(2M)}\,dx.
 \]
Since \(\|z\|^4\le q\sum_i z_i^4\), integrating the envelope
in~\eqref{sl:densitysecond} over the event that the absolute values of
three specified coordinates exceed \(b\) gives at most
a constant multiple of \(S\,b e^{-3b^2/(2M)}\), where the constant
depends only on \(q\).
The largest power is \(b^3 b^{-1}b^{-1}=b\);
if the fourth-degree factor belongs to an unrestricted coordinate,
its integral is a finite constant and the three tail powers are
\(b^{-3}\le b\). Thus the residual requires three tail exceedances,
whereas each pair deficit involves only two. Taylor's formula, the
zero value and derivative at independence, and the residual support
in Lemma~\ref{sl:pairidentity} therefore give
\begin{equation}
 |\slE_\Gamma R_w|
 \le C'_q S b e^{-3b^2/(2M)}.
 \label{sl:weakremainder}
\end{equation}

To compare this expected residual with the pair deficits, we obtain a lower
bound common to every pair of weights in \(w\).  Exchange \(u,v\) if
necessary so that \(u\le v\).  In~\eqref{sl:deficitintegral}, restrict the
integral to
\(z\in[ut/2,ut]\).  On this interval,
\[
 [c_{(u/2+v)t},c_{ut}]
       \subseteq[c_{z+vt},c_z].
 \]
By~\eqref{sl:quantile}, for all \(b\) above a threshold
depending only on \(m\), all endpoints in this display and
\(B=c_{(u+v)t}\) belong to \([b,b+K_m/b]\), where one may take
\(K_m=2\log(2m)+1\).
The displayed common interval has length at least
\(\ell/b\), where \(\ell=\tfrac12\log(3/2)>0\):
its asymptotic length times \(b\) is
\(\log((u/2+v)/u)\ge\log(3/2)\).
These bounds hold uniformly over the finite collection
\(u,v\in\mathbb N,\ u+v\le m\).

For \(0<s\le r_0\), the inequality \(\sinh x\ge x\)
in~\eqref{sl:foldedderivative} gives
\[
 Q_s'(y)\ge
 \frac{2s^2By}{\sqrt{2\pi}(1-s^2)^{3/2}}
       e^{-(B^2+s^2y^2)/(2(1-s^2))}.
 \]
Let
\(\kappa=(1+r_0^2)/[2(1-r_0^2)]\).
On the retained interval, \(By\ge b^2\) and
\[
 \frac{B^2+s^2y^2}{2(1-s^2)}
 \le\kappa(b+K_m/b)^2
 \le\kappa b^2+\kappa(2K_m+K_m^2).
 \]
Since the selected \(z\)-interval has length \(ut/2\ge\alpha/(2m)\),
\[
 D_s^{u,v}(\alpha/m)
 \ge
 \frac{\ell}{m\sqrt{2\pi}}\,
 e^{-\kappa(2K_m+K_m^2)}
       \alpha b\,s^2 e^{-\kappa b^2}.
 \]
For \(b\ge1\), the lower Mills bound implies
\(\alpha b=2b\slPhi(b)\ge\slphi(b)\).
Consequently
\begin{equation}
 D_s^{u,v}(\alpha/m)
 \ge c_{m,q}s^2 e^{-b^2/(1-r_0^2)},
 \qquad
 c_{m,q}=\frac{\ell}{2\pi m}
                   e^{-\kappa(2K_m+K_m^2)}>0.
 \label{sl:weakpair}
\end{equation}
At \(s=0\), both sides are zero, so the same bound holds.
Since
\[
 \frac{3}{2M}=\frac{15}{11}>
 \frac1{1-r_0^2},
 \]
the ratio of~\eqref{sl:weakremainder} to the sum of the pair lower
bounds converges to zero uniformly when \(S>0\). Choose the cutoff
so that this ratio is at most \(1/2\). The exact decomposition then gives
\[
 F_w(\Gamma,\alpha)
 \le\alpha-\tfrac12\sum_{g<h}D_{|\Gamma_{gh}|}^{w_g,w_h}(\alpha/m)
 <\alpha,
\]
because at least one correlation is nonzero and all are strictly below one.
When \(S=0\), equality follows from~\eqref{sl:independent}.
All lower bounds required for \(b\) depend only on \(m,q\).  We
choose their maximum and convert it to a positive common cutoff for
\(\alpha\).
\end{proof}

\subsection{Representatives separated at the critical scale}
\label{sl:separatedsubsection}

When every pair is separated at the critical scale, the pair deficits
control the residual. The preceding lemma already treats correlations
near zero. For the other cases, we will show that the deficit of a pair
with the largest absolute correlation exceeds the absolute expected residual.
We first bound simultaneous exceedances of three coordinates, since the
residual is supported on such events.

\begin{sllemma}[A trivariate orthant bound]
\label{sl:triple}
Let \(r\in[0,1)\), \(b>0\), and let \(V\) be any centered
trivariate Gaussian vector with unit marginal variances and
\(|\operatorname{Cov}(V_i,V_j)|\le r\).
For each choice of signs \(\varepsilon_i\in\{-1,1\}\),
\begin{equation}
 \slP(\varepsilon_1V_1\ge b,\varepsilon_2V_2\ge b,
                       \varepsilon_3V_3\ge b)
 \le
 \frac{(1+2r)^{5/2}}{(2\pi)^{3/2}(1-r)b^3}
           e^{-3b^2/[2(1+2r)]}.
 \label{sl:triplebound}
\end{equation}
\end{sllemma}

\begin{proof}
We compare the signed vector with an equicorrelated Gaussian vector
whose common correlation \(r\) is at least every signed correlation. We then bound the probability that all
three coordinates exceed the threshold.  First suppose the signed
covariance is positive definite.
For a Gaussian density, differentiating with respect to an
off-diagonal covariance entry \(\rho_{ij}\), while changing its
symmetric entry simultaneously, gives
\(\partial_{\rho_{ij}}f=\partial_{z_i}\partial_{z_j}f\).
This identity follows either by direct differentiation of the
density or of its Fourier transform.  Integrating the right-hand
side over the positive orthant above \(b\) gives
\[
 \frac{\partial}{\partial\rho_{ij}}
     \slP(V_\ell\ge b\ \forall\ell)
 =\int_b^\infty f(z_i=b,z_j=b,z_k)\,dz_k\ge0.
 \]
Interpolate to the equicorrelation matrix
\(E_r=(1-r)I_3+r\slone\slone^T\).
Its off-diagonal entries dominate all the signed entries.
The interpolating matrices are positive definite, so the
derivative identity proves the required orthant comparison.
For a singular initial covariance, use its mixture with
\(\tau I_3\), compare it with \(E_{(1-\tau)r}\), and let
\(\tau\downarrow0\).  Gaussian weak convergence is sufficient,
since a boundary coordinate equal to \(b\) has probability zero.

For \(E_r\), the determinant is
\((1-r)^2(1+2r)\), and
\(E_r^{-1}\slone=(1+2r)^{-1}\slone\).
At \(x=b\slone+y\), \(y_i\ge0\),
\[
 \tfrac12 x^T E_r^{-1}x
 =\frac{3b^2}{2(1+2r)}
       +\frac b{1+2r}\sum_i y_i
       +\tfrac12y^T E_r^{-1}y.
 \]
Dropping the last, nonnegative term and integrating the three linear
exponentials gives the factor
\(((1+2r)/b)^3\).  Multiplication by the determinant normalization
gives~\eqref{sl:triplebound}.
\end{proof}

\begin{sllemma}[Separated extremes]
\label{sl:separated}
Fix \(w\), and suppose that
\(\alpha_n\downarrow0\), \(\Gamma_n\in\slC_q\), \(b_n=c_{\alpha_n}\),
and
\[
 \Lambda_n=b_n^2(1-r_n)\longrightarrow\infty,
 \qquad r_n=\max_{g<h}|\Gamma_{gh,n}|.
 \]
For \(q=1\) the maximum need not be defined and equality holds.
For \(q\ge2\), the above conditions imply
\[
 F_w(\Gamma_n,\alpha_n)\le\alpha_n
 \quad\hbox{for all sufficiently large }n.
 \]
At every sufficiently large index for which \(\Gamma_n\ne I_q\),
the inequality is strict.
\end{sllemma}

\begin{proof}
For \(q\le2\), Lemma~\ref{sl:pairidentity} gives the assertion,
including strictness: the separation assumption excludes perfect
correlation at all sufficiently large indices.
For \(q\ge3\), let \(r_0=1/[10(q-1)]\).
Lemma~\ref{sl:weakcube} treats the indices with \(r_n\le r_0\),
with equality exactly when \(\Gamma_n=I_q\).
For the remaining indices, we select a pair with absolute correlation
\(r_n>r_0\), and denote its weights by \(u,v\).
Put $K_m=2\log(2m)+1$ and $\ell=\tfrac12\log(3/2)$.
The common intervals in the proof of \eqref{sl:weakpair} are available
for $b\ge B_m$, with a fixed $B_m\ge1$ depending only on $m$.
For $r\ge r_0$ and $\Lambda=(1-r)b^2\ge1$, their points $y$ and the
threshold $B=c_{(u+v)t}$ of~\eqref{sl:foldedderivative}, not to be
confused with the constant $B_m$, satisfy
\[
 |B-ry|\le(1-r)b+K_m/b,\qquad
 \frac{(B-ry)^2}{2(1-r^2)}
 \le\frac{(1-r)b^2}{2(1+r)}+K_m+K_m^2/2.
\]
Also $1-e^{-2rBy/(1-r^2)}\ge1-e^{-2r_0}$.
Use the first expression in \eqref{sl:foldedderivative}, integrate
over a $z$-interval of length at least $\alpha/(2m)$ and a $y$-interval
of length at least $\ell/b$, and use $\alpha\ge\slphi(b)/b$.
This gives a direct lower bound, without subtracting two tail approximations:
\begin{equation}
 D_r^{u,v}(\alpha/m)
 \ge\frac{c_w}{\sqrt{1-r}\,b^2}e^{-b^2/(1+r)},\qquad
 c_w=\frac{r_0\ell(1-e^{-2r_0})}{4\sqrt2\,\pi m}
                         e^{-K_m-K_m^2/2}>0.
 \label{sl:separatedpair}
\end{equation}
The bound is common to every allowed pair and every $b\ge B_m$,
$\Lambda\ge1$, and therefore applies eventually along the sequence.

We next compare this pair deficit with the entire residual. The support
condition requires three representative p-values to be small, and the
residual has absolute value at most \(C_q^R\). Taking a union over the
\(\binom q3\) triples and their eight sign patterns, and applying
Lemma~\ref{sl:triple} with \(r=r_n\) and \(b=b_n\) to each triple, gives,
with \(C_w=8\binom q3C_q^R3^{5/2}(2\pi)^{-3/2}\) and \(t_n=\alpha_n/m\),
\[
 |\slE_{\Gamma_n}R_w|
 \le
 \frac{C_w}{(1-r_n)b_n^3}
                  e^{-3b_n^2/[2(1+2r_n)]}.
 \]
Divide by~\eqref{sl:separatedpair}.  The difference of exponents is
\[
 \frac{3}{2(1+2r)}-\frac1{1+r}
       =\frac{1-r}{2(1+2r)(1+r)}.
 \]
It follows that
\begin{equation}
 \frac{|\slE_{\Gamma_n}R_w|}{D_{r_n}^{u,v}(t_n)}
 \le
 \frac{C_w}{c_w\sqrt{\Lambda_n}}
       e^{-\Lambda_n/[2(1+2r_n)(1+r_n)]}
 \le \frac{C_w}{c_w\sqrt{\Lambda_n}}e^{-\Lambda_n/12}
 \longrightarrow0.
 \label{sl:separatedratio}
\end{equation}
The remaining pair deficits are nonnegative. Once the last ratio
is at most \(1/2\), the exact decomposition gives
\[
 F_w(\Gamma_n,\alpha_n)
 \le\alpha_n-\tfrac12D_{r_n}^{u,v}(t_n)<\alpha_n
\]
on these indices, which proves the strict conclusion.
\end{proof}

\subsection{Signed contraction below the critical scale}
\label{sl:contractionsubsection}

The preceding argument does not apply to clusters whose internal
distances vanish at the critical scale.  For such clusters the
rejection probability divided by \(\alpha_n\) tends to one
(Lemma~\ref{sl:multicritical} below, with all limiting distances zero),
so a limit cannot decide on which side of \(\alpha_n\) it lies.  We
replace the members of such
a cluster by copies of its standardized mean.  To compare the two
rejection probabilities, we bound the chance that the original vector
rejects and the contracted vector does not, and then construct a more
likely event on which only the contracted vector rejects.  The first
bound uses the next lemma, which estimates the probability that a
coordinate exceeds a high threshold while differing substantially from
a nearly identical coordinate.

\begin{sllemma}[Exponential tilting with a marginal tail factor]
\label{sl:tilt}
Let \((X,D)\) be centered jointly Gaussian, possibly singular, with
\[
 \operatorname{Var}(X)=1,\qquad
 \operatorname{Var}(D)=2h,\qquad
 \operatorname{Cov}(X,D)=h>0.
 \]
If \(b>0\) and \(g>(b+b^{-1})h\), then
\begin{equation}
 \slP(X\ge b,D\ge g)
 \le \slPhi(b)
       \exp\left[-\frac{\{g-(b+b^{-1})h\}^2}{4h}\right].
 \label{sl:tiltbound}
\end{equation}
\end{sllemma}

\begin{proof}
We need a bound that retains the small probability of a marginal
exceedance as well as the rarity of a large difference. We obtain it by
weighting the Gaussian law exponentially in the difference. For
\(\theta\ge0\), this weighting shifts the mean of \(X\) by
\(\theta h\) and has normalizing factor \(e^{\theta^2h}\). Thus
\[
 \slE[e^{\theta D}\slone\{X\ge b\}]
     =e^{\theta^2h}\slPhi(b-\theta h).
 \]
The identity extends to a singular joint law through its Gaussian
moment-generating function or a nonsingular approximation.  Since the
normal hazard is increasing,
Lemma~\ref{sl:tails} gives, for all \(u\ge0\),
\[
 \log\frac{\slPhi(b-u)}{\slPhi(b)}
       =\int_{b-u}^b\frac{\slphi(v)}{\slPhi(v)}\,dv
       \le u\,\frac{\slphi(b)}{\slPhi(b)}\le u(b+b^{-1}).
 \]
Markov's inequality consequently bounds the required probability by
\[
 \slPhi(b)
     \exp\{-\theta g+\theta^2h+\theta(b+b^{-1})h\}.
 \]
Its minimum over \(\theta\ge0\) is attained at
\(\theta=\{g-(b+b^{-1})h\}/(2h)>0\), which proves the claim.
\end{proof}

\begin{sllemma}[Subcritical cluster contraction]
\label{sl:contraction}
Fix \(m\ge2\) and a partition \(B_1,\ldots,B_G\) of
\(\{1,\ldots,m\}\).
Let \(Z_n\sim N_m(0,\Sigma_n)\), \(\Sigma_n\in\slC_m\),
and let \(\alpha_n\downarrow0\), \(b_n=c_{\alpha_n}\).
Suppose that, after coordinate sign changes fixed along the
sequence within each block,
\begin{align}
 \delta_n&=\max_g\max_{i,j\in B_g}(1-\Sigma_{ij,n})>0,
 &\lambda_n&=b_n^2\delta_n\longrightarrow0,
 \label{sl:withincondition}\\
 b_n^2(1-|\Sigma_{ij,n}|)&\longrightarrow\infty,
 && i\in B_g,\ j\in B_h,\quad g\ne h.
 \label{sl:crosscondition}
\end{align}
Define the standardized means of the sign-changed coordinates
\[
 \overline Z_{g,n}=\frac1{|B_g|}\sum_{i\in B_g}Z_{i,n},
 \qquad
 \sigma_{g,n}^2=\operatorname{Var}(\overline Z_{g,n}),
 \qquad
 X_{g,n}=\overline Z_{g,n}/\sigma_{g,n},
 \]
and let \(\Gamma_n^c\) be the correlation matrix of
\((X_{1,n},\ldots,X_{G,n})\).
These means are well defined eventually and, with
\(w'=(|B_1|,\ldots,|B_G|)\),
\begin{equation}
 F_{\Sigma_n}^{(m)}(\alpha_n)
 \le F_{w'}(\Gamma_n^c,\alpha_n)
            -\frac{\alpha_n\sqrt{\lambda_n}}{160m^2}
 \quad\hbox{eventually}.
 \label{sl:contractionconclusion}
\end{equation}
More precisely,
\begin{equation}\label{sl:contraction_rate}
 \liminf_{n\to\infty}
 \frac{F_{w'}(\Gamma_n^c,\alpha_n)-F_{\Sigma_n}^{(m)}(\alpha_n)}
      {\alpha_n\sqrt{\lambda_n}}
 \ge\frac{2}{m\sqrt\pi}.
\end{equation}
On a subsequence on which a fixed block of size $w$ contains a pair
with $1-\Sigma_{ij,n}=\delta_n$, the lower bound is $w/(m\sqrt\pi)$. Thus every fixed constant strictly below
$2/(m\sqrt\pi)$ may replace $1/(160m^2)$ in
\eqref{sl:contractionconclusion}; the limiting constant itself need not
hold as an eventual finite-index bound. No positive-definiteness or lower
eigenvalue condition is imposed on either covariance matrix.  If instead
\(\delta_n=0\) at some index, a case excluded by~\eqref{sl:withincondition},
the aligned block members coincide and the original and contracted
indicators agree exactly.
\end{sllemma}

\begin{proof}
We work with the aligned vector, since coordinate sign changes leave
the two-sided p-values unchanged.  We suppress \(n\) when convenient.
Every within-block correlation lies in \([1-\delta,1]\), whence
\[
 1-\delta\le\sigma_g^2\le1.
 \]
In particular \(\sigma_g>0\) eventually.
Let \(H\) be the original rejection indicator and \(H^c\) the
indicator for the vector having \(|B_g|\) copies of \(X_g\).
We compare these indicators on the same probability space.  Their
difference is positive on a gained rejection, \(\{H=1,H^c=0\}\), and
negative on a lost rejection, \(\{H=0,H^c=1\}\).  Hence
\(F_{w'}(\Gamma^c,\alpha)-F^{(m)}_\Sigma(\alpha)
 =\slP(H=0,H^c=1)-\slP(H=1,H^c=0)\).  We show that a gain has
probability \(o(\alpha\sqrt\lambda)\), whereas a loss has probability at
least \(\{w/(m\sqrt\pi)+o(1)\}\alpha\sqrt\lambda\) along any subsequence
on which a fixed block of size \(w\ge2\) contains a pair attaining
\(\delta\).

\paragraph{Probability of gaining a rejection.}
Write
\[
 d_k=c_{k\alpha/m},\quad 1\le k\le m,\qquad d_m=b,
 \qquad
 g_\alpha=\min\left\{b,\min_{1\le k<m}(d_k-d_{k+1})\right\}.
 \]
The finite-set quantile expansion in Lemma~\ref{sl:tails} gives
\begin{equation}
 b g_\alpha\longrightarrow
       \log\frac{m}{m-1}>0.
 \label{sl:mingap}
\end{equation}
A gained rejection requires a difference within one block large enough
to disrupt the threshold comparison. Otherwise,
including all members of each selected block and replacing them by
copies of its standardized mean preserves a crossing. We make this observation precise through
the deterministic inclusion
\begin{equation}
 \{H=1,H^c=0\}
 \subseteq
 \bigcup_g\ \bigcup_{\substack{i,j\in B_g\\i\ne j}}\
 \bigcup_{\varepsilon=\pm1}
 \{\varepsilon Z_i\ge b,\
        \varepsilon(Z_i-Z_j)\ge g_\alpha\}.
 \label{sl:gaininclusion}
\end{equation}
If \(H=1\), choose a crossing rank \(k\) and \(k\) coordinates
with \(|Z_i|\ge d_k\).  Suppose no event on the right occurs.
For a selected coordinate \(i\), put
\(\varepsilon=\operatorname{sign}(Z_i)\).
Every other coordinate \(j\) in its block then satisfies
\[
 \varepsilon Z_j>\varepsilon Z_i-g_\alpha
                     \ge d_k-g_\alpha\ge0.
 \]
The coordinates in each selected block therefore have the same sign.
Complete these blocks, obtaining \(K\ge k\) coordinates.
If \(K=k\), every coordinate in each selected block was already selected;
their absolute block means are at least \(d_k\).
If \(K>k\), necessarily \(k<m\), and all such absolute block
means are at least
\[
 d_k-g_\alpha\ge d_{k+1}\ge d_K.
 \]
Dividing their means by \(\sigma_g\le1\) cannot decrease their
absolute values.  In both cases the \(K\) copies in the
contracted vector cross, contradicting \(H^c=0\).
This proves~\eqref{sl:gaininclusion}.

For any ordered pair in this inclusion, set
\(X=\varepsilon Z_i\), \(D=\varepsilon(Z_i-Z_j)\), and
\(h=1-\Sigma_{ij}\le\delta\).  If \(h=0\), then \(D=0\)
almost surely and the corresponding event is impossible.
Otherwise Lemma~\ref{sl:tilt} applies eventually.
Indeed, \(bg_\alpha\) has a positive limit while
\[
 b\delta(b+b^{-1})=\lambda(1+b^{-2})\longrightarrow0.
 \]
For \(0<h<g_\alpha/(b+b^{-1})\), differentiation shows that
\[
 \frac{d}{dh}
   \frac{\{g_\alpha-(b+b^{-1})h\}^2}{4h}
 =\frac14\left\{(b+b^{-1})^2-\frac{g_\alpha^2}{h^2}\right\}<0.
 \]
Replacing \(h\) by \(\delta\) therefore weakens the upper bound.
There are at most \(m(m-1)\) ordered pairs and two signs.  Using
\(2\slPhi(b)=\alpha\), we obtain from the union bound that
\begin{equation}
 \slP(H=1,H^c=0)
 \le m(m-1)\alpha
   \exp\left[
    -\frac{\{bg_\alpha-\lambda(1+b^{-2})\}^2}{4\lambda}
       \right]
 =o(\alpha\sqrt\lambda).
 \label{sl:gainbound}
\end{equation}
To obtain the final order bound, choose a positive constant smaller
than half the limit in~\eqref{sl:mingap}.  The numerator inside the
exponent is eventually bounded below by the square of this constant,
and \(e^{-c/\lambda}/\sqrt\lambda\to0\) for every fixed \(c>0\).  The bound retains the factor
\(\alpha\) for all relative rates at which \(\delta\) and \(\alpha\)
vanish.

\paragraph{Construction of a loss event.}
We construct a loss event within a block containing a pair \((i,j)\)
with \(1-\Sigma_{ij,n}=\delta_n\), the maximum in~\eqref{sl:withincondition}. There are only
finitely many choices of block and pair. We may therefore hold this
choice fixed along a subsequence: if the eventual bound failed, it
would fail along a subsequence with one fixed choice.
Fix a block \(B\) of size \(w\ge2\) and a pair \(i,j\in B\)
with \(1-\Sigma_{ij}=\delta\).  On the loss event, the standardized
block mean lies in a thin strip just above the full-block threshold
\(c_{w\alpha/m}\), some block coordinate stays below this threshold
while none reaches \(c_{(w-1)\alpha/m}\), and no coordinate outside the
block reaches \(b\) in absolute value.  Write
\[
 \overline Z=\frac1w\sum_{\ell\in B}Z_\ell,\qquad
 X=\overline Z/\sigma,\qquad
 W_\ell=\frac{Z_\ell-\overline Z}{\sqrt\delta},
 \quad \ell\in B.
 \]
The residuals describe the coordinates' deviations from the block mean.
Their sum is zero, so a sufficiently large difference between two
residuals forces at least one residual to be negative.  We use that coordinate to prevent the
original block from crossing its full-block threshold.
Let
\[
 a_\ell=\frac1w\sum_{j\in B}(1-\Sigma_{\ell j}),
 \qquad \overline a=\frac1w\sum_{\ell\in B}a_\ell.
 \]
The identities
\[
 \sigma^2=1-\overline a,\qquad
 \operatorname{Var}(Z_\ell-\overline Z)=2a_\ell-\overline a,
 \qquad
 \operatorname{Cov}(Z_\ell-\overline Z,\overline Z)
       =\overline a-a_\ell
 \]
and \(0\le a_\ell,\overline a\le\delta\) imply, for
\(\delta\le1/2\),
\begin{align}
 \sum_{\ell\in B}W_\ell&=0,\qquad
 \operatorname{Var}(W_\ell)\le2,\qquad
 |\operatorname{Cov}(W_\ell,X)|\le\sqrt{2\delta},
 \nonumber\\
 \operatorname{Var}(W_i-W_j)&=2,\qquad
 \operatorname{Var}(W_i-W_j\mid X)
       \ge2-8\delta.
 \label{sl:residualmoments}
\end{align}
These formulas remain valid for singular joint Gaussian laws because
we condition only on the single standard normal variable \(X\),
without inverting a block covariance matrix.

To separate the residual fluctuations from the block mean, subtract
their linear regressions on that mean. Put
$\kappa_\ell=\operatorname{Cov}(W_\ell,X)$ and
$R_\ell=W_\ell-\kappa_\ell X$. The centered Gaussian vector $R$ is
independent of $X$, may be singular, and has zero coordinate sum.
For $M=-\min_\ell R_\ell$ and $T=\max_\ell|R_\ell|$, symmetry gives, for
$\delta\le1/4$,
\begin{equation}\label{sl:residualcutoff}
 \Expect M=\tfrac12\Expect(\max R-\min R)
 \ge\tfrac12\Expect|R_i-R_j|
 \ge\sqrt{\frac{2-8\delta}{2\pi}},\qquad
 \Expect T^2\le2w.
\end{equation}
The equality holds only in expectation: pointwise, $M$ can be smaller
than half the range of $R$.  Because the pair $i,j$ attains $\delta$,
the lower bound tends to $1/\sqrt\pi$; for a pair with smaller
$1-\Sigma_{ij}$ it could tend to zero. We now put the block center just
above its full-block rejection threshold. A negative residual can then
pull one original coordinate below the threshold while every contracted
copy remains above it. Fix $A>0$, and consider
\[
 d_w=c_{w\alpha/m},\qquad X=x=d_w+u\sqrt\delta,
 \quad0<u<A,
 \qquad e_n=\max_\ell|\kappa_\ell|(d_w+A\sqrt\delta)=o(1).
\]
On $\{T\le A,M>u+e_n\}$, every $|W_\ell|\le A+e_n$ and some
$W_\ell<-u$. These inequalities will make the contracted block reject
while the original block fails.

\paragraph{Control of outside coordinates.}
We must also ensure that outside coordinates do not create another
rejection on the proposed loss event. The assumption separates them from
each block member; we first show that it also separates them from the
standardized block mean. This allows a conditional tail bound without
assuming independence from the block. Represent \(\Sigma\) by unit
Gram vectors \(v_\ell\).
For an anchor \(i\in B\), the vector
\(\overline v=w^{-1}\sum_{\ell\in B}v_\ell\) has norm \(\sigma\)
and
\[
 \|\overline v-v_i\|
       \le\frac1w\sum_{\ell\in B}\|v_\ell-v_i\|
       \le\sqrt{2\delta}.
 \]
Its normalized vector \(v_X=\overline v/\sigma\) satisfies
\(\|v_X-\overline v\|=1-\sigma\le\delta\); hence
\(\|v_X-v_i\|\le3\sqrt\delta\) for small \(\delta\).
For unit vectors define the projective distance
\[
 d_\pm(v,u)=\min\{\|v-u\|,\|v+u\|\}
                  =\sqrt{2(1-|\langle v,u\rangle|)}.
 \]
This distance treats a vector and its negative as the same point,
as required for two-sided p-values. On the resulting classes it is a
metric. The triangle inequality follows by combining the two minimizing
sign choices and applying the Euclidean triangle inequality.
For any outside coordinate \(j\notin B\), if
\(p_j=\operatorname{Cov}(Z_j,X)\), then
\[
 b\sqrt{2(1-|p_j|)}
 \ge b\sqrt{2(1-|\Sigma_{ij}|)}-3b\sqrt\delta
 \longrightarrow\infty
 \]
by~\eqref{sl:crosscondition} and \(b\sqrt\delta=\sqrt\lambda\to0\).
Thus \(b^2(1-|p_j|)\to\infty\).

Put \(h_j=1-|p_j|\).  Since \(x\le b+C_{m,A}/b\) on the strip,
eventually, uniformly there,
\[
 \frac{b-|p_j|x}{\sqrt{1-p_j^2}}
 \ge \frac{bh_j-C_{m,A}/b}{\sqrt{2h_j}}
 \ge \frac{b\sqrt{h_j}}{2\sqrt2}\longrightarrow\infty.
 \]
Here \(p_j^2<1\) eventually follows from the divergence just
proved.  The conditional law \(Z_j\mid X=x\) is normal with
mean \(p_jx\) and variance \(1-p_j^2\).  Its two tails are
bounded by twice the upper normal tail at the last display.
There are finitely many outside coordinates.  Therefore
\begin{equation}
 \sup_{d_w<x<d_w+A\sqrt\delta}
 \slP(\exists j\notin B:\ |Z_j|\ge b\mid X=x)
       \longrightarrow0.
 \label{sl:outside}
\end{equation}
When \(B\) contains the whole vector, this probability is zero.
The bound applies directly to the original outside coordinates.

\paragraph{Probability of losing a rejection.}
Write $\eta_n(A)$ for the conditional outside probability in
\eqref{sl:outside}. On the residual event above, provided there are no outside
exceedances, the contracted vector rejects through $w$ copies of
$X>d_w$. Since $Z_\ell=\sigma X+\sqrt\delta\,W_\ell$ with $\sqrt{1-\delta}\le\sigma\le1$, all coordinates in $B$ are positive eventually, while
\[
 \min_{\ell\in B}Z_\ell<d_w,\qquad
 \max_{\ell\in B}Z_\ell
 \le d_w+(2A+e_n)\sqrt\delta<d_{w-1}.
\]
The second inequality uses
$b(d_{w-1}-d_w)\to\log(w/(w-1))>0$ and
$b\sqrt\delta\to0$. Ranks below $w$ fail by the maximum bound;
rank $w$ fails by the minimum bound; larger ranks fail because the
outside has no absolute coordinate above $b$. Thus this event gives
$H=0,H^c=1$.

To estimate the loss probability, first average the possible width of
the strip over the residual fluctuations, and then subtract the chance
of an outside exceedance. The residual vector $R$ need not be independent
of the outside coordinates; a conditional union bound is enough.
Tonelli's theorem gives the following bound for the residual event:
\[
 \int_0^A\Prob(T\le A,M>u+e_n)\,du
 =\Expect[(M-e_n)_+\mathbf1\{T\le A\}]
 \ge\Expect M-\frac{2w}{A}-e_n.
\]
For the nonnegative residual contribution use the minimum density on
the strip; for the subtracted outside contribution use its maximum.
Reflecting the strip by a simultaneous sign change gives an equally
large disjoint contribution. Hence
\[
 \Prob(H=0,H^c=1)\ge2\sqrt\delta\left[
 \slphi(d_w+A\sqrt\delta)
       \left(\Expect M-\frac{2w}{A}-e_n\right)
 -\slphi(d_w)A\eta_n(A)\right].
\]
If the residual lower bound in parentheses is negative, replace it by
zero before using the minimum density; the displayed weaker inequality
still follows. For each fixed $A$,
$\slphi(d_w+A\sqrt\delta)/\slphi(b)\to w/m$ and
$2\slphi(b)/(\alpha b)\to1$. First take the lower limit along the
sequence, and then let $A\to\infty$. Equations
\eqref{sl:residualcutoff} and \eqref{sl:outside} yield
\begin{equation}\label{sl:lossbound}
 \liminf_n\frac{\Prob(H=0,H^c=1)}{\alpha_n\sqrt{\lambda_n}}
 \ge\frac{w}{m\sqrt\pi}.
\end{equation}
The gain estimate \eqref{sl:gainbound} is $o(\alpha\sqrt\lambda)$.
This proves the stronger fixed-block rate.
Only finitely many blocks and pairs can attain \(\delta_n\). Passing to a fixed choice on any potentially violating subsequence proves
\eqref{sl:contraction_rate}, since $w\ge2$, and hence
\eqref{sl:contractionconclusion}. At zero spread, aligned block members
coincide exactly and both indicators agree.
\end{proof}

\subsection{Finite critical clusters}
\label{sl:criticalsubsection}

We now consider clusters with finite limiting distances at the critical
scale.  Conditioning on one coordinate near the rejection threshold
gives a Gaussian limit for the relative fluctuations of the other
coordinates.  We call the conditioning coordinate the anchor.  The next
lemma expresses the limiting rejection ratio in terms of this Gaussian
vector and identifies when the ratio is strictly below one.  The argument
also covers singular limits.

\begin{sllemma}[One finite critical cluster]
\label{sl:critical}
Fix \(s\ge1\).  Let \(\alpha_n\downarrow0\), \(b_n=c_{\alpha_n}\),
and \(Z_n\sim N_s(0,\Sigma_n)\), where
\[
 b_n^2(1-\Sigma_{ij,n})\longrightarrow L_{ij}<\infty,
 \qquad 1\le i,j\le s.
 \]
There is a possibly singular Gaussian vector \(Y\) with
\[
 \slE Y_i=-L_{1i},\qquad
 \operatorname{Cov}(Y_i,Y_j)=L_{1i}+L_{1j}-L_{ij},
 \qquad Y_1=0.
 \]
Write \(Y_{[1]}\ge\cdots\ge Y_{[s]}\) for descending order
statistics.  Then
\begin{equation}
 \frac{F^{(s)}_{\Sigma_n}(\alpha_n)}{\alpha_n}
 \longrightarrow
 \Psi_s(L):=
 \slE\max_{1\le k\le s}\frac{k}{s}e^{Y_{[k]}}
 \le1.
 \label{sl:criticalfunctional}
\end{equation}
The inequality is strict if some \(L_{ij}>0\); if every
\(L_{ij}=0\), then \(\Psi_s(L)=1\).
\end{sllemma}

\begin{proof}
We proceed in three steps.  First, we condition on an anchor coordinate
to obtain the limit for positive exceedances.  Second, we justify taking
expectations and recover the two-sided rejection probability.  Finally,
a deterministic inequality gives the upper bound and its equality case.
Let \(U_n=Z_{1,n}\), \(a_{i,n}=\Sigma_{1i,n}\), and
\[
 V_{i,n}=b_n(Z_{i,n}-a_{i,n}U_n).
 \]
The vector \(V_n\) is independent of \(U_n\), and its covariance
matrix has entries
\[
 b_n^2(\Sigma_{ij,n}-a_{i,n}a_{j,n})
 \longrightarrow L_{1i}+L_{1j}-L_{ij}.
 \]
The limiting matrix is positive semidefinite, being a limit of
covariance matrices.  To compare the residual vectors on one probability
space, we couple \(V_n\) as \(K_n^{1/2}N\), where
\(N\) is a fixed standard Gaussian vector in \(\mathbb R^s\) and
\(K_n\) is the covariance matrix of the residual vector.  Continuity of
the positive semidefinite square root
gives \(V_n\to V\) almost surely.  Set \(Y_i=V_i-L_{1i}\).
The first row and column of \(K_n\) are zero, so \(Y_1=0\).
Moreover
\[
 \operatorname{Var}(Y_i)=2L_{1i},\qquad
 \slE e^{Y_i}=\exp(-L_{1i}+L_{1i})=1.
 \]

We first hold the residual fluctuations fixed and find how large the
anchor must be for rejection. Let \(H_n^+\) denote rejection using
positive exceedances alone, that is, at least \(k\) coordinates
\(Z_{i,n}\) exceed \(c_{k\alpha_n/s}\) for some \(k\).
All \(a_{i,n}>0\) eventually. Conditional on \(V_n\), the
event is \(U_n\ge A_n\), where
\[
 A_n=\min_{1\le k\le s}
       \operatorname{kth\,smallest}_{1\le i\le s}
       \left\{\frac{c_{k\alpha_n/s}-V_{i,n}/b_n}{a_{i,n}}\right\}.
 \]
Since \(U_n\) is standard normal and independent of \(V_n\),
\(\slP(H_n^+\mid V_n)=\slPhi(A_n)\), and hence
\(2\slP(H_n^+)/\alpha_n=\slE\{\slPhi(A_n)/\slPhi(b_n)\}\).  We first
find the almost sure limit of this ratio.  The quantile expansion and
\(b_n^2(1-a_{i,n})\to L_{1i}\) give, almost surely,
\[
 b_n\left\{
 \frac{c_{k\alpha_n/s}-V_{i,n}/b_n}{a_{i,n}}-b_n
 \right\}
 \longrightarrow \log(s/k)-Y_i.
 \]
For each rank, rejection requires the anchor to exceed the corresponding
ordered threshold. Taking the minimum over ranks gives the conditional
rejection threshold. The minimum and order-statistic maps are continuous
on finite-dimensional Euclidean space, including at ties. Consequently
\[
 b_n(A_n-b_n)
 \longrightarrow
       \min_{1\le k\le s}\{\log(s/k)-Y_{[k]}\}.
 \]
By~\eqref{sl:tailshift},
\begin{equation}
 \frac{\slPhi(A_n)}{\slPhi(b_n)}
 \longrightarrow
       \max_{1\le k\le s}\frac{k}{s}e^{Y_{[k]}}
 \quad\hbox{almost surely}.
 \label{sl:conditionalcritical}
\end{equation}

Almost sure convergence alone is not enough to pass to expectations.
We obtain uniform integrability from the exact marginal probabilities.
Define
\[
 T_{i,n}=
 \frac{\slPhi((b_n-V_{i,n}/b_n)/a_{i,n})}{\slPhi(b_n)}.
 \]
The numerator is the conditional exceedance probability of the standard
normal marginal \(Z_{i,n}\) above \(b_n\), so
\(\slE T_{i,n}=1\) exactly.
The same expansion gives \(T_{i,n}\to e^{Y_i}\) almost surely.
Since \(\slE e^{Y_i}=1\),
\[
 \slE|T_{i,n}-e^{Y_i}|
 =2-2\slE\min\{T_{i,n},e^{Y_i}\}\longrightarrow0
 \]
by dominated convergence, the minimum being bounded by
the integrable \(e^{Y_i}\).
Thus each sequence \(T_{i,n}\) is uniformly integrable,
as is their finite sum.
Every threshold \(c_{k\alpha_n/s}\) is at least \(b_n\),
so \(H_n^+\) implies at least one positive exceedance of
\(b_n\).  A conditional union bound yields
\[
 0\le\frac{\slPhi(A_n)}{\slPhi(b_n)}
       \le\sum_{i=1}^s T_{i,n}.
 \]
The left side of~\eqref{sl:conditionalcritical} is therefore uniformly
integrable, which gives convergence of its expectation.

We next compare positive rejection with two-sided rejection.
If a rejection involves both signs, there are indices \(i,j\)
with \(Z_{i,n}\ge b_n\) and \(Z_{j,n}\le-b_n\).
Their difference has variance \(2(1-\Sigma_{ij,n})=O(b_n^{-2})\).
Its probability of exceeding \(2b_n\) is bounded by
\(\slPhi(c b_n^2)=o(\slPhi(b_n))\) for a fixed positive
\(c\), provided its variance is nonzero; if its variance is
zero the event is impossible.  A common \(c\) applies to the finitely
many pairs after increasing their common variance upper bound if
necessary.  A union bound shows that all mixed-sign
events have probability \(o(\alpha_n)\).
By symmetry, positive and negative rejection have equal
probabilities, and their intersection is included in the same
mixed-sign bound.  Hence
\[
 F^{(s)}_{\Sigma_n}(\alpha_n)
       =2\slP(H_n^+)+o(\alpha_n).
 \]
Since \(\alpha_n=2\slPhi(b_n)\), expectation convergence
in~\eqref{sl:conditionalcritical} proves the limit.

To prove the upper bound and identify its equality case, take any vector
\(x\in(0,\infty)^s\) with descending order statistics \(x_{[k]}\).
Then
\[
 \max_k\frac{k}{s}x_{[k]}
       \le\frac1s\sum_{i=1}^s x_i.
 \]
For \(k<s\) the positive coordinates after rank \(k\) make
the inequality strict; for \(k=s\) equality holds only if
all coordinates are equal.  Taking \(x_i=e^{Y_i}\) proves
\(\Psi_s(L)\le1\).
If \(L_{ij}>0\), then
\(\operatorname{Var}(Y_i-Y_j)=2L_{ij}>0\), so, with probability one,
the coordinates are not all equal.  The resulting positive integrable difference has strictly positive
expectation.
If all \(L_{ij}=0\), \(Y=0\) almost surely and the functional
equals one.
\end{proof}

\begin{sllemma}[Several finite critical clusters]
\label{sl:multicritical}
Fix a partition \(B_1,\ldots,B_G\) of \(m\) physical coordinates
and put \(n_g=|B_g|\).  Let \(\alpha_n\downarrow0\),
\(b_n=c_{\alpha_n}\), and \(\Sigma_n\in\slC_m\).
Suppose that, after fixed within-block sign changes,
\[
 b_n^2(1-\Sigma_{ij,n})\longrightarrow L_{ij}^{(g)}<\infty
       \quad(i,j\in B_g),
 \]
and that
\(b_n^2(1-|\Sigma_{ij,n}|)\to\infty\) for coordinates in
different blocks.  Then
\begin{equation}
 \frac{F^{(m)}_{\Sigma_n}(\alpha_n)}{\alpha_n}
 \longrightarrow
       \sum_{g=1}^G\frac{n_g}{m}\Psi_{n_g}(L^{(g)})
 \le1.
 \label{sl:multicriticalformula}
\end{equation}
The limit is strictly below one if a within-block limiting
distance is positive.
\end{sllemma}

\begin{proof}
The blocks need not be independent. Instead, their separation makes it
unlikely that two blocks contain small p-values at the same time. We first
show that this probability is negligible relative to the nominal level.
Rejections can then be assigned to one block, up to that error.  For any bivariate
centered standard Gaussian vector \((Z_i,Z_j)\) with correlation \(\rho\), a union over four sign pairs and the
normal law of the corresponding signed sums gives
\[
 \slP(|Z_i|\ge b,|Z_j|\ge b)
 \le4\slPhi\left(b\sqrt{\frac2{1+|\rho|}}\right).
 \]
For \(c\ge1\), integration of the lower hazard bound
in~\eqref{sl:mills} gives
\[
 \frac{\slPhi(cb)}{\slPhi(b)}
 \le\exp\{-b^2(c^2-1)/2\}.
 \]
Therefore the cross-block pair probability divided by
\(\alpha_n=2\slPhi(b_n)\) is at most
\[
 2\exp\left\{-\frac{b_n^2(1-|\rho|)}{2(1+|\rho|)}\right\}
       \longrightarrow0.
 \]
A union over finitely many cross-block pairs shows that the
probability of p-values at most \(\alpha_n\) in two different
blocks is \(o(\alpha_n)\).
Outside this event, all coordinates contributing to a rejection belong
to one block.  The block rejection rule retains the original denominator
\(m\), and thus equals the ordinary \(n_g\)-coordinate Simes
rule at level \(n_g\alpha_n/m\).
Intersections of two such block events are included in the
same negligible cross-block event.  It follows that
\[
 F^{(m)}_{\Sigma_n}(\alpha_n)
 =\sum_{g=1}^G
     F^{(n_g)}_{\Sigma_n[B_g]}(n_g\alpha_n/m)
       +o(\alpha_n).
 \]
By~\eqref{sl:quantile},
\(c_{n_g\alpha_n/m}^2/b_n^2\to1\).
Lemma~\ref{sl:critical} applies to every summand with the
same limiting distances \(L^{(g)}\).  Any member of \(B_g\) may serve as
the first coordinate there: the left side
of~\eqref{sl:criticalfunctional} does not depend on how the
coordinates are ordered, so neither does its limit.
This gives~\eqref{sl:multicriticalformula}. Under the stated
positive-distance condition, strictness follows because all weights
\(n_g/m\) are positive.
\end{proof}

\subsection{Proof of the uniform small-level theorem}
\label{sl:inductionsubsection}

\begin{proof}[Proof of Theorem~\ref{sl:uniform}]
Fix the total number of tests \(m\).  We induct on the number \(q\)
of representatives, retaining their duplication multiplicities.  Our
claim is that each positive integer vector \(w\) of length \(q\)
and total \(m\) admits a positive cutoff common to all
\(\Gamma\in\slC_q\), below which
\(F_w(\Gamma,\alpha)\le\alpha\), with equality exactly when
\(\Gamma\in\mathcal E_q\). This strengthened claim includes every
multiplicity vector and the equality cases, because contraction
produces a shorter vector \(w'\) that is not known in advance, and exact
contraction needs strict inequality outside the equality class.  For
the finitely many vectors of smaller length, we use the minimum of
their cutoffs. Equality on
\(\mathcal E_q\) already holds at every level by
Lemma~\ref{sl:pairidentity}. For \(q=1\) all laws belong to this class.
For \(q=2\), the same lemma gives a strictly positive pair deficit
when \(0<|\Gamma_{12}|<1\); the remaining correlations, zero and
\(\pm1\), are exactly \(\mathcal E_2\). This proves the base cases
at all levels below one.

Suppose the assertion holds for fewer than \(q\) representatives
but fails for some fixed \(w\) of length \(q\).
Since equality on \(\mathcal E_q\) is known, failure means that no
common cutoff gives strictness outside this class. Choose a sequence
\(\alpha_n\downarrow0\) and
\(\Gamma_n\in\slC_q\setminus\mathcal E_q\) such that
\begin{equation}
 F_w(\Gamma_n,\alpha_n)\ge\alpha_n
       \quad\hbox{for every }n.
 \label{sl:countersequence}
\end{equation}
Write \(b_n=c_{\alpha_n}\).
By passing to a subsequence for the finitely many pairs, we may assume
that every
\(b_n^2(1-|\Gamma_{ij,n}|)\) has a limit in the compact
extended interval \([0,\infty]\).

We group representatives whose differences remain bounded at the
critical scale. Define the relation \(i\sim j\) if this limit is finite (including
\(i=j\)). To make these groups well defined, we must check transitivity.
Realize each \(\Gamma_n\) by unit Gram vectors \(v_{i,n}\).
The projective distance
\(\min\{\|v-u\|,\|v+u\|\}=\sqrt{2(1-|\langle v,u\rangle|)}\) between unit
vectors satisfies the triangle inequality (see the proof of
Lemma~\ref{sl:contraction}), which gives
\[
 \sqrt{2(1-|\Gamma_{ik,n}|)}
 \le \sqrt{2(1-|\Gamma_{ij,n}|)}
       +\sqrt{2(1-|\Gamma_{jk,n}|)}.
 \]
If \(i\sim j\) and \(j\sim k\), then \(b_n\) times the right side stays
bounded, and hence so does \(b_n^2(1-|\Gamma_{ik,n}|)\).
Its already selected limit is therefore finite.
Thus \(\sim\) partitions the representatives into classes.
The rescaled distances between different classes diverge.

Within each class, choose an anchor and signs aligning each
Gram vector with that anchor.
Each aligned vector is \(O(b_n^{-1})\) from the anchor;
therefore every pair of aligned vectors has inner product
tending to one.  In particular, for each pair, the aligned inner product
eventually equals the original absolute correlation.
The rescaled signed distances within each class have the finite limits
selected above.  Since there are finitely many
sign configurations, pass to one further subsequence so
that all signs are fixed.
Changing these signs does not affect any two-sided p-value.

Three cases exhaust the possibilities: some within-class limiting
distance is positive; all are zero and every class is a single
representative; or all are zero and some class has at least two
representatives.  Each case contradicts~\eqref{sl:countersequence}.

\paragraph{A positive finite limiting distance.}
If some finite within-class limiting distance is positive,
expand representative \(i\) into its \(w_i\) physical copies.
Apply Lemma~\ref{sl:multicritical} to the resulting physical
classes.  It gives a limiting rejection ratio strictly
below one, contradicting~\eqref{sl:countersequence}.

\paragraph{Singleton classes.}
It remains to consider the case where every finite within-class
limiting distance is zero.  If every class is a singleton representative,
\[
 b_n^2\left(1-\max_{i<j}|\Gamma_{ij,n}|\right)\longrightarrow\infty,
 \]
because there are finitely many pairs and each rescaled
distance diverges.
Since \(\Gamma_n\notin\mathcal E_q\) implies \(\Gamma_n\ne I_q\),
the strict conclusion of Lemma~\ref{sl:separated} contradicts
\eqref{sl:countersequence}.

\paragraph{Classes approaching exact duplication.}
Otherwise there are \(G<q\) classes.
Expand into physical copies and partition by those classes;
their physical sizes form a positive integer vector \(w'\)
with total \(m\) and length \(G\).
The largest aligned within-class value of one minus a correlation,
\(\delta_n\),
satisfies \(b_n^2\delta_n\to0\), and every cross-class
rescaled absolute distance diverges.
At indices where \(\delta_n=0\), exact contraction replaces
the classes by standard Gaussian centers. The contracted covariance
\(\Gamma_n^c\) is outside \(\mathcal E_G\). Indeed, if its centers
formed independent signed-copy groups, expanding each center into its
original aligned representatives would place \(\Gamma_n\) in
\(\mathcal E_q\), a contradiction. Conversely, contracting exact
copies of a law in \(\mathcal E_q\) preserves this structure.
Thus exact contraction and the strengthened induction hypothesis give
\[
 F_w(\Gamma_n,\alpha_n)
       =F_{w'}(\Gamma_n^c,\alpha_n)<\alpha_n
 \]
eventually.
If infinitely many indices have \(\delta_n>0\), restrict
to them and apply Lemma~\ref{sl:contraction}, obtaining
\[
 F_w(\Gamma_n,\alpha_n)
 <F_{w'}(\Gamma_n^c,\alpha_n)\le\alpha_n
 \quad\hbox{eventually}.
 \]
The covariance \(\Gamma_n^c\) may vary with \(n\) and may be singular.
The induction hypothesis applies because it is uniform over all
correlation matrices in the reduced dimension, including the singular
ones.
Both possibilities contradict~\eqref{sl:countersequence}.

This completes the induction.  Since there are finitely many ordered
positive integer multiplicity vectors with total \(m\), their cutoffs
have a positive minimum.  We reduce this minimum below one if necessary.
Applying the result to \(w=(1,\ldots,1)\) proves the upper bound
and the equality classification. Independence attains the upper bound,
so the supremum in~\eqref{sl:main} is exactly \(\alpha\).
The same finite minimum supplies the asserted common cutoff for every
multiplicity vector, completing the proof.
\end{proof}

\begin{slremark}[The role of the signed bounds]
\label{sl:scope}
A first-order approximation \(F_\Sigma(\alpha)=\alpha+o(\alpha)\)
alone does not determine the sign of the remainder, and pointwise
cutoffs alone do not give a uniform cutoff by compactness.  The proof
therefore distinguishes three regimes.  Nonzero finite critical
distances yield a strictly conservative limit.  Divergent distances
are controlled by Lemma~\ref{sl:weakcube} near independence and by the
relative bound~\eqref{sl:separatedratio} otherwise, whereas
zero finite distances require the signed
contraction~\eqref{sl:contractionconclusion}.  In the last regime, the
limiting ratio can equal one and the covariance may approach a singular
duplication law arbitrarily fast.  The number \(m\) remains fixed
throughout the induction.
\end{slremark}

\clearpage
\section{\texorpdfstring{Proof of the trivariate result at level \(1/20\)}{Proof of the trivariate result at level 1/20}}

\label{tri:proof}

For a positive semidefinite correlation matrix \(\Sigma\), let \(Z\sim
N_3(0,\Sigma)\), \(P_i=2\overline\Phi(|Z_i|)\), and
\[
 H=\mathbf1\{P_{(1)}\le\alpha/3\text{ or }P_{(2)}\le2\alpha/3
                     \text{ or }P_{(3)}\le\alpha\}.
\]
Until Section~\ref{tri:level-interval}, we write
\(F_\Sigma=F_\Sigma(\alpha)=\mathbb E_\Sigma H\) and fix
\(\alpha=1/20\). Write \(J_k\) for the \(k\times k\) all-ones matrix. A
\emph{perfect pair} means \(|\Sigma_{ij}|=1\) for some \(i\ne j\).

\begin{tritheorem}[Simes inequality for every trivariate covariance]
\label{tri:theorem}
For every positive semidefinite $3\times3$ correlation matrix,
\[
\Prob\{P_{(1)}\le1/60\text{ or }P_{(2)}\le1/30
                          \text{ or }P_{(3)}\le1/20\}\le1/20.
\]
Equality holds exactly at coordinate-sign changes and permutations of
\begin{equation}\label{tri:equalities}
I_3,\qquad J_3,\qquad \operatorname{diag}(J_2,1).
\end{equation}
The proof uses a finite computation under the arithmetic
assumptions stated in Section~\ref{tri:reproducibility}.
\end{tritheorem}

We divide the covariance space into three regions, following the outline
in Section~\ref{main:tri-proof-outline}. Lemma~\ref{tri:cube} establishes the
inequality near independence, where every absolute correlation is at most
\(7/16\), and Lemma~\ref{tri:pair} near perfect pairs, where some absolute
correlation is at least \(2047/2048\). Section~\ref{tri:completion} combines
these lemmas with a finite computational certificate,
Proposition~\ref{tri:finite-certificate}, which covers the
remaining feasible covariance matrices: it gives verified bounds on whole
regions, rather than only values at selected matrices. This argument is
independent of the covariance-uniform small-level theorem (Theorem~\ref{main:uniform}), and both
neighborhood sizes are specified. Together with Lemmas~\ref{tri:cube} and~\ref{tri:pair}, the certificate
establishes the result at level \(1/20\) in dimension three. Section~\ref{tri:level-interval}
then proves validity on a common interval around this level by retaining
the analytic bounds near equality and the certified margin elsewhere.

\subsection{Preliminaries and covariance geometry}

\label{tri:s2}

Write \(\Phi\) and \(\phi\) for the standard normal distribution and density, and \(\overline\Phi=1-\Phi\). Throughout the proof, set
\[
 \alpha=1/20,\qquad A=c_{\alpha/3},\quad B=c_{2\alpha/3},\quad C=c_\alpha,
 \quad c_u=\Phi^{-1}(1-u/2),
\]
and write
\begin{equation}
\Sigma(a,b,c)=\begin{pmatrix}1&a&b\\a&1&c\\b&c&1\end{pmatrix},\qquad
 \Delta(a,b,c)=1-a^2-b^2-c^2+2abc.                              \label{tri:eq:2-1}
\end{equation}
The three normal thresholds correspond, in order, to the three Simes
rejection thresholds. To use them in rigorous numerical bounds, we enclose
each true quantile between exact rational endpoints.
Let \(Q=2^{40}\). The quantile brackets used by the computation are
\[
\frac{2632208626561}{Q}<A<\frac{2632208626562}{Q},\qquad
\frac{2339810479419}{Q}<B<\frac{2339810479420}{Q},
\]
\[
\frac{2155003191023}{Q}<C<\frac{2155003191024}{Q}.
\]
These brackets can be checked with exact arithmetic or high-precision
interval arithmetic.
For reference, the decimal approximations of \(A,B,C\) are
\begin{equation}
A\approx2.393979799819,\quad B\approx2.128045234185,
 \quad C\approx1.959963984540.                                 \label{tri:eq:2-2}
\end{equation}
The brackets imply the following inequalities, which we use below:
\begin{equation}
\begin{gathered}
  19/10<C<49/25,\quad 8/5<C<13/5,\quad 9/5<B<14/5,\quad A<3,
  \quad A>2393/1000,\\
  A-B>1/5,\qquad 13/100<B-C<1/5.                              \label{tri:eq:2-3}
 \end{gathered}
\end{equation}

\subsubsection{One duplicated p-value: multiplicities (2,1)}

\label{tri:s2-1}

A perfect pair has two identical p-values, but both still count in the
three-test rule. We first calculate its rejection probability; this will
also serve as the comparison law for a nearly perfect pair. For a standard
bivariate normal pair \((U,W)\) with absolute correlation \(r\), let \(P_U=2\overline\Phi(|U|)\) and \(P_W=2\overline\Phi(|W|)\), and put
\(J_r(u,v)=\Prob(P_U\le u,P_W\le v)\). The three-test rule applied to
\((P_U,P_U,P_W)\) rejects when
\[
 P_U\le2\alpha/3\quad\text{or}\quad P_W\le\alpha/3
 \quad\text{or}\quad(P_U\le\alpha, P_W\le\alpha).
\]
Inclusion--exclusion and bivariate exchangeability give
\begin{equation}
F_{(2,1)}(r,\alpha)=\alpha-D_r,\qquad
 D_r=J_r(\alpha/3,\alpha)+J_r(2\alpha/3,\alpha)-J_r(\alpha,\alpha).
                                                                    \label{tri:eq:2-4}
\end{equation}
Writing \(f(z)=\Prob(P_W\le\alpha\mid P_U=z)\), the same quantity is
\begin{equation}
D_r=\int_0^{\alpha/3}\{f(z)-f(z+2\alpha/3)\}\,dz.              \label{tri:eq:2-5}
\end{equation}
For \(0<r<1\), with \(\sigma=\sqrt{1-r^2}\) and \(N\sim N(0,1)\),
\begin{equation}
q_r(y)=\Prob(|ry+\sigma N|\ge C),\qquad
 q_r'(y)=\frac r\sigma\phi\!\left(\frac{C-ry}{\sigma}\right)
       \left(1-e^{-2rCy/\sigma^2}\right)>0\quad(y>0).             \label{tri:eq:2-6}
\end{equation}
Since \(f(z)=q_r(c_z)\) and \(c_z\) decreases, \(D_r>0\) for \(0<r<1\).
At \(r=0\) or 1, \eqref{tri:eq:2-4} gives \(D_r=0\). Throughout this calculation,
the thresholds are those of the original three-test rule.

For the independence calculation, divide the p-value range into four bins.
The first three end at the successive Simes thresholds, and the fourth
contains all larger p-values. Under independence, these bins have probabilities
\((t,t,t,1-3t)\), \(t=\alpha/3\). We partition the rejection event according
to whether there is a first-bin coordinate, at least two second-bin
coordinates but no first-bin coordinate, or a remaining rejection
configuration in which all three p-values lie in the first three bins. Summing these probabilities gives
\begin{equation}
1-(1-t)^3+3t^2(1-2t)+t^3+4t^3=3t=\alpha.                      \label{tri:eq:2-7}
\end{equation}
With \(D_0=D_1=0\) in \eqref{tri:eq:2-4}, all matrices in \eqref{tri:equalities} therefore give equality, including arbitrary coordinate signs.

\subsubsection{Continuity and the canonical chambers}

\label{tri:s2-2}

Changing coordinate signs or permuting coordinates leaves the p-value multiset
unchanged. In dimension three, the two largest absolute correlations involve a common
coordinate. Sign changes can make both correlations nonnegative, and a
permutation can put them in decreasing order. Consequently every correlation matrix is
represented by
\begin{equation}
0\le d\le b\le a\le1,\qquad c=\varepsilon d,\quad
 \varepsilon\in\{-1,1\},\qquad \Delta(a,b,c)\ge0.                 \label{tri:eq:2-8}
\end{equation}
The two choices of sign give two closed parameter regions, which we call
canonical (or closed) chambers. We include both regions and retain their common
boundary, where the product of the three correlations is zero. In these
coordinates \(a=\max_{i<j}|\Sigma_{ij}|\), so Lemma~\ref{tri:cube} applies
when \(a\le7/16\) and Lemma~\ref{tri:pair} when \(a\ge2047/2048\);
Section~\ref{tri:s7-2} classifies analytic boxes by these two tests.

The rejection probability is continuous on the full space of positive
semidefinite correlation matrices. To see this, write each Gaussian vector as
\(\Sigma^{1/2}X\), with one common standard normal vector \(X\) and the continuous positive semidefinite square root \(\Sigma^{1/2}\). The discontinuity set of the rejection indicator is contained in the
nine marginal surfaces
\(|z_i|\in\{A,B,C\}\). Each surface has probability zero because the marginal distributions
remain standard normal, including when the joint law is singular.

The following geometric property will allow us to bound the error from varying
the covariance within a cell.

\begin{trilemma}[Feasible segments]\label{tri:segments}
A strict convex combination of two distinct positive semidefinite trivariate correlation
matrices is singular if and only if the endpoints share the same signed
perfect pair. In particular, it is positive definite if either endpoint
has no perfect pair.
\end{trilemma}
\begin{proof}
We first consider a singular matrix without a perfect pair. It has rank two,
and every nonzero kernel vector has three nonzero coordinates. Indeed, a
zero component would imply a linear relation between two unit-variance
Gaussian coordinates and hence a perfect pair. For a fixed vector \(u\) with all components nonzero, a symmetric
zero-diagonal difference matrix \(D\) with edges \(x,y,z\) and \(Du=0\) obeys
\[
 xu_1u_2+yu_1u_3=0,\quad xu_1u_2+zu_2u_3=0,\quad
 yu_1u_3+zu_2u_3=0.
\]
Adding the first two and subtracting the third gives \(x=0\), and similarly
\(y=z=0\). Thus two distinct trivariate correlation matrices cannot share a kernel
vector with three nonzero coordinates. For positive semidefinite matrices, the kernel of a
strict convex combination is the intersection of the two kernels. If that
intersection is nonzero, a vector in it must therefore have a zero
coordinate. Its other two coordinates are nonzero, because the diagonal
entries are one; the common relation between those two unit-variance
coordinates gives the same signed perfect pair at both endpoints.
Conversely, a common signed perfect pair supplies a common kernel vector.
This proves both directions.\end{proof}

When two distinct feasible endpoints have no common signed perfect pair,
this lemma makes the segment interior positive definite. We can then
integrate a bounded covariance gradient along that interior and use
continuity at possibly singular endpoints.
Thus a bound at one reference matrix can be extended to nearby feasible
matrices without assuming where an extremum occurs.
Section~\ref{tri:s7-1} uses the lemma in this way to prove the whole-cell
bound \eqref{tri:eq:7-1}; Lemmas~\ref{tri:cube} and~\ref{tri:pair} do not
use it.

\subsection{A neighborhood of independence}

\label{tri:s3}

\begin{trilemma}[Independence neighborhood]\label{tri:cube} At \(\alpha=1/20\),
\begin{equation}
\max(|a|,|b|,|c|)\le7/16
 \quad\Longrightarrow\quad
 F_\Sigma\le\alpha-\frac{a^2+b^2+c^2}{500}.                    \label{tri:eq:3-1}
\end{equation}
All signs are included. In particular, $F_\Sigma<\alpha$ away from independence.
\end{trilemma}

We expand the rejection probability around independence. Its quadratic term
is negative, and we bound all higher-order terms to show that they cannot
offset this decrease. High-precision integer arithmetic that rounds outward verifies the resulting
scalar inequalities.

\subsubsection{Pairwise decomposition}

\label{tri:s3-1}

We separate the changes attributable to pairs from the part involving
all three p-values. This makes the truncation error in
Section~\ref{tri:s3-2} small: under independence the three-way remainder
\(R\) defined below has second moment \(4\alpha^3/9\), compared with \(\alpha\)
for \(H\), and each pair term involves only one correlation. Recall that \(H\) is the three-test indicator. For \(K\subseteq\{1,2,3\}\), let \(H_K\)
be the same indicator after the p-values outside \(K\) are set to one,
without changing the three Simes thresholds.
Then
\begin{equation}
R=H-\sum_{i<j}H_{\{i,j\}}+\sum_iH_{\{i\}},\qquad
 F=\alpha-\sum_{i<j}D_{|\rho_{ij}|}(t)+\Expect_\Sigma R,              \label{tri:eq:3-2}
\end{equation}
where \(\rho_{ij}=\Sigma_{ij}\) and, in this subsection,
\(D_s(t)=2J_s(t,2t)-J_s(2t,2t)\).
The nonzero sorted patterns of \(R\), using the four bins of Section~\ref{tri:s2-1}, numbered 1 through 4,
are
\[
 111:1,\quad122:-1,\quad222:-2,\quad233:1,\quad333:1.
\]
For example, the pattern 122 means that one p-value is in the first bin and
two are in the second. The number after the colon is the value of \(R\).
The corresponding multiplicities are 1,3,1,3,1. It follows that
\begin{equation}
\Expect_I R=0,\qquad \Expect_I R^2=12t^3=4\alpha^3/9.                    \label{tri:eq:3-3}
\end{equation}
Checking all 64 bin patterns with exact rational arithmetic verifies
\eqref{tri:eq:2-7}, \eqref{tri:eq:3-2}, \eqref{tri:eq:3-3} and the boundary jumps used in Section~\ref{tri:s5}.

\subsubsection{Hermite coefficients and the tail bound}

\label{tri:s3-2}

We now expand the probability in powers of the correlations. Hermite
polynomials give the coefficients because they are orthogonal under the
independent Gaussian law. Let \(\operatorname{He}_j\) be the
probabilists' Hermite polynomials. We index the covariance expansion by
the powers of the three correlations. Each
correlation contributes its power at both of its endpoints. Thus powers
\((k,l,n)\), called edge counts, give Hermite degrees
\((k+l,k+n,l+n)\) at the three coordinates, called vertex degrees.
Since \(H\) is unchanged by any coordinate sign change, only triples with
all three edge counts of the same parity remain.
The coefficient of \(a^k b^l c^n\) is
\begin{equation}
C_{kln}=\frac{\Expect_I[H\operatorname{He}_{k+l}(Z_1)
                     \operatorname{He}_{k+n}(Z_2)
                     \operatorname{He}_{l+n}(Z_3)]}{k!l!n!}.   \label{tri:eq:3-4}
\end{equation}
This follows by expanding the Gaussian generating function
\(\exp(au_1u_2+bu_1u_3+cu_2u_3)\). In three dimensions the vertex degrees
uniquely determine the edge counts. Orthogonality of the Hermite polynomials
gives the following weights in the squared \(L^2(N(0,I))\) norm of the likelihood-ratio expansion used below:
\begin{equation}
W_{kln}=\binom{k+l}{k}\binom{k+n}{k}\binom{l+n}{l}.             \label{tri:eq:3-5}
\end{equation}
For positive even \(p\), put
\[
 A_p=2\phi(A)\operatorname{He}_{p-1}(A),\quad
 B_p=2\phi(B)\operatorname{He}_{p-1}(B),\quad
 C_p=2\phi(C)\operatorname{He}_{p-1}(C).
\]
The four bin integrals are \((A_p,B_p-A_p,C_p-B_p,-C_p)\). Degree zero uses
\((t,t,t,1-\alpha)\). Thus every coefficient can be calculated by summing over the 64 possible
bin configurations.
For three positive vertex degrees \(p,q,s\), the numerator of \(C_{kln}\) can also be written
\begin{equation}
\begin{split}
 &C_pC_qC_s-(A_pC_qC_s+C_pA_qC_s+C_pC_qA_s)\\
 &\quad-(B_pB_qC_s+B_pC_qB_s+C_pB_qB_s)\\
 &\quad+(A_pB_qC_s+A_pC_qB_s+B_pA_qC_s+B_pC_qA_s
             +C_pA_qB_s+C_pB_qA_s).
 \end{split}                                                     \label{tri:eq:3-6}
\end{equation}
When only one correlation appears, raised to degree \(p\), the numerator is
\begin{equation}
(2\alpha/3)C_p(C_p-A_p-B_p)+(1-\alpha)B_p(B_p-2A_p).            \label{tri:eq:3-7}
\end{equation}
The quadratic
coefficient is
\(C_{200}=-\ell\), where, with \(A_2,B_2,C_2\) as above,
\begin{equation}
\ell=(1-\alpha)B_2(2A_2-B_2)/2+
                 \alpha C_2(A_2+B_2-C_2)/3.                  \label{tri:eq:3-8}
\end{equation}

After computing finitely many coefficients, we need to bound the omitted
terms. We do this through the square-integrable likelihood ratio, which
converts expectations under independence to expectations under the target
covariance. For \(r=7/16<1/2\), \(\|\Sigma-I\|_{\mathrm{op}}\le2r<1\), so the Gaussian
likelihood ratio \(L_\Sigma=dN(0,\Sigma)/dN(0,I)\) belongs to \(L^2(N(0,I))\).
Because the event uses absolute values, we can average this likelihood
ratio over all coordinate sign changes:
\[
 \overline L_\Sigma(z)=2^{-3}\sum_{s\in\{-1,1\}^3}
              L_\Sigma(s_1z_1,s_2z_2,s_3z_3).
\]
Every indicator in the pairwise decomposition is coordinatewise even,
so its expectation against either likelihood is the same. The
parity-restricted Hermite series below is the series for
\(\overline L_\Sigma\). Averaging the Gaussian integrals over coordinate
sign changes gives
\begin{equation}
M(r)=\frac14\left\{[(1-3r^2)^2-4r^6]^{-1/2}
                   +3[(1+r^2)^2-4r^6]^{-1/2}\right\}.          \label{tri:eq:3-9}
\end{equation}
To derive this expression, integrate the product of two Gaussian likelihood
ratios, whose covariance perturbations are \(E,G\). The result is \(\det(I-EG)^{-1/2}\).
Averaging over coordinate signs for the equicorrelation perturbation with
edges \(r\) yields two determinant types and gives \eqref{tri:eq:3-9}.
Both determinants are positive because
\(2r<1\).

Define
\begin{equation}
N_j=\sum_{k+l+n=j,\ k,l,n\text{ same parity}}W_{kln},\qquad
 B_D=M(r)-1-\sum_{j=2}^D N_jr^{2j}.                            \label{tri:eq:3-10}
\end{equation}
Here \(M(r)=\sum_{j\ge0}N_jr^{2j}\) is the squared \(L^2(N(0,I))\) norm of
\(\overline L_\Sigma\) at the equicorrelation matrix with edges \(r\). Since
\(N_0=1\) and \(N_1=0\), \(B_D=\sum_{j>D}N_jr^{2j}\) is the tail of this
series; it is unrelated to the bin integrals \(B_p\).
With \(S=a^2+b^2+c^2\), each monomial in the omitted terms has absolute value at most
\(Sr^{j-2}\). This follows by bounding the product of two factors by \(S\)
and each remaining factor by the radius. Orthogonality
therefore bounds the norm of the omitted coordinatewise-even likelihood terms by
\((S/r^2)\sqrt{B_D}\).
We apply this bound to \(R\) in \eqref{tri:eq:3-2}. For each pair, the
pair-deficit indicator is the indicator that both p-values lie in the
first bin minus the indicator that both lie in the second bin. It equals
\(\mathbf1\{P_i\le t\}+\mathbf1\{P_j\le t\}-H_{\{i,j\}}\), and its expectation
under \(\Sigma\) is the pair deficit \(D_s(t)\) in \eqref{tri:eq:3-2}, with
\(s=|\rho_{ij}|\). Since each single-coordinate indicator has expectation
\(t\) under every \(\Sigma\), \(\mathbb E_\Sigma H_{\{i,j\}}\) depends on
\(\Sigma\) only through the expectation of this indicator. The indicator has
squared norm \(2t^2\) under independence. Beyond even degree \(D\), the tail of the pair's even-projected
(sign-averaged) bivariate likelihood ratio has norm \(s^{D+2}/\sqrt{1-s^4}\). Combining the two bounds gives
\begin{equation}
|F-P_D|\le S\left\{
 \frac{2\alpha^{3/2}}{3r^2}\sqrt{B_D}
 +\frac{\sqrt2\alpha}{3}\frac{r^D}{\sqrt{1-r^4}}\right\}
 =:S\varepsilon_D.                                           \label{tri:eq:3-11}
\end{equation}
Here \(P_D\) is the degree-\(D\) expansion of the full rejection probability.
The finite expansions agree term by term by the indicator identity
\eqref{tri:eq:3-2}. The error bound retains the squared-correlation
factor, so it can be compared with the negative quadratic term even
arbitrarily close to independence.

It remains to bound the computed terms above the quadratic order. We
keep negative contributions when useful and bound the positive ones by
powers of the largest allowed correlation. For even \(j\), define
\[
 A_j^*=C_{j00}+\frac13\sum_{\text{mixed},\ k+l+n=j}\max(C_{kln},0),
\]
and for odd \(j\), define
\(A_j^*=\frac13\sum_{k+l+n=j}|C_{kln}|\), always imposing the parity condition.
For even degree, all edge powers are even, so their monomials are
nonnegative. Dropping the negative coefficients of mixed triples, those
with at least two positive edge counts, therefore gives an upper bound
while retaining the sign of the pure-edge coefficient. Summing the weighted
arithmetic--geometric mean inequality
\(|a^kb^lc^n|\le(k|a|^j+l|b|^j+n|c|^j)/j\) over the permutations of
\((k,l,n)\), which share the coefficient \(C_{kln}\), gives the factor
\(1/3\) and bounds each full degree by
\(A_j^*(|a|^j+|b|^j+|c|^j)\). We may also discard a whole degree when this coefficient
is negative. Hence
\begin{equation}
F\le\alpha-S\left\{\ell-
      \sum_{j=3}^D\max(A_j^*,0)r^{j-2}-\varepsilon_D\right\}.  \label{tri:eq:3-12}
\end{equation}

\subsubsection{Finite coefficient bounds}

\label{tri:s3-3}

The finite calculation uses \(D=20\) and all 450 admissible nonconstant
coefficients. Each coefficient is a finite sum, over the 64 bin
configurations, of products of the explicit bin integrals above, so it can
be enclosed with high-precision integer arithmetic that rounds outward. The resulting bounds are
\begin{equation}
\begin{array}{ll}
 \ell>0.00366203765920,&
 \sum_{j=3}^{20}\max(A_j^*,0)r^{j-2}<0.00034671250901,\\
 \varepsilon_{20}<0.00095626812051,&
 \ell-\text{finite bound}-\varepsilon_{20}>0.00235905702969>1/500.
 \end{array}                                                   \label{tri:eq:3-13}
\end{equation}
The displayed bounds in \eqref{tri:eq:3-13} are rounded outward. The last
bound is rounded from the certified lower bound for the margin; subtracting
the three displayed bounds instead gives \(0.00235905702968\), which also
exceeds \(1/500\).
Equations \eqref{tri:eq:3-12}--\eqref{tri:eq:3-13} prove Lemma~\ref{tri:cube}.

\subsection{A neighborhood of perfect pairs}

\label{tri:s4}

\begin{trilemma}[Perfect-pair neighborhood]\label{tri:pair} At \(\alpha=1/20\), every positive semidefinite trivariate correlation matrix satisfying
\begin{equation}
\max_{i<j}|\Sigma_{ij}|\ge2047/2048                             \label{tri:eq:4-1}
\end{equation}
satisfies \(F\le\alpha\). Equality occurs only at coordinate-sign changes and permutations of
\(J_3\) and \(\operatorname{diag}(J_2,1)\).
\end{trilemma}

We first treat vectors close to full duplication. For the remaining cases,
we contract a nearly duplicated pair and compare the rejections gained and
lost by splitting it. The correlation between the pair's center and the
third coordinate determines which bound we use. These comparisons also
identify the singular equality configurations.

\subsubsection{Near full duplication}

\label{tri:s4-1}

After sign alignment let \(\delta=\max_{i<j}(1-\Sigma_{ij})\),
\(\lambda=C^2\delta\). We first establish the bound
\begin{equation}
C\ge8/5,\quad0<\lambda\le1/300
 \quad\Longrightarrow\quad
 F\le\alpha\{1-(7/100)\sqrt\lambda\}.                          \label{tri:eq:4-2}
\end{equation}
At \(\delta=0\), full duplication gives \(F=\alpha\), so the bound also
holds at \(\lambda=0\), with equality. The condition on \(C\) holds at \(\alpha=1/20\) by \eqref{tri:eq:2-3}; we state it because Section~\ref{tri:level-interval} reuses \eqref{tri:eq:4-2} at nearby levels.

The main decrease occurs because a pair that is not exactly
duplicated does not always cross the final threshold together. We must
subtract the possible rejections at earlier ranks. Choose a pair attaining
\(\delta\), with \(\rho=1-\delta\) and \(\sigma^2=\delta(2-\delta)\).
A Simes rejection either has all three p-values at most the nominal level
or occurs with fewer than three such p-values. Let \(N_\alpha\) be the
number of p-values at most \(\alpha\); we call a rejection with
\(N_\alpha<3\) a partial rejection. Therefore,
\begin{equation}
F\le J_\rho(\alpha,\alpha)+\Prob(H=1,N_\alpha<3).                \label{tri:eq:4-3}
\end{equation}
Under the hypotheses of \eqref{tri:eq:4-2}, \(\delta\le1/768\). In the positive integral
\[
 \alpha-J_\rho(\alpha,\alpha)
 =2\int_C^\infty\phi(x)\Prob(-C<\rho x+\sigma N<C)\,dx,
\]
retain \(C\le x\le C+\sqrt\delta/2\). The upper standardized endpoint is at
least \(-1/2\); the lower is less than \(-5\). Thus the conditional probability
exceeds \(1/4\), using \(\Phi(-1/2)>3/10\) and \(\Phi(-5)\le1/25\).
Mills' inequality \(\alpha\le2\phi(C)/C\) gives
\begin{equation}
\frac{\alpha-J_\rho(\alpha,\alpha)}\alpha
 \ge\frac{\sqrt\lambda}{8}e^{-\sqrt\lambda/2-\delta/8}
 >\frac3{25}\sqrt\lambda.                                    \label{tri:eq:4-4}
\end{equation}
The last inequality follows because the exponent is less than \(1/25\)
and \(e^{-x}>1-x\).

For a partial rejection there is an ordered pair with \(|Z_i|\ge B\) and
\(|Z_j|<C\). With \(g=B-C\), conditional monotonicity and a union over six
ordered pairs and two signs give
\begin{equation}
\Prob(H=1,N_\alpha<3)
 \le4\alpha\overline\Phi\left(\frac{g-\delta B}{\sqrt{\delta(2-\delta)}}\right).
                                                                    \label{tri:eq:4-5}
\end{equation}
To justify the bound for any pair with splitting \(h=1-\Sigma_{ij}\le\delta\), note that the endpoint
\((C-(1-h)B)/\sqrt{h(2-h)}\) increases with \(h\); its derivative has the sign
of \(B-(1-h)C>0\). Zero splitting yields no mismatch.

Let \(L=\log(3/2)\) and \(y=Cg\). The normal hazard bounds
\(x\le\phi(x)/\overline\Phi(x)\le x+x^{-1}\) imply
\[
 g\le L/C,\qquad L\le(1+C^{-2})\{y+y^2/(2C^2)\}.
\]
Using \(81/200<L<41/100\) and \(C\ge8/5\), substitution shows that
\(y\le27/100\) would contradict the latter inequality. Also
\(\delta CB<7\lambda/6\), because \(B/C<1+1/(5C)\le9/8\) by
\eqref{tri:eq:2-3}. Hence
\[
 C(g-\delta B)>27/100-7/1800=479/1800=:a_0.
\]
Mills' upper bound then gives
\begin{equation}
\frac{\Prob(H=1,N_\alpha<3)}{\alpha\sqrt\lambda}
 \le\frac4{a_0\sqrt\pi}e^{-75a_0^2}<1/20.                      \label{tri:eq:4-6}
\end{equation}
The scalar bound is less than 0.041862556359. Subtract \eqref{tri:eq:4-6} from \eqref{tri:eq:4-4} to obtain \eqref{tri:eq:4-2}. At our level, the condition
\(\delta\le1/2048\) is sufficient for this bound, as follows from
\(C<13/5\) and \(169/51200<1/300\).

\subsubsection{Rejections gained by splitting a pair}

\label{tri:s4-2}

We separate the pair's common movement from the difference between its
two members. Contracting the pair removes that difference and retains two
copies of the common standardized variable. Align a pair and write
\begin{equation}
Z_1=\sqrt{1-q^2}\,U+qV,\quad Z_2=\sqrt{1-q^2}\,U-qV,\qquad
 Z_3=rU+tV+eE,\quad r^2+t^2+e^2=1,                            \label{tri:eq:4-7}
\end{equation}
where \(U,V,E\) are independent standard normals. Change the sign of \(Z_3\) if necessary to ensure
\(r\ge0\); \(t,e\) have no sign restriction. Here \((Z_1,Z_2)\) is a pair
attaining \(\max_{i<j}|\Sigma_{ij}|\), aligned so that \(\Sigma_{12}\ge0\).
We call \(U\) the center of the pair and \(r\), the correlation of \(U\)
with \(Z_3\), the center correlation; the coefficients \(t,e\) are
unrelated to \(t=\alpha/3\) in Section~\ref{tri:s3} and to Euler's number. This representation also allows singular cases and zero coefficients. If the aligned pair correlation is at least \(2047/2048\), then
\begin{equation}
q=\sqrt{(1-\Sigma_{12})/2}\le1/64.                            \label{tri:eq:4-8}
\end{equation}
For \(q=0\), use \eqref{tri:eq:2-4} directly. For \(q>0\), let \(H^c\) replace the pair
by \((U,U)\), keeping three physical tests. Since \(U\) and \(Z_3\) have
correlation \(r\), \eqref{tri:eq:2-4} gives
\(\mathbb E H^c=F_{(2,1)}(r,\alpha)=\alpha-D_r\). Hence
\(F=\alpha-D_r-\Prob(H=0,H^c=1)+\Prob(H=1,H^c=0)\), where the loss
probability \(\Prob(H=0,H^c=1)\) counts rejections lost by splitting the
pair and the last term counts those gained.

If \(|qV|<\min(A-B,B-C,C)\), then \(H=1\) implies \(H^c=1\). We verify this
implication at each rank. A rank-one crossing by a split coordinate forces
\(|U|>B\), and hence a rank-two crossing by the duplicated pair. A rank-two
crossing involving a split coordinate and the third coordinate forces
\(|U|>C\), giving a full-rank crossing after contraction. If both split
coordinates cross rank two, equal signs allow averaging, whereas opposite
signs require \(|qV|\ge B\). The full-rank argument is the same with \(C\).
A rank-one crossing by the third coordinate is unchanged. It follows that
\begin{equation}
G:=\Prob(H=1,H^c=0)\le2\overline\Phi(13/(100q)).                 \label{tri:eq:4-9}
\end{equation}
The bound depends on the half-difference between the two coordinates,
\(qV\). It controls rejections present before contraction but absent
afterwards.
By \eqref{tri:eq:4-8}, \((13/(100q))^2/2\ge34.6112>34\). Mills' bound, \(e>2\), and
\(2/(x\sqrt{2\pi})<1\) for \(x\ge8.32\) give
\begin{equation}
G<2^{-34}<1/480000000=\alpha/24000000.                         \label{tri:eq:4-10}
\end{equation}

\subsubsection{Small center correlation}

\label{tri:s4-3}

For \(0\le r\le1/10\), we construct an event on which only the contracted
vector rejects. This gives a lower bound on the loss probability. Consider
\begin{equation}
B<U<B+q/4,\qquad -1<V<-1/2,\qquad |E|<3/10.                  \label{tri:eq:4-11}
\end{equation}
The original pair is positive, has minimum below \(B\), and maximum below
\(B+5q/4<A\). The third coordinate obeys
\[
 |Z_3|\le(1/10)(2.81)+1+3/10=1.581<C.
\]
Thus the contracted rule rejects at rank two, but the original rule does
not reject. On the displayed interval for the center, the normal density
divided by its value at the left endpoint is greater than \(49/50\) (see the end of this subsection), and
\(\phi(B)\ge B\alpha/3\ge3\alpha/5\). Also
\(\Prob(-1<V<-1/2)>1/10\) and \(\Prob(|E|<3/10)>1/5\). The loss probability is
therefore greater than
\begin{equation}
(q/4)(147\alpha/250)(1/10)(1/5)>q/40000.                      \label{tri:eq:4-12}
\end{equation}
Since \(q\le1/64<1/50\), Mills' bound in \eqref{tri:eq:4-9} gives
\[
 G/q<(100/13)e^{-169/8}<1/200000.
\]
Consequently
\begin{equation}
F\le F_{(2,1)}(r,\alpha)-q/50000<\alpha.                      \label{tri:eq:4-13}
\end{equation}
The density ratio is at least \(\exp(-Bq/4-q^2/32)\); for \(q\le1/40\) and
\(B<14/5\) this exceeds \(\exp[-(14/5)/160-1/51200]>49/50\). We use these bounds for their simple rational form.

\subsubsection{Intermediate center correlation}

\label{tri:s4-4}

Away from both independence and full duplication, the contracted pair
already has a strictly positive deficit. It is enough to show that this
deficit exceeds the possible gain from splitting; no new loss event is
needed. For \(1/20\le r\le19/20\),
\begin{equation}
D_r\ge\alpha/24000000.                                       \label{tri:eq:4-14}
\end{equation}
We prove this inequality by integrating over a rectangle on which the
integrand has a uniform lower bound. In \eqref{tri:eq:2-5}, retain
\(z\in[\alpha/6,\alpha/3]\). Since \(f(z)=q_r(c_z)\), the integrand there is
\(f(z)-f(z+2\alpha/3)=\int_{c_{z+2\alpha/3}}^{c_z}q_r'(y)\,dy\), with \(q_r'\)
from \eqref{tri:eq:2-6}. In this inner integral, the common interval
\[
 c_{5\alpha/6}\le y\le c_{3\alpha/4}
\]
is contained in \([c_{z+2\alpha/3},c_z]\). Its length is greater than
\(1/40\): \(\log(10/9)>1/10\), while the upper normal hazard on this interval
is below 4. To check the latter, the lower hazard bound gives
\(y\le C+\log(4/3)/C<C+3/16\), so
\(y+1/y<13/5+3/16+5/8<4\).

On this rectangle, \(-1<(C-ry)/\sigma<3\). The upper bound follows from
\(y\ge C\) and \(C(1-r)/\sigma\le C<3\); for the lower bound use
\(ry-C\le(19/20)(3/16)\) and \(\sigma>3/10\). Also
\(r/\sigma\ge1/20\), \(2rCy/\sigma^2>1/4\),
\(\phi(3)>1/1000\), and \(1-e^{-1/4}>1/5\). The derivative in \eqref{tri:eq:2-6} is therefore
greater than \(1/100000\). The retained area is at least
\((\alpha/6)(1/40)\), proving \eqref{tri:eq:4-14}.
For \(1/10\le r\le19/20\), \eqref{tri:eq:4-10}, \eqref{tri:eq:4-14}, and
\(F\le F_{(2,1)}+G\) give strict validity.

\subsubsection{Large center correlation}

\label{tri:s4-5}

For \(19/20\le r\le9999/10000\), we use the stronger bound
\begin{equation}
D_r>\alpha/4000                                             \label{tri:eq:4-15}
\end{equation}
to control the gain. Retain \(z\in[\alpha/3-\alpha/1000,\alpha/3]\) in \eqref{tri:eq:2-5}. By \eqref{tri:eq:2-3},
\[
 (19/20)(2393/1000)-49/25>313/1000>
       \sqrt{1-(19/20)^2}.
\]
The ratio \((rA-C)/\sqrt{1-r^2}\) increases in \(r\), its derivative having
numerator \(A-rC>0\). Thus \(f(z)>\Phi(1)>4/5\).

For \(w=z+2\alpha/3\in[999\alpha/1000,\alpha]\), let
\(y=c_{999\alpha/1000}\). The lower hazard bound gives
\[
 y-C\le\log(1000/999)/C<1/(999(19/10)).
\]
Since \(\sigma>7/500\),
\[
 (ry-C)/\sigma<1/25.
\]
The positive conditional tail is below
\(\Phi(1/25)<13/25\). The negative conditional tail is below
\(\overline\Phi(C/\sigma)<\overline\Phi(6)<1/100\), because
\(\sigma<313/1000\). Hence \(f(w)<53/100\), and
\(f(z)-f(w)>1/4\). Multiplying by the retained width \(\alpha/1000\) proves
\eqref{tri:eq:4-15}. The gain \eqref{tri:eq:4-10} is smaller than this deficit.

\subsubsection{Center correlation near one}

\label{tri:s4-6}

When the center and third coordinate are very highly correlated, the
small pair splitting forces all three original coordinates to be close
to full duplication. We verify this and use the earlier full-duplication
bound. If \(r\ge9999/10000\), \eqref{tri:eq:4-8} gives
\[
 \sqrt{1-q^2}>8190/8191,\qquad |t|\le\sqrt{1-r^2}<1/70.
\]
Both aligned correlations to the third coordinate are at least
\begin{equation}
\frac{9999}{10000}\frac{8190}{8191}-\frac1{64\cdot70}
 =\frac{916983497}{917392000}>\frac{2047}{2048}.                 \label{tri:eq:4-16}
\end{equation}
The original pair already satisfies the same lower bound. Thus
\(\delta\le1/2048\), and \eqref{tri:eq:4-2} gives strict validity unless \(\delta=0\).

These four closed ranges cover all center correlations. For \(q=0\), \eqref{tri:eq:2-4} gives strictness
unless \(r=0\) or 1. For \(q>0\), \eqref{tri:eq:4-13}, \eqref{tri:eq:4-14}--\eqref{tri:eq:4-15}, or \eqref{tri:eq:4-2} gives
strictness. This proves Lemma~\ref{tri:pair} and its equality statement.

Each scalar inequality used in this section compares explicit rational
numbers or values of elementary functions and of \(\Phi\) at rational points. Together
with the preceding
analytic argument, these inequalities establish the bound throughout the
perfect-pair neighborhood.

\subsection{Covariance-path identity}

\label{tri:s5}

For the remaining covariance matrices, we express the difference from the
independence probability as an integral along a covariance path. We first
derive the covariance derivatives by tracking changes in the rejection
indicator at its threshold boundaries. Both positive and negative
contributions are retained.

For \(\partial_cF\), condition \(Z_1\) on \(Z_2=x,Z_3=y\). The following forms of its
conditional mean and variance are convenient for numerical evaluation:
\begin{equation}
\mu=\frac{(a+b)(x+y)}{2(1+c)}+
      \frac{(a-b)(x-y)}{2(1-c)},\qquad
 v=1-\frac{(a+b)^2}{2(1+c)}-\frac{(a-b)^2}{2(1-c)}.             \label{tri:eq:5-1}
\end{equation}
These equal the usual regression expressions; in particular
\(v=\Delta/(1-c^2)\). The bivariate density is
\begin{equation}
\phi_2(x,y;c)=\frac{\exp[-(x+y)^2/(4(1+c))-(x-y)^2/(4(1-c))]}
                    {2\pi\sqrt{(1+c)(1-c)}}.                 \label{tri:eq:5-2}
\end{equation}
These expressions avoid the cancellation that would arise from expanding
the determinants in \eqref{tri:eq:5-1}--\eqref{tri:eq:5-2}.

The derivative depends on how the rejection decision changes when two
coordinates cross their bin boundaries together. Averaging those changes
over the remaining coordinate produces the following matrix. Let
\(Q_j(\mu,v)=\Prob(|N(\mu,v)|\ge d_j)\), with
\((d_1,d_2,d_3)=(A,B,C)\). Define
\begin{equation}
J(\mu,v)=\begin{pmatrix}
 0&-(1-Q_3)&-(Q_3-Q_2)\\
 -(1-Q_3)&1-Q_3&-(Q_2-Q_1)\\
 -(Q_3-Q_2)&-(Q_2-Q_1)&Q_3-Q_1
 \end{pmatrix}.                                               \label{tri:eq:5-3}
\end{equation}
For positive-definite \(\Sigma\),
\begin{equation}
\partial_cF=2\sum_{u,w=1}^3\sum_{\eta\in\{-1,1\}}
   \eta\,\phi_2(d_u,\eta d_w;c)
            J_{uw}(\mu(d_u,\eta d_w),v).                      \label{tri:eq:5-4}
\end{equation}
We evaluate the conditional mean separately in every term. The other
covariance derivatives follow by permuting coordinates.

\paragraph{Derivation of the covariance derivative.}
Varying an off-diagonal covariance changes both symmetric
entries, and direct differentiation of the Gaussian density gives
Price's identity \(\partial_cf_\Sigma=\partial_{z_2}\partial_{z_3}f_\Sigma\).
Integrate this identity against the step function \(H\); changes at its
bin boundaries determine the derivative. For this calculation, we renumber
the p-value bins starting from zero.
Using the bin immediately below each boundary to index that boundary, the
mixed jump at thresholds \(u,w\) for a remaining bin \(k\) is
\[
 H(k,u,w)-H(k,u+1,w)-H(k,u,w+1)+H(k,u+1,w+1).
\]
It is zero in remaining bin 0. The following list numbers the thresholds
from one to three, as in \eqref{tri:eq:5-3} and \eqref{tri:eq:5-4}:
threshold \(j\) lies at \(|z|=d_j\) and separates bins \(j-1\) and \(j\),
so the indices \(u\) and \(w\) in the formula above correspond to thresholds \(u+1\) and \(w+1\). In remaining bin 1 the nonzero jumps are
\((2,3),(3,2):-1\) and \((3,3):1\); in bin 2 they are
\((1,3),(3,1):-1\), \((3,3):1\); in bin 3 they are
\((1,2),(2,1):-1\), \((2,2):1\).
Conditional averaging gives \eqref{tri:eq:5-3}. The negative and positive boundary orientations
give \(\eta\), and global sign symmetry gives the factor 2. Thus distributional
integration by parts establishes \eqref{tri:eq:5-4}.

We now integrate the covariance derivatives along the straight path
from independence to the target matrix. Set \(\rho=(a,b,c)\). Since
\(F_I=\alpha\),
\begin{equation}
F_{\Sigma(\rho)}-\alpha
 =\int_0^1 h_\rho(t)\,dt,\qquad
 h_\rho(t)=\rho\cdot\nabla F_{\Sigma(t\rho)}.                   \label{tri:eq:5-5}
\end{equation}
The identity holds for every positive semidefinite endpoint, including
exact copies. To
justify a possibly singular endpoint, we bound the remaining integral
near the end of the path. For \(t<1\), the path covariance is positive
definite, and
\[
 \phi_2(x,y;t\rho_i)\le\frac1{2\pi\sqrt{1-t^2}}.
\]
There are at most eight nonzero entries in \eqref{tri:eq:5-3}, each
bounded by one, and two signs in \eqref{tri:eq:5-4}. Consequently
\begin{equation}\label{tri:path-terminal}
 |h_\rho(t)|\le\frac{16\|\rho\|_1}{\pi\sqrt{1-t^2}},\qquad
 \left|\int_{1-\delta}^1h_\rho(t)\,dt\right|
 \le\frac{16\|\rho\|_1}{\pi}\arccos(1-\delta),
 \quad 0\le\delta\le1.
\end{equation}
The first bound is integrable. We can therefore pass to the endpoint in
the integral and use covariance continuity to obtain \eqref{tri:eq:5-5}.
The value assigned to a conditional probability at the single terminal
point does not affect the integral. The numerical certificate applies
\eqref{tri:eq:5-5} only at the reference matrices of its covariance cells
(Section~\ref{tri:s7-1}). None of them has a perfect pair, so the bivariate
denominators remain positive. The terminal estimate above also justifies the path
identity at perfect-pair endpoints, which are treated analytically.

\subsection{Numerical integration with rigorous error bounds}

\label{tri:s6}

\subsubsection{Elementary-function bounds}

\label{tri:s6-1}

The numerical argument bounds the exact path integral
\eqref{tri:eq:5-5}, with every evaluation and integration error
included. Each computed quantity is an \emph{enclosure}: an interval
guaranteed to contain the exact value. We start with the two functions
in the integrand, the exponential function and the normal tail
\(\overline\Phi\). Standard floating-point library routines for these
functions do not guarantee error bounds, so we use tables instead. The tables give
lower and upper bounds for \(e^{-x}\) and \(\overline\Phi(x)\) at the
grid points \(x=ih_0\), where
\[
 h_0=1/4096.
\]
Between two neighboring grid points, the straight line through the
function values lies above the function, because both functions are
convex. The usual error bound for linear interpolation limits how far
above it can lie: if \(f\) is either function and \(L\) is its linear
interpolant on an interval of length \(h_0\), then
\begin{equation}
L(x)-\frac{h_0^2}{8}\sup f''\le f(x)\le L(x).                 \label{tri:eq:6-1}
\end{equation}
On the interval starting at \(ih_0\), \(\sup f''\le e^{-ih_0}\) for the
exponential and \(\sup f''\le(ih_0+h_0)\phi(ih_0)\) for the normal tail.
The tables stop at \(x=80\) and \(x=12\), respectively. Beyond them we
use
\begin{equation}
x\ge80:\ e^{-x}\le2^{-100},\qquad
 x\ge12:\ \overline\Phi(x)\le2^{-100},                         \label{tri:eq:6-2}
\end{equation}
which follow from \(e>8/3\), \((8/3)^{72}>2^{100}\) and Mills'
inequality \(\overline\Phi(12)<e^{-72}\).

Every table entry is verified in high-precision integer arithmetic
that rounds outward, so the tables can be trusted however they were
computed. The exponential values
follow from repeated multiplication by an enclosure of \(e^{-h_0}\).
The normal-tail values follow from \(\overline\Phi(0)=1/2\) by
integrating the normal density over one grid step at a time, using a
Taylor expansion with an explicit bound on the omitted terms. The
constants \(\pi\) and \(1/\sqrt{2\pi}\) and the three quantile brackets
of Section~\ref{tri:s2} are checked in the same way.

\subsubsection{Path derivatives and integration errors}

\label{tri:s6-2}

All other numerical calculations use interval arithmetic. Each quantity
is represented by an interval with floating-point endpoints, and each
arithmetic operation produces an interval that contains every exact
result obtainable from its inputs. A floating-point operation returns
the representable number nearest to the exact result. Moving one step
further out, to the adjacent representable number, therefore gives a
bound on the correct side. For example, the product of two intervals
uses the smallest and the largest of the four endpoint products, each
moved one step outward. Intervals for squares and variances are
intersected with their known ranges. Section~\ref{tri:reproducibility}
states the assumptions about the computer arithmetic.

We bound the integral in \eqref{tri:eq:5-5} by dividing \([0,1]\) into
short segments, whose endpoints are chosen so that every width and
midpoint is computed exactly. Let \(K\) enclose \(h_\rho\) on a segment
of width \(w\) and midpoint \(m\). Then \(wK\) encloses the integral
over the segment. When a finite bound \(M\ge\sup|h_\rho''|\) on the
segment is available, the midpoint rule gives a second enclosure,
\begin{equation}
w h_\rho(m)+[-Mw^3/24,Mw^3/24].                              \label{tri:eq:6-3}
\end{equation}
Both enclosures contain the exact integral, so we use their
intersection. The segments are chosen adaptively, but each segment
bound is valid however the segments were chosen. Adding the segment
bounds, again with outward rounding, encloses the whole integral.

The bound \(M\) comes from differentiating the analytic formulas
\eqref{tri:eq:5-1}--\eqref{tri:eq:5-4} twice along the path, again in
interval arithmetic: each intermediate quantity carries enclosures of
its value and of its first two derivatives, which are combined by the
sum, product and chain rules. The tables of Section~\ref{tri:s6-1} are
used only for function values. Their interpolants are never
differentiated, because the derivative of an interpolant need not bound
the derivative of the function.

A small conditional variance can make these derivative bounds large.
The following lower bound controls it along the path through its value
at the reference covariance \(\Sigma\), where the path ends. Write
\(v_i(\Sigma')\) for the conditional variance of the \(i\)th coordinate given
the other two under \(\Sigma'\). Along the path,
\begin{equation}
v_i((1-t)I+t\Sigma)\ge(1-t)+t v_i(\Sigma).                    \label{tri:eq:6-4}
\end{equation}
Indeed, a conditional variance is the smallest residual variance of a
linear predictor. Each residual variance is linear in the covariance
matrix, so their minimum is concave, which gives \eqref{tri:eq:6-4},
including at singular endpoints. When this lower bound is zero, as on
the last segment for a singular reference, we use only the enclosure
\(wK\).

\subsection{Covariance cells and classification by positive semidefiniteness}

\label{tri:s7}

\subsubsection{Uniform error within a covariance cell}

\label{tri:s7-1}

A covariance cell is a rectangular box in the three correlation coordinates.
Its feasible part consists of points that give positive semidefinite
correlation matrices. We bound the rejection probability throughout that
part using a bound at one reference matrix and bounds on the covariance
derivatives.
Let \(\mathcal K_3=\{\rho\in[-1,1]^3:\Sigma(\rho)\succeq0\}\) be the feasible set of correlation coordinates. Let \(\mathcal B=\prod_i[\ell_i,r_i]\) be a box in these coordinates, let
\(\rho_*\in\mathcal B\) be a feasible reference, and let \(G_i\) enclose the
\(i\)th component of \(\nabla F\) for every positive definite feasible
covariance in the box, with the same bounds valid as these matrices
approach a singular matrix. Thus each \(G_i\) is a uniform derivative
bound. We combine these intervals with the possible changes in the
correlation coordinates to bound the change in rejection probability.
Write \(F_\rho=F_{\Sigma(\rho)}\). If \(U_*\) is a certified upper bound for \(F_{\rho_*}-\alpha\), then
\begin{equation}
F_\rho-\alpha\le U_*+
 \sum_{i=1}^3\sup\bigl([\ell_i-\rho_{*,i},r_i-\rho_{*,i}]G_i\bigr)
 \qquad(\rho\in\mathcal B\cap\mathcal K_3).                    \label{tri:eq:7-1}
\end{equation}
We call \eqref{tri:eq:7-1} the whole-cell bound. Section~\ref{tri:s6}
supplies \(U_*\), and the sum, called the parameter error, bounds the
increase of \(F\) from the reference to any feasible point of the box.
The segment from the reference to the target stays feasible and in the box.
For the numerical boxes used below, every coordinate bound is strictly less
than one in absolute value, so neither endpoint has a perfect pair. If the endpoints are distinct,
singular or not, Lemma~\ref{tri:segments} makes the segment interior
positive definite; boundedness of \eqref{tri:eq:5-4} and endpoint continuity complete the
fundamental-theorem-of-calculus argument. For identical endpoints, the error
is zero. The bound therefore also covers singular targets.

We obtain each \(G_i\) by evaluating \eqref{tri:eq:5-4} in interval
arithmetic over the whole box, with conditional variances intersected
with \([0,1]\). This intersection removes values that no feasible
covariance matrix can produce, while allowing the box to contain
infeasible points. The terms of \eqref{tri:eq:7-1} keep the signs of
the derivative bounds. For example, if \(G_i\) contains only
nonnegative values, a reference with \(\rho_{*,i}=r_i\) makes the
\(i\)th term zero. No global sign of the derivatives is needed.

\subsubsection{Partitioning the canonical chambers}

\label{tri:s7-2}

To cover every feasible covariance, we begin with boxes covering each
closed chamber and repeatedly cut in half any box that is not yet
classified. All box endpoints are multiples of \(2^{-20}\). For each
sign we start with the following 32 root boxes:
\begin{equation}
a\in[j/32,(j+1)/32],\qquad b,d\in[0,1],\qquad j=0,\ldots,31.   \label{tri:eq:7-2}
\end{equation}
A box is cut at the exact midpoint of one coordinate. The boxes obtained
from a root box therefore form a binary tree, whose leaves, the boxes
that are never cut, partition the root box.

A box becomes a leaf once it is classified. It contains no point of its
chamber if tightening its bounds by the inequalities \(0\le d\le b\le
a\) leaves a lower bound above the corresponding upper bound, or if the
determinant test of Section~\ref{tri:s7-3} shows that it contains no
positive semidefinite matrix. It is an independence leaf if
\(a\le7/16\) throughout the box, and a perfect-pair leaf if
\(a\ge2047/2048\) throughout the box; Lemmas~\ref{tri:cube}
and~\ref{tri:pair} cover these leaves. Any other box becomes a
numerical leaf if it passes the criterion of Section~\ref{tri:s7-4},
and is cut again otherwise. A box touching \(a=1\) is cut in \(a\)
unless it is a perfect-pair leaf. For a numerical leaf, the box
\(\mathcal B\) in \eqref{tri:eq:7-1} is its tightened box, written in the
coordinates \((a,b,c)\) with \(c=\varepsilon d\). Its bounds are less than
one in absolute value, because tightening makes the upper bounds of
\(b\) and \(d\) at most that of \(a\), which is less than one.

\subsubsection{An exact corner test for positive semidefiniteness}

\label{tri:s7-3}

We need to distinguish boxes with no feasible covariance from boxes
that require a probability bound. The following determinant test either
excludes the whole target region or supplies a feasible reference within
it. Write \(A\) and \(B\) for the lower bounds of \(a\) and \(b\), and
\([D_-,D_+]\) for the range of \(d\), after the tightening of
Section~\ref{tri:s7-2}; the target region is this tightened box.
In this subsection \(A\) and \(B\) denote these corner coordinates, not the
normal quantiles of Section~\ref{tri:s2}.
For the negative chamber use the corner \((A,B,-D_-)\). Since
\(\Delta(a,b,-d)=1-a^2-b^2-d^2-2abd\) decreases in all nonnegative arguments,
a strictly negative determinant there excludes the entire target region as infeasible.

For the positive chamber use
\begin{equation}
c_* =\max(D_-,\min(D_+,AB)),\qquad \rho_*=(A,B,c_*).            \label{tri:eq:7-4}
\end{equation}
Here \(A\ge B\), \(D_-\le B\), and \(A,B\le1\).
If \(AB\le D_+\), then \(AB\le c_*\le B\). Since \(\Delta(A,B,\cdot)\) decreases on \([AB,B]\),
\begin{equation}
\Delta(A,B,c_*)\ge\Delta(A,B,B)
 =(1-A)(1+A-2B^2)\ge0.                                       \label{tri:eq:7-5}
\end{equation}
If \(AB>D_+\), then \(c_*=D_+<AB\le B\le A\). For all target points with
\(a\ge A,b\ge B,c\le D_+\), increasing \(c\) to \(D_+\) increases the
determinant, and decreasing \(a,b\) to \(A,B\) increases it. Thus
\begin{equation}
\Delta(a,b,c)\le\Delta(A,B,D_+).                              \label{tri:eq:7-6}
\end{equation}
Equations \eqref{tri:eq:7-5}--\eqref{tri:eq:7-6} prove that a strictly negative corner determinant
excludes only infeasible points. If it is nonnegative, \eqref{tri:eq:7-4} is a feasible
reference. Thus the test classifies the full target region rather than
only the grid points.

Any feasible point of \(\mathcal B\) may serve as the reference; the choice
affects only the size of the bound, not its validity. Every reference
is checked exactly. Its coordinates are less than one in absolute
value, so the leading \(2\times2\) principal minor is positive, and
positive semidefiniteness is equivalent to a nonnegative determinant.
The coordinates are fractions whose denominators are powers of two, so
the sign of the determinant is computed exactly with integers.

\subsubsection{The acceptance criterion}

\label{tri:s7-4}

A numerical leaf must show a strict margin throughout its feasible part,
after allowing for both integration error and variation in the
covariance. It is accepted only if the right side of
\eqref{tri:eq:7-1}, computed with outward rounding, is at most
\begin{equation}
-2^{-20}=-1/1048576.                                         \label{tri:eq:7-8}
\end{equation}
Otherwise the box is cut in half and both halves are classified
again.

\subsection{Finite certificate and verification}

\label{tri:s8}

\begin{triproposition}[Finite certificate for the remaining covariance cells]
\label{tri:finite-certificate}
Repeated halving of the $64$ root boxes \eqref{tri:eq:7-2} ends with
finitely many leaves, each classified as in Section~\ref{tri:s7-2}.
Every point of the two closed chambers lies in an independence leaf, a
perfect-pair leaf or a numerical leaf. The independence and
perfect-pair leaves satisfy the hypotheses \eqref{tri:eq:3-1}
and~\eqref{tri:eq:4-1}, respectively. For every numerical leaf,
with $\mathcal B$ its tightened box, a feasible reference and the
enclosures of Sections~\ref{tri:s6} and~\ref{tri:s7-1} make the right
side of the whole-cell bound \eqref{tri:eq:7-1} at most $-2^{-20}$.
\end{triproposition}

\begin{proof}[Verification of the certificate]
The certificate consists of the trees of boxes and, for each numerical
leaf, its reference, its derivative intervals \(G_i\), its integration
segments with their enclosures, and the resulting bound. Checking the
certificate has two parts: every local bound must
be valid, and the leaves must cover both chambers.

The local bounds are verified by computing them again. The function
tables are verified as in Section~\ref{tri:s6-1}. On every numerical
leaf, the derivative intervals and the segment enclosures are
recomputed with the interval arithmetic of Section~\ref{tri:s6-2}, and
each interval in the certificate must contain the recomputed one.
Since the recomputed interval contains the exact value, so does the
interval in the certificate.

The remaining checks use exact rational arithmetic, treating each
floating-point number as the exact fraction it represents. They rebuild
every tree from its root box, confirm each classification by the
criteria of Section~\ref{tri:s7}, and confirm that each reference lies
in the tightened box of its leaf and gives a positive semidefinite
matrix. They also confirm
that the segments of each numerical leaf cover \([0,1]\) without gaps
or overlaps and have exact widths and midpoints, that \(U_*\) is at least the sum of the upper bounds of
the segment enclosures, and that the right side of \eqref{tri:eq:7-1}
is at most \(-2^{-20}\). Hence the leaves cover both chambers, and
every numerical leaf satisfies the stated bound.
\end{proof}

\subsection{Completion of the theorem}\label{tri:completion}

\begin{proof}[Completion of Theorem~\ref{tri:theorem}]
We reduce an arbitrary covariance matrix to the canonical form
\eqref{tri:eq:2-8}. The verified partition places it in a root and a terminal
leaf. Since the matrix is canonical and feasible, neither the order test nor
the determinant test of Sections~\ref{tri:s7-2}--\ref{tri:s7-3} excludes it. On an independence leaf, we apply
\eqref{tri:eq:3-1}; on a perfect-pair leaf, we apply Lemma~\ref{tri:pair}; and on a
numerical leaf, Proposition~\ref{tri:finite-certificate} and the whole-cell
bound \eqref{tri:eq:7-1} give \(F_\Sigma-\alpha\le-2^{-20}\). These three cases establish the
probability inequality. The bound is strict on every numerical leaf.
Equation \eqref{tri:eq:3-1} is strict except at independence, while
Lemma~\ref{tri:pair} has precisely the stated singular equality cases.
Equations \eqref{tri:eq:2-7} and~\eqref{tri:eq:2-4}, together with full duplication, give equality
at all matrices in \eqref{tri:equalities}, completing the classification.
\end{proof}

The numerical part gives a strict bound throughout the feasible part
of every numerical cell between the two analytic neighborhoods. Both the testing level and
the neighborhood sizes are specified. The argument therefore requires
no unevaluated small-level cutoff or neighborhood radius.

\subsection{A common interval around the testing level}
\label{tri:level-interval}

Put \(\alpha_0=1/20\). Near an equality case, the gap between the
rejection probability and the nominal level tends to zero, so continuity
alone cannot show that its sign is preserved. Our bounds retain the covariance
factor that vanishes near independence and the splitting factor that
vanishes near exact copies. Elsewhere, the certified positive gap gives
an explicit interval of levels.

\begin{proof}[Proof of Corollary~\ref{main:tri-interval}]
For any positive semidefinite covariance and \(0<\alpha<\beta<1\), the rejection events
are nested. A change in the decision requires a marginal p-value to
cross one of the three thresholds. Marginal uniformity, also at
singular covariances, gives
\begin{equation}\label{tri:level-lipschitz}
 0\le F^{(3)}_\Sigma(\beta)-F^{(3)}_\Sigma(\alpha)
 \le\sum_{i=1}^3\sum_{k=1}^3\frac{k(\beta-\alpha)}3
 =6(\beta-\alpha).
\end{equation}
For \(d_\Sigma(\alpha)=\alpha-F^{(3)}_\Sigma(\alpha)\), this implies
\begin{equation}\label{tri:deficit-lipschitz}
 -5(\beta-\alpha)\le d_\Sigma(\beta)-d_\Sigma(\alpha)
 \le\beta-\alpha,\qquad
 |d_\Sigma(\beta)-d_\Sigma(\alpha)|\le5|\beta-\alpha|.
\end{equation}
We use the one-sided bounds separately on the two sides of \(\alpha_0\).
They preserve a fixed positive gap away from the equality cases. Near
those cases, where the gap can be arbitrarily small, we first need the
two more precise comparisons below.

\paragraph{The neighborhood of independence.}
Near independence, the error from changing the level must retain the
same squared-correlation factor as the existing deficit. Keep
\(r=7/16\), and put \(S=a^2+b^2+c^2\).
The even-projected likelihood introduced in Section~\ref{tri:s3-2}
has orthogonal Hermite terms indexed by edge counts \((k,l,n)\).
The vertex degrees \((k+l,k+n,l+n)\) determine these counts uniquely;
even projection keeps exactly the triples of common parity. Its only
term below total edge degree two is the constant. Bounding two factors
of every degree-\(j\) monomial by \(S\), and the rest by \(r\), gives
\[
 \|\overline L_\Sigma-1\|_2^2
 \le S^2r^{-4}\{M(r)-1\}=:K^2S^2.
\]
This is a bound for the projected likelihood. The unprojected
likelihood has degree-one terms and does not satisfy this estimate.
At \(r=7/16\), the two radicands in \eqref{tri:eq:3-9} are
\[
 \frac{642735}{4194304}>\frac{25}{169},\qquad
 \frac{5835951}{4194304}>\frac{400}{289}.
\]
Consequently \(M(r)<103/80\) and
\(K^2<94208/12005\).

Let \(D=H_\alpha-H_{\alpha_0}\). Under independence, the identity
\(\Expect_IH_\alpha=\alpha\) holds at every level, and nestedness
gives \(\Expect_ID=\alpha-\alpha_0\) and
\(\Expect_ID^2=|\alpha-\alpha_0|\). Because \(D\) is coordinatewise
even, \(L_\Sigma\) may be replaced by its even projection
\(\overline L_\Sigma\), and the level change \(\alpha-\alpha_0=\Expect_ID\)
cancels exactly:
\[
 d_\Sigma(\alpha)-d_\Sigma(\alpha_0)
 =-\Expect_I[(\overline L_\Sigma-1)D].
\]
Cauchy--Schwarz therefore bounds its absolute value by
\(KS\sqrt{|\alpha-\alpha_0|}\). The exact rational comparison
\[
 \frac{94208}{12005}\,2^{-21}
 =\frac{23}{6146560}
 <\left(\frac{39}{20000}\right)^2
\]
and the original \(S/500\) center-level bound imply, for
\(|\alpha-\alpha_0|\le2^{-21}\),
\begin{equation}\label{tri:interval-cube}
 \max(|a|,|b|,|c|)\le7/16
 \quad\Longrightarrow\quad
 F^{(3)}_\Sigma(\alpha)\le\alpha-S/20000,
\end{equation}
with strict inequality in this displayed bound when \(S>0\).
At \(S=0\), exact independence gives equality with \(\alpha\).
No coefficient computation at a nearby level is needed.

\paragraph{The neighborhood of perfect pairs.}
We keep the original covariance neighborhood
\(\max|\Sigma_{ij}|\ge2047/2048\). Write
\(A(\alpha)=c_{\alpha/3}\), \(B(\alpha)=c_{2\alpha/3}\), and
\(C(\alpha)=c_\alpha\). On the auxiliary interval
\(|\alpha-\alpha_0|\le2^{-21}\), these three quantiles move by less
than \(w=1/40000\). To prove this without assuming the desired
quantile bounds in advance, differentiate:
\[
 \left|\frac{d}{d\alpha}c_{k\alpha/3}\right|
 =\frac{k}{6\phi(c_{k\alpha/3})}<50
 \quad\hbox{when }0<c_{k\alpha/3}<5/2.
\]
Indeed, \(\phi(5/2)>7/648>1/100\), using
\(\sqrt{2\pi}<3\), \(e<3\), and
\(e^{1/8}<8/7\). At the center all three quantiles are below
\(12/5\). If a quantile first reached \(5/2\), it would have moved by at least
\(1/10\). Yet integrating the derivative bound up to that point gives
a movement of at most \(50\cdot2^{-21}<w\). This contradiction keeps all
three quantiles below the upper bound.

We now check the inputs to the original perfect-pair argument at every
nearby level, so that its gain and loss comparisons can be reused. Widen
both ends of each center bracket in Section~\ref{tri:s2} by
\(w\). Exact rational substitution preserves every inequality in
\eqref{tri:eq:2-3}; differences use an error of \(2w\).
For example, the tightest of these inequalities, \(C<49/25\), holds because
\[
 \frac{49}{25}-\left(\frac{2155003191024}{2^{40}}+w\right)
 =\frac{473110401}{42949672960000}>0.
\]
The additional level-dependent conditions are also retained:
\begin{equation}\label{tri:interval-pair-conditions}
 \frac1{40}<\alpha<\frac3{40},\qquad
 \frac{147\alpha}{50000}>\frac1{40000},\qquad
 2^{-34}<\frac{\alpha}{24000000},\qquad
 \frac{C(\alpha)^2}{2048}<\frac1{300}.
\end{equation}
For the middle two inequalities substitute
\(\alpha\ge1/20-2^{-21}\); for the last use \(C<49/25\).
We now check that these comparisons keep the factors that tend to zero
as the covariance approaches an equality case. Here \(q\) and \(r\) are the
parameters of \eqref{tri:eq:4-7}; in particular, \(r\) is now the center
correlation, not the radius \(7/16\).

For \(0<q\le1/64\) in \eqref{tri:eq:4-7}, the inequalities
\(A-B>1/5\), \(B-C>13/100\), and \(C>8/5\) retain the same
gain inclusion and give
\begin{equation}\label{tri:interval-gain}
 G\le2\barPhi\!\left(\frac{13}{100q}\right)<2^{-34},
 \qquad G/q<\frac{100}{13}e^{-169/8}<\frac1{200000}.
\end{equation}
For \(r\le1/10\), the lost event \eqref{tri:eq:4-11} still has
probability greater than \(147\alpha q/50000>q/40000\).
The original vector does not reject because the same quantile gaps
hold and the third coordinate remains below \(1.581<C\). Thus the strict
comparison retains its factor \(q\):
\[
 F^{(3)}_\Sigma(\alpha)
 <F_{(2,1)}(r,\alpha)-q/50000\le\alpha-q/50000.
\]
For \(1/10\le r\le19/20\), the rectangle in
Section~\ref{tri:s4-4} gives \(D_r\ge\alpha/24000000\).
Its length, derivative bound and p-value width use only the retained
quantile inequalities and fixed level ratios. For
\(19/20\le r\le9999/10000\), Section~\ref{tri:s4-5} likewise
gives \(D_r>\alpha/4000\). In particular, \(A>2393/1000\),
\(C<49/25\), and \(C>19/10\) preserve the positive- and
negative-tail bounds in that argument. Both deficits dominate
\eqref{tri:interval-gain} by \eqref{tri:interval-pair-conditions}.

If \(r\ge9999/10000\), the algebraic bound
\eqref{tri:eq:4-16} still gives \(\delta\le1/2048\).
Set \(\lambda=C(\alpha)^2\delta\). The conditions
\(C\ge8/5\), \(B<7C/6\), \(\lambda<1/300\), and the fixed tail ratio
\(\barPhi(B)/\barPhi(C)=2/3\) are precisely the level-dependent
inputs in \eqref{tri:eq:4-4}--\eqref{tri:eq:4-6}; the second follows
from the retained bounds \(B-C<1/5\) and \(C>8/5\). Therefore
\[
 F^{(3)}_\Sigma(\alpha)
 \le\alpha\{1-(7/100)\sqrt\lambda\}<\alpha
 \quad\hbox{if }\delta>0.
\]
At \(\delta=0\), full duplication gives equality. At \(q=0\),
the exact duplicated-pair identity holds at every level and is strict
for \(0<r<1\), with equality precisely at \(r=0,1\).
Thus singular laws and arbitrarily small positive splittings are covered
without dividing by a zero splitting parameter.

\paragraph{The remaining covariance matrices.}
Outside the two fixed original neighborhoods, the complete covariance
certificate gives \(d_\Sigma(\alpha_0)\ge m_0:=2^{-20}\).
Define
\begin{equation}\label{tri:interval-radius}
 I_0=\left[\alpha_0-\frac{m_0}{2},\,
                    \alpha_0+\frac{m_0}{6}\right]
 =\left[\frac{524283}{10485760},\,
                    \frac{1572869}{31457280}\right].
\end{equation}
This interval lies inside the auxiliary interval used for both analytic
estimates. If \(\alpha\le\alpha_0\), the upper one-sided bound in
\eqref{tri:deficit-lipschitz} gives
\(d_\Sigma(\alpha)\ge m_0-(\alpha_0-\alpha)\ge m_0/2\).
If \(\alpha\ge\alpha_0\), the lower bound gives
\(d_\Sigma(\alpha)\ge m_0-5(\alpha-\alpha_0)\ge m_0/6\).
Thus, for every matrix outside the two neighborhoods, singular or not, the
inequality is strict at every level of \(I_0\), including both endpoints. Together with the two analytic estimates this proves
the result. The matrices in \eqref{tri:equalities} give equality at
every level; every other covariance lies in one of the strict cases.
The argument uses the finite covariance certificate only at
\(\alpha_0\).
\end{proof}

\subsection{Arithmetic assumptions and verification scope}
\label{tri:reproducibility}

The computer-assisted part of the proof checks finitely many explicit
inequalities. Calculations with fractions are exact, and
high-precision integer calculations round outward. The other numerical
calculations use the interval arithmetic of Section~\ref{tri:s6-2},
with endpoints in the standard double-precision floating-point format
(IEEE 754). Their correctness rests on two assumptions about the
computer. First, each basic operation and each
square root returns the representable number nearest to the exact
result, also for very small numbers, which must not be replaced by
zero. Second, the operations are carried out in the stated order,
without regrouping them and without merging a multiplication and an
addition into a single rounded step. Correct execution and reliable
hardware are also assumed. These are the usual assumptions of
computer-assisted proofs based on interval arithmetic.

The checks cover the bin identities, coefficient bounds and scalar
inequalities for the two analytic neighborhoods. However, a few elementary bounds
used in Section~\ref{tri:s4}, such as \(\Phi(-5)\le1/25\) and
\(\overline\Phi(6)<1/100\), are not among the checks; each can be
checked by hand. For the remaining region, the checks cover the function tables, the reference matrices, the
integration segments and their error bounds, the derivative bounds,
the covering of both chambers by leaves, and the final comparison with
$-2^{-20}$. All of them succeeded. Together with the analytic
arguments, they give a computer-assisted proof at $1/20$ under the
stated assumptions.

Section~\ref{tri:level-interval} obtains the common interval of levels
analytically from the certified margin and the two analytic
neighborhoods. It uses no certificate computations at other levels.

\clearpage
\section{Pair-deficit expansion for sparse covariance graphs}\label{sec:lm:graph}

\label{sec:lm:g:1}

We prove the graph expansion by separating the pair deficits for the full
rejection probability from interactions involving more coordinates. The pair deficits may
approach zero, so an error bound is useful only if it is small relative
to their sum. We first express the rejection probability as a sum of
contributions from the original tests, counting every repeated copy.
We then bound the pair deficits, the interactions among larger connected
sets, and the contribution from ranks beyond a chosen cutoff. When every
correlation is at most $(\alpha/m)^4$ in absolute value,
Section~\ref{sec:lm:g:4} instead compares the full laws and keeps every
rank. Section~\ref{sec:lm:g:boundary} combines these bounds into the proof of
Theorem~\ref{thm:lm:graph}; Section~\ref{main:graph_roadmap} outlines the
argument in four steps. Throughout, the constants depend only on the fixed
parameters in the theorem. The normal tail estimates are the Mills and
hazard inequalities~\eqref{sl:mills} of Lemma~\ref{sl:tails}. Its quantile
expansion~\eqref{sl:quantile} is uniform only when the level multiplier
stays in a fixed compact set, so quantiles at the growing rank cutoff are
bounded directly from~\eqref{sl:mills}.

\begin{lmtheorem}[Pair-deficit expansion]
\label{thm:lm:graph}\label{thm:lm:graph_boundary}
Fix an integer \(W\ge1\), \(0<r_1<1\), and \(0<a_1<1\).
Fix \(0\le\beta<(1-r_1)/(1+r_1)\), and prescribe positive integers
\(\Delta_m\) satisfying
\begin{equation}
 \limsup_{m\to\infty}
 \frac{(\Delta_m+1)\log(\Delta_m+1)}{\log m}\le\beta.
 \label{eq:lm:g:boundary.1}
\end{equation}
For each \(m\), let \(X=(X_1,\ldots,X_n)\) be centered Gaussian
with positive semidefinite unit-diagonal representative covariance \(\Gamma\). Give
representative \(i\) an integer multiplicity \(1\le w_i\le W\), with
\(\sum_iw_i=m\), and assume
\begin{equation}
 |\Gamma_{ij}|\le r_1\quad(i\ne j),\qquad
 \max_i\#\{j\ne i:\Gamma_{ij}\ne0\}\le\Delta_m.
 \label{eq:lm:g:1.2}
\end{equation}
Set \(P_i=2\barPhi(|X_i|)\), repeat \(P_i\) exactly \(w_i\) times,
and apply the original physical thresholds \(k\alpha/m\). Denote the
rejection probability by \(F_m(\Gamma,w;\alpha)\). Let \(E\) be the set
of nonzero representative edges. For \(e=\{i,j\}\), let \(\Gamma^e\)
retain only its actual correlation, making all other representatives
independent without changing their weights. Define
\begin{equation}
 d_e=\alpha-F_m(\Gamma^e,w;\alpha),\qquad
 \mathcal D_m=\sum_{e\in E}d_e.
 \label{eq:lm:g:1.3}
\end{equation}
Put
\begin{equation}
 g=\frac{1-r_1}{(1+2r_1)(1+r_1)}>0,
 \label{eq:lm:g:1.4}
\end{equation}
and
\[
 \vartheta=\Theta_{r_1}(\beta)
 =\min\left\{g,\frac{1-\beta}{r_1}-\frac2{1+r_1}\right\}>0.
\]
Every nonzero edge has \(d_e>0\). For every
\(0<\varepsilon<\vartheta\), for all sufficiently large \(m\),
simultaneously for \(0<\alpha\le a_1\) and uniformly over
\eqref{eq:lm:g:1.2} with \(E\ne\varnothing\),
\begin{equation}
 \frac{|F_m(\Gamma,w;\alpha)-\alpha+\mathcal D_m|}{\mathcal D_m}
 \le(\alpha/m)^{\vartheta-\varepsilon}.
 \label{eq:lm:g:boundary.2}
\end{equation}
Consequently,
\begin{equation}
 F_m=\alpha-\mathcal D_m+o(\mathcal D_m),\qquad
 F_m\le\alpha-\tfrac12\mathcal D_m<\alpha
 \label{eq:lm:g:1.6}
\end{equation}
uniformly eventually when an edge is present. If \(E=\varnothing\),
weighted independence gives \(F_m=\alpha\) in every dimension.

The threshold depends only on \(W,r_1,a_1,\beta,\varepsilon\) and the
prescribed degree envelope, which is fixed before choosing the
threshold. No positive lower bound on a nonzero correlation is imposed.
The covariance may be singular, and its graph, signs, edge magnitudes
and multiplicities may vary jointly with \(m\) and \(\alpha\). Exact
signed-copy classes retain their integer multiplicities, bounded by
\(W\).
\end{lmtheorem}

At \(\beta=0\), the degree condition is equivalent to
\begin{equation}
 \Delta_m\log(\Delta_m+1)=o(\log m).
 \label{eq:lm:g:1.1}
\end{equation}
The theorem then gives, for every \(0<\varepsilon<g\),
\begin{equation}
 \sup_{\substack{(\Gamma,w)\text{ satisfying }\eqref{eq:lm:g:1.2}\\
 E\ne\varnothing}}
 \frac{|F_m(\Gamma,w;\alpha)-\alpha+\mathcal D_m|}{\mathcal D_m}
 \le(\alpha/m)^{g-\varepsilon},
 \label{eq:lm:g:1.5}
\end{equation}
for all sufficiently large dimensions, uniformly over \(0<\alpha\le a_1\). The same full exponent
holds throughout \(0\le\beta\le(1-r_1)/(1+2r_1)\).

\subsection{Physical ranks and pair deficits}\label{sec:lm:g:2}

We begin by expressing rejection as a sum of contributions from the original tests, including repeated p-values. Throughout this section, \(t=\alpha/m\) and \(c_a=\Phi^{-1}(1-a/2)\). For signed measures we use the full variation norm; in this convention, a probability measure has norm one. Let

\[
N_w(a)=\sum_iw_i\mathbf1\{P_i\le a\},\quad
R=\max(\{k\in\{1,\ldots,m\}:N_w(kt)\ge k\}\cup\{0\}).
\]

At the largest crossing, exactly \(R\) physical p-values are at most \(Rt\). Indeed, any larger count would give a crossing at that larger rank. This argument also applies when several p-values are tied.

Force every copy of group \(i\) to zero and call the new largest rank \(R^{[i]}\). Then group \(i\) is rejected exactly when \(P_i\le tR^{[i]}\), and on that event \(R=R^{[i]}\). Consequently

\begin{equation}
\mathbf1\{R>0\}=\sum_i\sum_{k=1}^m\frac{w_i}{k}
\mathbf1\{P_i\le kt,R^{[i]}=k\}. \label{eq:lm:g:2.1}
\end{equation}

Under independent uniform representatives, group \(i\) has expected contribution \(w_it\). Summing these contributions gives \(\alpha\) for arbitrary integer multiplicities. We will use this weighted equality as the reference throughout the proof.

For uniform, possibly dependent group p-values, define

\[
C_{i,k}(a)=\Prob(R^{[i]}\le k\mid P_i\le a).
\]

The forcing identity also identifies the source of a deficit. We compare the conditional forced-rank probabilities at consecutive physical thresholds. Their increments give the signed difference between the independent contribution and the actual expectation.

Since \(P_i\) is uniform, the \(k\)th term of \eqref{eq:lm:g:2.1} for group
\(i\) has expectation \(w_it\{C_{i,k}(kt)-C_{i,k-1}(kt)\}\). Finite
summation by parts over \(k\) then gives, for the expected contribution
\(B_i\) of group \(i\),

\begin{equation}
w_it-B_i=w_it\sum_{k=1}^{m-1}
\{C_{i,k}((k+1)t)-C_{i,k}(kt)\}. \label{eq:lm:g:2.2}
\end{equation}

The endpoint functions are \(C_{i,0}=0\) and \(C_{i,m}=1\); ranks below \(w_i\) have zero \(C_{i,k}\).

Consider a one-edge reference with selected weights \(u,v\le W\), total weight \(q_e=u+v\), and correlation magnitude \(s\in(0,r_1]\). For an edge \(e=\{i,j\}\in E\), we write \(s_e=|\Gamma_{ij}|\) and \(s_{\max}=\max_{e\in E}s_e\). Each background group, that is, each group outside the selected pair, has expected contribution equal to its weight times \(t\) in \eqref{eq:lm:g:2.1}. Conditional on the magnitude of either selected Gaussian variable, the magnitude of its partner is stochastically increasing in the conditioning value. Since \(R^{[i]}\le k\) is decreasing in the partner magnitude, \(a\mapsto C_{i,k}(a)\) is nondecreasing. The same argument applies after fixing the independent background p-values, so every conditional increment in \eqref{eq:lm:g:2.2} is nonnegative. Only the monotonicity of a single folded-normal conditional distribution is needed here; no PRDS assumption is made for a larger folded Gaussian vector.

Let \(G\) be the number of background representatives. On the event \(\mathcal E\) that forcing both selected groups to zero gives largest crossing rank exactly \(q_e\), every background p-value exceeds \(q_et\). Rescaling background p-values by \((P_j-q_et)/(1-q_et)\) identifies \(\mathcal E\) with weighted Simes acceptance at level \((m-q_e)t/(1-q_et)\). Thus

\begin{equation}
p_{\mathcal E}=(1-\alpha)(1-q_et)^{G-1}. \label{eq:lm:g:2.3}
\end{equation}

For \(G=0\) this expression is exactly one. If \(q_et\le1/2\), then

\begin{equation}
p_{\mathcal E}\ge p_0:=(1-a_1)e^{-4Wa_1}>0. \label{eq:lm:g:2.4}
\end{equation}

Indeed \(G\le m\), \(q_e\le2W\), and \(\log(1-q_et)\ge-2q_et\).

On any fixed background state in \(\mathcal E\), forcing the weight-\(u\) group gives rank \(u\) when the partner p-value exceeds \(q_et\), and rank \(q_e\) otherwise. Retaining the ranks \(u,\ldots,q_e-1\) in \eqref{eq:lm:g:2.2} telescopes to

\[
J_s(ut,q_et)-(u/q_e)J_s(q_et,q_et),
\]

where \(J_s(a,b)=\Prob(P_1\le a,P_2\le b)\) for the Gaussian pair. The other group contributes the analogous expression. Hence

\begin{equation}
d_e\ge p_{\mathcal E}D_{u,v,s}(t),\quad
D_{u,v,s}(t)=J_s(ut,q_et)+J_s(vt,q_et)-J_s(q_et,q_et). \label{eq:lm:g:2.5}
\end{equation}

The quantity \(D_{u,v,s}(t)\) depends only on the selected pair. Multiplying it by the probability of the background event gives the lower bound in \eqref{eq:lm:g:2.5}. The full pair deficit averages over every background configuration, so this bound retains only part of that deficit.

\subsection{Uniform pair estimates}\label{sec:lm:g:3}

To compare the remainder with the pair deficits, we need lower bounds for the deficits and upper bounds for changes in the pair laws. Both bounds must remain useful when the correlation approaches zero.

\subsubsection{A lower bound from the second derivative}\label{sec:lm:g:3.1}

We turn the pair difference into an integral with a positive integrand obtained from the second derivative. This proves its sign and gives a lower bound that keeps track of the correlation even near independence.

Put \(B=q_et\), \(z=c_B\), \(f(a)=J_s(a,B)\), \(y=c_a\), and \(d=1-s^2\). Conditioning gives

\[
f'(a)=\barPhi((z-sy)/\sqrt d)+\barPhi((z+sy)/\sqrt d).
\]

Since \(dy/da=-1/(2\phi(y))\), differentiation yields

\begin{equation}
-f''(a)=\frac{s}{2\sqrt d}
\exp\!\left(\frac{y^2}{2}-\frac{(z-sy)^2}{2d}\right)
\left(1-e^{-2szy/d}\right)>0\quad(0<a\le B). \label{eq:lm:g:3.1}
\end{equation}

As \(a\downarrow0\), we have \(f'(a)\to1\) for each fixed \(s>0\). Since \(-f''\) in \eqref{eq:lm:g:3.1} is nonnegative, integration, with an improper integral at the endpoint if needed, gives the identity

\begin{equation}
D_{u,v,s}(t)=\int_0^{ut}\int_0^{vt}[-f''(a+b)]\,db\,da. \label{eq:lm:g:3.2}
\end{equation}

For a lower bound, integrate only over \(a\in[ut/4,ut/2]\) and \(b\in[vt/4,vt/2]\). This rectangle stays away from the endpoint where the derivative may be unbounded. Its area is \(uvt^2/16\) and \(a+b\in[B/4,B/2]\). The normal hazard lower bound implies \(0\le y^2-z^2\le2\log4\). If \(L(y)\) denotes the exponent in \eqref{eq:lm:g:3.1},

\[
L(y)-L(z)=\frac{(1-2s^2)(y^2-z^2)+2sz(y-z)}{2(1-s^2)}
\ge-\frac{r_1^2}{1-r_1^2}\log4,
\qquad L(z)=\frac{sz^2}{1+s}.
\]

Also
\[
1-e^{-2szy/d}\ge(1-e^{-2})\min(1,sz^2).
\]

Therefore, with the explicit constant

\[
c_0=\frac{1-e^{-2}}{32}\,4^{-r_1^2/(1-r_1^2)},
\]

we obtain

\begin{equation}
d_e\ge p_0c_0\,uvt^2s\,
 e^{sz^2/(1+s)}\min(1,sz^2). \label{eq:lm:g:3.3}
\end{equation}

The quantitative bound holds when \(q_et\le1/2\). In addition, the strict positivity of the integrand in \eqref{eq:lm:g:3.2} gives \(d_e>0\) for every nonzero imperfect edge in every dimension. At \(s=0\), both the reference deficit and the pair difference are zero.

The bound is uniform in \(s\), including its limit at zero: it retains the factor \(s^2z^2\) whenever \(sz^2\le1\). Retaining this factor will let us compare errors with arbitrarily small pair deficits.

\subsubsection{Censored pair variation}\label{sec:lm:g:3.2}

For the first few ranks, the exact p-values above a given threshold do not matter. Replace every p-value above \(x\) by one symbol \(*\); we call this censoring. Let \(\nu\) be the resulting censored uniform marginal and \(\delta_s^x=\mu_s^x-\nu\otimes\nu\). This signed measure is the difference between the censored dependent-pair law and its independent counterpart. If \(c=c_x\ge1\), then

\begin{equation}
\|\delta_s^x\|\le C_{r_1}s^2c^2e^{-c^2/(1+s)},\quad
 C_{r_1}=\frac{56\cdot11940}{\pi(1-r_1^2)^{9/2}}. \label{eq:lm:g:3.4}
\end{equation}

We prove the bound by treating the continuous part and the atoms of the censored law. The variation norm controls changes in expectations of bounded functions of the censored values, including the function used to determine rejection. The key point is to retain the squared correlation in this bound, so that it can later be divided by the pair deficit. Let \(g_r(y,z)\) denote the bivariate Gaussian density and put \(d_0=1-r_1^2\). The folded density on the positive quadrant is \(q_s=2(g_s+g_{-s})\), with \(q'_0=0\). Direct differentiation gives

\[
\partial_r\log g_r=\frac r{1-r^2}
 +\frac{(1+r^2)yz-r(y^2+z^2)}{(1-r^2)^2}.
\]

Using \(2|yz|\le y^2+z^2\) gives

\[
|\partial_r\log g_r|\le2d_0^{-2}(1+y^2+z^2),\quad
|\partial_r^2\log g_r|\le10d_0^{-3}(1+y^2+z^2),
\]

and hence

\begin{equation}
|\partial_r^2g_r|\le14d_0^{-4}(1+y^2+z^2)^2g_r. \label{eq:lm:g:3.5}
\end{equation}

For \(y=c+a/c,z=c+b/c\), the density exponent is exactly

\[
\frac{c^2}{1+r}+\frac{a+b}{1+r}
 +\frac{a^2-2rab+b^2}{2(1-r^2)c^2}.
\]

The last term is nonnegative. Put \(h=a+b\). Since \(c\ge1\) and
\(a^2+b^2\le h^2\), the polynomial is bounded by
\[
 (1+y^2+z^2)^2\le c^4(3+2h+h^2)^2.
\]
The Jacobian is \(c^{-2}\). Thus, uniformly for \(|r|\le s\le r_1\),

\[
\int_c^\infty\!\int_c^\infty(1+y^2+z^2)^2g_r(y,z)\,dy\,dz
\le\frac{I}{2\pi\sqrt{d_0}}c^2e^{-c^2/(1+s)},
\]

where

\[
I=\int_0^\infty h(3+2h+h^2)^2e^{-h/2}\,dh
=9\cdot1!2^2+12\cdot2!2^3+10\cdot3!2^4
 +4\cdot4!2^5+5!2^6=11940.
\]

Taylor's integral remainder for \(q_s-q_0\), together with \eqref{eq:lm:g:3.5}, bounds the variation on the square where both p-values remain uncensored by \(14I s^2c^2e^{-c^2/(1+s)}/(\pi d_0^{9/2})\). To bound the other three parts of the censored law, we use the fact that both one-coordinate marginals of \(\delta_s^x\) vanish. Each part with exactly one censored coordinate is the negative of the corresponding marginal of the uncensored part. The mass with both coordinates censored equals the total mass of the uncensored part. Thus each of the four variations is bounded by the variation on the uncensored square, and their sum is at most four times that variation. This proves \eqref{eq:lm:g:3.4} without assuming that the density difference has a fixed sign.

\subsubsection{Consequences for the total deficit}\label{sec:lm:g:3.3}

We collect three consequences for later use. The first lower bound will match the sum of squared correlations in weak connected sets. The second will retain the stronger tail exponent needed for the other connected sets. The final estimate will compare the change in a pair law directly with that pair's deficit.

Let \(L=\log(1/t)=\log(m/\alpha)\). The elementary Mills bounds give the following quantile estimates, uniformly for \(2\le q_e\le2W\),

\begin{equation}
z_e^2=c_{q_e\alpha/m}^2=2L+O(\log L),\qquad
L\le z_e^2\le3L \label{eq:lm:g:3.6}
\end{equation}

eventually. For completeness, if \(z\ge1\) then

\[
\frac{\phi(z)}z\le2\barPhi(z)\le\frac{2\phi(z)}z.
\]

The upper inequality first bounds \(z^2\le2\log(1/B)+O(1)\); substitution in the lower inequality gives \(z^2\ge2\log(1/B)-\log\log(1/B)-O(1)\). Since \(B=q_et\) with \(2\le q_e\le2W\), \eqref{eq:lm:g:3.6} follows uniformly.

For \(a=sz_e^2\),

\[
e^{a/(1+s)}\min(1,a)\ge a/2.
\]

Together with \eqref{eq:lm:g:3.3} this gives a fixed \(b>0\) such that

\begin{equation}
d_e\ge b\,s_e^2Lt^2,\qquad
\mathcal D_m\ge bLt^2\sum_es_e^2. \label{eq:lm:g:3.7}
\end{equation}

One may take \(b=p_0c_0/2\) for all sufficiently large dimensions in the stated uniform range. No positive lower bound on the nonzero correlations is needed.

We also need a sharper lower bound that keeps the exponential dependence on the correlation. For \(z_e^2\ge1\),

\[
s\min(1,sz_e^2)\ge\frac{s^2z_e^2}{1+sz_e^2}\ge s^2/2.
\]

Using the more precise lower quantile bound in the preceding paragraph gives, for another fixed \(b'>0\),

\begin{equation}
d_e\ge b' s_e^2L^{-1}t^{\gamma(s_e)},\qquad
\gamma(s)=\frac2{1+s}. \label{eq:lm:g:3.8}
\end{equation}

Indeed, the lower quantile bound gives
\(e^{sz_e^2/(1+s)}\ge c\,(q_et)^{-2s/(1+s)}L^{-s/(1+s)}\) for a fixed
\(c>0\), and \(t^2(q_et)^{-2s/(1+s)}\ge(2W)^{-1}t^{\gamma(s)}\). Since
\(L^{-s/(1+s)}\ge L^{-1}\), the factor \(L^{-1}\) gives a lower bound
uniform over \(0<s\le r_1\).

Finally put \(x=Kt\), \(K\ge q_e\), \(c_x\ge1\), and \(z=z_e\). Hazard and Mills bounds give

\[
z^2-c_x^2\le2\log(K/q_e),\qquad e^{-z^2}/t^2\le2\pi q_e^2z^2.
\]

Dividing \eqref{eq:lm:g:3.4} by \eqref{eq:lm:g:3.3}, and using \(s/\min(1,sz^2)\le(1+sz^2)/z^2\), proves

\begin{equation}
\frac{\|\delta_s^{Kt}\|}{d_e}
\le C_* K^{\gamma(s)}z^2(1+sz^2)
\le C K^2L^2. \label{eq:lm:g:3.9}
\end{equation}

One may take \(C_*=4\pi C_{r_1}/(p_0c_0)\), since for
\(q_e=u+v\) and \(\gamma(s)\ge1\),
\(q_e^{2-\gamma(s)}/(uv)\le q_e/(uv)=1/u+1/v\le2\).
Dependence on the weight bound remains through \(p_0\). All constants
are independent of how rapidly \(s\) approaches zero and of the reference
background weight configuration.

\subsection{Very small correlations}\label{sec:lm:g:4}

When every correlation is very small, an absolute bound on the probability of a late crossing may exceed the pair deficits. Removing those ranks would then give too large an error for a relative comparison. We instead compare the full laws of the absolute Gaussian coordinates and retain all ranks. Let \(R\) be symmetric with zero diagonal, \(\|R\|_{\mathrm{op}}\le1/2\), and \(S=\operatorname{tr}(R^2)\). Denote by \(\mu_R\) the law of the coordinatewise absolute values of a
centered Gaussian vector with covariance \(I+R\), and let \(\mu_e\) retain
only edge \(e\). We establish

\begin{equation}
\left\|\mu_R-\mu_0-\sum_e(\mu_e-\mu_0)\right\|
\le6S^{3/2}. \label{eq:lm:g:4.1}
\end{equation}

We expand around independence through second order. After averaging over coordinate signs, the linear term disappears and the quadratic term is the sum of the individual edge quadratic terms. A dimension-independent bound on the third-order remainders will then allow us to subtract the complete one-edge laws.

To obtain a constant independent of dimension, let \(g_\lambda\) have covariance \(I+\lambda R\), $0\le\lambda\le1$. Represent a vector with this covariance as \(Z=(I+\lambda R)^{1/2}Y\), with \(Y\) standard normal, and put \(B=R(I+\lambda R)^{-1}\), \(b_B=\|B\|_F\le2\sqrt S\). Differentiating the log density at fixed \(Z\), then writing the result in terms of \(Y\), gives

\[
\mathcal L'=\tfrac12(Y^TBY-\operatorname{tr}B),\quad
\mathcal L''=\tfrac12\operatorname{tr}(B^2)-Y^TB^2Y,\quad
\mathcal L'''=-\operatorname{tr}(B^3)+3Y^TB^3Y.
\]

The centered Gaussian quadratic-form moment identities imply

\[
\|\mathcal L'\|_2=b_B/\sqrt2,\quad
\|\mathcal L''\|_2\le(1/2+\sqrt2)b_B^2,\quad
\Expect|\mathcal L'''|\le4b_B^3,
\]

and \(\Expect(Y^TBY-\operatorname{tr}B)^4\le60b_B^4\). Hence

\[
\|g'''_\lambda\|_1
\le\left[4+\frac3{\sqrt2}(1/2+\sqrt2)+\frac{60^{3/4}}8\right]b_B^3
<11b_B^3\le88S^{3/2}.
\]

Keeping the sharper bound $b_B\le\sqrt S/(1-\lambda/2)$ in the
Taylor integral, its remainder after second order has $L^1$ norm at most
\[
 11S^{3/2}\int_0^1\frac{(1-\lambda)^2}{2(1-\lambda/2)^3}\,d\lambda
 =11(4\log2-5/2)S^{3/2}<3S^{3/2}.
\]
The last inequality is equivalent to $\log2<61/88$, which follows by bounding the tail of
$2\sum_{j\ge0}3^{-2j-1}/(2j+1)$ by a geometric series. Averaging over coordinate signs cannot increase the \(L^1\) norm. This averaging removes \(g'_0\), and the averaged second derivative at zero is exactly

\begin{equation}
\operatorname{Even}(g''_0/g_0)=
\sum_{i<j}r_{ij}^2(z_i^2-1)(z_j^2-1). \label{eq:lm:g:4.2}
\end{equation}

Indeed, before sign averaging the expression is \((z^TRz)^2/4+\operatorname{tr}(R^2)/2-z^TR^2z\). Every product involving two distinct edges has an odd coordinate and therefore vanishes under the averaging. Formula \eqref{eq:lm:g:4.2} is precisely the sum of the individual edge second derivatives. Subtract the complete edge laws and apply the remainder bound to each. Since \(S_e=2r_e^2\) and \(\sum_e S_e^{3/2}\le S^{3/2}\), \eqref{eq:lm:g:4.1} follows. Taking coordinatewise absolute values is a measurable map and introduces no additional norm factor.

The signed remainder has total mass zero. The absolute value of its integral against the original physical rejection indicator is therefore bounded by half its full variation norm, giving

\begin{equation}
|F_m-\alpha+\mathcal D_m|\le3S^{3/2}. \label{eq:lm:g:4.3}
\end{equation}

If \(s_{\max}\le t^4\), then \(\|R\|_{\mathrm{op}}\le mt^4\le t^3\) and \(\sqrt S\le m s_{\max}\). By \eqref{eq:lm:g:3.7}, for \(S>0\),

\begin{equation}
\frac{|F_m-\alpha+\mathcal D_m|}{\mathcal D_m}
\le C\frac{\sqrt S}{Lt^2}
\le C\frac{mt^2}{L}=C\frac{\alpha t}{L}\le C\frac tL. \label{eq:lm:g:4.4}
\end{equation}

This proves the required relative bound when all correlations are very small, including for dense representative graphs. The comparison uses the full laws, so it remains valid when the deficit is smaller than every polynomial rate.

\subsection{Finite centered measures and late physical ranks}\label{sec:lm:g:5}

For the remaining correlations, we first restrict attention to ranks below a cutoff. The corresponding event depends only on the censored p-values. We expand their joint law into signed measures that isolate dependence within subsets of coordinates. The covariance graph will identify which subsets contribute and bound their number. A separate tail bound will then restore the omitted ranks.

Let \(V_i=P_i\) if \(P_i\le x\) and \(V_i=*\) otherwise, let \(\mu_T^x\) denote the joint law of \((V_i)_{i\in T}\), and write \(\mu^x\) for the joint law of all \(n\) censored values. For \(C\subseteq[n]\), define

\begin{equation}
A_C=\sum_{T\subseteq C}(-1)^{|C|-|T|}
\mu_T^x\otimes\nu^{\otimes(C\setminus T)},\qquad A_\varnothing=1. \label{eq:lm:g:5.1}
\end{equation}

Equivalently, \(A_C=\Expect[\bigotimes_{i\in C}(\delta_{V_i}-\nu)]\),
where the expectation is a finite signed measure, defined by integrating
against bounded measurable test functions and averaging over the joint
law of \(V_C\). Expanding the product gives \eqref{eq:lm:g:5.1}.
The tensor product is taken before expectation, so this identity does
not assume independence within \(C\). In the kernel expansion of Section~\ref{sec:lm:g:5.1},
each low-value factor has norm at most two and each fixed centering factor \(b_x\) defined there
has norm \(2x\), which gives \eqref{eq:lm:g:5.4} after taking
expectation and summing over the subsets.

Each measure is defined by a finite sum of marginal laws. The centered measures reconstruct the full censored joint law through the following sum:

\begin{equation}
\mu^x=\sum_{C\subseteq[n]}A_C\otimes\nu^{\otimes C^c}. \label{eq:lm:g:5.2}
\end{equation}

Singleton terms vanish because every coordinate has the same censored uniform marginal. If the induced covariance graph on \(C\) has components \(C_1,\ldots,C_b\), the zero cross-covariances make the corresponding Gaussian subvectors independent, giving

\begin{equation}
A_C=\bigotimes_{j=1}^b A_{C_j}. \label{eq:lm:g:5.3}
\end{equation}

This factorization uses marginal independence; it does not condition on variables outside \(C\). Consequently, \eqref{eq:lm:g:5.2} can be organized as a sum over collections of connected sets of size at least two. The sets in a collection must be disjoint, with no edge joining different sets; we call such a collection compatible. For an edge \(e\), \(A_e=\delta_{s_e}^x\), so the pair terms are exactly the changes in the censored one-edge laws bounded above.

\subsubsection{Bounds for connected terms}\label{sec:lm:g:5.1}

To bound a connected term, we separate each coordinate into its low-value contribution and a fixed centering measure. This reduces the variation norm to joint lower-tail probabilities. We will bound those probabilities either by a large independent subset or by a fixed number of correlated coordinates.

Writing \(b_x=\int_0^x(\delta_p-\delta_*)\,dp\), with \(\|b_x\|=2x\), expand

\[
\delta_{V_i}-\nu=\mathbf1\{P_i\le x\}(\delta_{P_i}-\delta_*)-b_x.
\]

Finite integration of these kernels yields

\begin{equation}
\|A_C\|\le2^k\sum_{T\subseteq C}x^{k-|T|}
\Prob(P_i\le x\text{ for every }i\in T),\quad k=|C|. \label{eq:lm:g:5.4}
\end{equation}

Let \(\chi=\Delta_m+1\). A proper coloring partitions the graph into classes with no internal edges. Gaussian coordinates within each class are independent, so coloring provides an independent set of size at least \(|T|/\chi\) in each \(T\). Keeping only that set in the joint tail event gives

\begin{equation}
\|A_C\|\le4^k x^{k/\chi}. \label{eq:lm:g:5.5}
\end{equation}

For smaller connected sets, a bound based on three selected coordinates is more useful. The following estimate is uniform even as the largest correlation within the set approaches zero:

\begin{equation}
\|A_C\|\le216\,4^k c_x^3x^{\lambda(r_C)},\quad
k\ge3,\quad \lambda(r)=\frac3{1+2r}. \label{eq:lm:g:5.6}
\end{equation}

Here \(r_C=\max_{i<j\in C}|\Gamma_{ij}|\), and \(c_x\ge1\). To prove it, a signed-sum Chernoff bound on any two or three coordinates gives respectively

\[
\Prob(P_i,P_j\le x)\le4e^{-\gamma(r_C)c_x^2/2},\quad
\Prob(P_i,P_j,P_k\le x)\le8e^{-\lambda(r_C)c_x^2/2}.
\]

Mills' inequality gives \(e^{-c_x^2/2}\le3c_xx\). Since \(\gamma\le2\) and \(\lambda\le3\), these are at most \(36c_x^2x^{\gamma(r_C)}\) and \(216c_x^3x^{\lambda(r_C)}\). In \eqref{eq:lm:g:5.4}, subsets of size at least three use any selected triple. Two-element subsets have the extra factor \(x^{k-2}\), and \(1+\gamma(r)\ge\lambda(r)\). Empty and singleton subsets give \(x^k\), with \(k\ge\lambda(r)\). There are \(2^k\) terms. This proves \eqref{eq:lm:g:5.6} with a constant independent of the size of the connected set. The argument also covers singular covariances, since a signed sum with zero variance cannot exceed the positive threshold.

The same two-coordinate bound in \eqref{eq:lm:g:5.4} gives

\begin{equation}
\|A_e\|\le156c_x^2x^{\gamma(s_e)}. \label{eq:lm:g:5.7}
\end{equation}

\subsubsection{A tail bound from graph coloring}\label{sec:lm:g:5.2}

We next bound crossings at the ranks omitted from the censored comparison. For \(Y_i=w_i\mathbf1\{P_i\le\alpha j/m\}\), apply H\"{o}lder's inequality across at most \(\chi\) color classes and use independence within each class. The secant bound for \(0\le w_i\le W\) gives

\[
\Expect e^{\theta\sum_iY_i}
\le\exp\left\{\frac{\alpha j}{W\chi}
(e^{\theta W\chi}-1)\right\}.
\]

Here the secant bound is
\(e^{\chi\theta w}-1\le(w/W)(e^{\chi\theta W}-1)\) for \(0\le w\le W\),
which holds by convexity in \(w\). Markov's inequality with
\(\theta=\log(1/\alpha)/(W\chi)\) turns this into a tail bound.
With \(h_\alpha=\log(1/\alpha)-1+\alpha>0\) and \(b_m=h_\alpha/(W\chi)\),

\begin{equation}
\Prob\{N_w(\alpha j/m)\ge j\}\le e^{-b_mj}. \label{eq:lm:g:5.8}
\end{equation}

This holds for the actual law, the independent weighted baseline, and every one-edge reference. Put
\[
 K=\max\left\{2W,\left\lceil\frac{20W\chi L}{h_\alpha}\right\rceil\right\},
 \qquad h_\alpha=\log(1/\alpha)-1+\alpha,
 \qquad x=Kt,\qquad T_m=\frac{e^{-b_m(K+1)}}{1-e^{-b_m}}.
\]
The following estimates also specify the uniformity used in the rest of the proof.
\begin{lemma}[Common cutoff and uniform factors]\label{mx:uniform_factors}
First assume the slow-degree condition \eqref{eq:lm:g:1.1}.
For all sufficiently large $m$, uniformly over $0<\alpha\le a_1$,
\[
 K\le C_{W,a_1}\chi(1+\log m)<m,
 \quad \log K=o(L),\quad
 \chi\log(\Delta_m+1)=o(L),
\]
and, for every fixed $p,C>0$,
\begin{equation}\label{mx:subpolynomial}
 C\Delta_m^pK^pL^p\exp\{C\chi\log(\Delta_m+1)\}
 =\exp\{o(L)\}=t^{-o(1)}.
\end{equation}
Moreover, $c_x^2=2(L-\log K)+O(\log L)$ and
$L\le c_x^2\le3L$ uniformly eventually.
Under the weaker condition
$\chi\log(\Delta_m+1)=O(\log m)$, the bounds on $K$ and the quantile conclusions
still hold, and every fixed polynomial $C\Delta_m^pK^pL^p$ is
$t^{-o(1)}$. Only the exponential factor in
\eqref{mx:subpolynomial} requires the little-o degree condition. More generally, the cutoff may be modified. Fix $C_0>0$ and let
$1\le\tilde\chi\le C_0\chi$. Define $K$ as above, but with $\tilde\chi$ in
place of $W\chi$ and with the numbers $20$ and $2W$ replaced by other
fixed positive numbers. Then, under either degree condition, $K\le C_{W,a_1}\tilde\chi(1+\log m)<m$
and $\log K=o(L)$, and the quantile and polynomial conclusions still hold. Every fixed polynomial in
$\tilde\chi$ is also $t^{-o(1)}$. The constants, and how large $m$ must
be, then depend also on $C_0$ and on these fixed numbers.
\end{lemma}
\begin{proof}
The cutoff must control late ranks without becoming too large when the testing level is very small. We first compare the logarithm of the level with the exponential tail rate. This will bound the cutoff in terms of dimension alone, while leaving the larger logarithmic factors in the subsequent estimates.

Write $u=\log(1/\alpha)\ge u_1=\log(1/a_1)>0$.
The ratio $(u-1+e^{-u})/u$ is positive and increasing, since
its derivative is $[1-(1+u)e^{-u}]/u^2>0$.
Thus $h_\alpha\ge c_1u$, while $h_\alpha\ge h_{a_1}>0$.
Consequently $L/h_\alpha\le C_{a_1}(1+\log m)$ and the
claimed cutoff bound follows. The degree assumption gives
$\chi\log(\Delta_m+1)=o(\log m)$, whereas $L\ge\log(m/a_1)$.
Also $\sup_{L\ge\log(m/a_1)}\log L/L\to0$.
Taking logarithms proves \eqref{mx:subpolynomial} for each fixed
finite collection of exponents. Mills' inequalities applied to $x=Kt$
give the quantile statement. In particular, this argument does not
require $L=O(\log m)$ when $\alpha$ is very small. Under the weaker condition, $\chi=O(\log m)$ and
$K=O((1+\log m)^2)$; the same logarithmic calculation proves the
polynomial assertion without discarding the exponential counting cost. The same calculation applies to the modified cutoff: since
$\tilde\chi\ge1$, the fixed numbers change only the constant in the
cutoff bound, and $\tilde\chi\le C_0\chi=O(\log m)$.
\end{proof}
Here and below, as in \eqref{mx:subpolynomial}, $t^{a-o(1)}$ denotes a
quantity at most $t^{a-\tau_m}$, where $\tau_m\ge0$ may change at each
occurrence, depends only on $m$ and the fixed parameters, and tends to
zero; the same $\tau_m$ serves every $0<\alpha\le a_1$ and every
$(\Gamma,w)$ in the class.
Since $(1-e^{-v})^{-1}\le1+v^{-1}$ and $h_\alpha\ge h_{a_1}$,
\begin{equation}
T_m\le(1+W\chi/h_\alpha)t^{20}\le C_{W,a_1}\chi t^{20}
=t^{20-o(1)}. \label{eq:lm:g:5.9}
\end{equation}

We retain the prefactor in the geometric tail bound because the degree bound may grow with dimension. Let \(H_K\) denote rejection at one of the physical ranks \(1,\ldots,K\). This indicator is a bounded function of \(V\), since larger p-values cannot affect these ranks. If rank \(K\) lies within a tied group, we still test that rank at \(Kt\). Using the count criterion \(N_w(jt)\ge j\) retains every such rank.

\subsubsection{The remainder of the centered expansion}\label{sec:lm:g:5.3}

We now collect the error into quantities that can be bounded separately: individual pair norms, larger connected terms, and the late ranks that were omitted. Let

\[
V_2=\sum_e\|A_e\|,\quad
B_3=\sum_{C\text{ connected},\ |C|\ge3}\|A_C\|,\quad V=V_2+B_3.
\]

After the independent baseline and the individual pair terms are removed, two types of terms remain. A collection either contains a connected set of at least three vertices or consists of at least two pairs. Selecting one of the larger connected sets bounds the total norm of the first type by \(B_3e^V\). The total norm of the second type is at most \(V_2^2e^V/2\). In both bounds we enlarge the sum by dropping the compatibility restriction and allowing repeated ordered choices. The remainder is therefore bounded by

\begin{equation}
(B_3+V_2^2/2)e^V. \label{eq:lm:g:5.10}
\end{equation}

Under the independent baseline the censored law is \(\nu^{\otimes n}\), and
under the one-edge reference \(\Gamma^e\) it is
\(\nu^{\otimes n}+A_e\otimes\nu^{\otimes e^c}\). Hence the integrals of
\(H_K\) against the baseline and pair terms give \(\alpha-\mathcal D_m\),
up to crossings above rank \(K\).
For \(E\ne\varnothing\), we integrate \(H_K\) against the signed remainder and use the rank-truncation bounds to recover the full rejection probabilities. This gives

\begin{equation}
|F_m-\alpha+\mathcal D_m|
\le(B_3+V_2^2/2)e^{V_2+B_3}+|E|T_m. \label{eq:lm:g:5.11}
\end{equation}

The sharper coefficient $|E|$ uses monotonicity of the rejection event
as the rank cutoff increases. Write each full-minus-truncated error as
$\delta\in[0,T_m]$. The restoration term is
$[\delta_{\rm actual}+(|E|-1)\delta_0]-\sum_e\delta_e$;
the bracketed sum and $\sum_e\delta_e$ are both nonnegative and at most $|E|T_m$. The bound thus
accounts for connected terms of every size and products of two or more
terms. Later estimates may harmlessly retain $2|E|T_m$. For the clipped
shares of one component used in Lemma~\ref{mx:lc}, the individual errors
need not have this sign, so
\eqref{mx:lc_remainder} retains the coefficient $2|E|$.

\subsection{Weak connected terms}\label{sec:lm:g:6}

The bound based on three coordinates is not enough when all correlations in a connected set are small. We now keep a squared-correlation factor in the bound for every such set, so that it can be compared with the pair deficit. For a set \(C\) of size \(k\ge2\), let $R_C=\Gamma_C-I_C$ and put \(S_C=\operatorname{tr}(R_C^2)\) and \(\eta_C=\max_i\sum_j|(R_C)_{ij}|\). If

\[
c=c_x\ge1,\qquad \eta_C\le\min(1/8,1/c^2),
\]

then\nopagebreak
\begin{equation}
\|A_C\|\le24S_C k^2c^4(8000x)^k. \label{eq:lm:g:6.1}
\end{equation}

We first compare each principal marginal with independence, retaining its squared-correlation factor throughout the Taylor estimate. We then substitute these marginal differences into the formula for the centered measure. Keeping track of the size of the set will allow the resulting bounds to be summed over all connected sets. For a principal marginal \(T\) of size \(j\), let \(g_\lambda\) have covariance \(I+\lambda R_T\) and write \(Z=\|z\|^2\), \(S_T=\operatorname{tr}(R_T^2)\le j\). With \(M=(I+\lambda R_T)^{-1}\), the zero trace gives \(\operatorname{tr}(MR_T)=-\lambda\operatorname{tr}(MR_T^2)\). Differentiating the log density and using \(\|M\|_{\mathrm{op}}\le2\), we obtain

\begin{equation}
|\mathcal L'|\le S_T+2\sqrt{S_T}Z,\quad
|\mathcal L''|\le2S_T+8S_TZ,\quad
|\mathcal L''+(\mathcal L')^2|\le16S_T(j+1+Z^2). \label{eq:lm:g:6.2}
\end{equation}

For the third bound, square the first inequality using \((a+b)^2\le2a^2+2b^2\), use \(S_T^2\le jS_T\), and \(8Z\le4+4Z^2\). The resulting polynomial is bounded by the final expression.

Expanding the determinant in terms of the eigenvalues and using the zero trace also gives

\begin{equation}
\frac{g_\lambda(z)}{g_0(z)}
\le\exp\{j\eta_C^2/[4(1-\eta_C)]\}
\exp\{\eta_CZ/[2(1-\eta_C)]\}
\le e^j e^{\eta_CZ}. \label{eq:lm:g:6.3}
\end{equation}

Averaging over coordinate signs makes the first derivative at zero vanish. If \(\lambda_T^x\) is the actual-minus-independent signed law restricted to all \(j\) p-values at most \(x\), second-order Taylor integration gives

\begin{equation}
\|\lambda_T^x\|\le8S_Te^j
\Expect_0[(j+1+\|Y\|^4)e^{\eta_C\|Y\|^2}
\mathbf1\{|Y_i|\ge c\ \forall i\}]. \label{eq:lm:g:6.4}
\end{equation}

For \(p=0,2,4\), substituting \(Y=c+a/c\) and using lower Mills gives

\[
I_p:=\Expect[|Y|^pe^{\eta_CY^2}\mathbf1\{|Y|\ge c\}]
\le2e\,c^px\int_0^\infty(1+a)^pe^{-3a/4}\,da
\le1300c^px.
\]

Here \(\eta_Cc^2\le1\), \(1-2\eta_C\ge3/4\), and the last integral for \(p=4\) is less than \(214\). Expanding \(\|Y\|^4\) under independence bounds \eqref{eq:lm:g:6.4} by

\begin{equation}
\|\lambda_T^x\|\le24S_Tj^2c^4(3900x)^j. \label{eq:lm:g:6.5}
\end{equation}

The empty and singleton differences are zero. For the kernel expansion in Section~\ref{sec:lm:g:5}, subtract the expectation under independence before taking the variation norm. This order keeps the cancellation at independence and therefore retains the squared-correlation factor. Each low-value kernel has norm two and each remaining fixed \(b_x\) factor has norm \(2x\). Therefore

\[
\|A_C\|\le2^k\sum_{T\subseteq C}x^{k-|T|}\|\lambda_T^x\|.
\]

Inserting \eqref{eq:lm:g:6.5}, using \(S_T\le S_C\) and \(j\le k\), proves \eqref{eq:lm:g:6.1}, since \(2(1+3900)<8000\). The same squared-correlation factor therefore controls the atomic parts of the censored law.

Call a connected set weak when

\begin{equation}
r_C\le\sigma:=\frac1{\Delta_m\max(8,c_x^2)}. \label{eq:lm:g:6.6}
\end{equation}

Each row of $R_C$ has at most $\Delta_m$ nonzero entries, each of
magnitude at most $r_C$, so a weak set has
$\eta_C\le\Delta_m r_C\le\min(1/8,c_x^{-2})$ and \eqref{eq:lm:g:6.1}
applies to it. A connected set that is not weak is called strong.
We first record a connected-set count that will also be used for strong
sets. If an edge $e=\{u,v\}$ is fixed, a connected $k$-set containing it
has a spanning tree containing $e$. Choose such a tree canonically, root
it at $u$, and place $v$ first among the children of $u$. An ordered tree
shape and the neighbor index for each of its other $k-2$ edges recover
the set. There are at most $4^{k-2}$ shapes: the Catalan numbers start at
one for $k=2$, and their successive ratios are less than four. Thus
\begin{equation}\label{eq:lm:edge_tree_count}
 N_k(e)\le(4\Delta_m)^{k-2}.
\end{equation}
Allowing repeated vertices or invalid embeddings only enlarges this
count. Neither covariance signs nor singularity affect it.

For weak sets it is better to retain the squared correlations along all edges rather than
assign the whole set to one edge. Interchanging nonnegative sums in
\eqref{eq:lm:g:6.1}, using $S_C=2\sum_{e\subset C}s_e^2$, gives, with
$q=32000\Delta_m x\le1/2$,
\begin{equation}
 B_{\mathrm w}\le\frac{3c_x^4}{\Delta_m^2}\sum_es_e^2
       \sum_{k\ge3}k^2q^k
 \le6912\cdot8000^3\Delta_m c_x^4x^3\sum_es_e^2.
 \label{eq:lm:g:6.7}
\end{equation}
Here $\sum_{k\ge3}k^2q^{k-3}\le36$ for $q\le1/2$, by differentiating the geometric series and summing the resulting terms, or by evaluating the nonnegative series at
$1/2$. Dividing by \eqref{eq:lm:g:3.7} gives
\begin{equation}
 B_{\mathrm w}/\mathcal D_m
       \le C\Delta_m K^3Lt=t^{1-o(1)}.
 \label{eq:lm:g:6.8}
\end{equation}
The last equality, which records the power up to a vanishing exponent, applies when the degree envelope satisfies \eqref{eq:lm:g:1.1} or, more
generally, \eqref{eq:lm:g:boundary.1}. The explicit polynomial bound
also applies under other envelopes when $32000\Delta_mKt\le1/2$.
The sum of squared correlations along edges cancels against the same factor in the deficit lower bound,
even when every nonzero correlation is arbitrarily small.

\subsection{Proof of the graph expansion}
\label{sec:lm:g:boundary}
We prove Theorem~\ref{thm:lm:graph} over its full degree range by keeping
the connected-set counting cost and comparing each strong set with the deficit of
its own largest-magnitude edge.

\begin{proof}
Write $\Delta=\Delta_m$, $\chi=\Delta+1$, $t=\alpha/m$, and $L=\log(1/t)$.
There is a nonnegative deterministic $\eta_m\to0$, independent of
$\alpha$, such that
\begin{equation}
 \chi\log\chi\le(\beta+\eta_m)\log m\le(\beta+\eta_m)L,
 \qquad \Delta=O(\log m),\quad\chi=o(\log m).
 \label{eq:lm:g:boundary.3}
\end{equation}
For the last relation in this display, if $\chi\ge c\log m$ along a subsequence for
some fixed $c>0$, then $\chi\log\chi\ge c\log m\,\log(c\log m)$, which
contradicts the first inequality for large $m$. Lemma~\ref{mx:uniform_factors}
gives $K\le C\chi(1+\log m)$, $\log K=o(L)$, and $L\le c_x^2\le3L$.
These conclusions, like the polynomial bound in the next sentence, need
only the lemma's weaker condition $\chi\log(\Delta_m+1)=O(\log m)$,
which follows from the first inequality in \eqref{eq:lm:g:boundary.3}.
Every fixed polynomial in $\Delta,K,L$ is $t^{-o(1)}$, uniformly for
arbitrarily small $\alpha$; the exponential counting cost will remain
explicit. The late-rank bound \eqref{eq:lm:g:5.9} is $T_m\le C\chi t^{20}$.

The full-distribution argument for very small correlations in \eqref{eq:lm:g:4.4} applies when
\(s_{\max}\le t^4\). In the remaining regime,
\begin{equation}
 s_{\max}>t^4,
 \label{eq:lm:g:7.1}
\end{equation}
\eqref{eq:lm:g:3.7} gives
\begin{equation}
 \mathcal D_m\ge bLt^{10}.
 \label{eq:lm:g:7.2}
\end{equation}
This coarse bound will be used only for late-rank restoration. The improved weak-set sum \eqref{eq:lm:g:6.8} gives
$B_{\mathrm w}/\mathcal D_m\le C\Delta K^3Lt=t^{1-o(1)}$.
Assign every remaining strong connected set to an edge of largest magnitude within that set, using a fixed rule for ties. By \eqref{eq:lm:g:6.6}, its assigned magnitude $s$ exceeds
$1/[\Delta\max(8,c_x^2)]$, so \eqref{eq:lm:g:3.8} and $\max(8,c_x^2)\le3L$ give
\begin{equation}
 d_e\ge c\Delta^{-2}L^{-3}t^{\gamma(s)},\qquad
 \gamma(s)=2/(1+s).
 \label{eq:lm:g:7.4}
\end{equation}

For larger connected sets, the three-coordinate bound alone does not give the desired degree range. We use more coordinates to improve the tail exponent, but keep their number fixed so that its constants do not grow with dimension. Coloring will control the remaining sum over arbitrarily large sets.

Fix an integer $J\ge3$, to be chosen after $\varepsilon$ and before $m$,
and put $a_j(s)=j/[1+(j-1)s]$. Any $j\le J$ selected coordinates with
correlations bounded in magnitude by $s$ satisfy
\[
 \Prob(P_{i_1},\ldots,P_{i_j}\le x)
 \le2^j e^{-a_j(s)c_x^2/2}
 \le6^Jc_x^Jx^{a_j(s)}.
\]
For each orthant, its signed sum has variance at most
$j[1+(j-1)s]$, and a zero-variance sum cannot exceed $jc_x>0$.
This proves the bound also for singular covariances.
In the kernel expansion \eqref{eq:lm:g:5.4}, retain $J$ coordinates when
$|T|\ge J$. For $1\le j=|T|<J$, keep the factor $x^{k-j}$:
$k-j+a_j(s)\ge a_J(s)$ for $k\ge J$, because the increments of $a_j$
are at most one. The empty term has exponent $k\ge J\ge a_J(s)$.
Consequently
\begin{equation}\label{lm:fixed_coordinate_tail}
 \|A_C\|\le6^J4^k c_x^J x^{a_J(r_C)}\quad(k\ge J).
\end{equation}
This estimate includes all atomic faces of the censored law.

For $3\le k<J$, the triple bound, \eqref{eq:lm:edge_tree_count}, and the
assigned-edge denominator give
\begin{equation}
 B_{\mathrm{s},<J}/\mathcal D_m\le t^{g-o(1)}.
 \label{eq:lm:g:boundary.4}
\end{equation}
Only fixed powers of $\Delta,K,L$ occur here. Indeed, dividing
\eqref{eq:lm:g:5.6} for a set assigned to an edge of magnitude $s$ by
\eqref{eq:lm:g:7.4} leaves $t^{\lambda(s)-\gamma(s)}$ times such powers,
and $\lambda(s)-\gamma(s)=(1-s)/[(1+2s)(1+s)]$ decreases in $s$, so it
is at least its value $g$ at $s=r_1$. For larger strong sets take
the minimum of \eqref{lm:fixed_coordinate_tail} and the coloring bound \eqref{eq:lm:g:5.5}. The fixed-coordinate bound supplies the needed correlation-dependent exponent. The coloring bound improves as the size of the connected set grows, making the sum over the largest sets convergent.
For an assigned edge of magnitude $s$, put $A=16\Delta$, $h=a_J(s)$,
$n=\lceil\chi h\rceil$, and $z=A x^{1/\chi}$. These four symbols are
local to this estimate: $A$ is unrelated to the centered measures $A_C$,
$h$ to $h_\alpha$, $n$ to the number of representatives, and $z$ to the
quantiles $z_e$. Consider the strong sets of size at least $J$ assigned
to this edge. The count
\eqref{eq:lm:edge_tree_count} bounds their total norm by
$C_Jc_x^J\Delta^{-2}\sum_{k\ge J}A^k\min(x^h,x^{k/\chi})$.
Here
\[
 \log z\le-[1-\beta-o(1)]L/\chi\longrightarrow-\infty.
\]
Once $z\le1/2$, enlarge the numerical sum to all $k\ge1$, sum through
$n$ using $x^h$ and above $n$ using $x^{k/\chi}$, to obtain
\[
 \sum_{k\ge J}A^k\min(x^h,x^{k/\chi})
 \le A^nx^h\{2+z/(1-z)\}\le3A^nx^h.
\]
This bound remains valid when $n<J$. Uniformly for $1\le h\le J$,
$\log A^n\le h\chi\log\chi+O_J(\chi+\log\chi)
\le\beta hL+o(L)$. For each assigned edge, divide the bound for its connected sets by that edge's deficit. Then multiply back by the nonnegative deficits and sum over edges. The powers $\Delta^{-2}$ cancel. The remaining fixed logarithmic factors contribute only a subpolynomial factor, giving
\begin{equation}
 \frac{B_{\mathrm{s},\ge J}}{\mathcal D_m}
 \le t^{\inf_{0\le s\le r_1}\{(1-\beta)a_J(s)-\gamma(s)\}-o(1)}.
 \label{eq:lm:g:boundary.5}
\end{equation}
There is no extra edge-count factor and no division by the deficit of a
different edge. The correlation-dependent index $n$ only sums a
numerical bound and does not change the rank cutoff $K$.

Let $f_J(s)=(1-\beta)a_J(s)-\gamma(s)$. Over the positive denominator
$[1+(J-1)s](1+s)$ its numerator is
$J(1-\beta)-2-(J-2+J\beta)s$. If $f_J(r_1)>0$, this numerator is positive
and decreasing throughout $[0,r_1]$, whereas its denominator increases;
hence its minimum is at $r_1$. Moreover
\[
 f_J(r_1)\uparrow(1-\beta)/r_1-2/(1+r_1)>0.
\]
Choose one finite $J$ making the gap to this limit smaller than
$\varepsilon/3$. Equations \eqref{eq:lm:g:boundary.4}--\eqref{eq:lm:g:boundary.5}
and the weak bound imply
$B_3/\mathcal D_m\le t^{\vartheta-\varepsilon/3-o(1)}$.
\paragraph{Pair products and rank restoration.}
\label{sec:lm:g:7.3}
We bound the remaining terms directly under
\eqref{eq:lm:g:boundary.1}. A product of pair terms will be small relative to the total deficit because one factor can be compared with that deficit, while the remaining pair norm tends to zero. Write \(\gamma_-=2/(1+r_1)>1\).
Equations \eqref{eq:lm:g:3.9} and \eqref{eq:lm:g:5.7}, together with
\(|E|\le m\Delta/2=\alpha\Delta/(2t)\), give
\[
 \frac{V_2}{\mathcal D_m}\le CK^2L^2=t^{-o(1)},\qquad
 V_2\le156|E|c_x^2x^{\gamma_-}
 \le C\alpha\Delta K^2L\,t^{\gamma_--1}
 =t^{\gamma_--1-o(1)}=o(1).
\]
Thus
\begin{equation}
 \frac{V_2^2}{\mathcal D_m}
 \le C\alpha\Delta K^4L^3t^{\gamma_--1}
 =t^{\gamma_--1-o(1)}.
 \label{eq:lm:g:7.9}
\end{equation}
For each one-edge reference, compare its truncated event with the
independent truncated event and then restore its late ranks. Summing
these individual comparisons proves, independently of the actual full
covariance law,
\begin{equation}
 \mathcal D_m\le V_2+2|E|T_m=o(1).
 \label{eq:lm:g:7.10}
\end{equation}
Indeed \(2|E|T_m\le C\alpha\Delta\chi t^{19}=o(1)\).
Now use the previously established bound
\(B_3/\mathcal D_m\le t^{\vartheta-\varepsilon/3-o(1)}\). Its exponent is positive for sufficiently large dimensions, so the bound \(\mathcal D_m=o(1)\) in \eqref{eq:lm:g:7.10}, which uses only
the one-edge laws, gives \(B_3=o(1)\). Together with the pair-norm estimate, this yields
\(e^{V_2+B_3}=1+o(1)\).
Using \eqref{eq:lm:g:7.2} only for late-rank restoration gives
\begin{equation}
 \frac{2|E|T_m}{\mathcal D_m}
 \le C\alpha\Delta\chi L^{-1}t^9=t^{9-o(1)}.
 \label{eq:lm:g:7.11}
\end{equation}
All fixed powers of \(\Delta,K,L\) are subpolynomial by
Lemma~\ref{mx:uniform_factors} and \eqref{eq:lm:g:boundary.3}; no
exponential counting factor has been discarded. The exponents
\(1,\gamma_--1,9\) exceed \(g\ge\vartheta\). The uncensored estimate
\eqref{eq:lm:g:4.4} also has a larger exponent. Absorb the remaining
uniform little-o losses and the finite sum into \(2\varepsilon/3\) in
\eqref{eq:lm:g:5.11}. This proves \eqref{eq:lm:g:boundary.2}, and
\eqref{eq:lm:g:1.6} follows from positivity of the deficits.

Finally, $(1-\beta)/r_1-2/(1+r_1)\ge g$ is equivalent to
$\beta\le(1-r_1)/(1+2r_1)$, proving the full-rate subrange.
\end{proof}

The method supplies no positive relative exponent at
$\beta=(1-r_1)/(1+r_1)$, and its constants are not uniform as $r_1$ tends
to one. A varying correlation cap therefore needs a separate argument.

\clearpage
\section{Uniform bounds for independent Gaussian blocks}
\label{mx:section}

We prove the mixed-block theorem, Theorem~\ref{main:mixed_theorem}, by bounding how much each original independent part contributes to the full rejection probability. This contribution is the expectation of its rejection share, defined below. We sum the bounds under the original full probability law; comparisons with other laws are used only to establish the bound for one part at a time. Throughout this section,
\[
 t=\alpha/m,\qquad L=\log(1/t),\qquad
 h_\alpha=\log(1/\alpha)-1+\alpha,\qquad 0<\alpha\le a_1<1.
\]
The constants may depend on the fixed upper bound for the level and the
stated size and correlation bounds, but not on the actual level,
covariance, or number of blocks. Upper bounds for the graph degrees are
prescribed in advance. All
statements about sufficiently large $m$ use one threshold for the whole
interval $(0,a_1]$. As in Appendix~\ref{sec:lm:graph},
$c_u=\Phi^{-1}(1-u/2)$, so a Gaussian p-value is at most $u$ exactly when
its absolute statistic is at least $c_u$.

For any set of original test coordinates $A$, let $N_A(u)$ be the number of its p-values at most $u$. If $R$ is
the largest BH crossing rank for the full vector, put
\begin{equation}\label{mx:shares}
 T_A=\frac{N_A(tR)}R\mathbf1\{R>0\},\qquad b_A=\Expect T_A,
 \qquad F_\Sigma(\alpha)=\sum_A b_A.
\end{equation}
When no hypotheses are rejected, define the random share to be zero. The sum is over a partition of the original coordinates, and every term
uses the same original law and the same $R$. Thus the random share is the
fraction of all rejections belonging to the selected part, and its
expectation is the part's contribution to the rejection probability.
The count at a largest positive crossing equals its rank, so $0\le T_A\le1$.
An independent uniform representative with $w$ identical physical copies
has expected share $wt$, even when the remaining p-values have an
arbitrary joint law. Indeed, force those $w$ copies to zero, condition on all remaining
p-values, and apply \eqref{eq:lm:g:2.1}. This observation supplies all
independent-reference equalities below.

The proofs compare either a full rejection probability or the expected
share of one part. The following table specifies both the quantity and
the probability law, since a bound for one does not automatically give
a bound for another.
\begin{center}
\small
\begin{tabular}{@{}>{\raggedright\arraybackslash}p{0.16\textwidth}>{\raggedright\arraybackslash}p{0.49\textwidth}>{\raggedright\arraybackslash}p{0.27\textwidth}@{}}
\toprule
Quantity & Definition and reference law & Use \\
\midrule
\(d_e\), graph expansion &
The deficit \(\alpha-F_m(\Gamma^e,w;\alpha)\) in the one-edge
Gaussian reference, with all other representatives independent and
all physical multiplicities retained. &
Expansion of the full rejection probability in~\eqref{eq:lm:g:boundary.2}. \\[0.5em]
\(d_e^V\), component share &
The deficit \(st-b_e\), where \(b_e\) is the expected share of
\(S\) with one edge retained inside \(S\) and the same actual
independent outside law \(V\). &
Original-part bound in Lemma~\ref{mx:lc}. Lemma~\ref{mx:capped} and~\eqref{gb:edge_deficit} denote this share deficit by \(d_e\), with a selected block \(A\) in place of \(S\); it is not the \(d_e\) of the first row. \\[0.5em]
\(c_NtG_t(d_C)\), terminal cluster &
A lower bound for \(|C|t-b_C\) for an unchanged isolated near-copy
cluster under the single terminal law; \(c_N\) and \(G_t\) are as in
Lemma~\ref{mx:near}, and \(d_C\) is the spread of \(C\). &
Lemma~\ref{fd:terminal} and the comparison in~\eqref{fd:margin}. \\
\bottomrule
\end{tabular}
\end{center}
The fixed-degree argument leading to \eqref{fd:margin} compares full
rejection probabilities under the original and terminal laws. After
clusters with larger rescaled spread are isolated and those with
sufficiently small rescaled spread are contracted, Lemma~\ref{fd:terminal} bounds all remaining cluster and
capped-graph shares under that one terminal law. Its graph deficit
\(D_G^{V,\mathrm{term}}\) retains the actual terminal independent
outside law; it is the deficit \(D_S^V\) of Lemma~\ref{mx:lc} for
\(S=G\), computed under the terminal law rather than the original law. In Theorem~\ref{main:mixed_theorem},
the bounds for the original independent parts are added under the
original law. Within this appendix, Lemma~\ref{mx:void} supplies a
background event of positive probability, Lemma~\ref{mx:lc} bounds the
expected share of the capped component,
Lemmas~\ref{mx:near}--\ref{mx:arbitrary_block} bound the expected share of
one bounded block, first for nearly identical coordinates and then for an
arbitrary covariance, and Section~\ref{mx:coordinate} shows that the
expected share of a single coordinate can exceed its nominal value.
Corollary~\ref{mx:prds} derives from this example the failure of both PRDS conditions. The
terminal law of the fixed-degree argument and the third row of the table
are needed only in Appendix~\ref{fd:section}.

\subsection{A positive-probability background event}
\label{mx:void_section}

We need a background event on which inserting a bounded number of
p-values cannot cause a crossing at a later rank. The following lemma
gives this event directly under the actual background distribution;
it does not assume a Simes inequality for that distribution.

\begin{lemma}[A uniform background event]\label{mx:void}
Fix integers $W,B,Q\ge1$, $0<\rho_0<1$, and $0<a_1<1$. Let a background
$V$, containing at most $m$ physical p-values, consist of a Gaussian
component with multiplicities at most $W$, representative correlations
bounded by $\rho_0$, and degree at most $\Delta_m$, together with mutually
independent blocks of size at most $B$ with superuniform marginals, $\Prob(V_i\le u)\le u$. The
Gaussian component is independent of all those blocks. Assume that a prescribed upper bound $M_m$ for its representative
edge count satisfies
\begin{equation}\label{mx:void_edge_envelope}
 M_m=o\!\left(m^\tau[\log(em)]^d\right),\qquad
 \tau=2/(1+\rho_0)\in(1,2),\quad d=(2-\tau)/2>0.
\end{equation}
There are a fixed integer $L_0\ge Q$,
a constant $\gamma>0$, and a common dimension threshold such that, simultaneously over
$0<\alpha\le a_1$, the event
\begin{equation}\label{mx:void_event}
 N_V(jt)=0\quad(1\le j\le L_0),\qquad
 N_V(jt)+Q<j\quad(L_0<j\le m)
\end{equation}
has probability at least $\gamma$.
\end{lemma}

\begin{proof}
We first construct the event under an independent-block reference law,
then transfer the event for early ranks to the actual law and bound crossings at the remaining
ranks separately. Replace only the Gaussian representative covariance by
the identity, keeping all multiplicities and every other block law unchanged.
The reference consists of independent blocks of size at most
$R_0=\max(W,B,Q,2)$. Write $\beta=(1+a_1)/2$. Choose $L_0\ge R_0$ so large
that
\begin{equation}\label{mx:void_reference_tail}
 \frac{\beta^{-1}e^{-h_\beta(L_0+1)/R_0}}
             {1-e^{-h_\beta/R_0}}\le\frac12.
\end{equation}
Let $A$ be the event that no p-value is at most $L_0t$. By a union bound within
each block and independence between blocks, for large $m$,
\[
 \Prob_0(A)\ge\prod_b(1-|b|L_0t)\ge e^{-2a_1L_0}.
\]
Conditioning each block on its own part of $A$ preserves block independence
and bounds each conditional marginal at $jt$ by
$jt/(1-R_0L_0t)$. Choose $m$ so that
$a_1/(1-R_0L_0a_1/m)\le\beta$.
For independent blocks of size at most $R_0$, convexity gives
$\Expect_0[e^{sN_V(jt)}\mid A]\le\exp\{\beta j(e^{sR_0}-1)/R_0\}$, because the
conditional mean of $N_V(jt)$ is at most $\beta j$. Taking
$s=\log(1/\beta)/R_0$, and allowing $Q\le R_0$ inserted zeros, gives
\[
 \Prob_0\{N_V(jt)+Q\ge j\mid A\}
 \le\beta^{-Q/R_0}e^{-h_\beta j/R_0}
 \le\beta^{-1}e^{-h_\beta j/R_0}.
\]
Thus the complete reference event in \eqref{mx:void_event} has probability
at least $\gamma_0=\tfrac12e^{-2a_1L_0}$.

To transfer the event to the actual law, we retain the labels and values of all uncensored p-values. Comparing only their total count would discard this information. For $n$ Gaussian representatives set
$\Xi=\sum_i\mathbf1\{P_i\le x\}\delta_{(i,P_i)}$, and let $\Xi^0$ be
its independent-uniform counterpart. In event total variation,
\begin{equation}\label{mx:absolute_tv}
 d_{\rm TV}(\mathcal L(\Xi),\mathcal L(\Xi^0))
 \le(2n+2|E|)x^2+2\sum_{\{i,j\}\in E}\Prob(P_i,P_j\le x).
\end{equation}
We compare both laws with an auxiliary Poisson process, so that the error can be expressed using only a label and its graph neighbors. Let $\Pi$ be the Poisson process with intensity
$dp$ on each label $i$, $0\le p\le x$. The process below adds points at this intensity and removes each existing point at rate one; the Poisson law is stationary under these two operations. Its immigration--death generator is
\[
 \mathcal A h(\xi)=\sum_i\int_0^x[h(\xi+\delta_{(i,p)})-h(\xi)]\,dp
 +\int[h(\xi-\delta_z)-h(\xi)]\,\xi(dz).
\]
The needed estimate concerns the change caused by adding two specified points, beyond the two separate one-point changes. Coupling the processes makes this second difference vanish as soon as either added point has disappeared.

For $0\le f\le1$, the integrated semigroup solution of
$\mathcal Ah=f-\Expect f(\Pi)$ has second difference at most
$\int_0^\infty2e^{-2s}\,ds=1$ in absolute value: couple the processes
by retaining the two added points with independent survival probability
$e^{-s}$ and use a common immigration process. The integral defining the
solution converges by the same coupling, with bound
$(|\xi|+nx)e^{-s}$; finite-time integration of the generator justifies
its limiting identity. This argument also allows multiple Poisson points
at a label.
Remove from $\Xi$ the points at label $i$ and at every graph neighbor of $i$ to obtain $\Xi_i$.
The mark at $i$ is independent of the entire $\Xi_i$, since all Gaussian
cross-covariances to that subvector vanish, also for singular laws.
Replacing the configuration in the immigration term by $\Xi_i$ costs
$(n+2|E|)x^2$. In the death term the point at $i$ is already removed, so
the cost is $2\sum_e\Prob(P_i,P_j\le x)$. The remaining terms cancel by
uniform marginal intensity and independence. Comparison of $\Xi^0$ with
$\Pi$ costs a further $nx^2$, proving \eqref{mx:absolute_tv}.

We next control the actual background at late ranks by its count variance. This requires a bound on the joint probability of small p-values for each pair joined by an edge in the covariance graph. We retain the logarithmic factor in that bound because it permits the stated edge-count envelope. For
$c=c_u\ge1$, write the two coordinates in the positive quadrant as $c+a,c+b$ with $a,b\ge0$.
There, for correlation $r$, the exponent of the bivariate normal density is
\[
 \frac{c^2}{1+r}+\frac{c(a+b)}{1+r}
 +\frac{a^2-2rab+b^2}{2(1-r^2)}.
\]
The last term is nonnegative. Integrating the product exponential and
summing the four quadrants, each of which reduces to this one with $r$
replaced by $\pm r$, yields, uniformly for $|r|\le\rho_0$,
\[
 \Prob(P_i,P_j\le u)\le
 \frac{2(1+\rho_0)^2}{\pi\sqrt{1-\rho_0^2}}
 c^{-2}e^{-c^2/(1+\rho_0)}
 \le\frac{C u^\tau}{[\log(e/u)]^d}.
\]
For the last step, the lower Mills bound gives
$e^{-c^2/2}\le\sqrt{2\pi}\,cu$ and $\log(e/u)\le3c^2$.
Increase the fixed constant for $c<1$. With physical multiplicities at
most $W$ and independent blocks of size at most $B$, this implies
\[
 \Expect N_V(\alpha v/m)\le a_1v,\qquad
 \operatorname{Var}N_V(\alpha v/m)
 \le C\left\{v+\frac{M_m(v/m)^\tau}{[\log(em/v)]^d}\right\}.
\]
Indeed $\operatorname{Var}Y\le B\,\Expect Y$ for a block count $Y\le B$,
and each covariance between the two threshold indicators on an edge is at most the probability that both indicators equal one, with
fixed multiplicity factors. Uniformity in the level follows directly
from $(\alpha y)^\tau/[\log(e/(\alpha y))]^d
\le y^\tau/[\log(e/y)]^d$ for $\alpha\le1$.
A union bound over every individual rank would waste the variance estimate. Instead, group ranks into geometric intervals. A crossing anywhere in one interval forces a count excess at its right endpoint, where Chebyshev's inequality applies.

Choose fixed $1<z<1/a_1$ and put $\eta=1-a_1z>0$. On a geometric rank
interval $(Kz^k,\min(Kz^{k+1},m)]$, a crossing after insertion of $Q$
zeros implies a count deviation at its right endpoint of at least
$\eta Kz^k/2$, once $K\ge2Q/\eta$. Chebyshev's inequality and the
geometric sum, using $1<\tau<2$, give
\begin{equation}\label{mx:void_large_connected}
 \Prob\{N_V(jt)+Q\ge j\text{ for some }j>K\}
 \le C_1/K+
 \frac{C_2M_mm^{-\tau}K^{\tau-2}}{[\log(em/K)]^d}.
\end{equation}
To justify the logarithmic geometric sum, put $\Lambda=\log(em/K)$
and $\kappa=\log z$. Its terms, after removing $K^{\tau-2}$, are
$e^{-(2-\tau)\kappa k}/(\Lambda-\kappa k)^d$ with $\Lambda-\kappa k\ge1$. For $\kappa k\le\Lambda/2$,
sum a geometric series with denominator at least $(\Lambda/2)^d$.
For $\kappa k>\Lambda/2$, the remaining sum is at most
$C e^{-(2-\tau)\Lambda/2}\le C'\Lambda^{-d}$. Replacing an interval endpoint
by its left endpoint costs only a fixed factor since their ratio is
at most $z$ and $\log(em/v)\ge1$.
Choose the fixed cutoff
\begin{equation}\label{mx:void_cutoff}
 K\ge\max\{L_0,2Q/\eta,8C_1/\gamma_0\}.
\end{equation}
The second term in \eqref{mx:void_large_connected} tends to zero by
\eqref{mx:void_edge_envelope}. At $x=Ka_1/m$, the joint-tail term in \eqref{mx:absolute_tv} is
$O(M_mm^{-\tau}/[\log m]^d)=o(1)$. Its remaining terms are
$O(m^{-1}+M_mm^{-2})=o(1)$, since $\tau<2$.
Censor further at $K\alpha/m$, repeat each representative according to its physical multiplicity, and adjoin the unchanged independent blocks. None of these operations increases total variation. The reference event through rank $K$ therefore transfers with error $o(1)$, uniformly in $\alpha$. Subtracting the probability of a later crossing under the actual law proves
\eqref{mx:void_event} with $\gamma=\gamma_0/2$, after increasing the common dimension threshold.

For later relative estimates we retain the stronger exponential count
tail supplied by coloring. With $c_m=\max\{W(\Delta_m+1),B,Q\}$ (in this proof only), color
representatives properly and copies separately, and use distinct colors
within each other independent block. H\"older and the Bernoulli moment
generating function give
\begin{equation}\label{mx:insertion_tail}
 \Prob\{N_V(jt)+Q\ge j\}\le\alpha^{-Q/c_m}e^{-h_\alpha j/c_m}.
\end{equation}
For $j\ge Q/(1-a_1)$ its right side increases in $\alpha$ and is at most
$a_1^{-1}e^{-h_{a_1}j/c_m}$. The polynomial tail used to construct the
positive-probability event is not substituted for this relative-error
bound.
\end{proof}

Under the degree condition \eqref{eq:lm:g:boundary.1} at the logarithmic boundary,
$|E|=O(m\log m)$, so \eqref{mx:void_edge_envelope} holds for every fixed correlation cap below one. More generally, the background lemma permits
$\Delta_m=o(m^{\tau-1}[\log(em)]^d)$. This larger range applies only to the positive-probability background event, not to the later relative expansion. When a bounded selected block is added to the background, we may therefore increase the fixed correlation cap of the background's Gaussian component to cover that block. The graph expansion is still applied only to the original component, at its original cap.

\subsection{The contribution of the capped component}
\label{mx:lc_section}

For the expected share of the capped component, with its independent outside law kept fixed, Lemma~\ref{mx:lc} gives a relative bound of the form in Theorem~\ref{thm:lm:graph_boundary}, with one-edge share deficits in place of deficits in the full rejection probability.

\begin{lemma}[A capped component with an independent background]
\label{mx:lc}
Let $S$ be a Gaussian component with physical size $s\le m$,
multiplicities at most $W$, correlation cap $0<\rho_0<1$, and
representative covariance graph of maximum degree at most $\Delta_m$,
where the prescribed $\Delta_m$ satisfy \eqref{eq:lm:g:boundary.1} for some
$0\le\beta<(1-\rho_0)/(1+\rho_0)$. Let its independent outside vector
consist of $m-s$ physical coordinates in arbitrary independent blocks of
superuniform p-values, each of size at most $B$. The parameters
$W,B,\rho_0,a_1,\beta$ and the sequence $(\Delta_m)$ are fixed. Keep that outside law fixed. Let $b_e$ be the expected
share of $S$ when its representative covariance retains just edge $e$, and put
\[
 d_e^V=st-b_e,\qquad D_S^V=\sum_e d_e^V,\qquad
 g_0=\frac{1-\rho_0}{(1+2\rho_0)(1+\rho_0)}.
\]
Put $\vartheta_0=\min\{g_0,(1-\beta)/\rho_0-2/(1+\rho_0)\}>0$.
Each nonzero edge has positive $d_e^V$ eventually. For every
$0<\varepsilon<\vartheta_0$, uniformly for $0<\alpha\le a_1$ and all allowed
laws, one common dimension threshold, depending only on
$W,B,\rho_0,a_1,\beta,\varepsilon$ and $(\Delta_m)$, gives
\begin{equation}\label{mx:lc_relative}
 |b_S-st+D_S^V|\le t^{\vartheta_0-\varepsilon}D_S^V.
\end{equation}
In particular $b_S<st$ when there is an edge. If there are no edges,
$b_S=st$ in every dimension and at every level.
\end{lemma}

\begin{proof}
We adapt the graph expansion to the expected share of this component,
keeping its independent outside law unchanged. The independent baseline
is exactly $st$ by the forcing argument following \eqref{mx:shares}. For a single edge, condition on every
p-value outside its two endpoints. Folded bivariate normal regression
implies that the conditional expectation of an increasing function of
the other endpoint is nondecreasing in the first p-value. In particular,
forcing that endpoint to zero makes each event $\{R^{[i]}\le k\}$
increasing in the other endpoint. Consequently every increment in
\eqref{eq:lm:g:2.2} is nonnegative. Independent representatives elsewhere
in $S$ have their exact nominal expected shares. This proves $d_e^V\ge0$
without any further restriction on dependence within the outside blocks.

For strictness and a quantitative bound, in a single-edge reference all coordinates
other than the copies of the two endpoints form independent blocks of size at most
$\max(W,B)$. The independent-reference part of the proof of
Lemma~\ref{mx:void}, with insertion allowance $2W$, gives a fixed positive
probability $\gamma_0$ of no outside crossing after inserting the pair.
For endpoint multiplicities $u,v$, $q_e=u+v$, and $r=|\Gamma_{ij}|$,
\begin{equation}\label{mx:edge_lower}
 d_e^V\ge\gamma_0\{J_r(ut,q_et)+J_r(vt,q_et)-J_r(q_et,q_et)\}.
\end{equation}
This lower bound comes from the increments up to the rank where all
copies of both pair members are counted. All other increments are
nonnegative. The curvature
calculation in \eqref{eq:lm:g:3.1}--\eqref{eq:lm:g:3.3}, with $p_0$
replaced by $\gamma_0$, together with the arguments for
\eqref{eq:lm:g:3.7}--\eqref{eq:lm:g:3.9}, which also use
\eqref{eq:lm:g:3.4} and the cutoff $K\ge2W\ge q_e$ defined below, therefore gives
\begin{equation}\label{mx:edge_estimates}
 d_e^V\ge c r^2Lt^2,\qquad
 d_e^V\ge c r^2L^{-1}t^{2/(1+r)},\qquad
 \|A_e^{Kt}\|/d_e^V\le CK^2L^2.
\end{equation}
These constants are uniform in the actual outside block laws. The
superuniform assumption is enough: the reference avoidance probability,
its conditional marginal bounds, and the coloring moment bound use only
$\Prob(V_i\le u)\le u$. Exact nominal shares belong to the selected independent
Gaussian representatives and hold for an arbitrary independent outside law.

Let $c_m=\max(\Delta_m+1,B)$ in this proof, and
\[
 K=\max\left(2W,\left\lceil20Wc_mL/h_\alpha\right\rceil\right),
 \qquad x=Kt.
\]
Color the representatives and use separate colors within each outside
block. Applying the weighted color bound of
\eqref{eq:lm:g:5.8} gives
\begin{equation}\label{mx:lc_late}
 T_m:=\sum_{j>K}\Prob\{N(jt)\ge j\}
       \le Cc_mt^{20}.
\end{equation}
This bound holds under the actual law, the independent baseline and every
single-edge reference, because the same coloring is proper for each of
their covariance graphs.
Lemma~\ref{mx:uniform_factors} applies to this cutoff, with $\tilde\chi=Wc_m\le WB(\Delta_m+1)$;
in particular $K<m$, $\log K=o(L)$, and $L\le c_x^2\le3L$ uniformly.

Let $R_K$ be the largest crossing among ranks at most $K$, and use the
clipped share $T_S^K=\min(1,N_S(tR_K)/R_K)$ for $R_K>0$, with value zero
otherwise. Clipping keeps this random share bounded by one even when the largest
crossing was omitted. It equals $T_S$ whenever there is no later crossing.
Average this clipped share over the fixed outside law, and expand only the censored Gaussian law of $S$. For each censored configuration of $S$, the average is a measurable function with values in $[0,1]$. Because $S$ is independent of the outside and that outside law is unchanged, the same function is integrated under the actual Gaussian law, the independent baseline, and each single-edge reference. Only the measure changes. We can therefore bound each signed-measure remainder by its variation norm.
With $V_2=\sum_e\|A_e\|$ and
$B_3=\sum_{C\subseteq S\text{ connected},\,|C|\ge3}\|A_C\|$ for
these measures censored at $x=Kt$,
\begin{equation}\label{mx:lc_remainder}
 |b_S-st+D_S^V|\le\tfrac12(B_3+V_2^2/2)e^{V_2+B_3}+2|E|T_m.
\end{equation}
The factor $1/2$ uses zero total mass of every signed remainder and the
bound $0\le f\le1$ for the common outside-averaged share function. The coefficient $2|E|$
is the sum of absolute truncation coefficients in
$b_S-\sum_e b_e+(|E|-1)st$.

We split according to the largest internal correlation for the same
reason as in the proof of Theorem~\ref{thm:lm:graph_boundary} at
\eqref{eq:lm:g:7.1}--\eqref{eq:lm:g:7.2}. If all correlations are too small, even a very small absolute truncation error could exceed the deficit, so the full-law comparison is needed. Otherwise a coarse positive lower bound for the deficit is enough to restore late ranks.

For $r_{\max}\le t^4$, apply the full folded-law cancellation
\eqref{eq:lm:g:4.1} directly to the outside-averaged full share, without
truncation. If $H=2\sum_e r_e^2$, the quantity $S$ of \eqref{eq:lm:g:4.1}, then $\sqrt H\le mt^4$ and
$\|\Gamma-I\|_{\rm op}\le mt^4\le t^3$. Equations
\eqref{eq:lm:g:4.1} and \eqref{mx:edge_estimates} bound the relative error
by $C\sqrt H/(Lt^2)\le C\alpha t/L$.

Otherwise $D_S^V\ge cLt^{10}$. Apply the centered-measure proof of
Theorem~\ref{thm:lm:graph_boundary} to the outside average of the clipped share $T_S^K$, the function of the
censored values of $S$ used in \eqref{mx:lc_remainder}. The measures concern only $S$, and each of their norms
bounds its integral against the common function in $[0,1]$.
The count of connected sets containing a fixed edge \eqref{eq:lm:edge_tree_count}, weak sum
\eqref{eq:lm:g:6.7}, and fixed-coordinate bound
\eqref{lm:fixed_coordinate_tail} are unchanged. Every edge denominator
is now \eqref{mx:edge_estimates}, under the same actual outside $V$.
In particular, the large strong-set sum is normalized by its own
$d_e^V$, giving \eqref{eq:lm:g:boundary.5} with $r_1=\rho_0$; it is
not divided by a one-edge deficit in the full rejection probability. The degree cost is retained
as $\beta a_J(s)$ with $J$ chosen after $\varepsilon$ and before $m$.

For clarity, the relative bounds for weak sets, smaller strong sets, larger strong sets, products of pairs, and
late ranks, the analogues of \eqref{eq:lm:g:6.8},
\eqref{eq:lm:g:boundary.4}, \eqref{eq:lm:g:boundary.5}, \eqref{eq:lm:g:7.9}
and \eqref{eq:lm:g:7.11}, are
\[
 t^{1-o(1)},\quad t^{g_0-o(1)},\quad
 t^{f_J(\rho_0)-o(1)},\quad
 t^{\gamma_--1-o(1)},\quad t^{9-o(1)},
 \qquad\gamma_-=2/(1+\rho_0).
\]
Here $f_J(s)=(1-\beta)J/[1+(J-1)s]-2/(1+s)$ and a fixed sufficiently
large $J$ approaches the positive limiting exponent in $\vartheta_0$.
All remaining factors are fixed polynomials in $\Delta_m,c_m,K,L$,
and hence subpolynomial uniformly in $\alpha$.
Separately, $V_2=o(1)$ and $D_S^V\le V_2+2|E|T_m=o(1)$, by summing
single-edge comparisons and their two truncation errors. Thus $B_3=o(1)$
and the exponential in \eqref{mx:lc_remainder} tends to one. The exponents $1$, $\gamma_--1$ and $9$ exceed $g_0$. For sufficiently large dimensions, the allowed loss in the exponent absorbs the remaining factors and proves
\eqref{mx:lc_relative}. No step assumes Simes validity for the outside law. The denominator is the component-share deficit with the actual outside distribution retained; it is generally different from a deficit for the full rejection probability.
\end{proof}

\subsection{Uniform estimates for a selected bounded block}
\label{mx:selected_section}

In this subsection, a Gaussian block $A$ of physical size $s\le Q$ is
independent of its entire outside vector. The latter consists of a capped
Gaussian component and of independent uniform-marginal blocks of size
at most $B$. The capped component has multiplicities at most $W$, a fixed
correlation cap below one and representative degree at most $\Delta_m$,
where the prescribed $\Delta_m$ satisfy the logarithmic-boundary
envelope \eqref{eq:lm:g:boundary.1}. The capped component may be
empty. The outside law is retained in every comparison of expected shares.
We first bound the expected share of a block whose coordinates are close to signed copies. We then show how replacing such a cluster by exact copies changes the share of the entire selected block. Finally, we partition an arbitrary block into finitely many clusters and use these comparisons to reduce it to representatives whose correlations are bounded away from perfect dependence. Set
\begin{equation}\label{mx:common_cutoff}
 c_m=\max\{W(\Delta_m+1),B,Q,2\},\quad
 K=\max\left(2Q,\left\lceil21c_mL/h_\alpha\right\rceil\right).
\end{equation}
Writing $z_*=c_{a_1}>0$, the following estimates hold for all sufficiently large dimensions, uniformly over the allowed levels:
\begin{equation}\label{mx:quantiles}
\begin{aligned}
 &K\le Cc_m(1+\log m)<m,\qquad \log K=o(L),\\
 &L\le c_{qt}^2,c_{Kt}^2\le3L\ (1\le q\le Q),\qquad
 c_t/c_{Kt}\longrightarrow1.
\end{aligned}
\end{equation}
Also $c_{jt}^2\ge L/2$ for $j\le\sqrt m$ and $c_{jt}\ge c_\alpha$
for all $j\le m$. These follow from Mills' inequalities and
Lemma~\ref{mx:uniform_factors}, with $c_m$ in place of
$W\chi$ and with the cutoff \eqref{mx:common_cutoff}; the assertion for
$j\le\sqrt m$ uses
$\log(1/(jt))\ge L/2$ and the logarithmic Mills correction.

Let $q_j=\min\{1,\alpha^{-Q/c_m}e^{-h_\alpha j/c_m}\}$. The coloring argument behind \eqref{mx:insertion_tail} allows empty color
classes, so \eqref{mx:insertion_tail} also applies to this outside vector
with the present $c_m$. Hence the numbers $q_j$ bound the probability that
the outside count can support a crossing after at most $Q$ zero p-values
are inserted. Differentiation as in \eqref{mx:insertion_tail} shows that, for
$j\ge J_0=\lceil Q/(1-a_1)\rceil$, they are at most
$a_1^{-1}e^{-h_{a_1}j/c_m}$. Hence
\begin{equation}\label{mx:count_sums}
 \sum_jq_j\le Cc_m,\qquad
 \sum_{j>J}q_j\le Cc_me^{-h_{a_1}J/c_m}\quad(J\ge J_0),
\end{equation}
whereas the sharper level-dependent estimate is
\begin{equation}\label{mx:rank_tail}
 \sum_{j>K}j q_j\le
 C\{c_m(K+1)+c_m^2\}\alpha^{-Q/c_m}e^{-h_\alpha K/c_m},
 \qquad
 \alpha^{-Q/c_m}e^{-h_\alpha K/c_m}\le t^{20}.
\end{equation}
The geometric prefactors use $h_\alpha\ge h_{a_1}$, and the last inequality
uses $Q/c_m\le1$ and $e^{-21L}=t^{21}$.

\subsubsection{Conditional shares and a monotone comparison}\label{sec:lm:n:2}

We now supply the conditional calculation used for a selected near-copy
block. It has three steps: express the share through one physical anchor,
compare with monotone one-sided functions, and quantify a loss in that
comparison. Lemma~\ref{mx:near} then shows that this loss, constructed in
Section~\ref{sec:lm:n:5}, exceeds the error from taking absolute values,
bounded in Section~\ref{sec:lm:n:3}, uniformly in the level and as the
spread tends to zero. All statements up to the final bounds are finite-dimensional
identities or inequalities. The testing level can be any
\(0<\alpha\le a_1\). Write \(N(u)=\#\{i:mP_i\le u\}\) and

\[
 R=\max(\{1\le r\le m:N(\alpha r)\ge r\}\cup\{0\}).
\]

When \(R>0\), \(N(\alpha R)=R\): a larger count would itself give a larger crossing. Thus the rejection indicator is the sum of the \(m\) random rejection shares: each rejected coordinate contributes the reciprocal of the total number rejected. If only physical coordinate \(i\) is forced to zero, with resulting largest crossing rank \(H_i\ge1\), restoration gives

\begin{equation}
 F_m=\sum_{i=1}^m \Expect\frac{\mathbf1\{mP_i\le\alpha H_i\}}{H_i}. \label{eq:lm:n:2.1}
\end{equation}

If a coordinate is rejected, it was already counted at every later threshold, so forcing it to zero cannot produce a later crossing. Conversely, restoring it below the forced threshold preserves that crossing and cannot increase the rank. This proves \eqref{eq:lm:n:2.1} with every physical tie retained. The identity holds pathwise, before expectations are taken, and therefore also for a vector in which one coordinate has already been forced to zero.

If \(H:(0,\infty)\to\{1,\ldots,m\}\) is nonincreasing, then

\begin{equation}
 I(H):=\int_0^\infty\frac{\mathbf1\{p\le\alpha H(p)\}}{H(p)}\,dp\le\alpha. \label{eq:lm:n:2.2}
\end{equation}

Indeed, the values of \(p\) selected by the indicator form an initial interval with endpoint \(a\). The nonincreasing integer-valued \(H\) has a terminal value \(h\) immediately to the left of \(a\). Then \(a\le\alpha h\) and \(H(p)\ge h\) before \(a\), so \(I(H)\le a/h\le\alpha\). Endpoint values do not change the integral.

The nonincreasing property also gives a strict loss when the forced rank drops before the nominal endpoint. We need to control both possible sources of this loss: selection may end early, or it may continue with a smaller terminal rank.

For the same nonincreasing $H$, a useful quantitative version applies when \(q\ge2\). If \(H\le q\), \(H=q\) on \((0,\ell)\), and \(H\le q-1\) on \((\ell,\infty)\), with \(0<\ell<\alpha q\), then

\begin{equation}
 \alpha-I(H)\ge\min\left\{\alpha-\ell/q,\ \ell/[q(q-1)]\right\}. \label{eq:lm:n:2.3}
\end{equation}

If the terminal value is \(h=q\), selection ends at \(\ell\) and \(I(H)=\ell/q\). Otherwise \(h\le q-1\); selection includes \((0,\ell)\), and splitting the integral there gives

\[
 I(H)\le\ell/q+(a-\ell)/h
 \le\alpha-\ell(1/h-1/q)
 \le\alpha-\ell/[q(q-1)].
\]

Fix a selected near-copy block \(C\) of size \(q\le Q\), independent of the entire vector of coordinates outside it, and fix an anchor \(i\in C\). Use sign-aligned coordinates within \(C\). Put \(r_j=\rho_{ij}\ge15/16\) for its other physical coordinates. Scalar Gaussian regression, valid also for singular covariances, represents the conditional block given \(mP_i=p\) as

\begin{equation}
 |Z_j|\overset d=|r_jz_m(p)+\varepsilon_j|,
 \qquad z_m(p)=c_{p/m},\quad0<p<m, \label{eq:lm:n:2.4}
\end{equation}

where \(\varepsilon\) is centered Gaussian, independent of \(p\) and of the entire outside vector, with \(\operatorname{Var}\varepsilon_j=1-r_j^2\). The two anchor signs are absorbed by central symmetry of the joint residual law. The density of the scaled anchor \(p\) is exactly \(1/m\); its two signs have already been included.

Keep the outside p-values fixed and force only the anchor to zero. Let \(H^a(p)\) use the actual sibling scaled values

\[
 q_j^a(p)=2m\barPhi(|r_jz_m(p)+\varepsilon_j|).
\]

The absolute value can make a sibling's conditional p-value change direction as the anchor varies. Removing that absolute value gives a monotone comparison, whose integral has the bound just proved. Let \(H^+(p)\) instead use the comparison values

\begin{equation}
 q_j^+(p)=2m\barPhi(r_jz_m(p)+\varepsilon_j). \label{eq:lm:n:2.5}
\end{equation}

The comparison values may exceed \(m\), which is allowed because we compare them only with thresholds at most \(\alpha m<m\). Each \(q_j^+\) is increasing in \(p\), and hence \(H^+\) is nonincreasing. Extend \(H^+\) constantly beyond \(m\). Its selected integrand is zero for \(p>\alpha m\), so \(I(H^+)\le\alpha\) for every realization of the residuals and the outside p-values.

Writing \(f^a(p)=\mathbf1\{p\le\alpha H^a(p)\}/H^a(p)\) and \(f^+(p)=\mathbf1\{p\le\alpha H^+(p)\}/H^+(p)\), exact conditional integration of the terms in \eqref{eq:lm:n:2.1} belonging to \(C\) gives

\begin{equation}
 b_C=\frac1m\sum_{i\in C} \Expect\int_0^m f_i^a(p)\,dp. \label{eq:lm:n:2.6}
\end{equation}

We next control the positive part of \(f^a-f^+\), which is the amount by which the actual share can exceed the comparison share. The comparison requires only the monotonicity of the one-sided conditional functions constructed above.

\subsubsection{The error from taking absolute values}\label{sec:lm:n:3}

The actual two-sided threshold indicator can differ from the one-sided comparison only when a sibling crosses the negative Gaussian threshold. To bound the resulting share error, we also require the outside count needed for a rejection at that rank. Keeping both events in the bound ensures that the error remains small as the block approaches exact copies. Put \(z_s=c_{\alpha s/m}\), and let \(N_o(u)\) count outside
physical p-values at scaled threshold \(u\). Write \(B_s\) for any bound
on \(\Prob\{N_o(\alpha s)\ge s-q\}\). For the mixed background in this
section, we take \(B_s=q_s\), with \(q_s\) defined just before
\eqref{mx:count_sums}. This bound is uniform in the actual testing level.

Since \(q_j^a\le q_j^+\), one has \(H^a\ge H^+\). A positive difference \(f^a-f^+\) requires \(H^a=s>H^+\), \(p\le\alpha s\), and at least one sibling satisfying

\begin{equation}
 r_jz_m(p)+\varepsilon_j\le-z_s. \label{eq:lm:n:3.2}
\end{equation}

When both indicators in the shares equal one, the larger denominator only reduces the share. If no opposite-sign sibling is counted at \(\alpha s\), the actual and comparison counts there agree, so the comparison also crosses at \(s\), a contradiction. Moreover a crossing at \(s\) requires \(N_o(\alpha s)\ge s-q\). The positive difference is at most \(1/s\).

The outside count event is independent of the anchor residuals. Thus, by summing over the crossing rank s and the sibling with the opposite sign, the expected positive integral error is at most

\begin{equation}
 \sum_{j\ne i}\sum_{s=1}^m\frac{B_s}{s}
 \int_0^{\alpha s}\Prob\{r_jz_m(p)+\varepsilon_j\le-z_s\}\,dp. \label{eq:lm:n:3.3}
\end{equation}

The integral is exactly

\begin{equation}
 2m\Prob\{Z_i\ge z_s,\widetilde Z_j\le-z_s\}, \label{eq:lm:n:3.4}
\end{equation}

for a standard Gaussian pair with correlation \(r_j\). The factor \(2\) is the Jacobian factor in \(dp=-2m\phi(z)\,dz\); the anchor signs were already included in the conditional representation. If \(r_j=1\) the event is impossible. Otherwise its occurrence forces \(Z_i-\widetilde Z_j\ge2z_s\); that difference has variance \(2(1-r_j)\). The normal hazard lower bound \(\phi(x)/\barPhi(x)\ge x\) yields, for \(a\ge z>0\),

\[
 \barPhi(a)/\barPhi(z)\le e^{-(a^2-z^2)/2}.
\]

Consequently, if \(\delta=\max_{k,l\in C}(1-\rho_{kl})>0\),

\begin{equation}
 2m\Prob\{Z_i\ge z_s,\widetilde Z_j\le-z_s\}
 \le\alpha s\exp\{-a(\delta)z_s^2\},\qquad
 a(\delta)=\delta^{-1}-\tfrac12. \label{eq:lm:n:3.5}
\end{equation}

Combining \eqref{eq:lm:n:3.3}--\eqref{eq:lm:n:3.5}, every anchor in that block has an expected integral of the positive error at most

\begin{equation}
 \mathcal E_m(\delta):=
 \alpha(Q-1)\sum_{s=1}^m B_s e^{-a(\delta)z_s^2}. \label{eq:lm:n:3.6}
\end{equation}

\subsubsection{A positive loss at two physical anchors}\label{sec:lm:n:5}

Let \(E\) be the background event in Lemma~\ref{mx:void}, with
insertion allowance \(Q\). It has probability at least \(\gamma>0\),
is independent of the anchor residuals, and makes the largest crossing rank equal to \(q\) when
\(q\le Q\) zeros are inserted. Inserting any other set of \(q\) values gives crossing counts no larger than inserting zeros. For a block that is not
an exact-copy group, choose a pair \((i,j)\) with \(1-\rho_{ij}=\delta>0\) and first take \(i\) as anchor. We construct a residual event on which one sibling's comparison p-value exceeds the common rejection threshold early enough to give a definite loss in the comparison integral. Put

\[
 b_0=\barPhi(1)>0,\quad A_Q=\sqrt{4Q/b_0},\quad
 G_i=\{\varepsilon_j\le-\sqrt\delta,
 \max_{\ell\ne i}|\varepsilon_\ell|\le A_Q\sqrt\delta\}.
\]

The residual variance of the selected partner is \(\delta(2-\delta)\ge\delta\), while every sibling residual variance is at most \(2\delta\). A normal tail bound for that partner, followed by a union bound and Markov's inequality for the condition on the maximum absolute residual, gives

\begin{equation}
 \Prob(G_i)\ge b_0-2Q/A_Q^2=b_0/2. \label{eq:lm:n:5.1}
\end{equation}

The bound uses only marginal residual variances, so it also applies when the joint residual covariance is singular.

Let \(z=c_{\alpha q/m}\), which is at least one for all sufficiently
large \(m\) by \eqref{mx:quantiles}, and put \(x=z\sqrt\delta\). On \(G_i\) define

\[
 \zeta=\max_{\ell\ne i}\frac{z-\varepsilon_\ell}{r_\ell},\qquad
 \ell_*=2m\barPhi(\zeta).
\]

The selected partner's negative residual gives \(\zeta\ge z+\sqrt\delta\). The bounded residuals and \(r_\ell\ge1-\delta\) give \(\zeta\le z+d\), with

\[
 d=(\delta z+A_Q\sqrt\delta)/(1-\delta)
 \le2\delta z+2A_Q\sqrt\delta.
\]

We need the rank drop to occur before the nominal endpoint, but we also need a lower bound on how long the initial rank persists. The two tail ratios below give these bounds. The hazard lower bound gives the upper ratio; integrating the upper hazard bound \(\phi(v)/\barPhi(v)\le v+1/v\) for \(v>0\) gives the lower ratio:

\begin{equation}
 \frac z{z+d}e^{-zd-d^2/2}
 \le\frac{\ell_*}{\alpha q}\le e^{-x}<1. \label{eq:lm:n:5.2}
\end{equation}

For \(\delta\le1/2\) and \(z\ge1\),

\[
 d/z\le1+\sqrt2A_Q,
 \quad zd+d^2/2\le3x^2+4A_Qx+A_Q^2\le4x^2+5A_Q^2.
\]

In the nontrivial case $q\ge2$, hence $Q\ge2$, define

\[
 c_Q=\frac{e^{-5A_Q^2}}{2+\sqrt2A_Q},\qquad
 b_Q=\min\{1/2,c_Q/(Q-1)\}>0.
\]

Then \eqref{eq:lm:n:5.2} implies

\begin{equation}
 c_Qe^{-4x^2}\le\ell_*/(\alpha q)\le e^{-x}. \label{eq:lm:n:5.3}
\end{equation}

On the outside event \(E\), \(H^+(p)\le q\) for every \(p\). It equals \(q\) exactly for \(p\le\ell_*\), since only then are all \(q-1\) siblings, as well as the forced anchor, counted at \(\alpha q\). For \(p>\ell_*\) at least one sibling is not counted there, and \(H^+(p)\le q-1\). The comparison rank can be extended beyond \(m\) while preserving these properties; its selected integral has already ended by \(\alpha q<m\).

Applying \eqref{eq:lm:n:2.3}, \eqref{eq:lm:n:5.3}, and \(1-e^{-x}\ge\frac12\min(x,1)\) gives

\begin{equation}
 \alpha-I(H^+)\ge\alpha b_Q\min\{x,1\}e^{-4x^2}
 \qquad\text{on }E\cap G_i. \label{eq:lm:n:5.4}
\end{equation}

Everywhere else the loss \(\alpha-I(H^+)\) is nonnegative. Independence of \(E\) and \(G_i\), the bound \(\Prob(E)\ge\gamma\), and \eqref{eq:lm:n:5.1} yield an expected loss at least

\begin{equation}
 K_0\min\{z\sqrt\delta,1\}e^{-4z^2\delta},\qquad
 K_0=\gamma(b_0/2)\alpha b_Q>0. \label{eq:lm:n:5.5}
\end{equation}

Here \(K_0\), which is unrelated to the rank cutoff \(K\), is
proportional to \(\alpha\), the scale of the comparison integral
in \eqref{eq:lm:n:2.2}. The same construction applies with \(j\) as anchor. These are two distinct physical coordinates even if there are other exact copies. Their losses can therefore be added in the physical-share identity. The two residual events need not occur together; each loss is averaged separately.

\subsubsection{The resulting bound for a near-copy block}

\begin{lemma}[The share of a near-copy block]\label{mx:near}
Suppose that, after sign changes, a selected block $C$ of size $w\le Q$
has all correlations at least $1-\delta$, where
$\delta=\max_{i,j\in C}(1-\rho_{ij})\le1/16$, and is independent of the
outside described above. There are $c_N>0$ and one common dimension threshold such that
\begin{equation}\label{mx:near_margin}
 b_C\le wt-c_NtG_t(\delta),\qquad
 G_t(\delta)=\min\{\sqrt{\delta L},1\}e^{-12\delta L},
 \quad 0<\delta\le1/16.
\end{equation}
For $\delta=0$, $b_C=wt$ exactly at every dimension and level.
\end{lemma}

\begin{proof}
We combine the preceding conditional estimates and make the comparison
uniform over the testing level. The remaining task is to show that the
folding error is smaller than the positive loss, even when the spread
tends to zero arbitrarily quickly. Use the
physical-anchor integral \eqref{eq:lm:n:2.6} for the anchors in
$C$, leaving the actual outside law fixed. The monotone one-sided comparison in
\eqref{eq:lm:n:2.5} has integral at most $\alpha$ at every anchor, by
\eqref{eq:lm:n:2.2}. On the event in Lemma~\ref{mx:void}, no outside
crossing can occur after inserting the $w$ cluster coordinates. Choose a
pair attaining $\delta$. For either of its endpoints, the residual event
and the Mills estimates in \eqref{eq:lm:n:5.1}--\eqref{eq:lm:n:5.5}
apply with $z=c_{wt}$. They give expected comparison loss at least
\[
 \gamma(b_0/2)\alpha b_Q\min\{z\sqrt\delta,1\}e^{-4z^2\delta}
 \ge c\alpha G_t(\delta),
\]
by \eqref{mx:quantiles}. The two anchors attaining the spread are distinct original
coordinates; all other comparison losses are nonnegative. This argument
uses only marginal residual variances and therefore permits singular
cluster covariance.

It remains to bound the error caused by taking absolute values,
uniformly when both $\alpha$ and $\delta$ tend to zero. The deterministic necessary condition for a folding error and the bivariate
opposite-sign bound \eqref{eq:lm:n:3.2}--\eqref{eq:lm:n:3.6}, with the
actual outside crossing probability replaced by $q_j$, give at each
anchor
\begin{equation}\label{mx:folding_start}
 E_{\rm fold}\le\alpha(Q-1)\sum_{j=1}^m
             q_j e^{-a(\delta)c_{jt}^2},\qquad
 a(\delta)=\delta^{-1}-\tfrac12.
\end{equation}
We use different bounds for early and late ranks. At early ranks the Gaussian quantiles have a useful lower bound; at late ranks the required outside count is unlikely. Applying the two ranges in \eqref{mx:count_sums} gives
\begin{equation}\label{mx:folding_split}
 E_{\rm fold}/\alpha\le Cc_m
 \left\{e^{-a(\delta)L/2}
       +e^{-h_{a_1}\sqrt m/c_m-a(\delta)c_\alpha^2}\right\}.
\end{equation}
We verify explicitly that both terms are $o(G_t(\delta))$ uniformly.
This relative estimate is needed because the loss may vanish arbitrarily
quickly as the block approaches exact duplication.
If $\delta L\ge1$, the early ratio is at most $t^7$, because
$-1/(2\delta)+1/4+12\delta\le-7$ for $\delta\le1/16$.
If $\delta L<1$, split $e^{-L/(2\delta)}$ into two equal factors.
One is at most $e^{-4L}$, and the other divided by $\sqrt{\delta L}$ is
at most $2/L$. The early ratio is thus at most
$2e^{12}t^{15/4}/L$.

For the late ratio when $\delta L\ge1$, $a(\delta)\ge31/2$ and
$e^{12\delta L}\le(m/\alpha)^{3/4}$. Mills' lower bound
\[
 \alpha\ge\sqrt{2/\pi}\,
      \frac{c_\alpha}{1+c_\alpha^2}e^{-c_\alpha^2/2}
\]
shows that
$\sup_{\alpha\le a_1}\alpha^{-3/4}e^{-(31/2)c_\alpha^2}
\le C_*<\infty$; for example one may take
$C_*=(\pi/2)^{3/8}\{2(1+z_*^{-1})\}^{3/4}$.
The late ratio before multiplication by the exponential factor from the outside-count bound is therefore at most
$C_*m^{3/4}$. When $\delta L<1$, use
$a(\delta)\ge1/(2\delta)+15/2$ and put $v=c_\alpha^2/\delta$.
Then
\[
 \frac{e^{-c_\alpha^2/(2\delta)}}{\sqrt{\delta L}}
 =\frac{\sqrt v e^{-v/2}}{c_\alpha\sqrt L}
 \le\frac1{z_*\sqrt L}.
\]
Combining all four bounds proves
\begin{equation}\label{mx:folding_relative}
 \frac{E_{\rm fold}}{\alpha G_t(\delta)}\le Cc_m
 \left[t^7+\frac{2e^{12}t^{15/4}}L+
 e^{-h_{a_1}\sqrt m/c_m}
 \left(C_*m^{3/4}+\frac{e^{12}}{z_*\sqrt L}\right)\right]=o(1).
\end{equation}
Indeed $c_m\le m^{1/8}$ eventually, so the last term is bounded by a
constant times $m^{7/8}e^{-h_{a_1}m^{3/8}}$, with a smaller second summand.
The earlier terms tend to zero by Lemma~\ref{mx:uniform_factors}.
Sum over the at most $Q$ anchors and divide their integrals by $m$.
The losses at the two selected anchors exceed the total error from
taking absolute values and give
\eqref{mx:near_margin}. If $\delta=0$, the cluster is one independent
uniform representative with $w$ copies, giving the exact equality.
\end{proof}

\subsubsection{Counting crossings with an additional event}

The next estimate counts threshold crossings together with an additional event, called a mark. We will use the mark to require that another coordinate has a small p-value. This allows the crossing bound to retain its dependence on the cluster spread, even when that spread is arbitrarily small.

\begin{lemma}\label{fd:marked_crossing}
Let $A,E$ be independent standard normals on a probability space
possibly carrying other dependent Gaussian variables. Put
$X_\theta=\cos\theta A+\sin\theta E$, $0\le\theta\le\Theta$, and
let $M_z$ be the event that the endpoint indicators of $|X|\ge z$
differ, where $z>0$. For any event $D_0$ on the same space,
\begin{equation}\label{fd:marked_crossing_bound}
 \Prob(M_z\cap D_0)
 \le\sum_{s=\pm1}\int_0^\Theta
 \phi(z)\Expect[|\dot X_\theta|\mathbf1_{D_0}\mid
                         X_\theta=sz]\,d\theta.
\end{equation}
In this integrand $\dot X_\theta$ has conditional second moment one.
\end{lemma}

\begin{proof}
We count crossings while retaining the additional event inside the expectation. Thus the estimate does not require that event to be independent of the path or its velocity.

An endpoint mismatch requires a crossing of one of the levels $z,-z$.
For a smooth rotation, count these crossings before taking expectation.
For fixed $\theta$, $(X_\theta,\dot X_\theta)$ is a pair of independent
standard normals. Change of variables at each crossing contributes
the Jacobian $|\dot X_\theta|$ and the density $\phi(z)$.
More explicitly, integrate first in the nondegenerate $(A,E)$ plane,
condition the other variables on that plane, and apply the
one-dimensional change-of-variables formula on the finitely many
monotone arcs of $\theta\mapsto X_\theta$. Tonelli's theorem then gives the displayed expression
for the expected marked crossing count. A tangency at either fixed
nonzero level would require the Gaussian radius
$\sqrt{A^2+E^2}$ to equal $z$, which has probability zero. Crossings
at endpoints also form a null event. The bound follows because the
marked mismatch indicator is at most the marked crossing count.
The conditional second moment is one by independence of the rotated
position and velocity. The velocity need not be independent of
$D_0$; Gaussian regression below controls its absolute moment on the marked tail.
\end{proof}

\subsubsection{Transporting the share of an entire block}
\label{sec:lm:c:6}

\begin{lemma}[Signed transport of a whole-block share]\label{mx:transport}
Let $C\subseteq A$ have size at least two and, after sign changes, spread
$0<\delta\le h^2/8$, where $0<h\le1/16$. Suppose that an actual anchor
$X\in C$ has absolute correlation at most $1-h$ with every coordinate in
$A\setminus C$. Let $H_1$ be the expected share of the whole block $A$
under the law to which the lemma is applied (in
Lemma~\ref{mx:arbitrary_block}, the law after the earlier contractions),
and let $H_0$ be its expected share after replacing
every coordinate in $C$ by the signed copy of that anchor, retaining all
other coordinates. There are $c_I>0$ and a dimension threshold $M_I$,
depending only on $h$ and the fixed parameters of this subsection, such
that for all $m\ge M_I$, all $0<\alpha\le a_1$, all blocks $A$, clusters $C$ and
anchors $X$ as above, and all outside laws of this subsection,
\begin{equation}\label{mx:transport_contraction}
 H_1\le H_0-c_I tG_t(\delta).
\end{equation}
No independence is assumed between $C$ and $A\setminus C$.
\end{lemma}

\begin{proof}
We compare the change in the whole-block share under two laws. The first is the law in the statement. In the second, which we call the independent law, the coordinates outside the cluster but inside the selected block are replaced by independent standard normals. The two changes will differ by less than the loss established for the second law. Sign-align $C$ and regress each coordinate on the actual
anchor $X$. For positive residual variance, write
$Z_i=\rho_{iX}X+\sigma_iE_i$, where
$\sigma_i^2=1-\rho_{iX}^2\le2\delta$ and $E_i$ is standard normal.
The joint residual vector is Gaussian and independent of $X$; it may be
singular, its coordinates may be dependent, and it may depend on the
actual coordinates in $A\setminus C$. Couple the endpoints by
\begin{equation}
 Z_i(v)=\sqrt{1-v^2\sigma_i^2}\,X+v\sigma_iE_i,
 \quad 0\le v\le1,\qquad
 \sigma_i^2=1-\rho_{iX}^2\le2\delta.
 \label{eq:lm:c:6.2}
\end{equation}
If $\sigma_i=0$, keep that coordinate identically equal to $X$.
Every coordinate in $A\setminus C$ and in the independent background
$V$ stays unchanged, together with their original joint law. This
construction is therefore a positive-semidefinite Gaussian path with
standard marginals, from exact copies of the actual anchor to the
original cluster. It does not require an inverse covariance matrix. Conditional on $Z_i(v)=z$, the derivative satisfies
$\Expect[\dot Z_i(v)^2\mid Z_i(v)=z]\le4\sigma_i^2$, since the path has constant
unit variance and hence $\dot Z_i(v)$ is independent of $Z_i(v)$.
Also $\sum_i\sigma_i\le Q\sqrt{2\delta}$. For any actual coordinate $Y\in A\setminus C$, the path's absolute
correlation with $Y$ is at most
$|\operatorname{Corr}(X,Y)|+\sigma_i\le1-h+\sqrt{2\delta}\le1-h/2$.
This uses the original joint residual law, including its dependence on $Y$.

The independent law makes the cluster independent of the other coordinates in the selected block. The near-copy bound will then give a signed decrease along the path. To transfer that sign back to the original law, we compare the two path increments and retain their common part before estimating the error.

Define $B_1,B_0$ to be the corresponding expected whole-block shares
when $A\setminus C$ is replaced by independent standard normals,
independent of $C$ and of the unchanged entire outside $V$, at the two
path endpoints. We bound
$(H_1-H_0)-(B_1-B_0)$ by comparing the two increments before taking an
absolute value. To do so, we split the whole-$A$ random share into a part
computed from ranks at most $K$ and a late-rank remainder; the early parts
are compared across the two couplings, and the late parts are bounded in
each coupling separately. Write $M_j$ for the event
that at least one coordinate path crosses a boundary $\pm c_{jt}$.
The probability of a path crossing is at most its expected crossing
count. The latter is obtained by integrating the conditional absolute
path derivative against the Gaussian density at each boundary, then
summing over the coordinate paths. For a positive residual scale, the
substitution $\theta=\arcsin(v\sigma_i)$ makes the path a rotation of
two independent standard normals. The change-of-variables argument
in Lemma~\ref{fd:marked_crossing} therefore applies, including with
an additional event in the integrand. Combining the crossing count
with the derivative bound gives
\begin{equation}\label{mx:marked_unconditional}
 \Prob(M_j)\le Ct\sqrt{\delta L}\,j\quad(1\le j\le m).
\end{equation}
Here $\phi(c_{jt})\le C_{a_1}jt\,c_t\le Cjt\sqrt L$. This inequality
also covers large ranks because $c_{jt}\ge c_{a_1}>0$.

For the early ranks, the additional event is that a coordinate of the selected block outside the cluster has a p-value below the censoring threshold. The fixed correlation gap makes this joint occurrence less likely than an unmarked cluster crossing.

For $j\le K$, use the conditional Gaussian moment with the additional event retained. Let
$r=\operatorname{Corr}(Z_i(v),Y)$, so $|r|\le1-h/2$. Conditional on
$Z_i(v)=x$, write $Y=rx+\sqrt{1-r^2}\,T$, and standardize a
positive-variance derivative as $U$. The pair $(U,T)$ is centered
jointly standard normal conditionally, possibly singular. Regression
therefore gives $\Expect(|U|\mid T)\le|T|+\sqrt{2/\pi}$.
For $|x|\le c_t$, the marked event implies $|T|\ge b$, where
\[
 b=\frac{c_{Kt}-|r|c_t}{\sqrt{1-r^2}}>0.
\]
This positivity holds for all sufficiently large $m$, uniformly in
$\alpha$ and $|r|\le1-h/2$, because $c_t/c_{Kt}\to1$ uniformly by
\eqref{mx:quantiles}. The weighted tail of $T$ is at most
$\Expect[(|T|+\sqrt{2/\pi})\mathbf1\{|T|\ge b\}]
\le C e^{-b^2/2}$, using $\Expect|T|\mathbf1\{|T|\ge b\}=2\phi(b)$.
Uniform quantile asymptotics give
\[
 \inf_{0\le |r|\le1-h/2}\frac{b^2}{2L}
 \longrightarrow\frac{h}{4-h}>h/4.
\]
Thus the early marked error has exponent $d=h/5$. Retaining the
original derivative scale $O(\sigma_i)$ and summing over at most $Q$
paths, $Q$ coordinates outside the cluster but inside the block, and $K$ thresholds yields
\begin{equation}\label{mx:marked_early}
 \Prob\{\text{an early path mismatch and an outsider below }Kt\}
       \le Ct\sqrt{\delta L}\,K^2t^d,\qquad d=h/5.
\end{equation}
The same estimate holds for the independent replacements for those coordinates.
Zero derivative variance gives a constant path and contributes nothing.

To see the cancellation, let $g$ be the whole-$A$ random share computed from ranks at most $K$ and clipped at one. After all coordinates in $A\setminus C$ are replaced by p-value one, the increment of $g$ along the path is the same function of the cluster path and $V$ in the original and independent couplings, so its expectation is the same in both. An increment can differ from this common value only when a cluster path changes an early threshold indicator and a coordinate in $A\setminus C$ has a small p-value. Equation~\eqref{mx:marked_early} bounds its expected magnitude. This cancellation is for the whole block and need not hold for the share of one coordinate.

For late ranks, each error event must still include a change in a cluster threshold indicator along the path. A bound based only on a late crossing would lose its dependence on the cluster spread. Put
$\ell=T_A-g$ and let subscripts $0,1$ denote path endpoints. A deterministic
case split gives
\begin{align}\label{mx:late_deterministic}
 |\ell_1-\ell_0|\le{}&
 2\sum_{j>K}\mathbf1_{M_j}\mathbf1\{N_V(jt)\ge j-s\}\notag\\
 &+2\sum_{k\le K}\mathbf1_{M_k}
       \mathbf1\{\text{some }j>K:\ N_V(jt)\ge j-s\}.
\end{align}
To see this, suppose first that the full ranks at the two endpoints differ.
Then a cluster threshold indicator changes at the larger of these ranks, and
the outside supplies all but at most $s$ of the counts there. If the full ranks agree, their whole-$A$ shares agree because
$N_A(tR)=R-N_V(tR)$. The difference of early truncated shares can still
remain on the late event; this is the second term in
\eqref{mx:late_deterministic}. These cases also include a late crossing at only one endpoint. The
factor two allows for both the full and truncated random shares.

The entire outside vector $V$ is independent of every marked path event. Apply
\eqref{mx:marked_unconditional} and \eqref{mx:rank_tail} to the first
term. For the second, sum \eqref{mx:marked_unconditional} over $k\le K$
and bound the late outside event by $\sum_{j>K}q_j$, which
\eqref{mx:rank_tail} also controls.
The result, for both couplings, is
\begin{equation}\label{mx:transport_error}
 |(H_1-H_0)-(B_1-B_0)|
 \le Ct\sqrt{\delta L}\,c_m^2(K^2+K+1)(t^d+t^{20}).
\end{equation}
The factors $c_m^2$ include the prefactors from both geometric sums
over late ranks. Because each late event was counted together with a
path mismatch, the bound retains the square-root dependence on the
cluster's spread.

Under the independent law, $B_0=st$. Independent replacement
singletons have their exact nominal expected shares, and Lemma~\ref{mx:near}
applies to $C$ with those singletons included in its background. Thus
$B_1-B_0\le-c_NtG_t(\delta)$. Divide
\eqref{mx:transport_error} by $tG_t(\delta)$. When $\delta L\le1$ the
ratio is at most
$Ce^{12}c_m^2(K^2+K+1)(t^d+t^{20})$.
When $\delta L>1$ it is at most
\[
 Cc_m^2(K^2+K+1)\sqrt L
       (t^{d-12\delta}+t^{20-12\delta}).
\]
Since $\delta\le h^2/8$ and $h\le1/16$,
$d-12\delta\ge h/5-3h^2/2\ge17h/160>d/2$;
both exponents are strictly positive. Lemma~\ref{mx:uniform_factors} makes both bounds tend to
zero uniformly. Taking $c_I=c_N/2$ proves the assertion.
\end{proof}

After every near-copy cluster has been contracted, the block consists of
exact-copy groups whose representative correlations are bounded away from
one in magnitude. The next lemma bounds the expected share of such a
block; Lemma~\ref{mx:arbitrary_block} combines it with
Lemma~\ref{mx:transport}.

\begin{lemma}[A selected capped block]\label{mx:capped}
Suppose the selected physical block $A$ has size at most $Q$. After exact
signed copies are compressed with their multiplicities retained, suppose
all distinct representative correlations have magnitude at most a fixed
$r<1$. Keep its actual mixed outside law unchanged. Define $d_e$ by
replacing only the covariance inside $A$ by its single-edge reference and
subtracting its expected whole-$A$ share from $st$. Put $D=\sum_e d_e$ and
$g_r=(1-r)/\{(1+2r)(1+r)\}$. For every $0<\varepsilon<g_r$, one common
dimension threshold gives
\begin{equation}\label{mx:capped_relative}
 |b_A-st+D|\le t^{g_r-\varepsilon}D.
\end{equation}
Each nonzero edge has positive $d_e$ eventually. Hence $b_A<st$ when there
is an edge; the edgeless case has expected share exactly $st$ at every level.
\end{lemma}

\begin{proof}
Here the centered expansion concerns only the selected bounded block. The actual outside component supplies the background event and late-rank bounds, but it does not enter the count of internal connected sets. This is why a fixed cap for the selected block is sufficient even when its outside graph grows.

The independent baseline is $st$. The representatives other than the two edge endpoints in a
single-edge reference are independent groups of at most $Q$ copies. Add
them to the background of Lemma~\ref{mx:void}; its parameters remain fixed.
The monotonicity and complete deficit argument proving
\eqref{mx:edge_lower}--\eqref{mx:edge_estimates} gives
\[
 d_e\ge c r_e^2Lt^2,\qquad
 d_e\ge c r_e^2L^{-1}t^{2/(1+r_e)},\qquad
 \|A_e^{Kt}\|\le CK^2L^2d_e.
\]
These deficits concern the expected share of the whole block. Every
reference retains the same actual outside law. For $r_{\max}\le t^4$, the full folded-law cancellation
\eqref{eq:lm:g:4.1} applies to the share averaged over the unchanged outside law. Since there are
at most $Q$ representatives, its relative error is at most $Ct^2/L$.
This step covers correlations smaller than every power of $t$.

Otherwise $D\ge cLt^{10}$. Use \eqref{mx:common_cutoff}; the full late
probability is at most $Cc_mt^{20}$ by physical coloring. The clipped
whole-block random share computed from the first $K$ ranks is bounded by one. In fact clipping is
unnecessary when $K\ge Q$: if its count exceeded $R_K$ while $R_K<K$,
there would be a crossing at $R_K+1$, whereas at $R_K=K$ the block count
is at most $Q\le K$. The centered expansion has exactly the remainder
\eqref{mx:lc_remainder}, with these internal edges.

We now bound each remainder relative to the sum of pair deficits. Call
an internal connected set weak when its largest absolute correlation is
at most $[Q\max(8,c_x^2)]^{-1}$; then $\eta_C\le\min(1/8,c_x^{-2})$. For weak connected sets use
\eqref{eq:lm:g:6.1} with the fixed internal size bound $Q$, and divide by
$cLt^2\sum_e r_e^2$. The result is at most $CK^3Lt$.
For each remaining connected set select an edge with largest correlation magnitude $r_e$.
It exceeds $c/L$, so that its own deficit is at least
$cL^{-3}t^{2/(1+r_e)}$. Equation~\eqref{eq:lm:g:5.6} and the finite number
of internal connected sets then give
$CK^3L^{9/2}t^{g_r}$. The block size is bounded, so only finitely many connected sets need to
be counted.
Finally \eqref{eq:lm:g:5.7} gives
$V_2\le CL K^2t^{2/(1+r)}$, while the pair norm-to-deficit bound gives
$V_2/D\le CK^2L^2$. Thus
\begin{equation}\label{mx:capped_errors}
 \frac{B_3}{D}\le C(K^3Lt+K^3L^{9/2}t^{g_r}),\quad
 \frac{V_2^2}{D}\le CK^4L^3t^{2/(1+r)},\quad
 \frac{2|E|T_m}{D}\le Cc_mt^{10}/L.
\end{equation}
As before, $D\le V_2+2|E|T_m=o(1)$ follows from the pair references
without using the conclusion. Hence $B_3=o(1)$. The uniform
subpolynomial bounds on $K,c_m,L$ prove \eqref{mx:capped_relative}.
\end{proof}

\begin{lemma}[An arbitrary selected bounded block]\label{mx:arbitrary_block}
Under the selected-block assumptions of this subsection, one common
dimension threshold gives $b_A\le st$ for every positive-semidefinite covariance of
$A$ and every $0<\alpha\le a_1$. The inequality is strict unless $A$
consists of independent groups of exact signed copies.
\end{lemma}

\begin{proof}
We partition the block into groups whose coordinates are nearly signed copies. To do so, choose an interval that contains none of the pairwise correlation distances. The resulting gap makes coordinates within a group close enough, and different groups far enough apart, to apply the preceding comparison. The gap may depend on the covariance, but it will be chosen from a fixed finite list. Taking the largest of the corresponding dimension thresholds will therefore give a uniform bound. Set
$N=Q(Q-1)/2$, $u_0=1/16$, and
$u_{k+1}=u_k^2/(65536Q^2)$ for $0\le k\le N$.
There are at most $N$ distinct distances $1-|\rho_{ij}|$, so one of the
$N+1$ intervals $(u_{k+1},u_k]$ is empty. Write $h=u_k$ and connect pairs
whose distance is at most $u_{k+1}$.

Realize the covariance by unit Gram vectors. In each connected component,
sign-align a spanning tree from an actual root. A tree edge then has
Euclidean length at most $\sqrt{2u_{k+1}}$. Every distance from the root is at
most $(Q-1)\sqrt{2u_{k+1}}$, and the tree path joining two vertices, which has
at most $Q-1$ edges, bounds every pair distance by
$Q\sqrt{2u_{k+1}}$. Thus each component has sign-aligned spread at most
$Q^2u_{k+1}=h^2/65536$. This also proves sign consistency: the aligned
inner products within the root ball are positive, so an edge of magnitude
at least $1-u_{k+1}$ cannot have the opposite aligned sign. Between
different components the distance exceeds $u_{k+1}$ and, by the empty
interval, exceeds $h$. Their absolute correlations are therefore less
than $1-h$.

Retain one actual representative in each component. Replace each
component by signed copies of its representative, one component at a time. Every retained representative is an original coordinate,
so its correlations to uncontracted components retain the same gap.
Lemma~\ref{mx:transport} applies at each step with $H_1=b_{\rm before}$ and
$H_0=b_{\rm after}$, and \eqref{mx:transport_contraction} gives
$b_{\rm before}\le b_{\rm after}-c_I tG_t(\delta)$ whenever the
component has positive spread. Thus the original expected share is bounded by the expected share under
the final reference law minus the sum of these positive losses. A zero-spread component is already an exact group of signed
copies and its compression changes no p-value. The terminal covariance
has cap $1-h$ and multiplicities at most $Q$, so
Lemma~\ref{mx:capped} bounds its expected share by $st$, strictly if a
nonzero edge remains. All comparisons keep the entire actual outside law fixed.
There are at most $Q$ steps and only the fixed finite collection
$h\in\{u_0,\ldots,u_N\}$, so the maximum of their dimension thresholds works for
every covariance of $A$. Equality can remain only when every collapsed
component had zero spread and the terminal representatives are
independent. This is exactly the asserted class.
\end{proof}

\subsection{Proof of Theorem~\ref{main:mixed_theorem}}

\begin{proof}
The proof bounds the expected share of each original part under the
original law and adds these bounds in \eqref{mx:shares}:
Lemma~\ref{mx:lc} treats $G$, and Lemma~\ref{mx:arbitrary_block} treats
each $A_j$. We first check the uses of the degree envelope in the selected-block
lemmas. The physical color count is $c_m=O(\log m)$ and
$K=O((1+\log m)^2)$, uniformly for $\alpha\le a_1$.
Lemma~\ref{mx:void} applies because the actual background has
$O(m\log m)$ edges, also after adding a bounded number of coordinates at any fixed cap
less than one. All subsequent near-copy and transport estimates use
only this positive probability, the exponential color tails, and fixed
polynomials in $c_m,K,L$. For example the folding error has an additional
factor $\exp\{-h_{a_1}\sqrt m/c_m\}$, and the transport exponents retain
$d-12\delta>d/2>0$. Only connected sets within the terminal selected block need to be counted, and that block has at most $Q$ coordinates. Its correlation cap may be $1-h$, but this larger cap is not used in a graph expansion for the outside component. The larger allowed degree range of the original graph therefore causes no additional loss in Lemmas~\ref{mx:near}, \ref{mx:transport}, or
\ref{mx:capped}, and hence in Lemma~\ref{mx:arbitrary_block}.

We now add the bounds for the original independent parts. Under the
original full law, Lemma~\ref{mx:lc} bounds the expected share of
$G$ by $|G|t$, where $|G|$ is its number of physical coordinates; the
actual independent Gaussian blocks $A_j$ form its outside. For each original $A_j$, retain the same original $G$ and every
other original $A_i$. Together they form an admissible mixed outside distribution for
Lemma~\ref{mx:arbitrary_block}, which gives $b_{A_j}\le |A_j|t$.
These dimension thresholds depend on the fixed size bounds, not on the
number of blocks or their combined size. Thus no error accumulates with
the number of blocks. Sum these bounds in \eqref{mx:shares}. All terms
refer to the original full distribution, so their sum is
$F_\Sigma(\alpha)\le mt=\alpha$. If any capped edge or any nontrivial
within-block dependence remains after exact-copy compression, its
corresponding original expected share is strictly below its nominal value. Conversely,
independent signed-copy groups have equality by forcing each group.
\end{proof}

Taking $G$ empty in Theorem~\ref{main:mixed_theorem} gives the
bounded independent-block conclusion in Section~\ref{main:mixed_blocks}.
All background colors are then bounded by $Q$, so the remaining
constants depend only on $Q,a_1$.

\subsection{An explicit obstruction to coordinatewise share bounds}
\label{mx:coordinate}

The preceding proof bounds the expected share of each whole independent block by its nominal contribution. The analogous bound for one coordinate need not hold: its expected share can exceed its nominal contribution even in a positive-definite block of size three. We now give an explicit example.

\begin{proposition}\label{mx:coordinate_theorem}
Fix $0<a_1<1$ and define
\begin{align*}
 h_0&=\sqrt{2/\pi},& a&=h_0/(1+h_0),&
 C_*&=\frac{2e^{-2}a^2}{1+\sqrt{2\log3}},\\
 \eta_0&=\min\{1,3c_{a_1}/4,\sqrt{2\pi}/8\},&&
 C_1=6+\eta_0^{-2}.&&
\end{align*}
For $m\ge3$, $0<\alpha\le a_1$, $t=\alpha/m$, and
$H_m=\sum_{j=1}^m j^{-1}$, put
\begin{equation}\label{mx:coordinate_parameters}
 q=\frac{C_*}{8C_1}\frac{t^4}{H_m},\qquad
 \varepsilon=\left(\frac{C_*^2}{384C_1}\right)^3
                  \frac{t^{21}}{H_m^3},\qquad
 \Gamma_q=\begin{pmatrix}
 1&1-2q^2&q\\1-2q^2&1&-q\\q&-q&1
 \end{pmatrix}.
\end{equation}
The correlation matrix $\Gamma=(1-\varepsilon)\Gamma_q+\varepsilon I_3$ is
positive definite. A centered Gaussian triple with this correlation,
independent of an arbitrary vector of $m-3$ outside p-values in $[0,1]$,
has third-coordinate expected share
\begin{equation}\label{mx:coordinate_lower}
 S_3\ge t+\frac{C_*^2}{32C_1}\frac{t^6}{H_m}>t.
\end{equation}
No assumptions on the outside marginals or dependence are needed.
\end{proposition}

\begin{proof}
We first construct a singular triple for which a small perturbation
increases the third coordinate's expected share above its independent
value. We quantify this increase and then make the covariance positive
definite while preserving it. Start with the singular triple
\[
 X_1=\sqrt{1-q^2}\,U+qV,\qquad
 X_2=\sqrt{1-q^2}\,U-qV,\qquad X_3=V,
\]
where $U,V$ are independent standard normals and are independent of the
outside p-values. Force $P_3$ to zero and write
$C_k(u;q)=\Prob\{R^{[3]}\le k\mid P_3\le u\}$. The exact forcing identity,
which is \eqref{eq:lm:g:2.2} with $w_3=1$ and uses only the uniformity of
$P_3$, is
\begin{equation}\label{mx:coordinate_forcing}
 S_3(q)=t-t\sum_{k=1}^{m-1}
       \{C_k((k+1)t;q)-C_k(kt;q)\}.
\end{equation}
At $q=0$, coordinate three is independent of the entire outside and
$S_3(0)=t$.

To see how the perturbation changes the forced rank, hold the outside p-values fixed. We compare the rank supported without the first two coordinates with the rank supported when both are counted. For a sufficiently small separation of that pair, a crossing supported by only one member will not determine the largest rank.

Fix an arbitrary outside profile. Put
$B_r=1+N_{\rm out}(rt)$ and define
\[
 R_0=\max\{r:B_r\ge r\},\qquad
 R_2=\max\{r:B_r+2\ge r\}.
\]
The forced coordinate gives $1\le R_0\le m-2$, and monotonicity of the
counts gives $R_2\ge R_0+2$. Set $d_r=c_{rt}$,
$b=q|V|$, and $x=|\sqrt{1-q^2}U|$. The absolute values of the two pair coordinates are
$x+b$ and $|x-b|$. Suppose
\begin{equation}\label{mx:coordinate_small_shift}
 b\le\eta<c_\alpha/4,\qquad
 2\eta<\min_{r<m}(d_r-d_{r+1}).
\end{equation}
A crossing at rank $r$ requiring only one pair member has $B_r=r-1$.
It cannot occur at $r=m$, because $B_m\le m-2$. If $x+b\ge d_r$, then
$x-b\ge d_r-2b>d_{r+1}>0$, and both pair members cross at rank $r+1$.
Such a one-member crossing is therefore redundant for the largest rank.
The alternative folded inequality $b-x\ge d_r$ is impossible under
\eqref{mx:coordinate_small_shift}.

It follows that, conditionally on $|V|=v$, the forced-rank distribution function is zero
for $k<R_0$, one for $k\ge R_2$, and, for every
$R_0\le k<R_2$, is the same function
\begin{equation}\label{mx:coordinate_proxy}
 f(q,z,v)=2\Phi\left(\frac{z+qv}{\sqrt{1-q^2}}\right)-1,
 \qquad z=c_{R_2t}.
\end{equation}
For $k$ in this range, all pair-dependent crossings require both members;
monotonicity of $B_r$ makes the least restrictive such threshold the one
at $R_2$. This also proves the formula at profiles with ties.
The derivative at zero is $2v\phi(z)$. It describes the first change in
the forced-rank probability as the pair departs from exact duplication.
Write $\mu(c)=\phi(c)/\barPhi(c)=\Expect(|V|\mid |V|\ge c)$. Differentiating
\eqref{mx:coordinate_forcing}, justified quantitatively below, and
telescoping the consecutive nonconstant distribution functions gives
\begin{equation}\label{mx:coordinate_derivative}
 D:=S_3'(0+)=2t\Expect_{\rm out}\left[
 \phi(c_{R_2t})\{\mu(c_{R_0t})-\mu(c_{R_2t})\}\right].
\end{equation}
Every integrand is positive. We need a lower bound uniform over outside
profiles and over the level.

For $c\ge0$, $\mu(c)\ge\max(c,h_0)\ge a(c+1)$. Restricting the integral
for $\mu(c)-c$ to
$[c+1/(c+1),c+2/(c+1)]$ gives
\begin{equation}\label{mx:hazard_lower}
 \mu(c)-c\ge\frac{e^{-2}\phi(c)}{(c+1)^2\barPhi(c)}
            \ge\frac{e^{-2}a}{c+1}.
\end{equation}
Indeed the density on that interval is at least $e^{-2}\phi(c)$ and the
excess over $c$ is at least $1/(c+1)$.
Until \eqref{mx:coordinate_D}, let $u=R_0t$ and $v=R_2t$; afterwards $u$
is again the cutoff in $C_k(u;q)$ and $v$ a value of $|V|$. The interval
$[v-2t,v]$ is contained in $[u,v]$,
and its left endpoint is at least $v/3$. For $s$ in this interval,
Mills' hazard lower bound yields $c_s^2-c_v^2\le2\log3$ and therefore
$(c_v+1)/(c_s+1)\ge(1+\sqrt{2\log3})^{-1}$. As
$\phi(c_v)/v=\mu(c_v)/2$, \eqref{mx:hazard_lower} gives
\[
 \frac{\phi(c_v)\{\mu(c_s)-c_s\}}s
 \ge\frac{e^{-2}a^2}{2(1+\sqrt{2\log3})}.
\]
Finally
\[
 \mu(c_u)-\mu(c_v)
       =\int_u^v\frac{\mu(c_s)-c_s}{s}\,ds.
\]
Substitution into \eqref{mx:coordinate_derivative} proves the following lower bound for every fixed
outside profile:
\begin{equation}\label{mx:coordinate_D}
 D\ge C_*t^2.
\end{equation}

The positive derivative shows an initial increase, but does not specify a perturbation that works uniformly over all outside profiles and levels. We therefore bound the finite-$q$ remainder. Since
$-c'_u=1/(2\phi(c_u))\ge\sqrt{2\pi}/2$,
$d_r-d_{r+1}\ge\sqrt{2\pi}\,t/2$. Set $\eta=\eta_0t$.
Because $t<1/3$, this choice satisfies
\eqref{mx:coordinate_small_shift} whenever $q|V|\le\eta$.
Extend the formula in \eqref{mx:coordinate_proxy} to all $v\ge0$.
For $0\le q\le1/2$, write $y=(z+qv)/\sqrt{1-q^2}$. Direct
differentiation gives
\[
 y'=\frac{v+qz}{(1-q^2)^{3/2}}\le2(v+z),\qquad
 y''=\frac{(1+2q^2)z+3qv}{(1-q^2)^{5/2}}\le4(v+z).
\]
Using $\phi(y)\le1/2$ and $y\phi(y)\le1/4$ gives
$|f''(q,z,v)|\le6(1+z^2+v^2)$.
Thus its Taylor remainder is at most $3q^2(1+z^2+v^2)$.
The probability of the excluded shift event, conditional on $P_3\le u$,
is at most
\[
 \Prob(q|V|>\eta\mid P_3\le u)
       \le\frac{q^2}{u\eta_0^2t^2},
\]
by the unconditional second moment of $V$. Also
$z^2\le c_t^2\le2L$, $\Expect(V^2\mid P_3\le u)\le1/u$, and
$t^2(2+2L)\le2$. Consequently, after averaging over any outside law,
\begin{equation}\label{mx:coordinate_CDF_remainder}
 |C_k(u;q)-C_k(u;0)-q A_k\mu(c_u)|
       \le\frac{C_1q^2}{ut^2},\qquad u\ge t,
\end{equation}
where $A_k$ is the outside average of
$2\phi(c_{R_2t})\mathbf1\{R_0\le k<R_2\}$.
This proves the differentiation used above, as well as its required
uniformity, without any tail moment assumption on the outside profile. The reciprocal threshold in this error bound sums to a harmonic factor. Controlling that factor determines how small the explicit perturbation must be.
Summing \eqref{mx:coordinate_CDF_remainder} in
\eqref{mx:coordinate_forcing} gives
\begin{equation}\label{mx:coordinate_share_remainder}
 |S_3(q)-t-qD|\le2C_1H_mq^2/t^2.
\end{equation}
With the value of $q$ in \eqref{mx:coordinate_parameters},
\eqref{mx:coordinate_D} and \eqref{mx:coordinate_share_remainder} imply
\begin{equation}\label{mx:coordinate_singular_margin}
 S_3(q)\ge t+\frac{C_*^2}{16C_1}\frac{t^6}{H_m}.
\end{equation}
The stated constants have $C_*<1$ and $C_1\ge6$, so $q<1/2$.

It remains to make the block positive definite without losing the
margin. Let $E_1,E_2,E_3$ be independent standard normals independent of
everything, and put $Z_i=\sqrt{1-\varepsilon}X_i+\sqrt\varepsilon E_i$.
The resulting vector has covariance $\Gamma$, with minimum eigenvalue at least $\varepsilon$.
Moreover $\Expect(Z_i-X_i)^2\le2\varepsilon$. For any threshold $d\ge0$ and
$h>0$, an absolute-value threshold indicator can change only if
$||X_i|-d|\le h$ or $|Z_i-X_i|>h$. Its probability is at most
$4\phi(0)h+2\varepsilon/h^2$. Taking $h=\varepsilon^{1/3}$ bounds it by
$4\varepsilon^{1/3}$. There are only $3m$ affected coordinate-threshold
indicators. If none changes, the entire BH outcome and the third coordinate's
random share are unchanged. Since random shares lie in $[0,1]$, their
expectations satisfy
\begin{equation}\label{mx:coordinate_PD_error}
 |S_3(Z)-S_3(X)|\le12m\varepsilon^{1/3}
  =\frac{C_*^2}{32C_1}(mt)\frac{t^6}{H_m}
  \le\frac{C_*^2}{32C_1}\frac{t^6}{H_m}.
\end{equation}
The parameters also have $0<\varepsilon<1/2$. Subtracting
\eqref{mx:coordinate_PD_error} from \eqref{mx:coordinate_singular_margin}
proves \eqref{mx:coordinate_lower} for every $m\ge3$ and every
$0<\alpha\le a_1$.
\end{proof}

\begin{corollary}[Failure of both PRDS conditions]\label{mx:prds}
With independent uniform outside p-values, the positive-definite Gaussian
example in Proposition~\ref{mx:coordinate_theorem} violates both
lower-tail PRDS and pointwise PRDS. For lower-tail PRDS the violation can
be witnessed by adjacent thresholds $kt,(k+1)t\le\alpha$ and the
increasing event $D_k=\{R^{[3]}\le k\}$.
\end{corollary}

\begin{proof}
The expected-share excess forces at least one conditional-rank increment
to be negative. We use that increment to exhibit both PRDS violations.
Raising any p-value cannot increase a BH rejection count, so $D_k$ is
increasing in the full p-value vector; its definition does not involve $P_3$, because the forced rank replaces that coordinate by zero.
Let $C_k(u)=\Prob(D_k\mid P_3\le u)$ for the positive-definite example.
The forcing identity and \eqref{mx:coordinate_lower} imply that, for
some $1\le k<m$,
\[
 C_k((k+1)t)-C_k(kt)
 \le-\frac{C_*^2}{32C_1}\frac{t^5}{(m-1)H_m}<0.
\]
This is a direct violation of lower-tail PRDS within the testing-level
range. Since $P_3$ is uniform,
$C_k(u)=u^{-1}\int_0^u \Prob(D_k\mid P_3=v)\,dv$.
The running average of a nondecreasing function is nondecreasing.
Therefore the conditional probability in this integral is not
nondecreasing on $(0,\alpha]$, proving failure of pointwise PRDS as well.
More explicitly, its average on $(kt,(k+1)t)$ is below its average on
$(0,kt)$, so there exist ordered values in these intervals exhibiting the
reversal. Under the positive-definite Gaussian law the conditional
probability has a version that is continuous in $v\in(0,1)$: the event
is specified by finitely many threshold inequalities and their
conditional boundaries have zero probability. Thus the violation is not
an artifact of a choice of conditional-probability version.
\end{proof}

Proposition~\ref{mx:coordinate_theorem} and Corollary~\ref{mx:prds}
together prove Proposition~\ref{main:coordinate_share}, with
$c_1=C_*/(8C_1)$, $c_2=(C_*^2/(384C_1))^3$ and $c_3=C_*^2/(32C_1)$; the
two PRDS conditions are those defined in Section~\ref{main:prds}.

These examples belong to the mixed-block class with $Q=3$ when the
outside p-values are independent Gaussian singleton tests. For large
$m$, Lemma~\ref{mx:arbitrary_block}, used in the proof of
Theorem~\ref{mx:theorem}, bounds the expected share of the entire
three-coordinate block, even though the third coordinate's expected
share exceeds $t$. The first two coordinates must then compensate for
that excess. Such compensation is asserted only where a
Simes inequality has been proved; the individual-share construction by
itself does not establish a finite-dimensional Simes inequality or a
counterexample to it.

\clearpage
\section{Increasing sizes of independent Gaussian blocks}
\label{gb:section}

We prove Theorem~\ref{main:growing_blocks} by tracking how the
bounds for an entire block depend on its size. The fixed-size result alone
does not give the rate in \eqref{main:growing_envelope}. We first divide
each block into groups of nearly identical coordinates and replace each
group by copies of one of its members. We then bound the share of the
remaining representatives. Two estimates make this possible. The gap between groups decreases
at most exponentially with block size. The lower bound for a pair
deficit has a constant independent of that gap. We also use the finite signed-measure estimates from
Sections~\ref{sec:lm:g:5} and~\ref{sec:lm:g:6}. In those estimates, the powers of the tail
quantile stay fixed as the size of a connected set increases.
Section~\ref{gb:outside} collects the outside estimates, and
Sections~\ref{gb:near_section}--\ref{gb:capped_section} rework
Lemmas~\ref{mx:near}, \ref{mx:transport} and~\ref{mx:capped}, keeping the
dependence on block size and on the correlation gap $h\ge h_Q$ of the
partition in Lemma~\ref{gb:partition}. The relative errors carry decaying
factors such as $t^{h/256}$ and $t^{h/6}$, while their prefactors grow at
most exponentially in $Q$ and polynomially in $\log m$ and $L$.
Section~\ref{gb:assembly} shows that \eqref{main:growing_envelope} makes
the decay dominate, uniformly in the level. Section~\ref{gb:other_bases}
replaces the base $31/10$ by any fixed base $B>B_*$.

Throughout this section, fix $0<a<1$ and a sequence $(Q_m)$ satisfying
\eqref{main:growing_envelope}. Write
\begin{equation}\label{gb:notation}
 \begin{gathered}
 Q=Q_m,\quad t=\alpha/m,\quad L=\log(1/t),\qquad 0<\alpha\le a,\\
 H_\alpha=\log(1/\alpha)-1+\alpha,\qquad c_u=\Phi^{-1}(1-u/2).
 \end{gathered}
\end{equation}
Here $H_\alpha$ is the quantity written $h_\alpha$ in Appendices~\ref{sec:lm:graph} and~\ref{mx:section}; below, $h$ denotes a correlation gap. Constants denoted by $C_a$ may depend on $a$ but on no varying parameter;
$C$ denotes an absolute constant. All eventual statements below hold
uniformly over $0<\alpha\le a$ and the allowed block laws. Their
dimension thresholds depend only on $a$ and $(Q_m)$: the size conditions
used inside the proofs follow from \eqref{main:growing_envelope} alone and
are checked in Section~\ref{gb:assembly}. We use the
full variation norm of a signed measure, so a probability measure has
norm one.

For an original block of tests $A$, let $V$ be all coordinates outside
that block. In every comparison for $A$, we keep the law of $V$
unchanged and independent of $A$. With the largest physical crossing
rank $R$, define $T_A$ and $b_A$ by \eqref{mx:shares}. These quantities
measure the block's contribution to rejection in the full testing
problem. The identity
$N(tR)=R$ at a positive largest crossing gives
\begin{equation}\label{gb:shares}
 0\le T_A\le1,\qquad F_\Sigma(\alpha)=\sum_A b_A.
\end{equation}
These statements include ties and exact duplicates. An independent
uniform representative $g$ with $w$ physical copies has expected share $wt$.
Indeed, after those copies are forced to zero, the resulting forced rank $R^{[g]}$
is at least $w$. Restoring the representative below $tR^{[g]}$ preserves
that rank. Conditional on the remaining p-values, its expected share is
therefore $(w/R^{[g]})tR^{[g]}=wt$. This argument allows arbitrary
outside dependence and always uses the original thresholds $kt$.

\subsection{Outside counts and a common cutoff}
\label{gb:outside}

We first bound the number of small p-values outside a selected block.
The outside consists of independent blocks of size at most $Q$, with
uniform one-coordinate marginals. A crossing at a large rank requires many small p-values from the
outside, even if every test in the selected block passes its threshold. We bound that probability without imposing any
condition on dependence within an outside block. A block count lies in $[0,Q]$.
We bound the exponential by the line joining its endpoint values.
Independence between blocks then gives
\[
 \Expect e^{\theta N_V(jt)}
 \le \exp\{\alpha j(e^{\theta Q}-1)/Q\}.
\]
Taking $\theta=\log(1/\alpha)/Q$ yields
\begin{equation}\label{gb:counts}
 \Prob\{N_V(jt)\ge j-Q\}\le q_j,
 \qquad q_j=\min\{1,\alpha^{-1}e^{-H_\alpha j/Q}\}.
\end{equation}
For $j\ge Q/(1-a)$, the expression before taking the minimum is
increasing in $\alpha\le a$: its logarithmic derivative has the sign of
$j(1-\alpha)-Q$. Summing geometric series consequently gives, for every $K\ge0$,
\begin{equation}\label{gb:count_sums}
 \sum_{j=1}^m q_j\le C_aQ,\qquad
 \sum_{j>K}j q_j\le C_a\{Q(K+1)+Q^2\}
                    \alpha^{-1}e^{-H_\alpha K/Q}.
\end{equation}

We also need an outside event on which a selected block can determine
the crossing rank by itself. The event excludes small outside p-values
at early ranks and leaves too few at later ranks for an additional
crossing. Its probability must not decrease too quickly as the allowed
block size grows. Put $\beta=(1+a)/2$, $r_a=(1+\beta)/(2\beta)$,
$b_a=\beta r_a=(1+\beta)/2$, and $\theta_a=\log(1/b_a)$.
Choose the following fixed positive constant:
\begin{equation}\label{gb:void_constants}
 \begin{gathered}
 D_a=\max\left\{1,\frac{\log2}{H_{b_a}(r_a-1)},
                 \frac{\theta_a+1+\log4}{H_{b_a}}\right\},\\
 \ell=\lceil D_aQ\rceil,\qquad
 \gamma_{Q,a}=\tfrac12e^{-2a\ell}.
 \end{gathered}
\end{equation}
Provided $Qa\ell/m\le1/2$ and
$a/(1-Qa\ell/m)\le\beta$, there is an event $\mathcal E_V$ of
probability at least $\gamma_{Q,a}$ on which
\begin{equation}\label{gb:void}
 N_V(\ell t)=0,\qquad N_V(jt)+Q<j\quad(\ell<j\le m).
\end{equation}
The event of no background p-values at or below \(\ell t\) has probability at least
$\prod_B(1-|B|\ell t)\ge e^{-2a\ell}$. Conditioning each independent
block on its part of this event preserves block independence and gives
conditional mean count at most $\beta j$ at threshold $jt$.
To control this probability cost as block size grows, we check
the later ranks in geometric intervals $(\ell r_a^k,\ell r_a^{k+1}]$.
Monotonicity of the count lets us bound all failures in one interval by
a single event at its upper endpoint. A failure at a
physical rank in such an interval implies
\[
 N_V(J_kt)>\ell r_a^k-Q,\qquad
 J_k=\min\{m,\lceil\ell r_a^{k+1}\rceil\}.
\]
Applying the same exponential bound conditionally, with parameter
$\theta_a/Q$, bounds the probability of a failure in this interval by
\[
 \exp\{\theta_a+\beta(e^{\theta_a}-1)/Q\}
 \exp\{-[\theta_a-\beta r_a(e^{\theta_a}-1)]\ell r_a^k/Q\}
 \le e^{\theta_a+1}e^{-H_{b_a}D_ar_a^k}.
\]
Here $\beta(e^{\theta_a}-1)=1/r_a-\beta<1$ and the bracket equals
$H_{b_a}$. Since $r_a^k\ge1+k(r_a-1)$, the sum over $k\ge0$ is at most
\[
 \frac{e^{\theta_a+1-H_{b_a}D_a}}
 {1-e^{-H_{b_a}D_a(r_a-1)}}\le\frac12.
\]
Including empty intervals only enlarges this bound. The clipped last
interval uses no rank above $m$. This proves \eqref{gb:void}. On this event,
inserting $w\le Q$ zero p-values makes the largest crossing rank
exactly $w$, including
when there are ties or exact copies. The cost is now
\begin{equation}\label{gb:void_cost}
 \ell\le D_aQ+1,\qquad
 \log(1/\gamma_{Q,a})\le2aD_aQ+2a+\log2.
\end{equation}

We next choose a rank cutoff that makes the total contribution from
later ranks negligible in the share comparisons. Below it, we use
explicit Gaussian threshold estimates. The cutoff counts physical
tests and is defined by
\begin{equation}\label{gb:cutoff}
 K=\max\{2Q,\lceil21QL/H_\alpha\rceil\},\qquad x=Kt.
\end{equation}
Uniformly in the eventual range considered here,
\begin{gather}
 K\le C_aQ(1+\log m)<m,\qquad
 L\le c_{wt}^2\le2L\quad(1\le w\le Q),\qquad
 \frac{1999}{1000}L\le c_x^2\le2L,\label{gb:quantiles}\\
 c_{jt}^2\ge L/2\quad(j\le\sqrt m),\qquad
 \alpha^{-1}e^{-H_\alpha K/Q}\le t^{20}.
 \label{gb:cutoff_tail}
\end{gather}
For the bound on $K$, write $u=\log(1/\alpha)$. The ratio
$(u-1+e^{-u})/u$ has a positive lower bound for
$u\ge\log(1/a)$, and $H_\alpha\ge H_a>0$. Hence
$L/H_\alpha\le C_a(1+\log m)$. Mills' inequalities give the quantile
bounds because $\log Q+\log\log m=o(\log m)$ under our envelope.
For the upper bound one can directly use
$2\barPhi(z)\le e^{-z^2/2}$, which gives $c_t^2\le2L$.
For the lower bounds, the lower Mills inequality in \eqref{sl:mills}, together with $c_u^2\le2\log(1/u)$, gives
$c_u^2\ge2\log(1/u)-\log\{2\log(1/u)\}-\log(2\pi)$ whenever $c_u\ge1$; no
bound on the multiplier of $t$ in $u$ is needed. Taking $u=Kt$ and using
$\log K=o(L)$ gives the displayed lower bound for $c_x^2$, uniformly for
$L\ge\log(m/a)$.
The last inequality uses $\alpha^{-1}\le t^{-1}$. None of these bounds
requires $L=O(\log m)$.

\subsection{A near-copy block with its size dependence retained}
\label{gb:near_section}

We first quantify how much a block of nearly identical statistics falls
short of its nominal share. The bound must shrink to zero when the
statistics become exact copies, and its constant must retain the cost
of increasing block size. For aligned coordinates in a block $C$, write
$\delta=\max_{i,j\in C}(1-\rho_{ij})$ and define
\begin{equation}\label{gb:G}
 G_t(\delta)=\min\{\sqrt{\delta L},1\}e^{-12\delta L},
 \qquad G_t(0)=0.
\end{equation}
For the outside vector, the required structure is independent blocks
with uniform one-coordinate marginals and size at most $Q$.

\begin{lemma}\label{gb:near}
Let $|C|=w\le Q$, let $C$ be independent of its outside, and suppose
$0<\delta\le1/16$. There are constants $a_{Q,a}>0$ with
$\log(1/a_{Q,a})\le C_aQ$ such that, uniformly eventually,
\begin{equation}\label{gb:near_bound}
 b_C\le wt-a_{Q,a}tG_t(\delta).
\end{equation}
At zero spread, $b_C=wt$ exactly. The estimate holds, in particular,
whenever $Q\log(Q+1)=o(\log m)$.
\end{lemma}

\begin{proof}
The proof has three parts. We first remove the absolute values to obtain
a monotone comparison, then find a positive decrease using a pair attaining
the maximum spread, and finally bound the error from restoring the absolute
values. These are the steps of
Sections~\ref{sec:lm:n:2}--\ref{sec:lm:n:5} and Lemma~\ref{mx:near}, where
the bound on block size is fixed; we now give the constants explicitly to
track block size. Put
\begin{equation}\label{gb:near_constants}
 b_0=\barPhi(1),\quad A_Q=\sqrt{4Q/b_0},\quad
 b_Q=\min\left\{\tfrac12,
 \frac{e^{-5A_Q^2}}{(2+\sqrt2 A_Q)(Q-1)}\right\},\qquad
 a_{Q,a}=\gamma_{Q,a}b_0b_Q/2.
\end{equation}
Fix an anchor $i$, set $p=mP_i$, and use scalar Gaussian regression.
Conditional on $p$, the other absolute coordinates have the law
\[
 |r_jc_{p/m}+\varepsilon_j|,\qquad
 r_j\ge1-\delta,\quad \Expect\varepsilon_j^2=1-r_j^2\le2\delta.
\]
The residual vector is independent of the anchor and the outside;
its covariance may be singular. Symmetry includes both anchor signs,
and $p$ has density $1/m$ on $[0,m]$.

Force the anchor's p-value to zero. First replace its siblings by the
unfolded scaled values $2m\barPhi(r_jc_{p/m}+\varepsilon_j)$. Here ``unfolded''
means that the Gaussian argument is used without its absolute value.
Their forced rank
$H^+(p)$ is positive and nonincreasing. Extend it constantly for $p\ge m$;
the selected integrand below is zero for $p>\alpha m$. For every
positive nonincreasing step function $H$,
\begin{equation}\label{gb:integral}
 I(H)=\int_0^\infty
       \frac{\mathbf1\{p\le\alpha H(p)\}}{H(p)}\,dp\le\alpha.
\end{equation}
Indeed, the selected set is an initial interval with endpoint $p_*$;
if $k$ is its last selected height, then $p_*\le\alpha k$ and
$H\ge k$ on that interval. If $H=w$ up to $v<\alpha w$ and
$H\le w-1$ thereafter, splitting at $v$ also gives
\begin{equation}\label{gb:integral_loss}
 \alpha-I(H)\ge
 \min\{\alpha-v/w,\ v/[w(w-1)]\}.
\end{equation}
These are the scalar integral bounds
\eqref{eq:lm:n:2.2}--\eqref{eq:lm:n:2.3}.

To obtain a strict decrease, take a pair attaining the maximum spread and use
either endpoint as anchor.
The event
\[
 \varepsilon_j\le-\sqrt\delta,\qquad
 \max_l|\varepsilon_l|\le A_Q\sqrt\delta
\]
has probability at least $b_0/2$: the attaining variance is
$\delta(2-\delta)\ge\delta$, and the union bound for the second
restriction costs at most $2Q/A_Q^2=b_0/2$. On this event and
$\mathcal E_V$, the forced rank equals $w$ until the last anchor
value $v$ that keeps all siblings below $wt$, and is smaller afterward.
Writing $z=c_{wt}$, the threshold for $v$ lies between
$z+\sqrt\delta$ and $z+2\delta z+2A_Q\sqrt\delta$, as shown before
\eqref{eq:lm:n:5.2}. Normal hazard and
Mills bounds give
\begin{equation}\label{gb:v_bounds}
 \frac{e^{-5A_Q^2}}{2+\sqrt2 A_Q}e^{-4z^2\delta}
 \le\frac v{\alpha w}\le e^{-z\sqrt\delta}.
\end{equation}
For the lower bound, the exponent from the upper threshold is bounded
by $4z^2\delta+5A_Q^2$, and its Mills prefactor is at least
$(2+\sqrt2 A_Q)^{-1}$. The upper bound follows by integrating the
hazard lower bound $u$ from $z$ to $z+\sqrt\delta$.
Equations \eqref{gb:integral_loss}--\eqref{gb:v_bounds},
$1-e^{-u}\ge\tfrac12\min(u,1)$, and \eqref{gb:quantiles} give an
integral loss at least
$\gamma_{Q,a}(b_0/2)\alpha b_QG_t(\delta)$ for each endpoint of that pair used as
anchor. Every other anchor has nonnegative comparison loss.
After division by $m$, the losses at the two endpoints sum to at least $2a_{Q,a}tG_t(\delta)$.

Restoring the absolute values can increase the share integrals. We
bound this positive error, retaining its dependence on the spread so
that it remains smaller than the loss just obtained. The deterministic condition
and the Gaussian difference calculation in
\eqref{eq:lm:n:3.2}--\eqref{eq:lm:n:3.5}, which lead to
\eqref{mx:folding_start}, show that a positive integral difference
at rank $j$ requires an opposite-sign sibling below $-c_{jt}$, an anchor
with $p\le\alpha j$, and at least $j-w$ outside p-values below the
rank's threshold.
For each anchor its expected positive integral error is at most
\begin{equation}\label{gb:folding}
 \alpha(Q-1)\sum_{j=1}^m q_j
       e^{-(\delta^{-1}-1/2)c_{jt}^2}.
\end{equation}
Here the Gaussian difference has variance at most $2\delta$, while
its required separation is at least $2c_{jt}$; the hazard comparison
with the anchor tail gives the displayed exponent. Independence of the
outside gives the factor $q_j$.

Split this sum into ranks below and above $\sqrt m$, and treat separately
the two spread ranges separated by $\delta L=1$. This bounds the total
block-share error by
\begin{equation}\label{gb:fold_error}
 C_aQ^3tG_t(\delta)\Omega,
 \quad \Omega=t^7+t^{15/4}/L+
       e^{-H_a\sqrt m/Q}(m^{3/4}+L^{-1/2}).
\end{equation}
In the early ranks and the case $\delta L\ge1$, the exponent after
normalization is at least seven because
$(\delta^{-1}-1/2)/2-12\delta\ge7$ for $\delta\le1/16$.
When $\delta L<1$, retain
$e^{-L/(4\delta)}/\sqrt{\delta L}\le C/L$ before using the remaining
factor $t^{15/4}$. For the late ranks, \eqref{gb:counts} and
\eqref{gb:count_sums} leave the displayed exponential in $\sqrt m/Q$.
The two spread ranges use, respectively,
\[
 \sup_{0<\alpha\le a}\alpha^{-3/4}e^{-(31/2)c_\alpha^2}<\infty,
 \qquad
 \frac{e^{-c_\alpha^2/(2\delta)}}{\sqrt{\delta L}}
 \le C_aL^{-1/2}.
\]
Thus \eqref{gb:fold_error} is uniform down to zero spread. If
$Q\log(Q+1)=o(\log m)$, then $Q^3\Omega=o(a_{Q,a})$ by
\eqref{gb:void_cost} and \eqref{gb:near_constants}. Subtracting this
error from the two anchor losses proves \eqref{gb:near_bound}. The
zero-spread statement follows from the exact independent-group share.
\end{proof}

The next step must transfer this loss to a cluster that is dependent on
the rest of its original block. In Lemma~\ref{gb:transport}, the
marked part of the transfer error carries a factor $t^{d_*h}$, where
$d_*>0.47$ is the exponent in \eqref{gb:marked}, while the loss decays like
$t^{\ell\delta}$ on spreads $\delta\le h/s^2$, with the correlation gap $h$ and gap ratio $s=21/10$ of Lemma~\ref{gb:partition}.
The transfer therefore requires $\ell/s^2<d_*$, which fails for the rate
$\ell=12$ of $G_t$, so a loss bound with a smaller rate is needed.
We obtain it on a smaller, fixed spread interval and record the constants
needed for the comparison. Write
\begin{equation}\label{gb:sharp_constants}
 \begin{gathered}
 \ell_*=\frac{103}{50},\quad D_*=\frac{25}{1764},\quad
 \varepsilon_*=\frac{1}{1000},\qquad
 A(D)=\frac{2-D}{2(1-D)^2},\\
 C(D,\varepsilon)=\frac{1}{4\varepsilon(1-D)^4}
                         +\frac{D}{2(1-D)^2},\qquad
 C_{\mathrm Y}=C(D_*,\varepsilon_*),\\
 M_Q=\frac{1+A_Q\sqrt{D_*}}{1-D_*},\qquad
 b_Q^\sharp=\min\left\{\tfrac12,
       \frac{e^{-C_{\mathrm Y}A_Q^2}}{2M_Q(Q-1)}\right\},\qquad
 a_{Q,a}^{\sharp}=\gamma_{Q,a}b_0b_Q^\sharp/2,\\
 G_{\ell,t}(\delta)=\min\{\sqrt{\delta L},1\}e^{-\ell\delta L},
 \qquad G_{\ell,t}(0)=0.
 \end{gathered}
\end{equation}
Here $A_Q$ and $b_0$ retain their definitions in
\eqref{gb:near_constants}; $A(D)$ is a scalar function, not a block.

\begin{lemma}[Sharper loss on a fixed spread interval]\label{gb:sharp_near}
Under the hypotheses of Lemma~\ref{gb:near}, if $0<\delta\le D_*$,
then, uniformly eventually,
\begin{equation}\label{gb:sharp_near_bound}
 b_C\le wt-a_{Q,a}^{\sharp} tG_{\ell_*,t}(\delta),\qquad
 \log(1/a_{Q,a}^{\sharp})
 \le\left(2aD_a+\frac{4C_{\mathrm Y}}{b_0}\right)Q
                                      +O_a(\log(Q+1)).
\end{equation}
More generally, replace $D_*,\varepsilon_*,\ell_*$ by any fixed
$0<D\le1/16$, $\varepsilon>0$, and $2<\ell\le12$ satisfying
$2[A(D)+\varepsilon]<\ell$. The same conclusion holds for
$\delta\le D$, with the corresponding constants. All zero-spread
cases remain exact equalities.
\end{lemma}

\begin{proof}
Retain the two anchors from a pair attaining the maximum spread,
and the residual event used above.
On that event and $\mathcal E_V$, the exact Gaussian threshold corresponding to the last anchor value defined above is
\[
 \zeta=\max_{j\ne i}\frac{z-\varepsilon_j}{r_j},\qquad
 z=c_{wt},\qquad v=2m\barPhi(\zeta).
\]
It satisfies
$z+\sqrt\delta\le\zeta\le T=(z+A_Q\sqrt\delta)/(1-\delta)$.
No condition on $A_Q\sqrt\delta$ is needed. For $z\ge1$ and
$\delta\le D$, Mills' inequalities and $T/z\le M_Q$ give
\[
 \frac{\barPhi(T)}{\barPhi(z)}
 \ge \frac{zT}{1+T^2}e^{-(T^2-z^2)/2}
 \ge\frac{1}{2M_Q}e^{-(T^2-z^2)/2}.
\]
Writing $u=z\sqrt\delta$, expansion of the exponent gives
\[
 \frac{T^2-z^2}{2}
 =A(\delta)u^2+\frac{A_Qu}{(1-\delta)^2}
             +\frac{A_Q^2\delta}{2(1-\delta)^2}
 \le[A(D)+\varepsilon]u^2+C(D,\varepsilon)A_Q^2.
\]
The inequality uses the monotonicity of $A$ and Young's inequality.
Thus $z^2\le2L$ gives
\[
 \frac{e^{-C(D,\varepsilon)A_Q^2}}{2M_Q}e^{-\ell\delta L}
 \le \frac{v}{\alpha w}\le e^{-z\sqrt\delta}.
\]
For the displayed numerical choices, even the weaker upper bound
$z^2\le2001L/1000$ suffices, because
\[
 \frac{103}{50}-\frac{2001}{1000}
             \left(A(25/1764)+\frac{1}{1000}\right)
 =\frac{41256347879}{3024121000000}>0.
\]
The first branch in \eqref{gb:integral_loss} is at least
$(\alpha/2)\min(\sqrt{\delta L},1)$; its second branch is at least
$\alpha e^{-C(D,\varepsilon)A_Q^2}e^{-\ell\delta L}/[2M_Q(Q-1)]$.
Their minimum is at least $\alpha b_Q^\sharp G_{\ell,t}(\delta)$.
The two anchor losses before folding therefore sum to at least
$2a_{Q,a}^{\sharp} tG_{\ell,t}(\delta)$.
Since $G_t=G_{12,t}\le G_{\ell,t}$, the absolute folding estimate
\eqref{gb:fold_error} transfers in the required direction. Its error is
$o(a_{Q,a}^{\sharp} tG_{\ell,t})$ under the same size condition as
Lemma~\ref{gb:near}. Finally the bound on $\log(1/\gamma_{Q,a})$ in \eqref{gb:void_cost},
$A_Q^2=4Q/b_0$, and $M_Q=O(\sqrt Q)$ give the stated logarithmic cost.
The argument only uses marginal normal tails of the residuals, and
therefore includes singular residual covariance and exact-copy siblings.
\end{proof}

\subsection{A quantitative partition and signed contraction}
\label{gb:partition_section}

We next partition an arbitrary positive semidefinite block. The
partition chooses its own correlation cutoff: coordinates within each
group are much closer to one another than coordinates in different
groups.

\begin{lemma}\label{gb:partition}
A correlation matrix on at most $Q$ coordinates admits a partition into
aligned components and a number $h\in[h_Q,1/16]$, where
$h_Q=(31/10)^{-2(Q-1)}/16$ and $s=21/10$, such that
\begin{equation}\label{gb:partition_bounds}
 \begin{gathered}
 \delta_C\le h/s^2\quad\text{within each component},\\
 |\rho_{iy}|\le1-h\quad\text{between different components}.
 \end{gathered}
\end{equation}
Selecting actual coordinates as component anchors preserves the latter
bound after any sequence of contractions.
\end{lemma}

\begin{proof}
The construction builds a forest from the shortest available edges.
We seek a next edge whose length exceeds a fixed multiple of the
total length already retained. That gap gives both a bound on distances within each tree
and a lower bound on distances between different trees. Represent the
correlation matrix by unit Gram vectors and use their projective
distance, which identifies a vector with its negative:
$d_{ij}=\min(\|v_i-v_j\|,\|v_i+v_j\|)=\sqrt{2(1-|\rho_{ij}|)}$.
Process edges in increasing order, accepting an edge when it creates no
cycle, and stop before length $\tau_0=1/\sqrt8$. Let the accepted lengths
be $e_1\le\cdots\le e_p<\tau_0$, where $p\le Q-1$, and set
$e_{p+1}=\tau_0$, $S_k=\sum_{j\le k}e_j$.
Choose the last index $k\in\{0,\ldots,p\}$ for which
\begin{equation}\label{gb:forest_gap}
 e_{k+1}>sS_k,
 \qquad \tau=e_{k+1},\qquad h=\tau^2/2.
\end{equation}
There is an eligible index after the initial zero edges, so $\tau>0$.
Retain the first $k$ accepted edges. Any edge joining distinct retained
components has length at least $\tau$: a shorter edge was processed
earlier and was either accepted or already joined vertices connected
by earlier edges. Align signs along each retained tree. Every aligned
within-component distance is at most $S_k<\tau/s$. This proves
\eqref{gb:partition_bounds}.

By maximality of $k$, $S_{j+1}\le(1+s)S_j$ for $j\ge k+1$,
whereas $S_{k+1}<(1+s)\tau/s$. If $k<p$, then
$\tau_0\le sS_p$, and hence $\tau_0\le(1+s)^{p-k}\tau$.
The same last inequality is immediate if $k=p$. Squaring gives
$h\ge(31/10)^{-2(Q-1)}/16$. The construction permits zero distances,
length ties, and singular Gram matrices. Anchors are original
coordinates, so their cross-component covariances retain the stated cap.
\end{proof}

For a component $C$ in this partition, the remaining coordinates of the original block
will generally depend on $C$. The next comparison therefore follows
the whole-block share while one component is replaced by anchor copies.

\begin{lemma}\label{gb:transport}
Let a block $A$, $|A|\le Q$, be independent of its outside $V$; here $A$
is an original block, possibly after earlier components have been
contracted to actual anchors.
Suppose a component $C\subset A$ has aligned spread
$0<\delta\le h/s^2$, $s=21/10$, and absolute cross-correlations
with $Y=A\setminus C$ bounded by $1-h$, where $h_Q\le h\le1/16$. Let $b_A^{\mathrm{con}}$ be its
share after contracting $C$ to an actual aligned anchor, with all
physical slots and the outside law retained. Uniformly eventually,
\begin{equation}\label{gb:transport_bound}
 b_A\le b_A^{\mathrm{con}}-\tfrac12a_{Q,a}^{\sharp}tG_{\ell_*,t}(\delta).
\end{equation}
For $\delta=0$, contraction changes no p-value and no share.
\end{lemma}

\begin{proof}
We compare the effect of contraction under the actual law with its
effect under a reference distribution. The reference replaces the rest
of the block by independent Gaussian singletons, retaining every test
slot and the unchanged outside law.
Under the reference, the near-copy estimate makes the original share
smaller than its contracted share. We bound the error in transferring
this comparison of the original and contracted shares to the actual law.
The path and crossing estimates below bound the difference between the
actual and reference contraction increments, first at ranks up to $K$ and
then at later ranks. We call the component $C$ the cluster. For a member of $C$ and its anchor $X$, scalar
regression gives a path between the contracted and original coordinates,
\begin{equation}\label{gb:path}
 Z_i(v)=\sqrt{1-v^2\sigma_i^2}\,X+v\sigma_i E_i,
 \quad 0\le v\le1,
 \quad \sigma_i^2=1-\rho_{iX}^2\le2\delta,
 \quad E_i\perp X.
\end{equation}
The residuals may depend on $Y$, but the entire path is independent of
$V$. The path has constant unit variance, so its velocity is independent of
its current coordinate, and $\Expect\dot Z_i(v)^2\le4\sigma_i^2$.
In the Gram representation, the path follows the shorter spherical arc
between $X$ and $Z_i$. Writing a point as a
nonnegative linear combination of these two endpoints, the sum of its
coefficients is at most $(1-\delta/2)^{-1/2}$. The cap at both physical
endpoints consequently gives
\[
 |\operatorname{Corr}(Z_i(v),Y_j)|
 \le\frac{1-h}{\sqrt{1-\delta/2}}
 \le1-h+\delta/4\le1-\kappa h,\qquad
 \kappa=1-\frac{1}{4s^2}=\frac{416}{441}.
\]
For the middle inequality, put $u=\sqrt{1-\delta/2}$ and use
$u(1+u)\ge2-\delta\ge2(1-h)$. This also proves the assertion after
previous components have been contracted to actual physical anchors.

Suppose for now that $\log K\le\eta_*hL$, where
$\eta_*=1/10000$. Section~\ref{gb:assembly} will verify this finite
condition and absorb the errors below. The hazard bound gives
$c_t^2-c_x^2\le2\log K$, so $c_t\le(1+\eta_*h)c_x$ because
$c_x^2\ge L$. For $j\le K$ and $0\le r\le1-\kappa h$,
\begin{equation}\label{gb:conditional_gap}
 \frac{c_x-r c_{jt}}{\sqrt{1-r^2}}
 \ge \frac{\kappa-\eta_*}{\sqrt{2\kappa}}\sqrt h\,c_x.
\end{equation}
Indeed, $[1-(1+\eta_*h)r]/\sqrt{1-r^2}$ decreases on this interval;
at its right endpoint its numerator is at least
$(\kappa-\eta_*)h$, while its denominator is at most
$\sqrt{2\kappa h}$.

The two contraction comparisons can differ at early ranks only when a
cluster threshold changes and another coordinate in the block has a
small p-value. We therefore need a joint bound for these two events.
Let $M_{ij}$ be the event that coordinate $i$'s threshold indicator at
rank $j$ differs between the endpoints of \eqref{gb:path}. Counting
crossings in the Gaussian rotation plane, as in
Lemma~\ref{fd:marked_crossing} after the substitution
$\theta=\arcsin(v\sigma_i)$ used in the proof of Lemma~\ref{mx:transport},
gives, for a fixed marked event $E$,
\[
 \Prob(M_{ij}\cap E)
 \le\sum_{z=\pm c_{jt}}\int_0^1
 \phi(z)\Expect[|\dot Z_i(v)|\mathbf1_E\mid Z_i(v)=z]\,dv.
\]
The crossing probability depends on the path's velocity as well as on
the marked event. We bound their joint contribution directly to retain
the stronger exponent. Given the current coordinate,
let $U$ be the standardized velocity and $W$ the standardized centered
marked coordinate. They are a possibly singular
Gaussian pair, so $U=bW+\sqrt{1-b^2}N$ for a standard normal $N$
independent of $W$. For $u\ge0$,
\[
 \Expect[|U|\mathbf1\{W\ge u\}]
 \le |b|\phi(u)+\sqrt{1-b^2}\sqrt{2/\pi}\barPhi(u)
 \le C e^{-u^2/2}.
\]
The estimate includes $|b|=1$; zero velocity contributes zero without
standardization. Apply it to both tails using \eqref{gb:conditional_gap},
$c_x^2\ge1999L/1000$, and $\phi(c_{jt})\le C_a jt\sqrt L$.
This gives
\begin{equation}\label{gb:marked}
 \Prob(M_{ij},|Y_l|\ge c_x)
 \le C_a jt\sqrt{\delta L}\,t^{d_*h}\qquad(j\le K),\qquad
 d_*=\frac{1999}{1000}\frac{(\kappa-\eta_*)^2}{4\kappa}.
\end{equation}
Numerically $d_*>0.47$.
The velocity and marked event need not be independent. The same
crossing calculation without the mark gives, with
$M_j=\bigcup_{i\in C}M_{ij}$,
\begin{equation}\label{gb:unmarked}
 \Prob(M_j)\le C_a Qjt\sqrt{\delta L}\qquad(1\le j\le m).
\end{equation}

We now specify the actual and reference expected shares. Write $H_1,H_0$ for the original and contracted expected whole-$A$ shares with
$(Y,V)$ unchanged; the subscripts label the two laws, not values of the rate $H_\alpha$. Write $B_1,B_0$ for the corresponding expected shares when
$Y$ is replaced by independent Gaussian singletons, independent of
$(C,V)$. The exact singleton and group shares, and
Lemma~\ref{gb:sharp_near}, give
\begin{equation}\label{gb:four_shares}
 B_0=|A|t,\qquad B_1-B_0\le-a_{Q,a}^{\sharp}tG_{\ell_*,t}(\delta).
\end{equation}
We now compare the actual contraction increment with the reference
increment. If the $Y$
p-values are set to one, with their slots retained, the changes in the
early share are the same function of $(C,V)$. A difference from this
common function requires both a changed threshold indicator in the
cluster and a p-value from $Y$
below $x$. Summing \eqref{gb:marked} over coordinates and ranks therefore
bounds the difference between the two changes in the early share.

For late ranks, there are two possibilities: a late crossing can
change, or a late crossing can remain present while the early share
changes. Let $g_K$ be the share defined using the largest
crossing through rank $K$, write $\Delta$ for the change between the two
endpoints of \eqref{gb:path}, and put
$\mathcal O_j=\{N_V(jt)\ge j-|A|\}$. The deterministic bound
\eqref{mx:late_deterministic} reads here
\begin{equation}\label{gb:late_deterministic}
 |\Delta(T_A-g_K)|
 \le2\sum_{j>K}\mathbf1_{M_j}\mathbf1_{\mathcal O_j}
  +2\sum_{k\le K}\mathbf1_{M_k}
                       \mathbf1_{\bigcup_{j>K}\mathcal O_j}.
\end{equation}
The first sum covers a changed late crossing; the second covers an
unchanged late crossing with a changed early share. If the largest rank
is unchanged, the full share is unchanged because
$N_A(tR)=R-N_V(tR)$. These observations prove the inclusion and show why
both terms are needed. Since the path and $V$ are independent,
\eqref{gb:count_sums}, \eqref{gb:cutoff_tail}, and
\eqref{gb:unmarked} apply to each term. The same bounds apply to the
comparison with independent $Y$.

Together, the early and late estimates give an error bound relative to
the near-copy loss:
\begin{equation}\label{gb:transport_error}
 \frac{|(H_1-H_0)-(B_1-B_0)|}{tG_{\ell_*,t}(\delta)}
 \le C_aQ^4K^2(1+\sqrt L)
                     \{t^{h/256}+t^{19}\}.
\end{equation}
Indeed, division of an early bound by $G_{\ell_*,t}(\delta)$ costs at most
$(1+\sqrt L)t^{-\ell_*\delta}$. The prefactor comes from counting the coordinates and ranks in
these bounds. The early
terms involve at most $Q^2$ pairs of a cluster path and a marked
coordinate, together with $\sum_{j\le K}j\le K^2$. The late terms use
\eqref{gb:unmarked}, \eqref{gb:count_sums} and the geometric bound
$\sum_{j>K}q_j\le C_aQ\alpha^{-1}e^{-H_\alpha K/Q}\le C_aQt^{20}$.
The exact margin is
\[
 d_*-\frac{\ell_*}{s^2}-\frac{1}{256}
 =\frac{7170072629173}{24460800000000000}>0,
\]
so the normalized early exponent is at least $h/256$.
The late bound contains $t^{20}$, and
$\ell_*\delta\le\ell_*D_*=103/3528<1$, leaving exponent at least
nineteen.
Section~\ref{gb:assembly} shows that the right side is at most
$a_{Q,a}^{\sharp}/2$. Equation~\eqref{gb:four_shares} then proves the lemma.
No division is made at $\delta=0$, where the coupled endpoints coincide.
\end{proof}

\subsection{A capped block with a varying correlation gap}
\label{gb:capped_section}

After these contractions, a selected original block still has at most
$Q$ physical coordinates represented by at most $Q$ distinct Gaussian
variables. Each representative retains its positive integer
multiplicity, their total is at most $Q$, and their pairwise absolute
correlations are at most $1-h$. We call such a block terminal.
Under the same outside law $V$, we prove \eqref{gb:capped_bound}: its
whole-block share is at most $|A|t-D/2$, where $D$ is the sum of the
one-edge deficits defined below, and it equals $|A|t$ exactly when the
representatives are independent. We track the dependence on the
correlation gap in the finite estimates behind Lemma~\ref{mx:capped},
because that gap can now decrease with dimension.

For an edge $e$ between representatives of multiplicities $u,v$, let
$b_A^e$ be the whole-$A$ share when this is the only retained internal
edge. All other representatives are independent in this reference, and
$V$ retains its actual law. Set
\begin{equation}\label{gb:edge_deficit}
 d_e=|A|t-b_A^e,\qquad D=\sum_e d_e.
\end{equation}
Edges with zero correlation may be omitted. Each deficit concerns the
whole block's share under its specified one-edge reference law. The
deficits in the full rejection probability used in the graph expansion
are different quantities.

\subsubsection{Pair curvature}
\label{gb:curvature_section}

Write $J_r(x,y)=\Prob(P_1\le x,P_2\le y)$ for a folded Gaussian pair
with absolute correlation $r$. Its two-group deficit is
\[
 D_{u,v,r}(t)=J_r(ut,(u+v)t)+J_r(vt,(u+v)t)
                         -J_r((u+v)t,(u+v)t).
\]
We first obtain a positive deficit for the two groups alone, then use
the background event to turn it into a deficit for their share in the
full problem. The constant in the first bound remains fixed as $h$ decreases.
This is needed because the partition's correlation gap can become small
as block size grows.
Let $0<r\le1-h$, $z=c_{(u+v)t}\ge1$, and assume $hz^2\ge1$.
For $f(p)=J_r(p,(u+v)t)$ and $y=c_p$, differentiation as in
\eqref{eq:lm:g:3.1} gives
\begin{equation}\label{gb:curvature}
 -f''(p)=\frac r{2\sqrt{1-r^2}}
 \exp\left\{\frac{y^2}{2}-\frac{(z-ry)^2}{2(1-r^2)}\right\}
 \left(1-e^{-2rzy/(1-r^2)}\right).
\end{equation}
The deficit is
$\int_0^{ut}\int_0^{vt}-f''(p+q)\,dq\,dp$.
Restrict to $p\in[ut/4,ut/2]$ and $q\in[vt/4,vt/2]$.
Then $c_{p+q}=z+e$ with $0\le e\le\log4/z$, by the hazard lower
bound. The exact exponent identity is
\[
 \frac{(z+e)^2}{2}-\frac{(z-r(z+e))^2}{2(1-r^2)}
 =\frac{rz^2}{1+r}+\frac{1+2r}{1+r}ze
   +\frac{1-2r^2}{2(1-r^2)}e^2
 \ge\frac{rz^2}{1+r}-\frac{(\log4)^2}{2hz^2}.
\]
The rectangle has area $uvt^2/16$. Moreover,
$1-e^{-2rzy/(1-r^2)}\ge(1-e^{-2})\min(1,rz^2)$.
Consequently
\begin{equation}\label{gb:bare_pair}
 D_{u,v,r}(t)\ge c_*uvt^2r\min(1,rz^2)e^{rz^2/(1+r)},
 \qquad c_*=\frac{(1-e^{-2})e^{-(\log4)^2/2}}{32}>0.
\end{equation}

To pass from the two-group deficit to $d_e$, augment the background $V$ by all the other
independent representative groups in the one-edge reference. Each added
group has at most $Q$ physical copies, so \eqref{gb:void} applies to
this augmented background with the same $\gamma_{Q,a}$. The original
outside law is unchanged. For a uniform representative of weight $w_i$,
the forcing identity \eqref{eq:lm:g:2.2} expresses its share deficit as a sum of differences between
conditional rank probabilities:
\[
 w_it\sum_{k<m}\left[
 \Prob(R^{[i]}\le k\mid P_i\le(k+1)t)
 -\Prob(R^{[i]}\le k\mid P_i\le kt)\right].
\]
Only the paired representative depends on this representative. A folded bivariate Gaussian makes the
partner p-value stochastically increasing with the conditioning p-value;
the displayed forced-rank events are increasing in that partner.
Thus all increments are nonnegative, conditionally on the augmented
background. On the event \eqref{gb:void}, the relevant endpoint ranks are exactly
$u,v,u+v$, and the two representative deficits sum to
$D_{u,v,r}(t)$. The other groups have their exact nominal shares. Hence
\begin{equation}\label{gb:pair_void}
 d_e\ge\gamma_{Q,a}D_{u,v,r_e}(t).
\end{equation}
Using \eqref{gb:bare_pair}, Mills' inequalities, $u+v\le Q$, and
\eqref{gb:quantiles}, we obtain
\begin{equation}\label{gb:pair_lower}
 d_e\ge b_{Q,a}r_e^2Lt^2,\qquad
 d_e\ge b'_{Q,a}r_e^2L^{-1}t^{2/(1+r_e)},
 \qquad
 b_{Q,a}=\gamma_{Q,a}c_*/2,\quad
 b'_{Q,a}=b_{Q,a}/(10Q).
\end{equation}
For the first inequality use
$e^{rz^2/(1+r)}\min(1,rz^2)\ge rz^2/2$.
For the second, $r\min(1,rz^2)\ge r^2/2$. Put $w=u+v$ and
$\kappa_r=2r/(1+r)\le1$. The lower Mills inequality gives
$e^{z^2/2}\ge1/(\sqrt{2\pi}zwt)$ for $z\ge1$. Since
$\sqrt{2\pi}zQ\ge1$ and $z^2\le3L$, it follows that
\[
 e^{rz^2/(1+r)}
 \ge\frac{t^{-\kappa_r}}{\sqrt{2\pi}zQ}
 \ge\frac{t^{-\kappa_r}}{10QL}.
\]
This proves the displayed choice of $b'_{Q,a}$. The logarithms of the
reciprocals of both constants are $2aD_aQ+O_a(\log(Q+1))$.
The first lower bound will control errors when correlations are very
small; the second retains the tail exponent needed to compare larger
connected sets with a single edge. Both lower bounds hold for arbitrarily small
positive $r_e$; no absolute error has been substituted for a vanishing
deficit.

\subsubsection{Finite signed measures}
\label{gb:finite_measures}

To study crossings through the cutoff rank, we need the small p-values
and their ranks, but not the exact values of larger p-values. Replace
every representative p-value above $x=Kt$ by a single symbol
$*$, let $\nu$ be the resulting censored uniform law, and let $\mu_T^x$ be the joint censored law of the representatives in a set $T$. We then separate
the departure from independence into contributions from representative
sets. For a representative set $C$ define
\begin{equation}\label{gb:centered}
 A_C=\sum_{T\subseteq C}(-1)^{|C|-|T|}
             \mu_T^x\otimes\nu^{\otimes(C\setminus T)}.
\end{equation}
The finite identities \eqref{eq:lm:g:5.2}--\eqref{eq:lm:g:5.3}
reconstruct the full censored law from these terms. Singletons vanish,
and a disconnected set factors over its Gaussian covariance components.
The expansion is therefore over compatible collections of connected
sets of at least two vertices. The factorization uses independence of
Gaussian subvectors with no covariance edges between them; it does not
condition on an outside vertex.

For the block share, the single-edge terms supply the negative correction. Larger
connected sets and products of edge terms form the error that must be
smaller than that correction. The following estimates retain all dependence
relevant to $Q,h$:
\begin{align}
 \|A_e\|&\le C r_e^2c_x^2
                          e^{-c_x^2/(1+r_e)},\label{gb:pair_derivative}\\
 \|A_e\|&\le156c_x^2x^{2/(1+r_e)},\label{gb:pair_absolute}\\
 \|A_C\|&\le216\,4^{|C|}c_x^3
                     x^{3/(1+2r_C)},\quad |C|\ge3,\label{gb:triple}\\
 \|A_C\|&\le24 S_C|C|^2c_x^4(8000x)^{|C|}
 \quad\text{if }\eta_C\le\min(1/8,c_x^{-2}).\label{gb:weak}
\end{align}
Here $r_C$ is the largest internal absolute correlation,
$S_C=\operatorname{tr}(\Gamma_C-I)^2$, and
$\eta_C=\max_i\sum_{j\ne i}|\rho_{ij}|$.
For \eqref{gb:pair_derivative}, when $r_e\le1/2$ use
\eqref{eq:lm:g:3.4} with the fixed cap $1/2$: twice differentiate the
folded pair density, cancel its first derivative at zero by sign
averaging, and integrate \eqref{eq:lm:g:3.5}. When $r_e\ge1/2$, the
centered kernel and the signed-sum tail bound instead give
\[
 \|A_e\|\le4\Prob(P_i,P_j\le x)+12x^2
 \le28e^{-c_x^2/(1+r_e)}
 \le112r_e^2c_x^2e^{-c_x^2/(1+r_e)}.
\]
Here $x\le e^{-c_x^2/2}$ and $c_x\ge1$. Thus the pair derivative
bound has an absolute constant even when the terminal gap tends to zero.

For a product of pair terms, one factor will use an absolute bound
that does not retain the squared correlation. Equation~\eqref{gb:pair_absolute}
is the separate absolute estimate already given in \eqref{eq:lm:g:5.7}.
Expanding the two centered censored coordinates gives
$\|A_e\|\le4\Prob(P_i,P_j\le x)+12x^2$.
A union over four signed sums bounds the joint probability by
$4e^{-c_x^2/(1+r_e)}$. Since
$e^{-c_x^2/2}\le3c_xx$ and $c_x\ge1$, this probability is at most
$36c_x^2x^{2/(1+r_e)}$, which proves \eqref{gb:pair_absolute}.

For \eqref{gb:triple}, use the finite kernel expansion
\eqref{eq:lm:g:5.4}. Any three retained low coordinates have probability
at most $8e^{-3c_x^2/(2(1+2r_C))}$ by their eight signed sums, and hence
at most $216c_x^3x^{3/(1+2r_C)}$. Terms with two retained lows have an
extra factor $x$, and
$1+2/(1+r)\ge3/(1+2r)$. Terms with fewer lows also have sufficient
powers of $x$. Summing the subsets gives the factor $4^{|C|}$.
This proof uses no nonsingularity assumption.

Finally, the proof of \eqref{eq:lm:g:6.1} gives \eqref{gb:weak} with
its displayed constants for every set size. To make the dependence on set size
explicit, fix a $j$-coordinate set $T\subseteq C$ and put $S_T=\operatorname{tr}(\Gamma_T-I)^2$.
Along $I+v(\Gamma_C-I)$, the log density $\mathcal L$ of the marginal of $T$ satisfies
\[
 |\mathcal L''+(\mathcal L')^2|
 \le16S_T(j+1+\|z\|^4),\qquad
 g_v(z)/g_0(z)\le e^j e^{\eta_C\|z\|^2}.
\]
The eigenvalues lie in $[7/8,9/8]$. For $p=0,2,4$ and a standard normal $N$, substituting
$z=c_x+u/c_x$ and using $\eta_Cc_x^2\le1$ gives
\[
 \Expect[|N|^p e^{\eta_CN^2}\mathbf1\{|N|\ge c_x\}]
 \le1300c_x^p x.
\]
The remaining integral is bounded using
$\int_0^\infty(1+u)^4e^{-3u/4}\,du<214$.
Sign averaging removes the first derivative at independence. Taylor's
remainder gives $24S_Tj^2c_x^4(3900x)^j$ on the low marginal, and the
finite centered expansion gives \eqref{gb:weak}, since
$2(1+3900)<8000$. In particular, the power $c_x^4$ is independent of
$|C|$.

When every correlation is extremely small, we compare the uncensored
probability distributions directly. This avoids a separate truncation error that might exceed the
vanishing pair deficit. We use the full-distribution estimate for very small correlations in \eqref{eq:lm:g:4.1}.
For zero-diagonal $R_0$ with $\|R_0\|_{\mathrm{op}}\le1/2$, it implies the weaker bound
\begin{equation}\label{gb:microscopic}
 \left\|\mu_{R_0}-\mu_0-\sum_e(\mu_e-\mu_0)\right\|
 \le30\{\operatorname{tr}(R_0^2)\}^{3/2}.
\end{equation}
Its constant is independent of dimension. Indeed, along $I+vR_0$, put
$B=R_0(I+vR_0)^{-1}$. The Gaussian density log derivatives are
\[
 \mathcal L'=\tfrac12(Y^TBY-\operatorname{tr}B),\quad
 \mathcal L''=\tfrac12\operatorname{tr}B^2-Y^TB^2Y,\quad
 \mathcal L'''=-\operatorname{tr}B^3+3Y^TB^3Y.
\]
The quadratic-form moment bound used in Section~\ref{sec:lm:g:4}
gives $\|g_v'''\|_1\le11\|B\|_F^3\le88\{\operatorname{tr}(R_0^2)\}^{3/2}$.
The third-order remainder is at most $15\{\operatorname{tr}(R_0^2)\}^{3/2}$.
Averaging over coordinate signs removes the first-order and mixed distinct-edge
second-order terms; the remaining quadratic term is the sum of the
individual edge terms. Subtracting their exact laws costs at most another
15, using $\sum_e(2r_e^2)^{3/2}\le\{\operatorname{tr}(R_0^2)\}^{3/2}$.
This proves \eqref{gb:microscopic} with no assertion about orthogonality
of higher-order terms.

\subsubsection{Relative errors for the selected block}
\label{gb:capped_errors_section}

We first compare a single-edge measure with its positive share deficit.
Dividing \eqref{gb:pair_derivative} by
\eqref{gb:bare_pair}--\eqref{gb:pair_void}, and using
$c_{(u+v)t}^2-c_x^2\le2\log K$, gives
\begin{equation}\label{gb:pair_ratio}
 \frac{\|A_e\|}{d_e}\le\mathcal C_{Q,h}K^2L^2,
 \qquad
 \mathcal C_{Q,h}\le \frac{CQ^2}{\gamma_{Q,a}c_*},\qquad
 \log\mathcal C_{Q,h}\le2aD_aQ+O_a(\log(Q+1)).
\end{equation}
The retained subscript $h$ records the bound $1-h$ on absolute correlations
after contraction, but the displayed
constant does not depend on this gap. To check the quotient directly, put
$z=c_{(u+v)t}$ and use
$e^{-z^2}/t^2\le2\pi(u+v)^2z^2$ and
$r/\min(1,rz^2)\le(1+rz^2)/z^2$. Together with
$z^2-c_x^2\le2\log K$, these are precisely the cancellation in
\eqref{eq:lm:g:3.9}. No inverse-gap or inverse-correlation factor is
introduced.

We compare uncensored distributions directly when all correlations
are sufficiently small, and use the censored expansion otherwise. In the latter case a
nonzero lower bound for the total deficit also controls the cost of
restoring late ranks. There are two cases. If every $r_e\le t^4$, then
$\|\Gamma-I\|_{\mathrm{op}}\le Qt^4<1/2$ eventually.
Average the whole-block share over the unchanged outside law. The
result is a function bounded by one, against which we integrate
\eqref{gb:microscopic}.
Writing $S=2\sum_e r_e^2$, the error relative to $D$ is at most
\begin{equation}\label{gb:micro_relative}
 \frac{30S^{3/2}}{(b_{Q,a}/2)Lt^2S}
 \le\frac{60Q}{b_{Q,a}L}t^2.
\end{equation}
Because this comparison uses the uncensored law, its error stays relative to the
deficit however fast all correlations tend to zero. If there is no
edge, the independent representatives instead give $b_A=|A|t$ exactly.

In the second case, an edge exceeds $t^4$, so
$D\ge b_{Q,a}Lt^{10}$. Let
\[
 V_2=\sum_e\|A_e\|,\qquad
 B_3=\sum_{C\text{ connected},\ |C|\ge3}\|A_C\|.
\]
Integrate the early whole-block share over the same outside law in every
term of \eqref{gb:centered}. The result is a bounded function of the
censored representatives, and its individual edge terms are exactly the
corresponding early-share differences. The early share is bounded by
one: at a largest crossing below $K$ the count equals its rank, and at
rank $K$ its numerator is at most $Q\le K$. Dropping compatibility
among connected terms bounds the remainder beyond single edges by
$(B_3+V_2^2/2)e^{B_3+V_2}$, as in \eqref{eq:lm:g:5.10}.
The full share and the share through the cutoff can differ only if a
late crossing occurs,
which requires at least $j-Q$ outside lows. Equations
\eqref{gb:count_sums} and \eqref{gb:cutoff_tail} bound each such error
by $C_aQt^{20}$ under the original, independent, and one-edge laws.
Restoring the late terms therefore gives
\begin{equation}\label{gb:expansion}
 |b_A-|A|t+D|
 \le(B_3+V_2^2/2)e^{B_3+V_2}+C_aQ^3t^{20}.
\end{equation}

We bound the remainder from larger connected sets, products of pairs,
and late ranks in turn. Call a connected set weak when
$r_C\le[Q\max(8,c_x^2)]^{-1}$. Summing \eqref{gb:weak}, retaining
$S_C\le2\sum_e r_e^2$, and using $x^{|C|}\le x^3$ before summing
subsets, gives
\begin{equation}\label{gb:weak_relative}
 \frac{B_{3,\mathrm{weak}}}{D}
 \le C b_{Q,a}^{-1}Q^2 16000^QK^3Lt.
\end{equation}
Indeed, $|C|^2\le Q^2$, $c_x^4\le9L^2$, and the subset sum of the
remaining constants is bounded by $16000^Q$. Thus no power of $L$
depending on $Q$ is introduced.

In a strong set, $r_C>1/(3QL)$ eventually. Compare its bound
\eqref{gb:triple} with the deficit of an edge attaining the largest
absolute correlation, using the second inequality in
\eqref{gb:pair_lower}. The exponent difference satisfies
\[
 \frac3{1+2r_C}-\frac2{1+r_C}
 =\frac{1-r_C}{(1+2r_C)(1+r_C)}\ge h/6.
\]
The factor $r_C^{-2}$ costs at most $9Q^2L^2$, the quantile factor
costs $C L^{3/2}$, and the inverse deficit contributes another $L$.
The sum of $4^{|C|}$ over subsets is at most $5^Q$. Consequently
\begin{equation}\label{gb:strong_relative}
 \frac{B_{3,\mathrm{strong}}}{D}
 \le C(b'_{Q,a})^{-1}Q^2 5^QK^3L^{9/2}t^{h/6}.
\end{equation}
For the product of pair terms, use \eqref{gb:pair_ratio} on one factor
and the separate absolute bound \eqref{gb:pair_absolute} on the other:
\begin{equation}\label{gb:pair_product}
 V_2\le C Q^2LK^2t^{2/(2-h)},\qquad
 \frac{V_2^2}{D}
 \le C\mathcal C_{Q,h}Q^2K^4L^3t^{2/(2-h)}.
\end{equation}
The first factor is normalized by its own positive deficit; the second
uses an absolute tail bound. The split at correlation $1/2$ in
\eqref{gb:pair_derivative} removes the inverse-gap prefactor. Finally,
\begin{equation}\label{gb:late_relative}
 \frac{C_aQ^3t^{20}}D
 \le C_aQ^3b_{Q,a}^{-1}t^{10}/L.
\end{equation}

Section~\ref{gb:assembly} shows that all the relative errors
\eqref{gb:micro_relative} and
\eqref{gb:weak_relative}--\eqref{gb:late_relative} tend to zero.
The exponential factor in \eqref{gb:expansion} tends to one as well:
$d_e\le |A|t$ implies $D\le Q^3t$, so $B_3=o(D)=o(1)$, while the
first bound in \eqref{gb:pair_product} gives $V_2=o(1)$ directly.
It follows that every terminal capped block satisfies, uniformly
eventually,
\begin{equation}\label{gb:capped_bound}
 b_A\le |A|t-D/2\quad(D>0),\qquad
 b_A=|A|t\quad(D=0).
\end{equation}
The zero case consists of independent representatives by
\eqref{gb:pair_lower}.

\subsection{Uniform absorption and summation under the original law}
\label{gb:assembly}

We now show that the prescribed size envelope makes each error smaller
than the corresponding positive deficit, uniformly in the level.
The envelope \eqref{main:growing_envelope} implies
\begin{equation}\label{gb:absorption_budget}
 Q=O(\log\log m),\qquad \log(1/h_Q)=O(Q),\qquad
 \frac{Q+\log\log m}{h_Q\log m}\longrightarrow0.
\end{equation}
For example, the first statement follows by taking logarithms in
$h_Q\log m\to\infty$. In particular, the finite conditions for
\eqref{gb:void} hold eventually, and
$Q\log(Q+1)=o(\log m)$ as needed for Lemma~\ref{gb:near}.

The logarithm of every prefactor in
\eqref{gb:transport_error}, \eqref{gb:pair_ratio}, and
\eqref{gb:micro_relative}--\eqref{gb:late_relative}, including the
reciprocals of $a_{Q,a}^{\sharp},b_{Q,a},b'_{Q,a}$, is bounded by
\begin{equation}\label{gb:prefactors}
 C_a\{Q+\log K+\log L\}.
\end{equation}
Here the powers of $K,L$ are fixed numerical powers.
The relative errors in \eqref{gb:transport_error}, \eqref{gb:micro_relative}
and \eqref{gb:weak_relative}--\eqref{gb:late_relative} each contain a factor
$t^{ch}$ with a fixed $c>0$: the slowest are $t^{h/256}$ in
\eqref{gb:transport_error} and $t^{h/6}$ in \eqref{gb:strong_relative},
and the others contain fixed powers $t^p$ with $p\ge1$.
We have $\log K\le C_a+\log Q+\log\log m$. Moreover,
$L\ge L_0:=\log(m/a)$, and $\log L/L$ decreases for $L\ge L_0$
eventually. Dividing \eqref{gb:prefactors} by $h_QL$ therefore gives
$o(1)$ uniformly over all $0<\alpha\le a$. For every fixed $c>0$,
each displayed prefactor multiplied by $t^{ch}$, $h\ge h_Q$, tends to
zero. This also gives $h c_{wt}^2\ge1$ and
$\log K=o(hL)$, the conditions needed for the pair-curvature and transport estimates.
The errors involving positive powers in \eqref{gb:fold_error} are absorbed the same way;
its late factor $e^{-H_a\sqrt m/Q}$ dominates
$e^{C_aQ}m^{3/4}$. Thus, above one common dimension threshold, the near-copy error is
smaller than $a_{Q,a}^{\sharp} tG_{\ell_*,t}(\delta)$, and the right side
of \eqref{gb:transport_error} is smaller than $a_{Q,a}^{\sharp}/2$.
All calculations remain uniform as $\alpha$ tends to zero.

It remains to combine the comparisons. For each original independent block $A$, apply
Lemma~\ref{gb:partition}. Contract its components successively to
actual anchors, preserving all physical copies. The block's original outside
$V_A$ remains unchanged throughout. Previously contracted components
preserve every cross-component cap, while the uncontracted components
retain their original spreads. Lemma~\ref{gb:transport} bounds the share
before contraction strictly below the share after contraction whenever
the spread is positive; at zero spread the shares agree. The terminal
representatives have cap $1-h$ and total
multiplicity at most $Q$, so \eqref{gb:capped_bound} applies. Hence
\begin{equation}\label{gb:original_bound}
 b_A\le |A|t-\frac{a_{Q,a}^{\sharp}}2t
               \sum_{C\subset A}G_{\ell_*,t}(\delta_C)
                    -\frac12D_A^{\mathrm{term}},
\end{equation}
where the sum is over the chosen partition components and $D_A^{\mathrm{term}}$ is the sum $D$ in \eqref{gb:edge_deficit} for the terminal block of $A$; these one-edge deficits all retain the same unchanged $V_A$.

The left side of \eqref{gb:original_bound} always refers to the original
Gaussian law, although different blocks use different auxiliary
references. Summing by \eqref{gb:shares} proves
$F_\Sigma(\alpha)\le\alpha$, without an additional error depending on the number
of blocks. If a component has positive spread or a terminal edge has
nonzero correlation, the bound is strict. Otherwise every component is
an exact signed-copy group and all terminal representatives are
independent, including across original blocks. Conversely, the exact
independent-group shares prove equality for every such law at every
dimension and level. This proves the claimed eventual classification.

The proof includes singular block matrices. The partition uses Gram
vectors; transport uses scalar Gaussian regression and rotation planes;
and the strong-set estimates use signed sums. Density arguments are
used only for capped pairs or in the automatically positive-definite
neighborhoods of independence used for the weak-correlation and
full-distribution estimates. Endpoint limits, if needed, follow
by positive semidefinite square-root coupling: all discontinuities of the physical
rejection event lie on finitely many one-coordinate threshold surfaces,
which have probability zero. No absolute approximation error is divided
by a vanishing spread or deficit. This completes the proof of Theorem~\ref{main:growing_blocks}.

\subsection{Other fixed bases}
\label{gb:other_bases}

The following result permits every fixed base $B>B_*$ and shows
that condition \eqref{main:growing_simple}, with its strict inequality, suffices.
The threshold $B_*$ comes from the transport step: with gap ratio
$s=B-1$, the marked-crossing error decays like $t^{dh}$ for exponents $d$
below, but arbitrarily close to, $\kappa/2$ with $\kappa=1-1/(4s^2)$,
whereas the sharper near-copy loss needs a rate $\ell>2$ on spreads
$\delta\le h/s^2$; the requirement $\kappa/2>2/s^2$ is equivalent to
$s^2>17/4$, that is, to $B>B_*$.

\begin{corollary}[A family of fixed bases]\label{gb:general_bases}
Fix $B>B_*:=1+\sqrt{17}/2$, put $s=B-1$, and define
\begin{equation}\label{gb:general_parameters}
 \begin{gathered}
 \kappa=1-\frac{1}{4s^2},\qquad
 \ell_B=2+\min\left\{\tfrac12,\frac{s^2-17/4}{4}\right\},\qquad
 M_B=\frac{\kappa}{2}-\frac{\ell_B}{s^2}>0,\\
 h_0(B)=\min\left\{\frac{1}{16},\frac{s^2(\ell_B-2)}{32}\right\},
 \qquad h_Q=h_0(B)B^{-2(Q-1)}.
 \end{gathered}
\end{equation}
The validity and equality conclusions of
Theorem~\ref{main:growing_blocks} hold if
\[
 \frac{Q_m+\log\log m}{h_{Q_m}\log m}\longrightarrow0.
\]
The dimension threshold may also depend on the fixed $B$.
In particular, condition~\eqref{main:growing_simple} suffices.
\end{corollary}

\begin{proof}
The same proof works for another fixed base if the marked-crossing error
still decays faster than the near-copy loss. The choices in the statement
leave a positive margin between those exponents. Here
$M_B\ge(s^2-17/4)/(4s^2)>0$, because $\ell_B\le s^2/4+15/16$,
$\kappa=1-1/(4s^2)$ and $s^2>17/4$ when $B>B_*$. Use the forest
construction with gap ratio $s$ and stopping threshold
$\sqrt{2h_0(B)}$. It gives within-component spread at most $h/s^2$ and
$h\ge h_Q$. Put $D_B=h_0(B)/s^2$ and
$\varepsilon_B=(\ell_B-2)/16$. For $0\le D\le1/16$,
$A(D)-1\le2D$, so
\[
 2[A(D_B)+\varepsilon_B]
 \le2+4D_B+2\varepsilon_B
 \le2+\frac{\ell_B-2}{4}<\ell_B.
\]
Thus Lemma~\ref{gb:sharp_near} applies with these fixed constants.
Also $\ell_B\le5/2<12$ and $\ell_BD_B<1$, preserving the folding
comparison and the late exponent nineteen.

Choose the two fixed tolerances
$\eta_B=\xi_B=M_B/8$. The lower quantile bound
$c_x^2\ge(2-\xi_B)L$ and $\log K\le\eta_BhL$ hold uniformly under
the envelope. The same two-endpoint path bound and marked-moment
calculation give exponent $d_Bh$, where
\[
 d_B=\frac{(2-\xi_B)(\kappa-\eta_B)^2}{4\kappa},\qquad
 \frac{\kappa}{2}-d_B\le\eta_B+\frac{\xi_B\kappa}{4}
 \le \frac{5M_B}{32}.
\]
Consequently $d_B-\ell_B/s^2>M_B/2=:\mu_B>0$.
The normalized transport error has $t^{\mu_Bh}$ in place of
$t^{h/256}$. The near-loss reciprocal still has logarithm $O_{a,B}(Q)$;
all capped estimates retain their previous form. The absorption and
original-law summation therefore apply without change. All choices are
fixed before $m$, $\alpha$, or the covariance law varies.

If the limsup in \eqref{main:growing_simple} is strictly smaller than
$1/(2\log B_*)$, choose a single $B>B_*$ sufficiently close to $B_*$
that the same limsup is smaller than $1/(2\log B)$. The growth condition
for this fixed base then follows, proving the asserted consequence.
This argument supplies neither the full envelope at $B_*$ nor the
endpoint case where this limsup equals the stated bound.
\end{proof}

\clearpage
\section{Graph comparisons and arbitrary correlations at fixed degree}
\label{fd:section}

This appendix proves Theorem~\ref{main:fixed_degree}, which allows
arbitrary correlations subject to a fixed covariance degree bound. The main obstacle
is a group of nearly identical coordinates that remains dependent on
its neighbors. We select clusters whose internal correlations are
closer to one than their boundary correlations, and use their rescaled spreads
to choose between two comparisons.

For larger rescaled spreads, we make a cluster independent of its complement
while preserving both marginal laws; we call this operation isolation.
For smaller rescaled spreads, we replace the cluster by copies of one of its
actual members, retaining that member's dependence on the outside;
we call this operation contraction. Isolation has an error small
relative to the cluster's deficit, and contraction gives a deficit of
square-root order. These operations change the joint law of the full
vector, so we compare full rejection probabilities along the way. After these comparisons, we
apply Lemmas~\ref{mx:lc} and~\ref{mx:near} to shares under one final
law and sum them. We state the theorem first, then
develop its covariance comparison and geometric tools.
Lemma~\ref{fd:hierarchy} selects the clusters, Lemma~\ref{fd:isolation}
and \eqref{fd:isolation_relative} bound the isolation error,
Lemma~\ref{fd:dependent_micro} gives the square-root deficit from
contraction, and Lemma~\ref{fd:terminal} bounds the shares under the final
law; the proof of Theorem~\ref{fd:theorem} chains these comparisons.
The appendix ends
with an example having common neighbors and the consequences across
dimensions.

Throughout this appendix, put
\[
 t=\alpha/m,\qquad L=\log(1/t),\qquad
 c_u=\Phi^{-1}(1-u/2),\qquad I_\alpha=\log(1/\alpha)-1+\alpha,
\]
and fix $0<a_1<1$. Here $I_\alpha$ is the quantity written $h_\alpha$ in Appendices~\ref{sec:lm:graph} and~\ref{mx:section}. A covariance
degree bound refers to physical coordinates unless representatives
are expressly specified. Constants
may depend on fixed degree and size bounds, fixed correlation caps and $a_1$, but not on
the actual covariance, $\alpha\in(0,a_1]$, or a positive spread.
Write $H$ for the full rejection indicator and $H_K$ for rejection at
some rank at most $K$. All comparisons retain the $m$ physical slots
and thresholds $jt$, including slots assigned p-value one.

For a physical set $A$, let $N_A(u)$ be the number of its p-values at most $u$, and let $R$
be the largest full crossing rank, including zero. If $R>0$, then
$N(tR)=R$: if the count were larger, that count would itself be a
larger crossing. Thus, including ties and exact duplicates,
\begin{equation}\label{fd:shares}
 b_A=\Expect\left[\frac{N_A(tR)}R\mathbf1\{R>0\}\right],
 \qquad F_\Sigma(\alpha)=\sum_A b_A
\end{equation}
for every physical partition, with every share evaluated under the
same full law and the fraction defined as zero when $R=0$. A set $C$ is
aligned when its coordinates are signed so that all their correlations
are positive; its spread is then $d_C=\max_{i,j\in C}(1-|\Sigma_{ij}|)$,
and its rescaled spread is $d_CL$. Define
\[
 G_t(d)=\min\{\sqrt{dL},1\}e^{-12dL}\quad(d>0),
 \qquad G_t(0)=0.
\]

\begin{theorem}\label{fd:theorem}
Theorem~\ref{main:fixed_degree} holds. More precisely, the disjoint
aligned clusters $\mathcal C$ constructed in Lemma~\ref{fd:hierarchy}
have spreads $d_C$ and sizes at most $D+1$. There are constants
$c,\lambda_*>0$ and $M$, depending only on $D,a_1$, such that
\begin{equation}\label{fd:margin}
 F_\Sigma(\alpha)\le\alpha-\tfrac12D_G^{V,\mathrm{term}}-ct\left\{
   \sum_{0<d_CL\le\lambda_*}\sqrt{d_CL}
   +\sum_{d_CL>\lambda_*}G_t(d_C)\right\}
 \quad(m\ge M,\ 0<\alpha\le a_1).
\end{equation}
Here $D_G^{V,\mathrm{term}}$ is the sum of one-edge \emph{share}
deficits of the actual terminal graph, retaining its physical
multiplicities and its actual terminal independent outside law; see
Lemma~\ref{fd:terminal}. It is not the original full-event pair deficit.
In this range equality holds exactly for independent groups of signed
exact copies.
\end{theorem}

Lemma~\ref{fd:hierarchy} selects the clusters in \eqref{fd:margin} using a
spanning forest built from projective distances between Gaussian Gram vectors. We call a positive spread
microscopic when the rescaled spread is at most the fixed cutoff in the
first sum. The proof of Theorem~\ref{fd:theorem} in Section~\ref{fd:assembly}
specifies the terminal graph and its independent outside law.

\subsection{Covariance comparison and cluster geometry}

We need two shared tools: a derivative formula for comparing full
rejection probabilities, and a geometric bound on high-correlation
components. The comparison will retain the total number of physical
tests and every crossing rank. Its error must be small relative to a
cluster deficit that may tend to zero arbitrarily quickly.

The covariance derivative used below is a version of Price's
identity; see \citet{Voigtlaender2021Price}. It expresses the change in
an event probability through jumps across coordinate thresholds. We
give the form for a finite grid, including its proof, so that we can
apply it to both early and late crossing ranks.

\begin{lemma}\label{ug:grid}
Let \(f\) be a bounded function constant on the cells determined by
finitely many coordinate hyperplanes. Let \(\Sigma(u)\) be a
differentiable positive definite covariance path with fixed unit
diagonal. If \(T_i\) is the set of breakpoints in coordinate \(i\),
then
\begin{equation}\label{ug:grid_identity}
 \frac{d}{du}\Expect_u f(Z)
 =\sum_{i<j}\dot\Sigma_{ij}(u)
   \sum_{a\in T_i,\ b\in T_j}
   \phi_{\Sigma_{ij}(u)}(a,b)
   \Expect_u[J_{ij}^{a,b}f\mid Z_i=a,Z_j=b].
\end{equation}
Here \(\phi_\rho\) is the standard bivariate normal density and
\(J_{ij}^{a,b}f=f_{++}-f_{+-}-f_{-+}+f_{--}\) is the mixed jump,
with the other coordinates fixed. This is the change in the jump across
one threshold when the other threshold is crossed. The signs indicate
which side of the two thresholds is used.
\end{lemma}

\begin{proof}
The Fourier transform of a Gaussian density is
\(\exp(-\xi^T\Sigma\xi/2)\). Differentiating an off-diagonal
entry of the symmetric matrix therefore differentiates the
density by \(\partial_i\partial_j\), with coefficient one.
On any compact part of the path where the covariance is positive definite,
the density and its derivatives have an integrable Gaussian bound.
This bound justifies differentiation under the integral. Integrate \(\partial_i\partial_j\)
over each cell; in the coordinates \((z_i,z_j)\) this gives a signed sum
of the density at the four corners of the cell. Grouping these corner
terms by breakpoint pairs \((a,b)\) leaves exactly the four adjacent
values in the mixed jump. Integration over the
remaining coordinates gives the conditional expectation in
\eqref{ug:grid_identity}. Unbounded cells are obtained by taking
limits, with the Gaussian boundary terms tending to zero.
\end{proof}

The next estimate treats isolation when the boundary correlations are
small. Later, in Lemma~\ref{fd:isolation}, we use the same argument with an
arbitrary fixed cap below one, as needed for the general theorem; the
proof of Theorem~\ref{fd:theorem} applies that lemma rather than Lemma~\ref{ug:comparison}.

\begin{lemma}\label{ug:comparison}
Suppose that the coordinates in \(S\) have covariance degree at
most \(D\), and the remaining coordinates are independent
Gaussian blocks of size at most \(Q\), independent of \(S\).
Let \(C\subset S\), and suppose
\(\lvert\Sigma_{ij}\rvert\le1/4\) for every
\(i\in C,j\notin C\). Write
\[
 B_C=\sum_{i\in C,j\notin C}|\Sigma_{ij}|,\qquad
 \chi=\max\{D+1,Q,2\},\qquad s_0=2(D+1).
\]
Let \(\Sigma^0=\Sigma_C\oplus\Sigma_{C^c}\), and set
\begin{equation}\label{ug:cutoff}
 K=\max\left\{s_0,2,
   \left\lceil\frac{(s_0+10)\chi L}{I_\alpha}\right\rceil\right\}.
\end{equation}
For every fixed \(a_1<1\), uniformly over \(0<\alpha\le a_1\)
and all sufficiently large \(m\),
\begin{equation}\label{ug:comparison_bound}
 |F_\Sigma(\alpha)-F_{\Sigma^0}(\alpha)|
 \le C B_C\{K^4L t^{8/5}+\chi t^8\}.
\end{equation}
The conclusion includes singular covariance matrices.
\end{lemma}

\begin{proof}
We bound early and late ranks separately along the same covariance
path. Use the positive semidefinite path
\[
 \Sigma(u)=u\Sigma+(1-u)\Sigma^0,\qquad 0\le u\le1.
\]
It changes only the cross-covariance block and preserves every
original zero covariance. First suppose the matrices are positive
definite. As at the start of this appendix, \(H\) is the full rejection
indicator and \(H_K\) the indicator of a crossing at some rank at most
\(K\). For \(H_K\), each coordinate has \(2K\) breakpoints
\(\{\pm c_{jt}:j\le K\}\). A mixed jump has absolute value at most
four. Since every changed correlation has absolute value at most
\(1/4\), its density at any such pair of breakpoints is at most
\[
 C\exp\{-4c_{Kt}^2/5\}.
\]
The Mills inequalities give
\(e^{-c_x^2/2}\le Cx\sqrt{\log(1/x)}\) as \(x\downarrow0\).
Thus Lemma~\ref{ug:grid}, integrated over \(u\), bounds the
early-rank contribution by
\[
 C B_C K^2 e^{-4c_{Kt}^2/5}
 \le C B_C K^{18/5}L^{4/5}t^{8/5}
 \le C B_C K^4L t^{8/5}.
\]

For the late ranks, consider the remaining indicator
\(g=H-H_K\) on the full grid. The key is to remove both endpoints of a changed edge and all their
original covariance neighbors. The remaining
vector is independent of the two endpoints and must supply most of
the small p-values needed for a late crossing. Specifically, for a changed edge \(\{i,j\}\), let \(W\) contain all coordinates
except \(i\), \(j\), and their original covariance neighbors within
\(S\). Thus the omitted coordinates form the union of the two closed
neighborhoods. All external blocks
are included in \(W\). At most \(s_0\)
coordinates are omitted, and \(Z_W\) is jointly independent of
\((Z_i,Z_j)\) throughout the path. Put
\[
 E_{ij}=\{\text{for some }\ell>K,\ 
              N_W(\ell t)\ge\ell-s_0\}.
\]
If \(E_{ij}\) fails, a late crossing is impossible for every
assignment of the omitted coordinates. Consequently
\[
 |J_{ij}^{a,b}g|\le4\mathbf1_{E_{ij}},\qquad
 \Expect_u[|J_{ij}^{a,b}g|\mid Z_i=a,Z_j=b]
 \le4\Prob_u(E_{ij}).
\]
The conditioning fixes only the two endpoints. Thus joint independence
of \(W\) and those endpoints gives the displayed bound, even though
coordinates inside the remaining vector may depend on one another.

Color the graph on \(S\) with \(D+1\) colors, and reuse at most
\(Q\) colors within each independent external block. Every color
class consists of jointly independent coordinates. H\"older's
inequality and the Bernoulli moment generating function yield
\[
 \Expect e^{vN_W(x)}
 \le \exp\{mx(e^{v\chi}-1)/\chi\}.
\]
Taking \(v=\log(1/\alpha)/\chi\) gives
\[
 \Prob\{N_W(\ell t)\ge\ell-s_0\}
 \le \alpha^{-s_0/\chi}e^{-I_\alpha\ell/\chi}.
\]
Summing the geometric series, and then using
\eqref{ug:cutoff}, shows that
\begin{equation}\label{ug:late_probability}
 \Prob(E_{ij})
 \le C_{a_1}\chi\alpha^{-s_0/\chi}
               e^{-I_\alpha K/\chi}
 \le C_{a_1}\chi t^{10}.
\end{equation}
Indeed, the exponential is at most \(t^{s_0+10}\), and
\(\alpha\ge t\). The cutoff makes this probability small enough that the bound can
be summed over every pair of thresholds on the full grid. There are at most \(4m^2\)
pairs of breakpoints for \(g\), and their bivariate densities
are uniformly bounded. The late-rank contribution is therefore at most
\(CB_C\chi m^2t^{10}\le CB_C\chi t^8\).

We also need \(K\le C_{a_1,D,Q}(1+\log m)<m\) uniformly
in \(\alpha\le a_1\) for large \(m\). To see this, put
\(u=\log(1/\alpha)\). The ratio \(u/(u-1+e^{-u})\) is bounded
for \(u\ge\log(1/a_1)>0\), and
\[
 \frac{L}{I_\alpha}
 =\frac{\log m+u}{u-1+e^{-u}}
 \le C_{a_1}(1+\log m).
\]
This remains true when \(L=O(\log m)\) fails because the level tends
to zero very quickly.

For a singular path, replace it by
\((1-\varepsilon)\Sigma(u)+\varepsilon I\). Its changed entries
are \((1-\varepsilon)\Sigma_{ij}\), so the preceding bound uses
the original \(B_C\). The graph zeros and coloring persist.
Let \(\varepsilon\downarrow0\). Each event boundary is contained
in finitely many coordinate hyperplanes, all null under the
continuous standard normal margins. Weak convergence therefore
passes the bound to the singular law.
\end{proof}

The geometric tool bounds the sizes of high-correlation components
and shows how to align their coordinate signs. We use it in Lemma~\ref{fd:hierarchy} to select the clusters
for the fixed-degree theorem.

\begin{lemma}\label{ug:geometry}
Let a positive semidefinite correlation matrix have covariance
degree at most \(D\), and let \(r<1\) satisfy
\[
 r\ge1-\{256(D+1)^2\}^{-1}.
\]
Every component of the graph with edges
\(|\Sigma_{ij}|>r\) has at most \(D+1\) vertices. Its coordinates
can be signed so that all pairwise correlations are positive
and its spread is at most \(D^2(1-r)<1/256\).

\end{lemma}

\begin{proof}
A high-correlation component must be small because its members are
all correlated with one another. Represent the covariance by unit vectors in a real Hilbert space. Along a spanning tree of a high-correlation component,
choose signs successively so that each tree edge has positive
inner product. Its endpoint distance is less than
\(\sqrt{2(1-r)}\). If the component contained \(D+2\) vertices,
take a connected subset of exactly \(D+2\) vertices and a
spanning tree. Every pair of its signed vectors is at distance
less than \((D+1)\sqrt{2(1-r)}\). Every pair of signed vectors therefore has inner product greater than
\(1-(D+1)^2(1-r)\ge255/256>0\). These \(D+2\) vertices would
form a clique in the original graph, contradicting the degree bound.
Within a component of size at most \(D+1\), the same argument
uses path length at most \(D\) and proves the spread bound.
\end{proof}

\subsection{Counts, thresholds and the final share bounds}

We first collect three tools: a bound on counts at late ranks, estimates
for the spacing of Gaussian thresholds, and share bounds for the final
law. Each is uniform as the level tends to zero.
The count bounds also hold for superuniform p-values with a
\emph{dependency graph}: disjoint vertex sets with no edge between them
carry independent random subvectors. For Gaussian coordinates their
covariance graph has this property, including at singular covariance.
For a vector of length at most $m$ and dependency degree at most $d$,
color its graph with $\chi=\max(d+1,2)$ colors. Coordinates within one
color are jointly independent, and their marginal low probabilities are
at most the testing threshold. H\"older's inequality and the Bernoulli moment generating
function give
\[
 \Expect e^{vN(jt)}
 \le\exp\{\alpha j(e^{v\chi}-1)/\chi\}.
\]
Taking $v=\log(1/\alpha)/\chi$, for every fixed insertion allowance
$s\ge0$ we obtain the following bound. This allowance counts the coordinates
that may be added to the count:
\begin{equation}\label{fd:count}
 \Prob\{N(jt)\ge j-s\}
 \le\min\{1,\alpha^{-s/\chi}e^{-I_\alpha j/\chi}\}.
\end{equation}
For $j\ge s/(1-a_1)$, the logarithmic derivative of the second bound
with respect to $\alpha$ is
$\{j(1-\alpha)-s\}/(\chi\alpha)\ge0$. Evaluating at $a_1$ and
enlarging a constant for the finitely many smaller ranks yields
\begin{equation}\label{fd:count_uniform}
 \Prob\{N(jt)\ge j-s\}\le C e^{-\kappa j},
 \qquad \kappa>0,
\end{equation}
uniformly over $\alpha\le a_1$. In geometric tail sums below,
$(1-e^{-I_\alpha/\chi})^{-1}$ is bounded by a constant because
$I_\alpha\ge I_{a_1}>0$.

The bound on gained rejections in the proof of
Lemma~\ref{fd:independent_micro} compares within-cluster
differences with gaps between successive Gaussian thresholds. We therefore need both lower bounds on the gaps and bounds on the
threshold sizes.
Set
$B=c_t+c_t^{-1}$ and $z_j=c_{jt}$. Mills' inequalities
give
\[
 c_u^2=2\log(1/u)+O\{\log\log(e/u)\},\qquad
 x\le h(x):=\phi(x)/\barPhi(x)\le x+x^{-1}\quad(x>0).
\]
The ratio $I_\alpha/\log(1/\alpha)$ is bounded below by a positive
constant on $(0,a_1]$. Therefore a cutoff
$K=\max\{K_0,\lceil C_0L/I_\alpha\rceil\}$, with fixed constants,
satisfies, uniformly for all sufficiently large $m$,
\begin{equation}\label{fd:quantiles}
 \begin{gathered}
 K\le C(1+\log m)<m,\qquad \log K=o(L),\qquad
 L\le c_{Kt}^2\le3L,\qquad L\le B^2\le3L,\\
 c_t/c_{Kt}\longrightarrow1,\qquad c_{wt}/B\ge1/2
 \quad\text{for each fixed bounded }w.
 \end{gathered}
\end{equation}
Indeed, $L/I_\alpha=(\log m+\log(1/\alpha))/I_\alpha
\le C(1+\log m)$, and Mills' logarithmic correction is $o(L)$.
The bounds hold even when $L=O(\log m)$ fails because the level tends
to zero quickly.
For $j<m$, monotonicity of $h$ gives
\begin{equation}\label{fd:gaps}
 g_j:=z_j-z_{j+1}
 =\int_{jt}^{(j+1)t}\frac{du}{u h(c_u)}
 \ge\frac1{(j+1)B}.
\end{equation}
Set $g_m=z_m$. The same lower bound holds eventually because
$z_m\ge c_{a_1}>0$.

Finally, we record the share bounds under the law that will result from
isolation and contraction. Specifying this law matters: all shares in
the final sum must refer to the same distribution.

\begin{lemma}\label{fd:terminal}
Fix a correlation cap $0<r<1$, representative multiplicity and degree
bounds $W,d$, a block-size bound $Q$, and $a_1<1$. Suppose the full
law consists of a capped Gaussian component $G$ with these
representative bounds, independent of mutually independent Gaussian
blocks $C$ of physical size at most $Q$. After sign changes, every
such block has spread $d_C\le1/16$. Under this one law, there are
$c_N>0$ and a common minimum dimension such that
\begin{equation}\label{fd:terminal_bounds}
 b_G\le |G|t-(1-t^{g_r/2})D_G^V,\qquad
 b_C\le |C|t-c_NtG_t(d_C),
 \quad g_r=\frac{1-r}{(1+2r)(1+r)}.
\end{equation}
Here $D_G^V=\sum_e d_e^V$ under the actual independent outside law.
If a nonzero representative edge remains, then $D_G^V>0$, and the
first bound gives $b_G<|G|t$. A block with zero spread has exactly its nominal share.
\end{lemma}

\begin{proof}
Apply Lemma~\ref{mx:lc} with multiplicity bound $W$, correlation cap
$r$, constant degree envelope $\Delta_m=d$ and outside block-size
bound $Q$. A constant envelope satisfies its required condition
\eqref{eq:lm:g:boundary.1} with $\beta=0$. The one-edge deficits in that application are
\[
 d_e^V=|G|t-b_G^{e,V},
\]
where $V$ is this actual independent outside and only edge $e$ is
retained inside $G$. Apply \eqref{mx:lc_relative} with $\beta=0$ and $\varepsilon=g_r/2$. With $\rho_0=r$ and $\beta=0$, Lemma~\ref{mx:lc} has
$\vartheta_0=\min\{g_r,1/r-2/(1+r)\}=g_r$, because
$1/r-2/(1+r)=(1-r)/\{r(1+r)\}$, which exceeds $g_r$ since $r<1+2r$.
Hence $\varepsilon<\vartheta_0$ and $\vartheta_0-\varepsilon=g_r/2$.
Its relative remainder is at most $t^{g_r/2}$ times the sum of these deficits,
and each nonzero edge has positive deficit by
\eqref{mx:edge_lower}--\eqref{mx:edge_estimates}. This proves the first
bound and, when a nonzero edge remains, the strict inequality $b_G<|G|t$. The quantity $d_e^V$ is the share
deficit with the actual outside law retained; it differs from a
deficit in the full rejection probability or in a separate bivariate test.

For one block $C$, its outside consists of the same capped component
$G$ and all the other independent blocks. This is exactly the outside
specified at the start of Section~\ref{mx:selected_section} for
Lemma~\ref{mx:near}, with bounded multiplicities,
constant representative degree and bounded individual block size.
That lemma gives the second bound, including its uniform control
\eqref{mx:folding_relative} as $d_C$ tends to zero. If $d_C=0$, forcing
and restoring the whole independent group gives share $|C|t$, as in
\eqref{mx:shares}. The number of blocks does not enter either lemma's
constants. Empty parts have zero share.
\end{proof}

\subsection{A forest partition of nearly identical coordinates}

We need clusters whose internal variation is much smaller than their
separation from the rest of the vector. The next construction gives this
gap and ensures that retaining one actual member of each cluster leaves
a graph with a fixed correlation cap. Both properties are needed for the
later comparisons.

\begin{lemma}\label{fd:hierarchy}
Let $\Sigma$ have physical covariance degree at most $D\ge1$, and put
\[
 q=D+1,\qquad\varepsilon_0=(256q^2)^{-1},\qquad
 h_{\min}=\varepsilon_0\,49^{-D}.
\]
There are disjoint aligned clusters $C$, each of size between two and
$q$, with numbers $h_C\in[h_{\min},\varepsilon_0]$ such that
\begin{equation}\label{fd:hierarchy_bounds}
 d_C<h_C/36,\qquad
 |\Sigma_{iy}|\le1-h_C\quad(i\in C,y\notin C).
\end{equation}
Retaining one actual anchor per cluster and every remaining singleton
leaves representatives of degree at most $D$, multiplicity at most $q$
and correlation cap $r_T=1-h_{\min}<1$.
\end{lemma}
\begin{proof}
Represent $\Sigma$ by unit Gram vectors $v_i$, also when it is singular.
The projective chord distance is
$\ell(i,j)=\min(\|v_i-v_j\|,\|v_i+v_j\|)
=\sqrt{2(1-|\Sigma_{ij}|)}$.
It satisfies the triangle inequality: choose minimizing signs for two
successive differences and multiply them. A tree can consequently be
sign-aligned with signed distances along any path bounded by its length.

Form base components with $\ell<\tau_0:=\sqrt{2\varepsilon_0}$.
Lemma~\ref{ug:geometry} makes each a covariance clique of size $n\le q$.
Inside a nontrivial base component, process the edges of the
\emph{complete} graph of projective distances from shortest to longest,
accepting an edge unless it creates a cycle. This is Kruskal's algorithm;
break ties by labels and stop before $\tau_0$. The accepted edges form a spanning tree. The selection below is the one in the proof of
Lemma~\ref{gb:partition}, with factor six in place of $21/10$ and the
threshold $\tau_0$ above.
Write their lengths as $e_1\le\cdots\le e_p<\tau_0$, with $p=n-1\le D$,
and put $e_{p+1}=\tau_0$, $S_0=0$, $S_j=\sum_{i\le j}e_i$.
Choose the largest $k\in\{0,\ldots,p\}$ with $e_{k+1}>6S_k$.
Such an index exists: $k$ equal to the number of zero accepted lengths
qualifies, since then $S_k=0<e_{k+1}$ (with $e_{p+1}=\tau_0$ when all
accepted lengths are zero). Put $\tau=e_{k+1}$ and $h_C=\tau^2/2$ for every
nontrivial component of the first $k$ edges.

Any edge between two distinct retained components, singletons included,
has length at least
$\tau$. Indeed an earlier shorter edge either would have been accepted,
joining them, or was rejected only because an earlier path already
joined its endpoints. Ties at $\tau$ give the required weak boundary
inequality. Along a retained tree the total signed path length is at
most $S_k<\tau/6$, so its spread is less than $h_C/36$.
For $j>k$, maximality gives $e_{j+1}\le6S_j$. If $k<p$, then
$S_{k+1}<7\tau/6$, $S_p\le7^{p-k-1}S_{k+1}$, and
$\tau_0\le6S_p<7^{p-k}\tau$. If $k=p$, $\tau=\tau_0$.
Thus $\tau^2/2\ge\varepsilon_0 7^{-2D}=h_{\min}$ in all cases, including
$k=0$, when the base component forms no cluster. In particular,
$h_C\ge h_{\min}$.
Different base components have correlations at most
$1-\varepsilon_0\le1-h_C$ in magnitude. This proves the boundary
bound against the whole complement, not just the base component.
Retained anchors are actual original coordinates, so their graph is an
induced original subgraph. Two distinct representatives lie in different
retained components. If both are in one base component, their projective
distance is at least $\tau$ by the edge-length bound above, so their
correlation is at most $1-\tau^2/2\le1-h_{\min}$ in magnitude. If they are in different base
components, it is at most $1-\varepsilon_0\le1-h_{\min}$. Hence the common
cap is $1-h_{\min}$.
\end{proof}

As the covariance matrix varies, the values of $h_C$ need not belong
to one fixed finite list. Their common positive
lower bound supplies the uniform constants below. More generally, the
same proof works with factor $c>\sqrt{24}$ in place of six, giving
$h_{\min}=\varepsilon_0(c+1)^{-2D}$ and $d_C<h_C/c^2$.
The resulting conservative isolation exponent, the lower bound in
\eqref{fd:isolation_exponent} with $c$ in place of six, is
$(1/2-12/c^2)\varepsilon_0(c+1)^{-2D}$.
For a fixed $D\ge1$, differentiating its logarithm shows that it is
maximized at the unique root above $\sqrt{24}$ of
$Dc^3-24(D+1)c-24=0$. We use six to keep the constants simple; optimizing
this particular exponent does not evaluate the full dimension threshold.

The preceding degree bound counts representatives. The successive
contractions by Lemma~\ref{fd:dependent_micro} also need neighbor and
degree bounds that count every repeated test, because the crossing ranks
use all of them. During contraction, every group of copied coordinates remains a
subset of an original clique. If the group has $w$ members, its anchor has at
most $D-w+1$ original outside neighbors. A neighboring group has at
most $D$ members, since its anchor already has that outside neighbor.
Counting every current copy, the physical degree is therefore at most $w-1+D(D-w+1)\le D^2$, and we may use
\begin{equation}\label{fd:degree_accounting}
 d_1=D^2.
\end{equation}
An uncontracted cluster of size $k$ has at most $k(D-k+1)$ original
outside neighbors, each with at most $D$ current copies. Thus a common
neighbor bound is
\begin{equation}\label{fd:neighbor_accounting}
 n_0=D\lfloor(D+1)^2/4\rfloor.
\end{equation}
For $D=1$ an actual selected pair has no outside neighbor, and the stated
bound is harmless. Isolation only deletes edges. The internal covariance
of a cluster and its boundary cap persist until its own contraction.

\subsection{Isolation when the rescaled spread exceeds the cutoff}

We use isolation for clusters whose rescaled spread exceeds the
microscopic cutoff. We extend the full-grid comparison of
Lemma~\ref{ug:comparison} to any fixed boundary cap below one. On a straight-line path between positive semidefinite
covariance matrices, suppose every changed entry has absolute value at
most $\rho<1$, and let $B_\partial$ be the sum of absolute changes
over unordered pairs. For a function $f\in[0,1]$ on the grid
$\{\pm c_{jt}:j\le K\}$, Lemma~\ref{ug:grid} gives, for positive
definite endpoints and then by the regularization below,
\begin{equation}\label{fd:finite_grid}
 |\Delta\Expect f|
 \le C_\rho B_\partial K^2e^{-c_{Kt}^2/(1+\rho)}
 \le C_\rho B_\partial K^4L t^{2/(1+\rho)}.
\end{equation}
Indeed, every mixed jump is at most four in absolute value and there
are at most $4K^2$ boundary pairs per changed edge. Their bivariate
density is bounded by
$C_\rho\exp\{-(x^2+y^2)/(2(1+\rho))\}$. Mills' inequalities give
the second bound. If $K$ is fixed, the sharper resulting estimate is
\begin{equation}\label{fd:fixed_grid}
 |\Delta\Expect f|\le C B_\partial t^pL^{p/2},
 \qquad p=2/(1+\rho)>1.
\end{equation}
For singular endpoints, apply the same regularization as in
Lemma~\ref{ug:comparison}: use
$(1-\varepsilon)\Sigma(u)+\varepsilon I$ and let $\varepsilon$ decrease to
zero. The cap and zero entries persist, the changed coefficients are
multiplied by $1-\varepsilon$, and the null marginal threshold boundaries
allow passage to the limit with the original $B_\partial$.
For the isolation path $u\Sigma+(1-u)\Sigma^0$ of a cluster $C$ in the next
lemma, $B_\partial=\sum_{i\in C,\,j\notin C}|\Sigma_{ij}|$, which is the
quantity $B_C$ of Lemma~\ref{ug:comparison}.

\begin{lemma}\label{fd:isolation}
Suppose the physical degree is bounded by $d$, and all correlations
across a selected cluster boundary are at most $\rho<1$ in absolute
value. Let $\Sigma^0=\Sigma_C\oplus\Sigma_{C^c}$, set
$s_0=2(d+1)$, $\chi=\max(d+1,2)$, and take
\[
 K=\max\left\{s_0,2,
 \left\lceil (s_0+10)\chi L/I_\alpha\right\rceil\right\}.
\]
Then, uniformly eventually over $0<\alpha\le a_1$,
\begin{equation}\label{fd:isolation_bound}
 |F_\Sigma(\alpha)-F_{\Sigma^0}(\alpha)|
 \le C_\rho B_\partial\{K^4Lt^{2/(1+\rho)}+t^8\},
\end{equation}
where $C_\rho$ depends only on $\rho$, $d$ and $a_1$.
\end{lemma}

\begin{proof}
Follow the early- and late-rank decomposition in the proof of
Lemma~\ref{ug:comparison}, now allowing the boundary cap $\rho<1$.
Interpolate by $u\Sigma+(1-u)\Sigma^0$ and apply
\eqref{fd:finite_grid} to $H_K$. For $H-H_K$, remove the original
closed neighborhoods of the endpoints $\{i,j\}$ of a changed edge.
The remaining vector $W$ is jointly independent of $(Z_i,Z_j)$ along
the path, and at most $s_0$ coordinates were removed. The same
mixed-jump argument on the complete $m$-rank grid requires
\[
 E_{ij}=\{N_W(\ell t)\ge\ell-s_0
                   \text{ for some }\ell>K\}.
\]
If this event fails, no assignment of the removed coordinates can
produce a late crossing. The indicator of a late crossing without an early crossing is
therefore zero for every such assignment, at all four corners of each
mixed jump and whatever happens at early ranks, so
$J_{ij}^{a,b}(H-H_K)=0$ off $E_{ij}$. Hence the conditional absolute mixed jump is
bounded by $4\Prob(E_{ij})$. By \eqref{fd:count},
\[
 \Prob(E_{ij})\le C\alpha^{-s_0/\chi}
                          e^{-I_\alpha K/\chi}\le Ct^{10}.
\]
The last step uses $e^{-I_\alpha K/\chi}\le t^{s_0+10}$ and
$\alpha\ge t$. On the full grid there are at most $4m^2$ boundary
pairs, whose densities are uniformly bounded by the fixed cap.
The late contribution is at most
$CB_\partial m^2t^{10}\le CB_\partial t^8$. The singular case follows
by the zero-preserving regularization described above.
\end{proof}

We now compare the isolation error with the lower bound $c_NtG_t(d_C)$
on the cluster's share deficit that Lemma~\ref{fd:terminal} gives for the
unchanged isolated cluster under the terminal law.
For a selected cluster $C$, isolation occurs before any contraction.
Thus $B_\partial\le Dq$, while its local boundary cap is $1-h_C$.
First retain this local cap in the exponent. We will bound the
constants uniformly afterward:
\begin{equation}\label{fd:isolation_exponent}
 \frac2{2-h_C}-1-12d_C
 >\left(\frac12-\frac{12}{36}\right)h_C
 \ge h_{\min}/6>0.
\end{equation}
The density constants in \eqref{fd:isolation_bound} can be bounded by
the worst cap $r_T$, but using the worst cap in the exponent before
subtracting $12d_C$ would lose this estimate.
Fix $\lambda_*>0$. If $d_CL>\lambda_*$, divide the isolation error by
$tG_t(d_C)$ and use
$\min(\sqrt{d_CL},1)\ge\min(\sqrt{\lambda_*},1)$. Uniformly,
\begin{equation}\label{fd:isolation_relative}
 \frac{|F_\Sigma-F_{\Sigma^0}|}{tG_t(d_C)}
 \le C\{K^4Lt^{h_{\min}/6}+t^{7-12\varepsilon_0/36}\}=o(1).
\end{equation}
Here $K=O(1+\log m)$, while $L\ge\log(m/a_1)$, so fixed positive
powers of $e^{-L}$ dominate every displayed polynomial uniformly on that
half-line. The bound is uniform over all realized local scales
$h_C\in[h_{\min},\varepsilon_0]$. The positive rescaled-spread cutoff
is needed here; the next comparison handles arbitrarily small spreads.

\subsection{A background event without a correlation cap}

For the microscopic comparison, we need a background event on which
adding a bounded cluster cannot create a late crossing. We prove that
this event has positive probability directly, without assuming a Simes
inequality for the background.

\begin{lemma}\label{fd:void}
Fix a degree bound $d$, an insertion allowance $q_0\ge2$ and
$0<a_1<1$. Let $V$ consist of at most $m$ superuniform p-values with a
dependency graph of maximum degree at most $d$. There are an integer
$J_0\ge q_0$, a constant $\gamma>0$ and a common minimum dimension such
that
\begin{equation}\label{fd:void_event}
 \Prob\{N_V(jt)=0\ (j\le J_0),\quad
       N_V(jt)+q_0<j\ (J_0<j\le m)\}\ge\gamma
\end{equation}
for every $0<\alpha\le a_1$. In particular this holds for every
standard-marginal Gaussian background of covariance degree at most $d$.
\end{lemma}
\begin{proof}
We first bound the probability that every background p-value exceeds
the early cutoff. Conditional on that event, we bound the probability
of a later crossing, even after the allowed insertions, by one half.
We choose the following cutoff to leave that probability margin.
Put $Q_0=\max(d+1,q_0,2)$ and $\beta=(1+a_1)/2$. Choose $J_0\ge Q_0$
so that
\[
 \frac{\beta^{-1}e^{-I_\beta(J_0+1)/Q_0}}
      {1-e^{-I_\beta/Q_0}}\le\frac12.
\]
Set $A_i=\{V_i\le J_0t\}$, $p=J_0t$, $x=2p$, and
$A=\bigcap_iA_i^c$. For large $m$, uniformly in $\alpha\le a_1$,
\[
 x\le1/2,\qquad p\le x(1-x)^d,\qquad
 a_1(1-x)^{-(d+1)}\le\beta.
\]
Conditioning many background p-values to avoid the early thresholds
can increase the chance that another is small. We first show that this
increase remains controlled by the degree bound. Specifically, we prove
$\Prob(A_i\mid\bigcap_{j\in T}A_j^c)\le x$ for every $i\notin T$,
with every conditioning event of positive probability. For $T$ empty it follows from
$\Prob(A_i)\le p\le x$. If the bound holds through conditioning sets of
size $k-1$, successive exclusion proves
$\Prob(\bigcap_{j\in T}A_j^c)\ge(1-x)^{|T|}>0$ for $|T|\le k$.
Now split a size-$k$ set into neighbors $T_1$ of $i$ and nonneighbors
$T_2$. If $T_1$ is empty, whole-subvector independence gives the result.
Otherwise write the conditional probability as a ratio of two
probabilities conditional on $\bigcap_{T_2}A_j^c$. The numerator is at
most $\Prob(A_i\mid\bigcap_{T_2}A_j^c)=\Prob(A_i)$. Expanding the denominator over
the vertices of $T_1$ one at a time gives factors that each condition on
at most $k-1$ vertices, so their product is at least $(1-x)^{|T_1|}$.
Thus the conditional probability is at most $p/(1-x)^d\le x$.
The induction proves positivity of the conditioning events together
with the bound on the conditional probabilities.
In particular,
\[
 \Prob(A)\ge(1-x)^{|V|}\ge e^{-2|V|x}\ge e^{-4a_1J_0}.
\]

Conditioning on $A$ need not preserve color-class independence. Instead,
fix a subset $S$ in one color and let $B_S=\bigcap_{i\in S}\{V_i\le jt\}$.
Its closed neighborhood has size at most $(d+1)|S|$. The subvector on
$S$ is independent of the whole complement of that neighborhood.
Conditioning on the absence of events in that neighborhood can increase
the probability by a factor of at most $(1-x)^{-(d+1)|S|}$, by the
bound on conditional probabilities just proved. Independence within the unconditioned
color and superuniformity therefore give
\[
 \Prob(B_S\mid A)\le[jt(1-x)^{-(d+1)}]^{|S|}.
\]
Let $c=(1-x)^{-(d+1)}$, and use $Q_0$ colors, allowing empty classes.
For a class $I_\ell$, expand the moment generating function over all its
subsets, with $y=e^{vQ_0}-1\ge0$:
\[
 \Expect(e^{vQ_0N_\ell(jt)}\mid A)
 =\sum_{S\subset I_\ell}y^{|S|}\Prob(B_S\mid A)
 \le(1+yjtc)^{|I_\ell|}.
\]
H\"older under the conditional measure yields
\[
 \Expect(e^{vN_V(jt)}\mid A)
 \le\exp\{\beta j(e^{vQ_0}-1)/Q_0\}.
\]
Taking $v=\log(1/\beta)/Q_0$ gives conditional insertion probability
at most $\beta^{-q_0/Q_0}e^{-I_\beta j/Q_0}$, whose sum after $J_0$ is
at most $1/2$. Intersecting its complement with $A$ proves the lemma with
$\gamma=\tfrac12e^{-4a_1J_0}$.

An explicit sufficient lower bound on $m$ is
\[
 m\ge\max\{J_0,\lceil2J_0a_1/x_*\rceil\},\qquad
 x_*=\min\{1/2,1-(a_1/\beta)^{1/(d+1)},1-2^{-1/d}\},
\]
where the last entry is omitted for $d=0$. These inequalities ensure
exactly the three conditions on $x$ used above.
The dependency-graph assumption is stronger than pairwise independence
or zero covariances for a non-Gaussian background; it must not be replaced
by either. Gaussian covariance graphs satisfy it, because jointly Gaussian
subvectors with zero cross-covariance are independent. Atoms, exact copies within allowed graph cliques and empty
backgrounds require no change.
\end{proof}

\subsection{Contraction with an independent background}

When the within-cluster variation is too small for the isolation bound,
we compare rejection events before and after replacing the cluster by
exact copies. We first do this with an independent outside vector. The
next subsection will retain dependence on the cluster's neighbors.
Here $F(C,V)$ is the full rejection probability of the vector formed by
the cluster $C$ and its outside $V$, $\operatorname{dup}_wA$ puts $w$
exact copies of $A$ in the slots of $C$, and $F(C,Y,V)$ below is defined
in the same way. The proofs keep the notation $z_k=c_{kt}$,
$B=c_t+c_t^{-1}$, $h$, $g_k$ and $\kappa$ of
\eqref{fd:count_uniform}--\eqref{fd:gaps}.

\begin{lemma}\label{fd:independent_micro}
Fix $q_0\ge2$, a degree bound $d$ and $a_1<1$. Let an aligned
Gaussian cluster $C$, with $2\le w=|C|\le q_0$ and spread
$\delta=\max_{i,j\in C}(1-\rho_{ij})$, be independent of its
entire outside $V$. Assume $V$ has superuniform marginals and a
dependency graph of degree at most $d$. There are $c_0,\lambda_0>0$ and a common minimum dimension such that,
for any actual aligned anchor $A$ in $C$,
\begin{equation}\label{fd:independent_micro_bound}
 F(C,V)\le F(\operatorname{dup}_w A,V)
             -c_0t\sqrt{\delta L}
 \quad(0<\delta L\le\lambda_0, 0<\alpha\le a_1).
\end{equation}
For $\delta=0$ the comparison is an identity. Every probability uses
the original total number $m$ of tests.
\end{lemma}

\begin{proof}
We compare the original cluster with a vector of exact copies. A gained
rejection occurs only for the original vector; a lost rejection occurs
only for the copy reference. A gain requires two cluster statistics
to differ enough to span adjacent thresholds. We construct a loss event in which the standardized cluster mean passes a
threshold but one cluster member does not. We first contract to the standardized physical mean. Put
\[
 \overline Z=w^{-1}\sum_{i\in C}Z_i,\qquad
 \sigma^2=\operatorname{Var}(\overline Z)\in[1-\delta,1],\qquad
 X=\overline Z/\sigma,\qquad \mu=\delta B^2.
\]
We choose a small fixed upper bound on $\mu$; by
\eqref{fd:quantiles} it suffices to choose $\lambda_0$ small.

First consider a rejection gained by the original vector relative to
the contracted vector. At an original crossing rank $k$, the outside
count is at least $k-w$. Unless, for some ordered pair $i\ne j$
in $C$ and sign $s$, one has
\begin{equation}\label{fd:gain_inclusion}
 N_V(kt)\ge k-w,\qquad sZ_i\ge z_k,\qquad
 s(Z_i-Z_j)\ge g_k,
\end{equation}
the contracted vector still rejects, possibly at a larger rank. To see
this, suppose that some, but not all, cluster members are selected and that no
ordered pair and sign satisfy all the conditions
in~\eqref{fd:gain_inclusion}. Then $k<m$, because a crossing at rank $m$
selects every slot. The selected members share one sign, since two of
opposite signs would satisfy \eqref{fd:gain_inclusion} with
$s(Z_i-Z_j)\ge2z_k\ge g_k$; every member has that sign and exceeds
$z_{k+1}$ in it.
Completing the cluster adds at least one physical slot. If $n_C<w$
cluster members were selected, the completed rank is
$k'=k+w-n_C\ge k+1$. Also $k'\le m$, since
$k\le N(kt)=n_C+N_V(kt)\le n_C+m-w$; the outside has only $m-w$ slots. All cluster members exceed $z_{k+1}\ge z_{k'}$ in
the selected sign, so their mean also passes that rank, and
standardization can only increase its signed magnitude. The previously
selected outside coordinates exceed $z_k\ge z_{k'}$ in absolute value. If the
whole cluster is selected instead and no such pair and sign exist,
all its members have the same sign and its mean passes $z_k$. At $k=m$, opposite signs
give a signed difference of at least $g_m=z_m$, while common signs preserve
rejection. If no cluster member was selected, the outside crossing
is unchanged. Hence every gained rejection satisfies
\eqref{fd:gain_inclusion} for some original crossing rank, ordered pair
and sign, including ties.

We bound the probability of such a large difference using a Gaussian tail
calculation. For a centered Gaussian pair $(U,T)$ with variances $1,2\eta$ and
covariance $\eta>0$, exponential tilting gives
\[
 \Expect[e^{vT}\mathbf1\{U\ge z\}]
 =e^{v^2\eta}\barPhi(z-v\eta).
\]
Since the normal hazard is increasing,
$\barPhi(z-v\eta)\le\barPhi(z)e^{v\eta\,h(z)}$ for $v\ge0$.
Markov's inequality, optimized over $v\ge0$, therefore yields
\begin{equation}\label{fd:tilted_pair}
 \Prob(U\ge z,T\ge g)
 \le\barPhi(z)
       \exp\left\{-\frac{[g-\eta\,h(z)]_+^2}{4\eta}\right\}.
\end{equation}
For $U=sZ_i$, $T=s(Z_i-Z_j)$, we have $\eta=1-\rho_{ij}\le\delta$.
At $\eta=0$ a positive difference is impossible, so no division by zero
is made. For $\eta>0$ the pair is nonsingular, since its
covariance determinant is $\eta(2-\eta)>0$, even when the cluster
covariance is singular. Using $h(z_k)\le B$,
\eqref{fd:gaps}, outside independence and
\eqref{fd:count_uniform}, then summing the finitely many cluster
pairs, gives
\begin{equation}\label{fd:micro_gain}
 \Prob(\text{gained rejection})
 \le Ct\sum_{k\ge1}k e^{-\kappa k}
 \exp\left\{-\frac{[1/(k+1)-\mu]_+^2}{4\mu}\right\}
 =o(t\sqrt\mu)
 \quad(\mu\downarrow0).
\end{equation}
To verify the last order bound, split the ranks according to which of
the two sources of rarity applies. At early ranks the required within-cluster difference is
unlikely; at later ranks the outside count is unlikely. Put $u=\mu^{-1/3}\ge3$. For $k\le\lfloor u\rfloor$,
$1/(k+1)-\mu\ge1/(2u)$, giving an exponential at most
$e^{-u/16}$. For later ranks discard the exponential factor for the difference and bound the
geometric sum by $C(u+1)e^{-\kappa u}$. Both bounds divided by
$\sqrt\mu=u^{-3/2}$ tend to zero. This estimate is uniform over
all smaller positive spreads and retains every physical rank.

We next construct a lost rejection of order $t\sqrt\mu$. The standardized cluster mean will lie just above the cluster
threshold, while a negative residual keeps one member below it. Let
\[
 a_i=w^{-1}\sum_j(1-\rho_{ij}),\qquad
 \overline a=w^{-1}\sum_i a_i,\qquad
 W_i=(Z_i-\overline Z)/\sqrt\delta.
\]
Then $\sum_iW_i=0$, $\operatorname{Var}(W_i)\le2$,
$\sigma^2=1-\overline a$, and
\[
 \operatorname{Cov}(W_i,X)
 =\frac{\overline a-a_i}{\sigma\sqrt\delta}.
\]
For a pair attaining $\delta$, scalar conditioning gives
\[
 \operatorname{Var}(W_i-W_j\mid X)
 =2-\frac{(a_i-a_j)^2}{\delta(1-\overline a)}\ge1
\]
when $\delta$ is small, because $0\le a_i,\overline a\le\delta$.
Set $a_0=1/(4q_0)$ and $T_0=\sqrt{60q_0}$, and condition on
\begin{equation}\label{fd:loss_strip}
 z_w<X<z_w+a_0\sqrt\delta.
\end{equation}
For small fixed $\mu$, the conditional means of $W_i$ have absolute
value at most one and their second moments are at most three.
The difference $W_i-W_j$ for a pair attaining the spread is conditionally normal
with variance at least one. It has probability at least $2\barPhi(1)>1/5$ of absolute value
at least one; a centered normal minimizes that two-sided probability
among its translates. A union and Markov bound gives
\[
 \Prob(\max_i|W_i|>T_0\mid X)\le3q_0/T_0^2=1/20.
\]
A zero-sum vector of length $w$ whose largest and smallest entries
differ by at least one has
minimum at most $-1/w$. Thus, uniformly on \eqref{fd:loss_strip},
with conditional probability at least $1/10$,
\begin{equation}\label{fd:loss_residual}
 \max_i|W_i|\le T_0,\qquad \min_iW_i\le-1/q_0.
\end{equation}

Decrease the fixed microscopic cutoff so that
$(a_0+T_0)\sqrt\delta<1/(wB)$. On these events every original
cluster coordinate is positive for large $m$, its minimum is below
$z_w$, and its maximum is below $z_{w-1}$ by \eqref{fd:gaps}.
Intersect with the independent background event of Lemma~\ref{fd:void},
using allowance $q_0$. For ranks $j\le w-1$ all cluster coordinates
miss the threshold and the outside count is zero; at rank $w$ the
cluster minimum misses it. For $w<j\le J_0$ there are at most $w$
selected coordinates, and for $j>J_0$ the insertion condition excludes
a crossing. Hence the original vector has no crossing at any rank, while the
reference vector with $w$ copies of $X$ rejects at rank $w$.

For small fixed $\mu$,
$\phi(z_w+a_0\sqrt\delta)\ge\phi(z_w)/2$.
Also $\phi(z_w)\ge wt z_w/2$ and $z_w/B\ge1/2$ by
\eqref{fd:quantiles}. Integration over \eqref{fd:loss_strip}, including
the two fixed positive event probabilities, yields
\begin{equation}\label{fd:micro_loss}
 \Prob(\text{lost rejection})\ge c t\sqrt\mu.
\end{equation}
The loss dominates the gain after we choose the microscopic cutoff so
that \eqref{fd:micro_gain} is at
most half of \eqref{fd:micro_loss}. Since $L\le B^2\le3L$, this
proves \eqref{fd:independent_micro_bound} for contraction to $X$.
Both $X$ and any actual aligned anchor $A$ are standard normal and
independent of the entire $V$, so their duplicated vectors with $V$
have the same law. This last replacement is an equality in distribution
only under that independence assumption. All conditioning above is
on the scalar $X$; singular residual covariance is permitted. For
$\delta=0$, the aligned coordinates are exact copies almost surely.
\end{proof}

\subsection{Preserving the dependence on common neighbors}

We now retain the dependence between a cluster and its outside
neighbors, including common neighbors.
The marked-crossing estimate in Lemma~\ref{fd:marked_crossing}, proved
in Appendix~\ref{mx:section} for the transport of block shares, bounds a
threshold crossing jointly
with the event that a neighbor's absolute Gaussian value exceeds its
threshold. This latter event is the mark, and it may depend on the
crossing derivative.

\begin{lemma}\label{fd:dependent_micro}
Fix $q,n_0,d$, a boundary cap $0\le r_0<1$ and $a_1<1$. Let a full
Gaussian vector be partitioned as $(C,Y,V)$, where $2\le|C|\le q$,
$|Y|\le n_0$, the covariance degree of $(Y,V)$ is at most $d$, and
$C$ is independent of the entire $V$. Suppose $C$ is aligned with
spread $\delta$ and every correlation between $C$ and $Y$ is bounded
in absolute value by $r_0$. There are $c_*,\lambda_*>0$ and one
common minimum dimension such that contraction to any actual aligned anchor $A$
satisfies
\begin{equation}\label{fd:dependent_micro_bound}
 F(C,Y,V)\le F(\operatorname{dup}_{|C|}A,Y,V)
                    -c_*t\sqrt{\delta L}
 \quad(0<\delta L\le\lambda_*,\ 0<\alpha\le a_1).
\end{equation}
The actual outside distribution and the anchor's outside covariances
are retained. For $\delta=0$ the comparison is an identity.
\end{lemma}

\begin{proof}
We compare the change caused by contraction with the corresponding
change under an independent background. The neighbors may themselves
depend on the remaining outside coordinates, so replacing them by
independent singletons would alter the background law. Instead, the
reference removes their dependence on the cluster while preserving the
joint law of all outside coordinates. To construct it,
orthogonally project $Y$ onto the Gaussian linear span of $V$, writing
$Y=TV+E_0$ for a deterministic matrix $T$. This projection remains
defined on the span when $V$ is singular. Since $C$ is independent of
$V$, the pair $(C,E_0)$ is independent of $V$. Introduce $E'_0$ with
the same law as $E_0$, independent of the original variables, and put
$Y'=TV+E'_0$. Then
\[
 (Y',V)\stackrel{d}{=}(Y,V),\qquad C\ \text{is independent of }(Y',V).
\]
This construction gives an actual pair and a reference pair of
testing problems. Write $H_1,H_0$ for the actual rejection probabilities
before and after anchor contraction, and $B_1,B_0$ for those using
$Y'$. Lemma~\ref{fd:independent_micro}, applied to the preserved
outside marginal, gives
\begin{equation}\label{fd:surrogate_increment}
 B_1-B_0\le-c_0t\sqrt{\delta L}
\end{equation}
throughout a sufficiently small fixed microscopic range. We use only
the difference of the two reference probabilities. Neither is assumed
to equal $\alpha$.

To compare the actual and reference contraction increments, we couple each coordinate to the retained
anchor by a Gaussian rotation. For each member of $C$, write
$Z_i=\cos\theta_i A+\sin\theta_i E_i$, where $E_i$ is standard
normal and independent of $A$, and $0\le\theta_i\le C\sqrt\delta$.
Let $M_{ij}$ be the event that the indicators of $|Z_i|\ge z_j$ and
$|A|\ge z_j$ differ, where $z_j=c_{jt}$, and write $\Delta$ for the change
of a quantity under contraction.
If $Z_i=A$ almost surely, its mismatch events are empty and it is
omitted. Otherwise use this scalar rotation between $A$ and $Z_i$,
holding the outside fixed. The whole cluster, and hence each rotation plane, is independent of $V$.
At intermediate angles the rotation is a positive linear combination
of its endpoints:
\[
 X_i(\theta)=
 \frac{\sin(\theta_i-\theta)A+\sin\theta Z_i}{\sin\theta_i}.
\]
The coefficients sum to at most $1/\cos(\theta_i/2)$.
Because both endpoint correlations with every neighbor are bounded by
$r_0$, its path correlation is at most
$R_\delta=r_0/\sqrt{1-\delta/2}=r_0+O(\delta)$.
The entire path is positive semidefinite, with standard marginals;
the independent reference neighbor has zero correlation with it.
We compare the two increments at ranks up to a cutoff $K$ and at later
ranks. Up to $K$, they can differ only if contraction changes a cluster
threshold indicator while a neighbor has a p-value below $Kt$, and the
next calculation bounds this joint event; later crossings need many small
p-values in $V$, which \eqref{fd:count} makes unlikely.

Put $\chi=\max(d+1,2)$, $s_0=q+n_0$ and
\[
 K=\max\left\{s_0,2,
     \left\lceil(s_0+3\chi)L/I_\alpha\right\rceil\right\},
 \qquad 0<\beta_0<\frac{1-r_0}{1+r_0}.
\]
The crossing probability is weighted by the rotation velocity, which
may depend on the neighbor. We therefore bound the expectation of the absolute velocity on the
marked event directly.
Conditional on a crossing
$X_i(\theta)=x$, write $Y_\ell=rx+\sqrt{1-r^2}\,\xi$, with
$|r|\le R_\delta$. The unit rotation velocity $U$ and $\xi$ are jointly
standard normal conditionally, possibly singular. Regression gives
$\Expect(|U|\mid \xi)\le|\xi|+\sqrt{2/\pi}$.
For $|x|\le c_t$, the mark $|Y_\ell|\ge c_{Kt}$ implies
$|\xi|\ge u=(c_{Kt}-|r|c_t)/\sqrt{1-r^2}>0$ eventually. Hence
\[
 \Expect[|U|\mathbf1\{|Y_\ell|\ge c_{Kt}\}\mid X_i(\theta)=x]
 \le2\phi(u)+2\sqrt{2/\pi}\,\barPhi(u)\le C e^{-u^2/2}.
\]
For $0\le r\le R_\delta$ the quotient
$(c_{Kt}-r c_t)^2/(1-r^2)$ decreases in $r$ while
$R_\delta c_t<c_{Kt}$. Uniformly on $\delta L\le\lambda_*$,
\[
 \frac{(c_{Kt}-R_\delta c_t)^2}{2(1-R_\delta^2)L}
 \longrightarrow\frac{1-r_0}{1+r_0}>\beta_0,
\]
since $\delta\to0$, $c_t/c_{Kt}\to1$ and $c_{Kt}^2/(2L)\to1$.
Integrating Lemma~\ref{fd:marked_crossing}, with path length
$O(\sqrt\delta)$ and $\phi(z_j)\le Cjt\sqrt L$, gives
\begin{equation}\label{fd:early_marked}
 \Prob(M_{ij},|Y_\ell|\ge c_{Kt})
 \le Cjt\sqrt{\delta L}\,t^{\beta_0}.
\end{equation}
The same calculation applies to $Y'$. It permits singular joint
conditional laws and a velocity fully determined by the neighbor.

We first compare the actual and reference increments at early ranks. Set every neighbor p-value
to one without deleting its slot. The changes caused by contraction
then agree as functions of $(C,V)$ for the two coupled backgrounds.
Before this replacement, a difference requires both a cluster mismatch
at some $j\le K$ and a p-value below the cutoff in $Y$ or $Y'$. Summing
\eqref{fd:early_marked} over the bounded cluster and neighbor sets
and over $j\le K$ bounds the difference of early expected increments
by $Ct\sqrt{\delta L}K^2t^{\beta_0}$.

For late ranks, either the late-crossing indicator can change or the
early indicator subtracted from it can change. Let $M_j$ be the union of the
cluster endpoint mismatches at threshold $jt$, and put
\[
 \mathcal O_j=\{N_V(jt)\ge j-s_0\},\qquad \mathcal O_{>K}=\bigcup_{j>K}\mathcal O_j.
\]
For either coupled background, the deterministic bound is
\begin{equation}\label{fd:late_increment}
 |\Delta(H-H_K)|\le
 2\sum_{j>K}\mathbf1_{M_j}\mathbf1_{\mathcal O_j}
 +2\sum_{k\le K}\mathbf1_{M_k}\mathbf1_{\mathcal O_{>K}}.
\end{equation}
To verify it, write $H-H_K=\Lambda_K(1-H_K)$, where $\Lambda_K$ indicates a
crossing at some rank after $K$. A change of $\Lambda_K$ requires a late
mismatch and enough outside coordinates for that crossing. If an
early crossing changes while a late crossing remains present, the
second sum bounds the change. This second possibility is retained
even when every late threshold indicator itself is unchanged.

The unmarked crossing bound gives
$\Prob(M_j)\le Cjt\sqrt{\delta L}$ for all $j\le m$.
Each $M_j$ is a function only of $C$ and is independent of $V$.
Using \eqref{fd:count} in both terms of
\eqref{fd:late_increment}, their total expectation is at most
\[
 Ct\sqrt{\delta L}(K^2+K+1)
       \alpha^{-s_0/\chi}e^{-I_\alpha K/\chi}
 \le Ct\sqrt{\delta L}(K^2+K+1)t^3.
\]
For the last inequality,
$e^{-I_\alpha K/\chi}\le t^{s_0/\chi+3}$ and $\alpha\ge t$.
The geometric prefactors are uniform because $I_\alpha\ge I_{a_1}$.
Combining the early and late estimates for the two couplings yields
\begin{equation}\label{fd:increment_error}
 |(H_1-H_0)-(B_1-B_0)|
 \le Ct\sqrt{\delta L}(K^2+K+1)(t^{\beta_0}+t^3)
 =o(t\sqrt{\delta L}).
\end{equation}
The error contains the same spread-dependent factor as the decrease
under the independent reference. Thus the error remains small relative to that decrease even when
the positive spread tends to zero arbitrarily quickly. The error
coefficient is uniform over the stated class and level interval. Choose one minimum
dimension so that this error is at most half the deficit in
\eqref{fd:surrogate_increment}.
This proves \eqref{fd:dependent_micro_bound}. The projection and
rotation arguments are defined directly for singular Gaussian laws.
For zero spread all aligned coordinates already equal the anchor.
\end{proof}

\subsection{Assembly and equality}
\label{fd:assembly}

\begin{proof}[Proof of Theorem~\ref{fd:theorem}]
We combine the comparisons in three stages: isolate the clusters with
larger rescaled spread, contract the microscopic clusters, and sum the
shares under the final law. Fix $D,a_1$, and construct the clusters of
Lemma~\ref{fd:hierarchy}. Their common gap, cap and size bounds are fixed
functions of $D$, although the local scales may vary with the covariance.
Lemma~\ref{fd:dependent_micro} applies with $q=D+1$,
$n_0=D\lfloor(D+1)^2/4\rfloor$, $d=d_1=D^2$ and $r_0=r_T$.
Choose its fixed microscopic cutoff $\lambda_*>0$, and partition
the selected clusters as
\[
 \mathcal L=\{C:d_CL>\lambda_*\},\qquad
 \mathcal S=\{C:0<d_CL\le\lambda_*\},\qquad
 \mathcal Z=\{C:d_C=0\}.
\]

First isolate every cluster in $\mathcal L$, one at a time. Doing all
isolations before any contraction keeps the original degree bound
available for these comparisons. Each intermediate covariance is a
direct sum of already isolated original principal blocks and the
remaining original principal submatrix.
It is positive semidefinite, retains degree at most $D$, and its
remaining boundary entries are original entries or zero. Internal
spreads are unchanged. Hence Lemma~\ref{fd:isolation} applies at each
step with $d=D$, boundary cap $\rho=1-h_C$ and $B_\partial\le Dq$, and the
uniform relative estimate
\eqref{fd:isolation_relative}, over the common local-scale interval, gives a deterministic $\omega_m\to0$ such that
\begin{equation}\label{fd:assembly_isolation}
 F_{\rm original}\le F_{\rm isolated}
       +\omega_m t\sum_{C\in\mathcal L}G_t(d_C).
\end{equation}

Next contract the clusters in $\mathcal S$ to actual aligned anchors.
For each cluster, let $Y$ contain all its current physical outside
neighbors, and let $V$ be the remaining coordinates. Joint Gaussianity
and the zero cross-covariance imply independence of the entire cluster
and $V$. Equations \eqref{fd:degree_accounting} and
\eqref{fd:neighbor_accounting} give the required common degree and
neighbor bounds. Correlations between the cluster and current outside coordinates are
either original correlations with uncontracted coordinates or retained
anchors, or are zero; hence they are capped by $r_T$. The cluster's spread remains unchanged until that cluster is
contracted. Lemma~\ref{fd:dependent_micro} therefore gives
\begin{equation}\label{fd:assembly_contraction}
 F_{\rm isolated}\le F_{\rm terminal}
                -c_*t\sum_{C\in\mathcal S}\sqrt{d_CL}.
\end{equation}
Contract the clusters in $\mathcal Z$ by exact identities; no bound
dividing by their spread is used.

The terminal law consists of the mutually independent unchanged
clusters in $\mathcal L$, independent also of a capped weighted
Gaussian graph $G$. Its representatives are actual original anchors
and unselected singletons, with degree at most $D$, multiplicities
at most $q$, and cap $r_T$. Every cluster in $\mathcal L$ has size
at most $q$ and spread less than $\varepsilon_0/36<1/16$.
Thus Lemma~\ref{fd:terminal} applies with $W=Q=q$, $d=D$ and
$r=r_T$. Sum \eqref{fd:terminal_bounds} using
\eqref{fd:shares}, under this one terminal law, to obtain
\begin{equation}\label{fd:assembly_terminal}
 F_{\rm terminal}\le\alpha-(1-t^{g_{r_T}/2})D_G^{V,\mathrm{term}}
                  -c_Nt\sum_{C\in\mathcal L}G_t(d_C).
\end{equation}
Choose one $M(D,a_1)$ satisfying all earlier minimum-dimension requirements and
$\omega_m\le c_N/2$ and $t^{g_{r_T}/2}\le1/2$. Equations
\eqref{fd:assembly_isolation}--\eqref{fd:assembly_terminal} prove
\eqref{fd:margin} with $c=\min(c_*,c_N/2)$. Each cluster's comparison
error is bounded by a fixed fraction of its own deficit. Thus the
number of clusters may grow without introducing an error independent
of their spreads.

Any positive selected spread gives strict inequality. If all selected
spreads vanish, no isolation occurred and all contractions were
identities. A nonzero terminal representative edge gives strictness
by Lemma~\ref{fd:terminal}. Otherwise the original vector is exactly
independent groups of signed exact copies. Conversely, let an
independent uniform p-value $U_g$ occupy $w_g$ physical copies and
force that entire group to zero, obtaining rank $R^{[g]}\ge w_g$.
Restoration gives the same rank whenever $U_g\le tR^{[g]}$, and
otherwise that group has no rejection share. Conditional on all
other groups,
\[
 \Expect\left[\frac{w_g}{R^{[g]}}
        \mathbf1\{U_g\le tR^{[g]}\}
        \,\middle|\,\text{other groups}\right]=w_gt.
\]
Here $tR^{[g]}\le\alpha<1$. Summing over groups proves exact
equality at every dimension and level. Every singular limit in the
proof was treated by Gram geometry, Gaussian projection or
zero-preserving covariance interpolation with null marginal
boundaries, so no additional equality case is introduced.
\end{proof}

\subsection{Examples and consequences}
\label{fd:scope}

\subsubsection{A shared obstruction to PRDS}

We verify the two PRDS conditions separately. In each case we find an
increasing event whose conditional probability decreases as the
conditioning p-value increases. One calculation conditions on a fixed
p-value; the other conditions on a p-value being at most a threshold.
The required monotonicity can differ between these two forms of
conditioning.

\begin{lemma}\label{ug:path_prds}
Let \(a=1/4\) and
\[
 P=\begin{pmatrix}1&a&0\\a&1&a\\0&a&1\end{pmatrix}.
\]
The two-sided p-values from \(N(0,P)\) fail both pointwise PRDS
and lower-tail PRDS.
\end{lemma}

\begin{proof}
We first check pointwise PRDS. The eigenvalues of \(P\) are \(1,1\pm\sqrt2a>0\).
Write \(Z_2=aZ_1+aZ_3+\sigma E\), where
\(\sigma^2=1-2a^2=7/8\) and \(Z_1,Z_3,E\) are independent
standard normals. The event
\[
 A=\{|Z_2|\le1\}\ \cup\ \{|Z_3|\le9\}
\]
is increasing in p-value coordinates. Its conditional
probability \(q(x)=\Prob(A\mid Z_1=x)\) is even and smooth.
Differentiation under the Gaussian integral gives
\[
 q''(0)=a^2\int_{|y|>9}\phi(y)\int_{-1}^1
 \left\{\frac{(z-ay)^2}{\sigma^4}-\frac1{\sigma^2}\right\}
 \phi_\sigma(z-ay)\,dz\,dy>0,
\]
since \(|z-ay|>5/4>\sigma\) on the integration region.
As \(q\) is even, \(q(x)\) therefore increases for small positive \(x\).
Conditioning on p-values near one corresponds to absolute Gaussian values
near zero, so \(\Prob(A\mid P_1=u)=q(c_u)\) decreases for \(u\) sufficiently
close to one. This contradicts pointwise PRDS.

The pointwise calculation does not by itself establish failure of the
lower-tail condition, which averages over conditioning values. For that
condition, put \(\tau^2=1-a^2=15/16\).
Given \(Z_3=y\), the distribution of \(Z_2\) has mean \(ay\)
and variance \(\tau^2\). Conditioning also on \(Z_1=0\) replaces
\(\tau\) by \(\sigma\).
For \(\mu>9/4\) and \(s\in[\sigma,\tau]\),
\[
 \frac{d}{ds}\Prob\{|\mu+sE|\le1\}
 =\frac{(\mu-1)\phi((\mu-1)/s)
       -(\mu+1)\phi((\mu+1)/s)}{s^2}>0.
\]
Both arguments exceed one, where \(v\phi(v)\) is strictly
decreasing. The case of negative \(\mu\) follows by symmetry. Integration
over \(|y|>9\) therefore proves \(\Prob(A)>q(0)\).
Let \(X=|Z_1|\), with density \(f_X\) and survival function
\(\overline F_X\), and put
\[
 H(z)=\Expect[q(X)\mid X\ge z].
\]
Direct differentiation yields
\[
 H'(z)=\frac{f_X(z)}{\overline F_X(z)}\{H(z)-q(z)\}>0
\]
for all sufficiently small positive \(z\), by continuity and
\(H(0)>q(0)\). Hence
\(\Prob(A\mid P_1\le u)=H(c_u)\) decreases for \(u\) near one,
contradicting lower-tail PRDS.
\end{proof}

\subsubsection{A graph with common neighbors}

\begin{proof}[Proof of Proposition~\ref{main:fixed_degree_example}]
This example shows that common neighbors need not prevent the
fixed-degree conclusion. For $m=n+3$, use coordinates $X_1,X_2,Y,V_1,\ldots,V_n$ with
\[
 \rho_{X_1X_2}=1-\delta_m,\qquad
 \rho_{X_1Y}=\rho_{X_2Y}=1/4,\qquad
 \delta_m=\{8192\sqrt{\log(m+1)}\}^{-1}.
\]
Give consecutive coordinates on the path $Y,V_1,\ldots,V_n$
correlation $1/8$, and set every other off-diagonal entry to zero.
The graph is connected and has degree at most three. In the
symmetric and antisymmetric coordinates $(X_1\pm X_2)/\sqrt2$ of the
first pair, the antisymmetric variance is $\delta_m>0$ and that
coordinate is
independent of the rest. The symmetric variance is $2-\delta_m$.
Taking the Schur complement of the symmetric coordinate changes only
the diagonal entry at $Y$, to
\[
 1-\frac{2(1/4)^2}{2-\delta_m}\ge\frac{29}{31}
 \qquad(\delta_m\le1/16).
\]
The remaining symmetric path matrix is strictly diagonally dominant
with positive diagonal, so it is positive definite. Hence the
original matrix is positive definite.

The pair's two correlations with its common neighbor $Y$ have total
absolute value $1/2$. Its spread tends to zero but
$\delta_m\log(m/\alpha)\to\infty$ at every fixed positive level.
Thus the pair's boundary correlations are not bounded by a fixed multiple
of $\sqrt{\delta_m}$, and its rescaled spread $\delta_m\log(m/\alpha)$ leaves
every bounded interval; Theorem~\ref{main:fixed_degree} needs neither
property. Connectedness and positive
definiteness prevent a decomposition into independent bounded blocks
or compression into exact-copy representatives with a fixed
correlation cap.

For $n\ge3$, change the path correlations $V_1V_2,V_2V_3$ to $1/4$.
Strict diagonal dominance is preserved, so positive definiteness
still holds. The marginal on $V_1,V_2,V_3$ is exactly the path law
in Lemma~\ref{ug:path_prds}. Its increasing event, ignoring all other
coordinates, retains the same conditional probabilities and violates
both pointwise and lower-tail PRDS near conditioning p-value one.
This proves the final assertion without changing the degree bound.
\end{proof}

The coordinate-share construction in
Proposition~\ref{main:coordinate_share}, with independent test coordinates added,
instead has covariance degree at most two and witnesses both PRDS
failures below the actual testing level. Its covariance and witness
may depend on $(m,\alpha)$; this is distinct from the near-one
marginal witness above.

\subsubsection{Consequences across dimensions}

We next obtain a single small level that works in every dimension for
each fixed degree bound. To combine the finitely many smaller dimensions with the
large-dimension result, fix one ceiling $a_1$ and its
minimum dimension $M(D,a_1)$. For $2\le j<M(D,a_1)$, let $a_j>0$ be the
covariance-uniform cutoff of Theorem~\ref{main:uniform} in dimension $j$;
the symbol $a_1$ remains the chosen ceiling. Then
\[
 a_*(D)=\min\bigl(\{a_1\}\cup\{a_j:2\le j<M(D,a_1)\}\bigr)>0
\]
gives validity in every dimension of the fixed-degree class for
$0<\alpha\le a_*(D)$; dimension one is immediate, and if
$M(D,a_1)\le2$ the second set is empty and $a_*(D)=a_1$. Throughout this interval, the small-level theorem gives the signed-copy
equality cases in dimensions $2\le m<M(D,a_1)$ (in dimension one,
equality always holds and the single coordinate forms one group), and the
fixed-degree theorem gives the same characterization in the remaining
dimensions.
The cutoff remains unevaluated and is not asserted uniformly over $D$.

There is a second, qualitative consequence: the allowed degree can
increase with dimension at some unspecified rate. We postpone each
increase until the dimension is large enough for the fixed-degree
theorem at that new degree. For example, set $N_1=1$, choose increasing integers
$N_q\ge\max\{q+1,M(q,1-1/q),N_{q-1}+1\}$, $q\ge2$, and set
$D_m=\max\{q:N_q\le m\}$ for $m\ge N_2$. For every fixed ceiling
$a_1<1$, eventually $1-1/D_m\ge a_1$ and
$m\ge M(D_m,1-1/D_m)$, so the theorem applies to degree at most
$D_m$. This nonconstructive consequence provides no prescribed
growth rate and does not justify substituting $D=m-1$ into the
fixed-degree theorem.

\clearpage
\setlength{\bibsep}{0pt}
\begingroup
\interlinepenalty=10000

\endgroup
\end{document}